\documentclass[12pt]{article}

\usepackage{amsmath,amssymb,amsthm}
\usepackage{authblk}
\usepackage{enumitem}
\usepackage{graphicx}
\usepackage{mathtools}
\usepackage{overpic}
\usepackage{subcaption}
\usepackage{tikz}
\usepackage{vmargin}

\usetikzlibrary{decorations.markings}
\usetikzlibrary{hobby}
\usetikzlibrary{positioning, arrows.meta, decorations.markings, calc, intersections}
\tikzset{arrowmark/.style={postaction={decorate, decoration={markings, mark=at position #1 with {\arrow{>};}}}},arrowmark/.default={.5}}

\newtheorem{theorem}{Theorem}
\newtheorem{proposition}[theorem]{Proposition}
\newtheorem{corollary}[theorem]{Corollary}
\newtheorem{lemma}[theorem]{Lemma}

\numberwithin{theorem}{section}

\theoremstyle{definition}
\newtheorem{definition}[theorem]{Definition}
\newtheorem{remark}[theorem]{Remark}

\numberwithin{equation}{section}

\DeclareMathOperator{\diag}{diag}
\DeclareMathOperator{\Prob}{Prob}
\DeclareMathOperator{\Res}{Res}
\DeclareMathOperator{\supp}{supp}
\DeclareMathOperator{\Tr}{Tr}
\DeclareMathOperator{\Log}{Log}

\let\Re\undefined
\let\Im\undefined
\let\mod\undefined
\DeclareMathOperator{\Re}{Re}
\DeclareMathOperator{\Im}{Im}
\DeclareMathOperator{\mod}{mod}

\makeatletter
\renewcommand{\@makefntext}[1]{\noindent #1}
\makeatother

\title{Asymptotic analysis of a $3\times3$-doubly periodic random lozenge tiling model of a hexagon}

\author{Max J. van Horssen \& Arno B.J. Kuijlaars {\small \authorcr \ \authorcr Department of Mathematics, KU Leuven, \authorcr Celestijnenlaan 200 B bus 2400, 3001 Leuven, Belgium \authorcr max.vanhorssen@kuleuven.be, arno.kuijlaars@kuleuven.be}}

\begin{document}

\maketitle

\begin{abstract}
    We study a random lozenge tiling model of a regular hexagon with weights that are doubly periodic with period $3$. This model features a smooth phase with six cusp points. We obtain explicit descriptions of the limit shape and the local scaling limits of the correlation kernel. Both limits are expressed as contour integrals on the Harnack curve that underlies the model, or on its double cover that is obtained by gluing together two copies along the support of an equilibrium measure in an external field.

    Our asymptotic analysis is based on the double contour integral representation of the correlation kernel due to Duits and Kuijlaars. This double contour integral contains a factor that is given in terms of the solution of a Riemann--Hilbert problem for matrix-valued orthogonal polynomials. The asymptotic behavior of this Riemann--Hilbert problem was recently studied by one of the authors using the Deift--Zhou method of steepest descent. In this work, we complement this approach with the asymptotic analysis of the double contour integral. The integrals are moved to the double cover of the Harnack curve. We analyze the effect of the first two transformations of the Deift--Zhou method, and deform the contours to the relevant saddle points. We visualize our asymptotic analysis on the amoeba of the Harnack curve, in the spirit of the work of Berggren and Borodin on the doubly periodic Aztec diamond.
\end{abstract}

\begingroup
\renewcommand{\thefootnote}{}
\endgroup


\section{Introduction}


\subsection{Overview}
This paper is devoted to the asymptotic analysis of a random lozenge tiling model of a regular hexagon with weights that are doubly periodic with period $3$. Such models have attracted considerable attention in the last decade because of the fascinating coexistence of three asymptotic phases, called frozen, rough, and smooth (or solid, liquid, and gaseous). In particular, the appearance of the smooth phase is a relative new phenomenon.

Significant progress has been made on the closely related model of domino tilings of the Aztec diamond with doubly periodic weights \cite{Bai23, BCJ18, Ber21, BB25, BB26+, BD23, CDF26+, CJ16, DK21, JM23}, and extensions thereof \cite{BBG26+, BB26} that also involve the more general Fock's weights \cite{BCdT23,BdT27}.

In models with doubly periodic weights, an essential role is played by a spectral curve that is known to be a Harnack curve \cite{KOS06}. The genus of the Harnack curve is equal to the number of smooth phases in the case of the Aztec diamond \cite{BB25}.

An important feature in the hexagon case is that the Harnack curve is not enough to describe the geometry of the model, but one has to go to a double cover of the associated Riemann surface $\mathcal{R}$ \cite{BB26}. The double cover $\widehat{\mathcal{R}}$ is formed from two copies of $\mathcal{R}$ by cross-gluing them along a suitable cut (or cuts). A central question is where to place this cut.

Our answer to this question is that the cut should be the support of the equilibrium measure in a certain external field on $\mathcal{R}$. This insight originates from earlier works \cite{BGK23,DK21,Kui25}, where the random tiling model was related to matrix-valued orthogonal polynomials (MVOPs). The equilibrium measure will also be related to the asymptotic behavior of the MVOPs as their degrees tend to infinity.

The equilibrium problem was solved for a restricted class of doubly periodic weights in \cite{Kui25}. In this paper, we build on this work and consider the same class of weights. For this class, the Harnack curve has genus one, and the equilibrium measure is supported on a non-contractible closed contour (a cycle). The double cover is unramified and remains of genus one. Other choices of weights or variations in the dimensions of the hexagon may lead to a support that becomes an arc, or the union of a cycle with an arc. Then, the double cover is ramified as discussed in \cite{BB26}. We intend to study such situations and the transitions leading to them in future work, see also the outlook in Section~\ref{sec:outlook}.

The equilibrium measure is the minimizer of an equilibrium problem on a contour on $\mathcal{R}$. The minimization problem gives Euler--Lagrange variational conditions that are satisfied by the minimizer. These variational conditions consist of an equality on the support, and a corresponding inequality on the remaining part of the contour. The equality provides the necessary algebraic and geometric structure on the double cover $\widehat{\mathcal{R}}$, and we make use of it to construct a differential $\omega$ that is meromorphic on $\widehat{\mathcal{R}}$. We emphasize that $\omega$ is not meromorphic on the original Riemann surface $\mathcal{R}$. The inequality is essential for the success of the asymptotic analysis.

The equilibrium problem is not needed for the doubly periodic Aztec diamond. Here, one can instead use Wiener--Hopf techniques for the factorization of the weight matrix \cite{Ber21, BD19, BD23}. These techniques fail for the hexagon as explained in \cite{KP25}. 

\begin{figure}[t]
	\centering
    \hspace{-0.2cm}
	\raisebox{-0.19cm}{\includegraphics[width=0.32\linewidth]{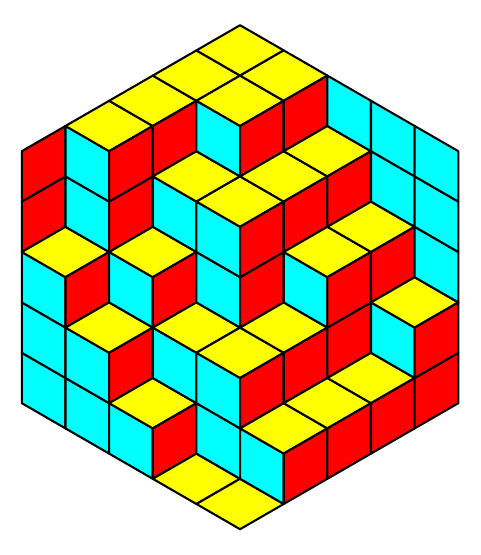}}
    \hspace{0.1cm}
    \includegraphics[width=0.3\linewidth]{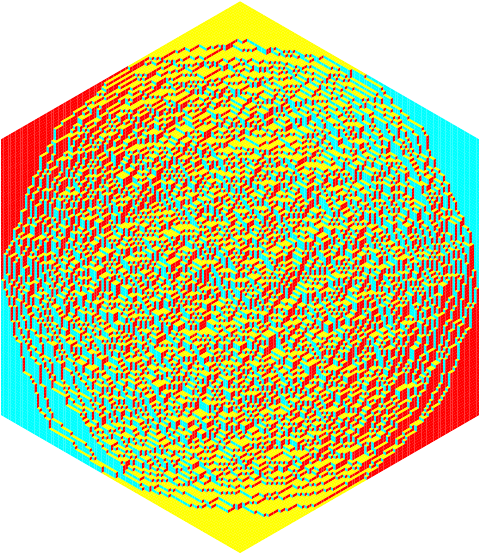}
    \hspace{0.3cm}
    \includegraphics[width=0.3\linewidth]{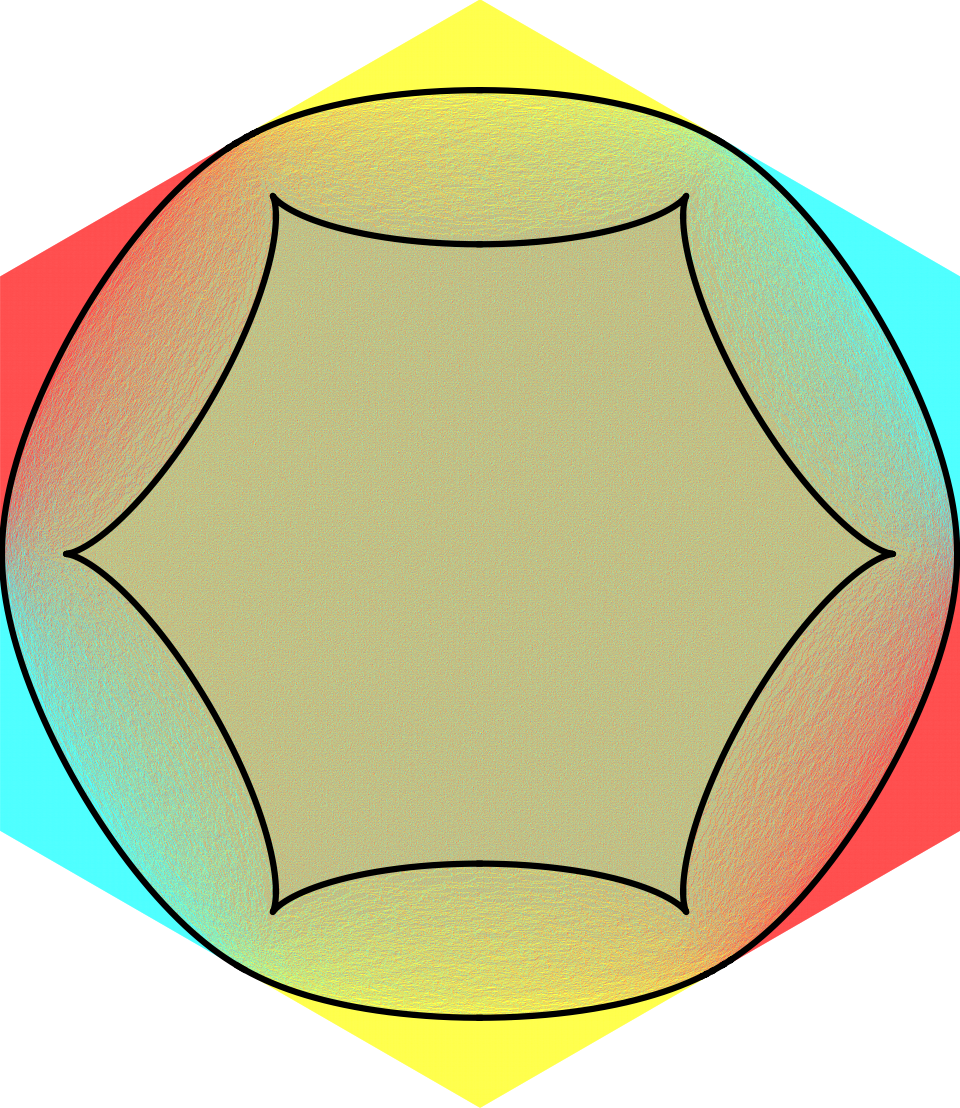}
    \hspace{0.2cm}
    \caption{Random lozenge tilings of a regular hexagon. The left subfigure shows a small-size example where the individual lozenges are still visible. The middle subfigure shows a large-size random tiling with uniform probability
    in which one can clearly observe the arctic circle phenomenon. The right subfigure shows a large-size random tiling with the $3\times3$-doubly periodic weights under consideration. In addition to the frozen and rough phases, there is a smooth phase appearing in the middle. The boundaries between the phases are called the arctic curves, which are also shown in the right subfigure.}
    \label{fig:tiling_simulations_1}
\end{figure}


\subsection{Random tiling model}
We consider large random tilings of a regular hexagon with three types of lozenges, see the left subfigure of Figure~\ref{fig:tiling_simulations_1}. As the size of the system increases, one observes the arctic circle phenomenon: near each of the corners there is a frozen phase with only one type of lozenge (with very high probability), while in the middle all three lozenges appear.

The middle subfigure is for the uniform model, which is by now very well-studied, see e.g.~\cite{Agg23, Gor21, Hua26+, Joh02, Pet14, Pet15}, and the many references cited therein. The right subfigure shows a random tiling with $3\times3$-doubly periodic weights that we consider in this paper. In addition to the frozen and rough phases, there is a smooth phase with six cusp points that appears in the middle. The subfigure also contains the arctic curves separating the three asymptotic phases.

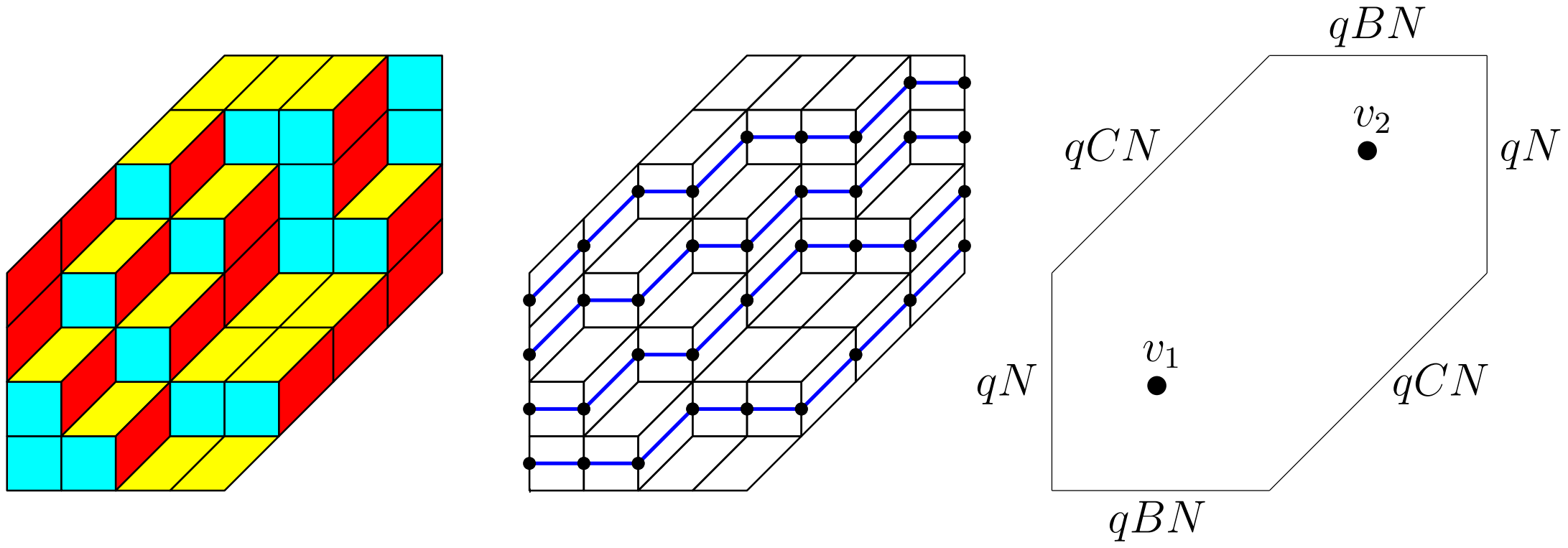
\begin{figure}[t]
    \centering
    \hspace{0.25cm}
	\includegraphics[width=0.3\linewidth]{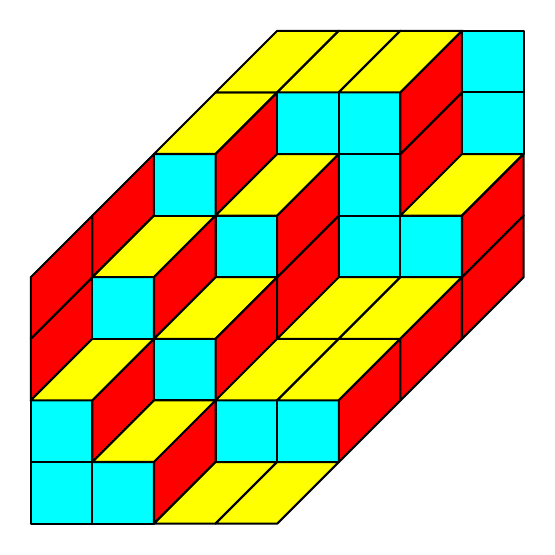}
    \hfill
	\includegraphics[width=0.3\linewidth]{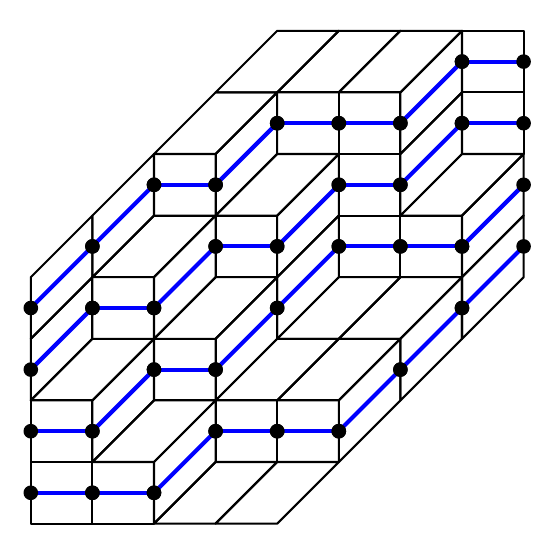}
    \hfill
    \begin{overpic}[width=0.3\linewidth]{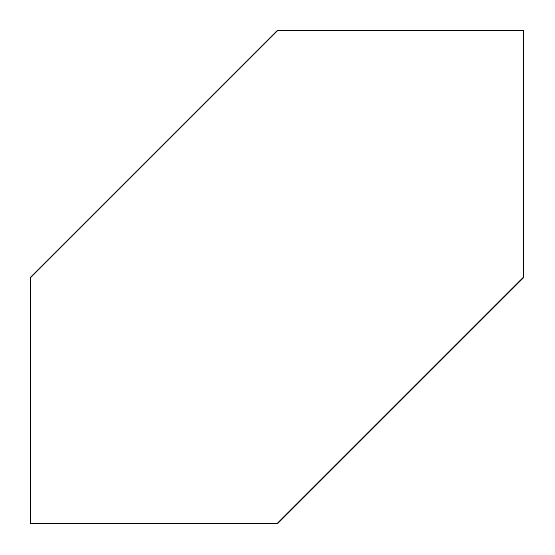}
        \put(27,27){\circle*{4}}
        \put(24,32){$v_1$}
        \put(70,75){\circle*{4}}
		\put(67,80){$v_2$}
        \put(-10,25){$qN$}
        \put(8,73){$qCN$}
        \put(62,98){$qBN$}
        \put(97,73){$qN$}
        \put(75,25){$qCN$}
        \put(17,-3){$qBN$}
	\end{overpic}
    \hfill
    \hspace{0.25cm}
    \caption{A lozenge tiling of a hexagon (left subfigure) and the corresponding system of non-intersecting paths (middle subfigure) together with the random particles. The right subfigure shows the dimensions of the hexagon as in \eqref{eq:hexagon_dimensions}--\eqref{eq:BC_condition}, and it illustrates the vertices $v_1$ and $v_2$ from \eqref{eq:v1_v2_pq}.}
    \label{fig:tiles_paths_dimensions}
\end{figure}

Our periodic weights can be described in terms of an equivalent dimer model as in \cite{KOS06}. However, we present and study it from the point of view of non-intersecting paths. Any tiling of a hexagon is equivalent to a system of non-intersecting paths as indicated in Figure~\ref{fig:tiles_paths_dimensions}. It is convenient to apply an affine transformation to the hexagon, and let the paths go from left to right, with horizontal steps and steps that go diagonally up.

Paths are then on the directed graph $(V,E)$ with vertex set $V = \mathbb{Z}^2$ and edge set $E=E_1\cup E_2$, where $E_1=\{(v,v+e_1)\mid v\in V\}$ and $E_2=\{(v,v+e_1+e_2) \mid v\in V\}$ are the sets of horizontal edges and diagonal edges, respectively. We use $e_1=(1,0)$ and $e_2=(0,1)$. Putting particles on the paths as shown in the middle subfigure of Figure~\ref{fig:tiles_paths_dimensions}, one also obtains a random particle system (in case of a random tiling).

A weighting $w \colon E \to \mathbb{R}_{>0}$ of the edges is described entirely by two infinite matrices $(a_{j,k})_{j,k\in\mathbb{Z}}$ and $(b_{j,k})_{j,k\in\mathbb{Z}}$ of positive numbers, namely 
\begin{equation}\label{eq:ab_weights}
	\begin{aligned}
        a_{j,k} &= w((v,v+e_1)), \\
        b_{j,k} &= w((v,v+e_1+e_2)),
	\end{aligned}
    \qquad \text{if } v=(j,k)\in\mathbb{Z}^2.
\end{equation}
The weighting is called $p,q$-doubly periodic (with $p,q\geq1$) if
\begin{equation}\label{eq:ab_periodic}
    a_{j+np,k+mq} = a_{j,k}, \quad b_{j+np,k+mq} = b_{j,k}
\end{equation}
for all $j,k,n,m\in\mathbb{Z}$. Such a weighting has period $p$ in the horizontal direction and period $q$ in the vertical direction. 

Any weighting induces a non-uniform probability measure on tilings in the following way. A tiling $\mathcal{T}$ gives rise to the non-intersecting path system $\mathcal{P}$ on the graph $(V,E)$ and then by definition, we assign the probability
\begin{equation}\label{eq:T_prob}
    \Prob\left[\mathcal{T}\right] = \frac{1}{Z}\prod_{e\in\mathcal{P}} w(e)
\end{equation}
to the tiling $\mathcal{T}$. Here, $Z$ is the normalizing constant (or partition function).

The weighting also gives rise to a random particle system, a subset of the vertices lying within the hexagon, which is known to be a determinantal point process. This means that there exists a kernel $K \colon V \times V \to \mathbb{R}$ such that for any $k\geq1$ and any collection of distinct vertices $v_1,\ldots,v_k$, one has 
\begin{equation}\label{eq:K_determinantal}
    \Prob\!\left[
    \begin{array}{l}
        \text{there are particles} \\
    	\text{at each $v_1, \ldots, v_k$}
    \end{array}
   \right] = \det\bigl[K(v_i,v_j)\bigr]_{i,j=1}^k.
\end{equation}
This is a remarkable result, the proof of which relies on a combination of the Lindstr\"om--Gessel--Viennot lemma \cite{GV85,Lin73} from combinatorics and the Eynard--Mehta theorem \cite{EM98} from random matrix theory.


\subsection{Periodic weights}
For the case of periodic weights \eqref{eq:ab_weights}--\eqref{eq:ab_periodic}, Duits and Kuijlaars \cite{DK21} gave a double contour integral representation for the correlation kernel $K$ of the determinantal point process. We now describe this formula in detail.

We assume, for convenience, that the hexagon is of size
\begin{equation}\label{eq:hexagon_dimensions} 
	qN \times qBN \times qCN, \quad N\in\mathbb{N},
\end{equation}
see also the right subfigure of Figure~\ref{fig:tiles_paths_dimensions}. A regular hexagon corresponds to $B=C=1$. We further assume the parameters $B$ and $C$ are such that 
\begin{equation}\label{eq:BC_condition}
	CN\in\mathbb{N}, \quad \frac{q}{p}(B+C)N\in\mathbb{N}.
\end{equation}
The assumptions \eqref{eq:hexagon_dimensions}--\eqref{eq:BC_condition} mean that the dimensions of the hexagon are well adjusted to the periodicity parameters $p$ and $q$.

We consider $N$ as the large-size parameter of the system and write $K_N$ for the correlation kernel.

In order to state the formula from \cite{DK21}, we need the following $q \times q$ matrices, called transition matrices,
\begin{equation}\label{eq:Tjz}
    T_j(z) = \begin{pmatrix} a_{j,0} & b_{j,0} & 0 & \cdots & \cdots & 0 \\ 0 & a_{j,1} & b_{j,1} & 0 & \ddots & \vdots \\ \vdots & \ddots & \ddots & \ddots & \ddots & \vdots \\ \vdots & \ddots & \ddots & \ddots & \ddots & 0 \\ 0 & \cdots & \cdots & 0 & a_{j,q-2} & b_{j,q-2} \\ b_{j,q-1} z & 0 & \cdots & \cdots & 0 &   a_{j,q-1} \end{pmatrix}, \quad j=0,1, \ldots, p-1,
\end{equation}
that contain the weights \eqref{eq:ab_weights}. We also need their products
\begin{align}
    T_{0,j}(z) &= T_0(z)\cdots T_{j-1}(z), \quad j=0, \ldots, p-1, \label{eq:T0jz} \\
    W(z) &= T_0(z)\cdots T_{p-1}(z), \label{eq:Wz}
\end{align}
with the understanding that $T_{0,0}(z) = I_q$.

Then, for two vertices
\begin{equation}\label{eq:v1_v2_pq}
    v_1=(px_1+j_1,qy_1+k_1), \quad v_2=(px_2+j_2,qy_2+k_2),
\end{equation}
with $x_1,x_2,y_1,y_2\in\mathbb{Z}$, $j_1,j_2\in\{0,1,\ldots,p-1\}$, and $k_1,k_2\in\{0,1,\ldots,q-1\}$, we have the following formula (see \cite[Theorem~4.7]{DK21}):
\begin{multline}\label{eq:KN_v1_v2}
    K_N(v_1,v_2) = \left[-\frac{\chi_{px_1+j_1>px_2+j_2}}{2\pi i}\oint_\mathbb{T} T_{0,j_2}^{-1}(z) W(z)^{x_1-x_2} T_{0,j_1}(z) \frac{dz}{z^{y_1-y_2+1}}\right]_{k_2,k_1} \\
    + \left[\frac{1}{(2\pi i)^2}\oint_\mathbb{T}\oint_\mathbb{T}T_{0,j_2}^{-1}(z_2) \frac{W(z_2)^{\frac{q}{p}(B+C)N-x_2}}{z_2^{CN+N-y_2}}R_N(z_2,z_1)\frac{W(z_1)^{x_1}}{z_1^{y_1}}T_{0,j_1}(z_1) \frac{dz_1dz_2}{z_1}\right]_{k_2,k_1},
\end{multline}
where $\mathbb{T} = \{z \mid |z| = 1\}$ denotes the unit circle in the complex plane. All integrals in \eqref{eq:KN_v1_v2} are taken entrywise and $[\cdots]_{k_2,k_1}$ means that $(k_2,k_1)$-entry is taken. We use the convention that the rows and columns of $q \times q$ matrices are indexed by $\{0,1, \ldots, q-1\}$.

The double contour integral in \eqref{eq:KN_v1_v2} contains the factor $R_N(z_2,z_1)$ that we have not explained yet. It is given in terms of the solution of a Riemann--Hilbert (RH) problem. The RH problem asks for a $2q \times 2q$ matrix-valued function $\mathbf{Y}_N$ such that
\begin{itemize}
	\item $\mathbf{Y}_N \colon \mathbb{C}\setminus\mathbb{T} \to \mathbb{C}^{2q \times 2q}$ is analytic,
	\item $\mathbf{Y}_N$ has limiting values on $\mathbb{T}$, denoted by $\mathbf{Y}_{N,+}$ and $\mathbf{Y}_{N,-}$ ($\mathbf{Y}_{N,+}$ is the limiting value from the interior of the unit circle, and $\mathbf{Y}_{N,-}$ the limiting value from the exterior), that satisfy
    \begin{equation}\label{eq:YN_jump}
        \mathbf{Y}_{N,+}(z) = \mathbf{Y}_{N,-}(z) \left(\begin{matrix} \text{\raisebox{3pt}{$I_q$}} & \dfrac{W(z)^{\frac{q}{p}(B+C)N}}{z^{(1+C)N}} \\ 0_q & I_q \end{matrix}\right)\!, \quad z\in\mathbb{T},
    \end{equation}
	\item as $z\to\infty$,
    \begin{equation}\label{eq:YN_at_infty}
        \mathbf{Y}_N(z) = \left(I_{2q} + \mathcal{O}(z^{-1})\right) \begin{pmatrix} z^N I_q & 0_q \\ 0_q & z^{-N} I_q \end{pmatrix}.
    \end{equation}
\end{itemize}
There is a unique solution to this RH problem, and
\begin{equation}\label{eq:RN_formula}
    R_N(z_2,z_1) = \frac{1}{z_1-z_2} \begin{pmatrix} 0_q & I_q \end{pmatrix} \mathbf{Y}_N^{-1}(z_2) \mathbf{Y}_N(z_1) \begin{pmatrix} I_q \\ 0_q \end{pmatrix}.
\end{equation}
The function $R_N$ is $q \times q$ matrix-valued and its entries are polynomial in $z_1$ and $z_2$ of degree at most $N-1$ in both variables.

The solution $\mathbf{Y}_N$ of the RH problem is related to the matrix-valued orthogonal polynomial (MVOP) $P_N(z)$ of degree~$N$ that satisfies the orthogonality relations
\begin{equation}\label{eq:PN_orthogonal} 
	\frac{1}{2\pi i}\oint_\mathbb{T} P_N(z) \frac{W(z)^{\frac{q}{p}(B+C)N}}{z^{(1+C)N}} z^k dz = 0_q, \qquad k = 0,1, \ldots, N-1,
\end{equation}
through the relation
\begin{equation}\label{eq:PNz}
    P_N(z) = \begin{pmatrix} I_q & 0_q \end{pmatrix} \mathbf{Y}_N(z) \begin{pmatrix} I_q \\ 0_q \end{pmatrix}.
\end{equation}
The orthogonality relations \eqref{eq:PN_orthogonal} are of non-Hermitian type with a matrix-valued weight on the unit circle that varies with $N$. The matrix-valued weight is a rational function in $z$ with only a pole at $z=0$. This is due to the assumptions \eqref{eq:hexagon_dimensions}--\eqref{eq:BC_condition} on the parameters. The integral \eqref{eq:PN_orthogonal} is thus determined by its residue at $z=0$. By Cauchy's theorem, $\mathbb{T}$ can be replaced by any contour that goes around $0$ once in the positive direction.

The RH problem is the matrix-valued analogue of the well-known RH problem for orthogonal polynomials due to Fokas, Its, and Kitaev \cite{FIK92}. The extension to MVOPs was given in \cite{CM12,GIM11}.


\subsection{Strategy for the asymptotic analysis}
The formulas \eqref{eq:KN_v1_v2}--\eqref{eq:RN_formula} suggest a research program for an asymptotic analysis of the kernel \eqref{eq:KN_v1_v2} that proceeds in two steps. In the first step, one analyzes the RH problem in the limit $N\to\infty$ by means of the Deift--Zhou method of steepest descent \cite{DZ93}. This method has been very successful for the asymptotic analysis of scalar orthogonal polynomials, see \cite{Dei99, DKMVZ99}, where the fundamental steps of the method were first worked out in the context of orthogonal polynomials.

If the Deift--Zhou method can be made to work, then one may have control over the large-$N$ behavior of $R_N$ due to \eqref{eq:RN_formula}. In the second step, one then analyzes the double contour integral, which is the main term in \eqref{eq:KN_v1_v2}, by means of a steepest descent analysis of integrals.

We would like to carry out this program for the lozenge tilings of a hexagon for all periodicities $p$ and $q$, for all possible sizes of the hexagon (encoded by parameters $B$ and $C$), and all possible $p \times q$-doubly periodic weights.

For the periodicities $p=q=2$, the MVOPs reduce to scalar orthogonal polynomials \cite{Cha21,CDKL20,GK21}, and then the RH problem can be analyzed in the large-$N$ limit. However, this case does not give rise to a smooth phase in the random tiling model.

For $p=q=3$, $B=C=1$, and a restricted class of weights, the RH analysis was done successfully in the recent paper \cite{Kui25}. In this paper, we continue with the second step, namely the steepest descent analysis of the double contour integral for the special cases of \cite{Kui25}. We actually do not use the outcome of the RH analysis in a direct way. That is, we do not derive an asymptotic formula for $R_N$ that we then use in \eqref{eq:KN_v1_v2} to find the large-$N$ behavior of the double contour integral in \eqref{eq:KN_v1_v2}. Rather, we follow the first two steps of the RH analysis in \cite{Kui25} and study their effect on the double contour integral. The steepest descent analysis of \cite{Kui25} is ultimately used to show in Lemma~\ref{lem:TN_bounded} that certain factors remain bounded as $N\to\infty$, and as a consequence do not contribute to the large-$N$ behavior of the correlation kernel.

Before we state our results on the limiting behavior of $K_N$, we first discuss the geometry that underlies the problem.


\section{Geometry of the problem}


\subsection{Assumptions}
As already mentioned, we consider the case $p=q=3$, and the same positive weights $a_{j,k}$ and $b_{j,k}$ as in \cite{Kui25}. That is, the weights satisfy the periodicity condition \eqref{eq:ab_periodic},
\begin{equation}\label{eq:ab_assumption_1}
    \begin{aligned}
        a_{0,k} a_{1,k} a_{2,k} &= 1, \quad k=0,1,2, \\
        b_{0,k} b_{1,k+1}, b_{2,k+2} &= 1, \quad k=0,1,2, \\
        \frac{a_{j,0} a_{j,1} a_{j,2}}{b_{j,0} b_{j,1} b_{j,2}} &= 1, \quad j=0,1,2,
    \end{aligned}
\end{equation}
and
\begin{equation}\label{eq:ab_assumption_2}
    a_{0,0} b_{1,0} \neq a_{1,1} b_{0,0} \quad \text{or } \quad a_{0,1} b_{1,1} \neq a_{1,2} b_{0,1},
\end{equation}
see \cite[(1.4)--(1.7)]{Kui25}. Note that \cite{Kui25} uses the numbering $1$, $2$, and $3$ for the indices $j$ and $k$, while in this paper we are using $0$, $1$, and $2$.

 
\subsection{Harnack curve}
For any choice of parameters, the characteristic equation
\begin{equation}\label{eq:spectral_curve0}
    P(z,\lambda) = \det\left(\lambda I_q - W(z)\right) = 0 
\end{equation}
defines a Harnack curve \cite{KO06,KOS06}. Recall that $W(z)$ is defined in \eqref{eq:Wz}. This is a real algebraic curve with the special property that the amoeba map
\begin{equation}\label{eq:amoeba_map}
    \Log \colon (z,\lambda) \mapsto (\log|z|,\log|\lambda|) 
\end{equation}
is at most two-to-one on the curve. This map is precisely one-to-one on the real locus of \eqref{eq:amoeba_map}. The Harnack curve property implies that the eigenvalues $\lambda_j(z)$ of $W(z)$ can be ordered such that
\[
    |\lambda_1(z)|>|\lambda_2(z)|>\cdots>|\lambda_q(z)|>0 \quad \text{for } z\in\mathbb{C}\setminus\mathbb{R}.
\] 
This ordering provides a natural sheet structure on the Riemann surface associated with \eqref{eq:spectral_curve0}, where we take $\lambda = \lambda_k$ on the $k$th sheet.

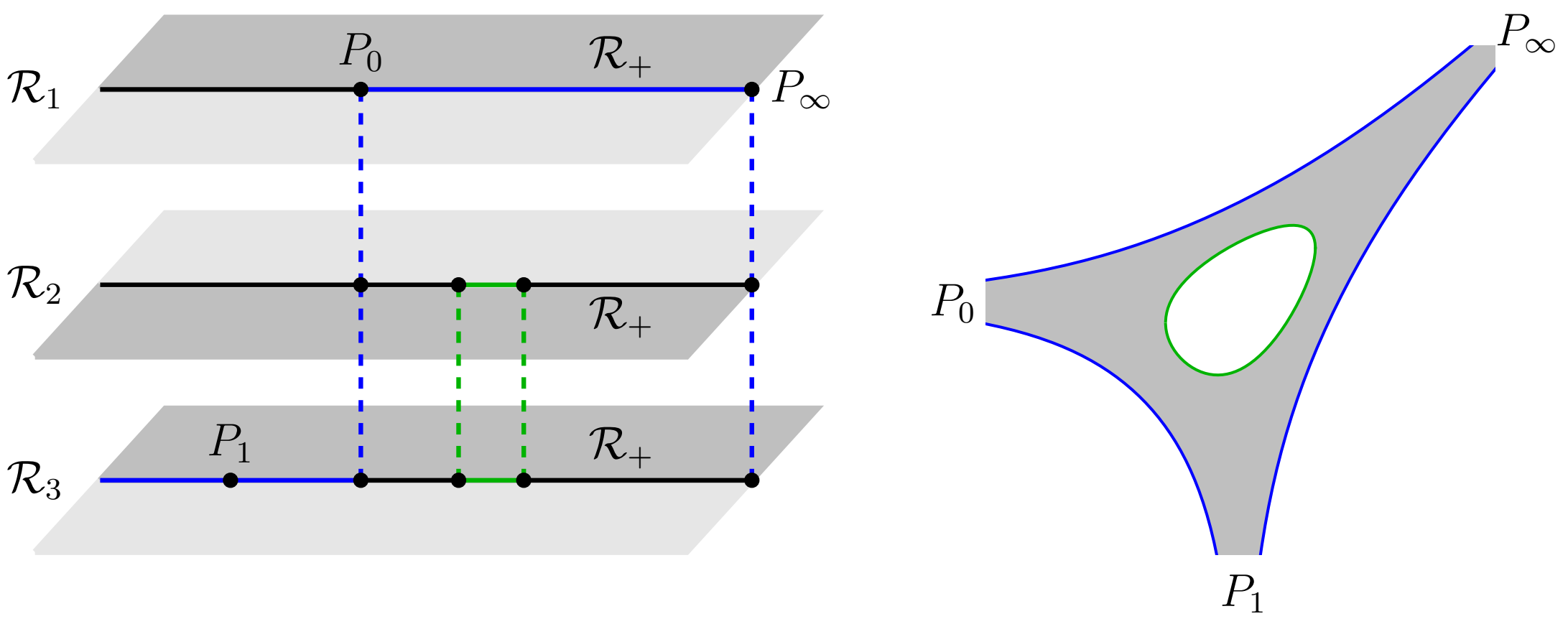
\begin{figure}[t]
    \centering
    \hspace{0.1cm}
    \begin{subfigure}[t]{0.6\textwidth}
        \begin{tikzpicture}[scale=0.6](15,10)(0,0)
        \begin{scope}[ultra thick,decoration={markings,mark=at position 0.5 with {\arrow{>}}}]
        \filldraw[gray!20!white] (6,2.9) --++(10,0) --++(1,1.1) --++(-10,0) --++(-1,-1.1);
        \filldraw[gray!50!white] (7,4.0) --++(10,0) --++(1,1.1) --++(-10,0) --++(-1,-1.1);
        \filldraw[gray!50!white] (6,-0.1) --++(10,0) --++(1,1.1) --++(-10,0) --++(-1,-1.1);
        \filldraw[gray!20!white] (7,1.0) --++(10,0) --++(1,1.1) --++(-10,0) --++(-1,-1.1);
        \filldraw[gray!20!white] (6,-3.1) --++(10,0) --++(1,1.1) --++(-10,0) --++(-1,-1.1);
        \filldraw[gray!50!white] (7,-2.0) --++(10,0) --++(1,1.1) --++(-10,0) --++(-1,-1.1);

        \draw[dashed,help lines,very thick,blue] (11,1)--(11,4);
        \draw[dashed,help lines,very thick,blue] (17,1)--(17,4);

        \draw[dashed,help lines,very thick,blue] (11,-2)--(11,1);
        \draw[dashed,help lines,very thick,green!70!black] (12.5,-2)--(12.5,1);
        \draw[dashed,help lines,very thick,green!70!black] (13.5,-2)--(13.5,1);
        \draw[dashed,help lines,very thick,blue] (17,-2)--(17,1);

        \draw[very thick] (7,4)--++(4,0);
        \draw[very thick,blue] (11,4)--++(6,0);
        \draw[very thick] (7,1)--++(5.5,0);
        \draw[very thick] (13.5,1)--++(3.5,0);
        \draw[very thick,green!70!black] (12.5,1)--++(1,0);
        \draw (9,-2) node[above] {$P_1$};
        \draw[very thick] (11,-2)--++(1.5,0);
        \draw[very thick] (13.5,-2)--++(3.5,0);
        \draw[very thick,blue] (7,-2)--++(4,0);
        \draw[very thick,green!70!black] (12.5,-2)--++(1,0);

        \filldraw (11,4) circle (2pt);
        \filldraw (17,4) circle (2pt);
        \draw (11,4) node[above] {$P_0$};
        \draw (17,4) node[right] {$P_\infty$};
        \filldraw (11,1)  circle (2pt);
        \filldraw (12.5,1)  circle (2pt);
        \filldraw (13.5,1)  circle (2pt);
        \filldraw (17,1)  circle (2pt);
        \filldraw (9,-2) circle (2pt);
        \filldraw (11,-2) circle (2pt);
        \filldraw (12.5,-2) circle (2pt);
        \filldraw (13.5,-2) circle (2pt);
        \filldraw (17,-2) circle (2pt);

        \draw (6,4) node {$\mathcal{R}_1$};
        \draw (6,1) node {$\mathcal{R}_2$};
        \draw (6,-2) node {$\mathcal{R}_3$};
        \draw (15,4.5) node {$\mathcal{R}_+$};
        \draw (15,0.5) node {$\mathcal{R}_+$};
        \draw (15,-1.5) node {$\mathcal{R}_+$};
        \end{scope}
        \end{tikzpicture}
    \end{subfigure}
    \begin{subfigure}[t]{0.35\textwidth} 	
        \begin{overpic}[scale=0.5]{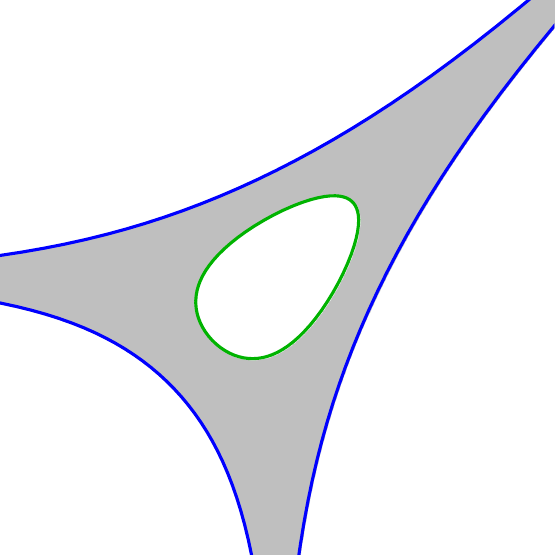}	
		\put(-11,47){$P_0$}
		\put(46,-10){$P_1$}
		\put(100,100){$P_\infty$}
        \end{overpic}
    \end{subfigure}
    \hspace{0.1cm}
    \vspace{0.25cm}
    \caption{The Riemann surface $\mathcal{R}$ as a three-sheeted cover of the complex plane (left subfigure), together with the bounded (green) and unbounded (blue) ovals, and the distinguished points $P_\infty$, $P_1$, and $P_0$. The amoeba map \eqref{eq:amoeba_map} maps the real ovals onto the boundary of the amoeba (right subfigure), and it maps $\mathcal{R}_+$ (shown in darker gray in the left subfigure) bijectively onto the interior of the amoeba (indicated by the same shade of gray in the right subfigure).}
    \label{fig:Riemann_surface_amoeba}
\end{figure}

The assumptions \eqref{eq:ab_assumption_1}--\eqref{eq:ab_assumption_2} lead to the explicit form
\begin{equation}\label{eq:spectral_curve}
    P(z,\lambda) = (\lambda-z-1)^3-27(1+\beta) \lambda z
\end{equation}
for a certain $\beta>0$, see \cite[Lemma~2.1]{Kui25}. Let $\mathcal{R}$ be the compact Riemann surface associated with \eqref{eq:spectral_curve}. The genus is one. It has three sheets $\mathcal{R}_1$, $\mathcal{R}_2$, and $\mathcal{R}_3$ as shown in the left subfigure of Figure~\ref{fig:Riemann_surface_amoeba}. The image of $\mathcal{R}$ under the amoeba map \eqref{eq:amoeba_map} is called the amoeba $\mathcal{A}$, and it is shown in the right subfigure.

The real locus consists of two closed contours that are called the real ovals, as shown in Figure~\ref{fig:Riemann_surface_amoeba}. The unbounded oval $\mathcal{O}_u$ is on the first and third sheets, and it contains the distinguished points
\begin{equation}\label{eq:Pinfty_P1_P0}
	P_\infty = (\infty,\infty), \quad P_1 = (-1,0), \quad P_0 = (0,1),
\end{equation}
that are mapped by the amoeba map \eqref{eq:amoeba_map} to infinity. We collect these points into the set $\mathcal{P}=\{P_\infty,P_1,P_0\}$. There is also a bounded oval $\mathcal{O}_b$ that is on the second and third sheets. The branching along this oval occurs at the points where $z$ attains its extremal values, which are given explicitly
\begin{equation}\label{eq:B1_B3}
	B_1 = (z_{\max},\lambda_{\max}), \quad B_3 = (z_{\min},\lambda_{\min}),
\end{equation}
with
\begin{equation}\label{eq:zmax_zmin}
	\begin{aligned}
		z_{\max} &= 1+2\beta+2\sqrt{\beta(1+\beta)}, \quad\, z_{\min} = 1+2\beta-2\sqrt{\beta(1+\beta)}, \\
		\lambda_{\max} &= -1-\beta-\sqrt{\beta(1+\beta)}, \quad\; \lambda_{\min} = -1-\beta+\sqrt{\beta(1+\beta)},
	\end{aligned}
\end{equation}
see \cite[Remark~2.2]{Kui25}. The points $B_1$ and $B_3$ play a role in Section~\ref{sec:KO_map_a}. There, we also use $B_0$ and $B_2$ to denote the points on the bounded oval where $\lambda$ attains its extremal values.  

The complement of the real locus has two connected components, which we denote by $\mathcal{R}_+$ and $\mathcal{R}_-$, where $\mathcal{R}_+$ contains the upper half-plane of the first sheet, see also Figure~\ref{fig:Riemann_surface_amoeba}. The amoeba map \eqref{eq:amoeba_map} maps $\mathcal{R}_+$ bijectively onto the interior of the amoeba.


\subsection{Double cover}
The double cover $\widehat{\mathcal{R}}$ is obtained by taking two copies of $\mathcal{R}$ and gluing them together along a certain cut. The cut is the support of an equilibrium measure in a certain external field. The equilibrium measure depends on the weights, but also on the relative dimensions of the hexagon. As in \cite{Kui25}, we consider the regular hexagon, that is
\begin{equation}\label{eq:BC_assumption}
    B = C = 1,
\end{equation}
Then, for the weights \eqref{eq:ab_assumption_1}--\eqref{eq:ab_assumption_2}, the equilibrium measure $\mu_\mathrm{eq}$ was determined in \cite{Kui25}. Its support is the union
\begin{equation}\label{eq:mu_support}
    \supp(\mu_\mathrm{eq}) = \mathbb{T}_1\cup\mathbb{T}_2
\end{equation}
of the two unit circles on the first and second sheets. 

The double cover is easily visualized if one uses the complex torus. Being of genus one, $\mathcal{R}$ is conformally equivalent to a complex torus $\mathbb{C}\slash(\mathbb{Z}+\tau\mathbb{Z})$ by means of the Abel map. Because the Harnack curve is a real curve, we may take $\tau\in i\mathbb{R}_{>0}$. As in \cite[Section~2.2]{Kui25}, we take the Abel map such that the points $P_\infty$, $P_1$, and $P_0$ are mapped to $1/6$, $1/2$, and $5/6$, respectively. The unbounded oval $\mathcal{O}_u$ is then mapped to $[0,1]$ and the bounded oval $\mathcal{O}_b$ to $[0,1]+\frac{1}{2}\tau$. Additionally, $\mathbb{T}_1\cup\mathbb{T}_2$ is mapped to the vertical segments $[-\frac{\tau}{2},\frac{\tau}{2}]$ and $\mathbb{T}_3$ to $[-\frac{\tau}{2},\frac{\tau}{2}]+\frac{1}{2}$, see \cite[Proposition~2.4]{Kui25}. All the coordinates are taken modulo the lattice $\mathbb{Z} + \tau \mathbb{Z}$. 

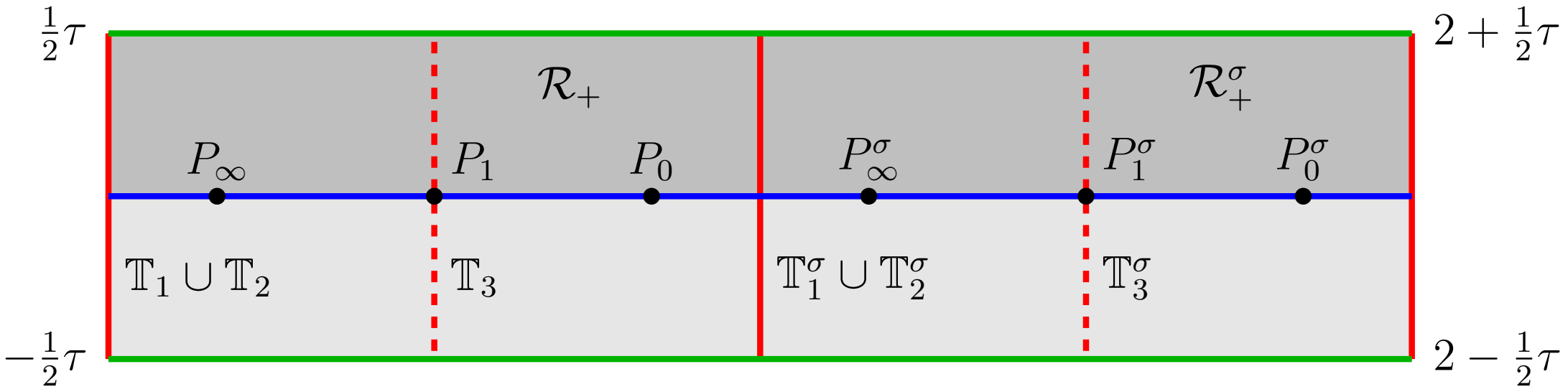
\begin{figure}[t]
    \centering
    \begin{tikzpicture}[scale=0.5](15,10)(0,0)
        \filldraw[gray!50!white] (-6,0) --(18,0) --(18,3) --(-6,3) --(-6,0);
        \filldraw[gray!20!white] (-6,0) --(18,0) --(18,-3) --(-6,-3) --(-6,0);

        \draw (-6.1,3) node[left] {$\frac{1}{2}\tau$};
        \draw (-6.1,-3) node[left] {$-\frac{1}{2}\tau$};
        \draw (18.1,3) node[right] {$2+\frac{1}{2}\tau$};
        \draw (18.1,-3) node[right] {$2-\frac{1}{2}\tau$};

		\draw[ultra thick,red] (-6,-3)--(-6,3) (6,-3)--(6,3) (18,-3)--(18,3);
        \draw[ultra thick,dashed,red] (0,-2.85)--(0,2.9)  (12,-2.85)--(12,2.9);
        \draw (-6,-1.5) node[right] {$\mathbb{T}_1\cup\mathbb{T}_2$};
        \draw (6,-1.5) node[right] {$\mathbb{T}_1^\sigma\cup\mathbb{T}_2^\sigma$};
        \draw (0,-1.5) node[right] {$\mathbb{T}_3$};
        \draw (12,-1.5) node[right] {$\mathbb{T}_3^\sigma$};

		\draw[ultra thick,blue] (-6,0)--(18,0);
		\draw[ultra thick,green!70!black] (-6,3)--(18,3) (-6,-3)--(18,-3);

        \draw (2.5,2) node {$\mathcal{R}_+$};
        \draw (14.5,2) node {$\mathcal{R}^\sigma_+$};

		\filldraw (-4,0) circle (4pt);
		\filldraw (0,0) circle (4pt);
		\filldraw (4,0) circle (4pt);
        \draw (-4,0) node[above] {$P_\infty$};
        \draw (0,0) node[above right] {$P_1$};
        \draw (4,0) node[above] {$P_0$};

		\filldraw (8,0) circle (4pt);
		\filldraw (12,0) circle (4pt);
		\filldraw (16,0) circle (4pt);
		\draw (8,0) node[above] {$P_\infty^\sigma$};
		\draw (12,0) node[above right] {$P_1^\sigma$};
		\draw (16,0) node[above] {$P_0^\sigma$};
    \end{tikzpicture}
    \caption{The double complex torus with the cycles $\mathbb{T}_1\cup\mathbb{T}_2$ and $\mathbb{T}_3$ (red). The figure also shows the bounded (green) and unbounded ovals (blue), and the points $P_\infty$, $P_1$, and $P_0$ on $\mathcal{R}$, as well as their copies $P_\infty^\sigma$, $P_1^\sigma$, and $P_0^\sigma$ on $\mathcal{R}^\sigma$.}
    \label{fig:double_torus}
\end{figure}

Since we cut along $\mathbb{T}_1\cup\mathbb{T}_2$, which corresponds to the vertical segment $[-\frac{1}{2}\tau,\frac{1}{2}\tau$] as in Figure~\ref{fig:double_torus}, the double cover is conformally equivalent to the double complex torus $\mathbb{C}\slash(2\mathbb{Z}+\tau\mathbb{Z})$, as shown in Figure~\ref{fig:double_torus} with $\mathcal{R}$ on the left and its copy on the right.

The double cover comes with a holomorphic involution
\begin{equation}\label{eq:sigma}
    \sigma \colon \widehat{\mathcal{R}} \to \widehat{\mathcal{R}}
\end{equation}
that interchanges the two copies. We use $\mathcal{R}^\sigma = \sigma(\mathcal{R})$ to denote the second copy. More generally, we adopt the convention that a superscript $\sigma$ indicates the image of some notion in $\mathcal{R}$ to $\mathcal{R}^\sigma$. For example, we write $\mathcal{P}^\sigma=\{P_\infty^\sigma,P_1^\sigma,P_0^\sigma\}$. We also put
\begin{equation}\label{eq:calP_hat}
    \widehat{\mathcal{P}} = \mathcal{P}\cup\mathcal{P}^\sigma = \{P_\infty,P_1,P_0,P_\infty^\sigma,P_1^\sigma,P_0^\sigma\}.
\end{equation} 
Likewise, we use $\mathcal{O}_u^\sigma$ and $\mathcal{O}_b^\sigma$ to denote the unbounded and bounded ovals on $\mathcal{R}^\sigma$. We also use
\[
    \widehat{\mathcal{O}}_u = \mathcal{O}_u\cup\mathcal{O}_u^\sigma, \quad \widehat{\mathcal{O}}_b = \mathcal{O}_b\cup\mathcal{O}_b^\sigma	
\]
to denote the unbounded and bounded ovals on the double cover.

There is also an anti-holomorphic involution
\begin{equation}\label{eq:tau}
	p = (z,\lambda) \mapsto \bar{p} \coloneqq (\bar{z},\bar{\lambda})
\end{equation}
that interchanges $\widehat{\mathcal{R}}_+$ and $\widehat{\mathcal{R}}_-$, where $\widehat{\mathcal{R}}_\pm = \mathcal{R}_\pm\cup\mathcal{R}^\sigma_\pm$.


\section{Statement of results} 


\subsection{Classification of asymptotic phases}
Our first results concern the classification of asymptotic phases (frozen, rough, and smooth) in the hexagon. These phases are determined by the zeros of a certain meromorphic differential.

\begin{definition}\label{def:omega}
    The differential $\omega$ is the unique meromorphic differential of the third kind, i.e., with only simple poles, on $\widehat{\mathcal{R}}$ with purely imaginary periods, and with poles and residues as listed in Table~\ref{table:poles_omega}.
\end{definition}

\begin{table}[ht]
    \centering
    \renewcommand{\arraystretch}{1.25}
    \setlength{\tabcolsep}{5pt}
    \begin{tabular}{c | c c c c c c}
        \hline
        Pole & $P_\infty$ & $P_1$ & $P_0$ & $P_\infty^\sigma$ & $P_1^\sigma$ & $P_0^\sigma$ \\ \hline
        Residue & $-3$ & $3$ & $-3$ & $3$ & $-3$ & $3$ \\
        \hline
    \end{tabular}
    \renewcommand{\arraystretch}{1}
    \caption{Poles and residues of $\omega$.}
    \label{table:poles_omega}
\end{table}

Since the residues listed in Table~\ref{table:poles_omega} add up to zero, there exist meromorphic differentials of the third kind with the prescribed poles and residues. The additional requirement of having purely imaginary periods then determines $\omega$ uniquely.

Besides $\omega$, we use the meromorphic differentials $d\log z$ and $d\log\lambda$ that come from the coordinate functions on $\mathcal{R}$. They extend to meromorphic differentials on~$\widehat{\mathcal{R}}$. All these differentials combine to give us the family of meromorphic differentials
\begin{equation}\label{eq:omega_xi1_xi2_1}
    \omega_{(\xi_1,\xi_2)} \coloneqq \omega - \xi_1d\log\lambda + \xi_2d\log z, \quad (\xi_1,\xi_2) \in \mathbb{R}^2.
\end{equation}
The differential \eqref{eq:omega_xi1_xi2_1} has the same poles as $\omega$, but with residues that depend on $\xi_1$ and $\xi_2$ as listed in Table~\ref{table:poles_omega_xi1_xi2}.

\begin{table}[ht]
    \centering
    \renewcommand{\arraystretch}{1.25}
    \setlength{\tabcolsep}{5pt}
    \begin{tabular}{c | c c c c c c}
        \hline
        Pole & $P_\infty$ & $P_1$ & $P_0$ & $P_\infty^\sigma$ & $P_1^\sigma$ & $P_0^\sigma$ \\
        \hline
        Residue & $-3+3\xi_1-3\xi_2$ & $3-3\xi_1$ & $-3+3\xi_2$ & $3+3\xi_1-3\xi_2$ & $-3-3\xi_1$ & $3+3\xi_2$ \\
        \hline
    \end{tabular}
    \renewcommand{\arraystretch}{1}
    \caption{Poles and residues of $\omega_{(\xi_1,\xi_2)}$.}
    \label{table:poles_omega_xi1_xi2}
\end{table}

We are particularly interested in $(\xi_1,\xi_2)$ that belong to the (open) hexagon 
\begin{equation}\label{eq:calH}
	\mathcal{H} = \{(\xi_1,\xi_2)\in\mathbb{R}^2 \mid -1<\xi_1<1, -1<\xi_2<1, -1<\xi_1-\xi_2<1\}. 
\end{equation}
Note that $(\xi_1,\xi_2)\in\mathcal{H}$ if and only if the residues listed in Table~\ref{table:poles_omega_xi1_xi2} alternate in sign.
 
\begin{remark}\label{rem:double_zero}
    Generically, $\omega_{(\xi_1,\xi_2)}$ has simple poles at the points in $\widehat{\mathcal{P}}$, see \eqref{eq:calP_hat}. For $(\xi_1,\xi_2) \notin \mathcal{H}$, it may happen that one of the residues in Table~\ref{table:poles_omega_xi1_xi2} vanishes. Then, we consider $\omega_{(\xi_1,\xi_2)}$ to have both a zero and a pole at the corresponding point $P \in\widehat{\mathcal{P}}$, even though there is a zero-pole cancellation. In general, there is no genuine zero there. However, if $\omega_{(\xi_1,\xi_2)}$ does have a genuine zero at $P$, then we consider it a double zero in the following theorem.
\end{remark}

The following theorem is stated for general $(\xi_1,\xi_2)\in\mathbb{R}^2$, although the main interest is for $(\xi_1,\xi_2)\in\mathcal{H}$.

\begin{theorem}\label{thm:four_zeros_on_bounded_oval}
    Let $(\xi_1,\xi_2)\in\mathbb{R}^2$. Then, the following hold.
    \begin{enumerate}
        \item[\rm{(a)}] The meromorphic differential $\omega_{(\xi_1,\xi_2)}$ has six zeros on $\widehat{\mathcal{R}}$, counted with multiplicity. There are at least four distinct zeros with odd multiplicity on the bounded oval $\widehat{\mathcal{O}}_b$.
        \item[\rm{(b)}] The remaining two zeros can be either both on the bounded oval, both on the unbounded oval, or form a complex-conjugate pair away from the real ovals.
        \item[\rm{(c)}] There is at most one zero of $\omega_{(\xi_1,\xi_2)}$ with higher multiplicity, which can only occur on the real ovals and the multiplicity is at most three.
        \item[\rm{(d)}] For $(\xi_1,\xi_2)\notin\mathcal{H}$, there are two simple zeros on the unbounded oval, except at
        \begin{equation}\label{eq:turning_points}
            \mathfrak{T} = \left\{(\tfrac{1}{2},-\tfrac{1}{2}),(1,\tfrac{1}{2}),(\tfrac{1}{2},1),(-\tfrac{1}{2},\tfrac{1}{2}),(-1,-\tfrac{1}{2}),(-\tfrac{1}{2},-1)\right\}.
        \end{equation}
        For $(\xi_1,\xi_2) \in \mathfrak{T}$,
       a double zero occurs at one of the points in $\widehat{\mathcal{P}}$, in the sense of Remark~\ref{rem:double_zero}.
    \end{enumerate}
\end{theorem}

The points in \eqref{eq:turning_points} are called turning points because, at each of these points, the tangent to the boundary $\partial\mathfrak{L}$ of the rough phase is also tangent to the boundary $\partial\mathcal{H}$ of the hexagon, as we will show in Lemma~\ref{lem:parallelity_property}(d).

As a consequence of Theorem~\ref{thm:four_zeros_on_bounded_oval}, we have a disjoint union
\[
    \overline{\mathcal{H}} = \mathfrak{S}\cup\mathfrak{L}\cup\mathfrak{G}\cup\mathfrak{A}_\mathrm{out}\cup\mathfrak{A}_\mathrm{in},
\]
where $\mathfrak{S}$, $\mathfrak{L}$, and $\mathfrak{G}$ are the three asymptotic phases defined as follows, and $\mathfrak{A}_\mathrm{out}$ and $\mathfrak{A}_\mathrm{in}$ are the boundaries between them.

\begin{definition}\label{def:phases}
    Let $(\xi_1,\xi_2)\in\overline{\mathcal{H}}$. Then, one says that $(\xi_1,\xi_2)$ belongs to
    \begin{itemize}
    	\item the frozen (solid) phase $\mathfrak{S}$ if $\omega_{(\xi_1,\xi_2)}$ has two simple zeros on the unbounded oval $\widehat{\mathcal{O}}_u$,
    	\item the rough (liquid) phase $\mathfrak{L}$ if $\omega_{(\xi_1,\xi_2)}$ has two 
    	(necessarily simple) zeros that are away from the real ovals,
    	\item the smooth (gaseous) phase $\mathfrak{G}$ if all six zeros of $\omega_{(\xi_1,\xi_2)}$ are simple and lie on the bounded oval $\widehat{\mathcal{O}}_b$,
        \item the outer arctic curve $\mathfrak{A}_\mathrm{out}$ if $\omega_{(\xi_1,\xi_2)}$ has a zero with higher multiplicity on the unbounded oval $\widehat{\mathcal{O}}_u$,
        \item the inner arctic curve $\mathfrak{A}_\mathrm{in}$ if $\omega_{(\xi_1,\xi_2)}$ has a zero with higher multiplicity on the bounded oval $\widehat{\mathcal{O}}_b$.
    \end{itemize}
\end{definition}

For $(\xi_1,\xi_2)\in\mathfrak{A}_\mathrm{out}$, the zero with higher multiplicity is always a double zero, and the remaining four zeros are all simple. For $(\xi_1,\xi_2)\in\mathfrak{A}_\mathrm{out}$, the zero with higher multiplicity is either a double zero or triple zero, and the remaining zeros are simple. This follows from Theorem~\ref{thm:four_zeros_on_bounded_oval}.

For $(\xi_1,\xi_2)\in\mathfrak{L}$, the two zeros of $\omega_{(\xi_1,\xi_2)}$ that are away from the real ovals form a complex-conjugate pair. One of these zeros lies in $\widehat{\mathcal{R}}_+$ and the other is its complex conjugate in $\widehat{\mathcal{R}}_-$. Thus, there exists a well-defined map
\begin{equation}\label{eq:KO_map}
	\Omega \colon \mathfrak{L} \to \widehat{\mathcal{R}}_+, \quad (\xi_1,\xi_2) \mapsto s
\end{equation}
that maps $(\xi_1,\xi_2)\in\mathfrak{L}$ to the unique zero $s$ of $\omega_{(\xi_1,\xi_2)}$ that lies in $\widehat{\mathcal{R}}_+$. Following \cite{BB26}, we call this map the KO map (after Kenyon and Okounkov \cite{KO07,KOS06}). This map also appeared prominently in the recent works \cite{BB25,BB26}. 

\begin{theorem}\label{thm:KO_map}\leavevmode
	\begin{enumerate}
        \item[\rm{(a)}] The KO map  \eqref{eq:KO_map} is a diffeomorphism between $\mathfrak{L}$ and $\widehat{\mathcal{R}}_+$. 
        \item[\rm{(b)}] The KO map extends continuously to the boundary, leading to a homeomorphism between the closures that is also denoted by $\Omega$,
        \begin{equation}\label{eq:KO_map_ext_1}
            \Omega \colon \overline{\mathfrak{L}} \to \widehat{\mathcal{R}}_+ \cup\widehat{\mathcal{O}}_b \cup\widehat{\mathcal{O}}_u.
        \end{equation}
        Furthermore, it maps the six turning points $\mathfrak{T}=\overline{\mathfrak{L}}\setminus\mathcal{H}$ as in \eqref{eq:turning_points} onto the six distinguished points $\widehat{\mathcal{P}}$ as in \eqref{eq:calP_hat} (in the same order).
	\end{enumerate}
\end{theorem}

The statement (a) is analogous to \cite[Theorem~4.9]{BB25}, where it was proved in the context of periodic tilings of the Aztec diamond, see also \cite[Proposition~10]{BB26}. We give an independent proof below. The statement (b) is the analogue of \cite[Theorem~4.11]{BB25}.

After composition with the amoeba map \eqref{eq:amoeba_map}, one obtains the homeomorphism
\begin{equation}\label{eq:Log_KO_map}
    \Log\circ\,\Omega \colon \overline{\mathfrak{L}}\cap\mathcal{H}=\overline{\mathfrak{L}}\setminus\mathfrak{T} \to \mathcal{A},
\end{equation}
from the closure of the rough phase within $\mathcal{H}$ to the amoeba $\mathcal{A}$. The map \eqref{eq:Log_KO_map} is two-to-one due to the appearance of the double cover. Taking two copies of the amoeba, we can consider \eqref{eq:Log_KO_map} as a homeomorphism to a double amoeba $\widehat{\mathcal{A}} = \mathcal{A}\cup\mathcal{A}^\sigma$. The double amoeba is instrumental in simplifying the steepest descent analysis, see Section~\ref{sec:local_scaling_limit}. This insight is due to \cite{BB25}, where the authors use the amoeba to simplify the steepest descent analysis for periodic tilings of the Aztec diamond.

It follows from Theorem~\ref{thm:KO_map} that the rough phase $\mathfrak{L}$ is doubly connected, while the smooth phase $\mathfrak{G}$ is simply connected. Further properties are collected in the following theorem.

\begin{theorem}\label{thm:six_cusps}\leavevmode
    \begin{enumerate}
    \item[\rm{(a)}] The boundary of the rough phase $\mathfrak{L}$ is given by
    \begin{equation}\label{eq:boundary_rough_phase_1}
        \partial\mathfrak{L}=\mathfrak{A}_\mathrm{out}\cup\mathfrak{A}_\mathrm{in}.
    \end{equation}
    \item[\rm{(b)}] Let $(\eta_1,\eta_2) \in \partial\mathfrak{L}$. If $\omega_{(\eta_1,\eta_2)}$ has a double zero, then $\partial\mathfrak{L}$ is smooth near $(\eta_1,\eta_2)$.
    \item[\rm{(c)}] Let $(\eta_1,\eta_2) \in \partial\mathfrak{L}$. If $\omega_{(\eta_1,\eta_2)}$ has a triple zero, then $(\eta_1,\eta_2) \in \mathfrak{A}_\mathrm{in}$ and $\partial\mathfrak{L}$ has an inward cusp at $(\eta_1,\eta_2)$.
    \item[\rm{(d)}] The inner arctic curve $\mathfrak{A}_\mathrm{in}$ has six cusp points. 
    \end{enumerate}
\end{theorem}

Since $\mathfrak{A}_\mathrm{out}$ has no cusps, it follows from Theorem~\ref{thm:six_cusps}(b) that $\mathfrak{A}_\mathrm{out}$ is a smooth closed contour. In Lemma~\ref{lem:parallelity_property}(b), we also prove that it is convex. The inner arctic curve has six cusp points. The six parts of $\mathfrak{A}_\mathrm{in}$ in between the cusp points are smooth and convex  when viewed from the rough region, see again Lemma~\ref{lem:parallelity_property}(b). They are concave when viewed from the smooth phase $\mathfrak{G}$, and then the cusps are outward pointing, which agrees with the Figure~\ref{fig:tiling_simulations_1}.

Theorem~\ref{thm:six_cusps}(d) is a special case of \cite[Corollary~48]{BB26}, and the proof relies on the remarkable parallelity property that was observed in \cite[Corollary~4.14]{BB25}, see also Lemma~\ref{lem:parallelity_property}. We give a somewhat different proof in Section~\ref{sec:six_cusps}.

In \cite{Pio26}, algebraic equations describing the arctic curves of periodic tilings of the Aztec diamond and certain periodic tilings of the hexagon were studied. In particular, \cite[Theorem~5.7]{Pio26} provides an explicit formula for the degree of the algebraic equation that applies to our model, and it follows that the algebraic equation describing both $\mathfrak{A}_\mathrm{out}$ and $\mathfrak{A}_\mathrm{in}$ has degree $12$.


\subsection{Height function and the limit shape}
The height function $h_N \colon V \to \mathbb{Z}_{\geq0}$ is the random variable that assigns to a vertex $v\in V$ the number of paths that are below $v$. More precisely, if $v=(x,y)$ with $0\leq x\leq 6N$, then any tiling gives us a particle system with $3N$ points $(x,y_j)$, $j= 0,\ldots, 3N-1$ at level $x$, see also Figure~\ref{fig:tiles_paths_dimensions}. Thus,
\begin{equation}\label{eq:hNv}
    h_N(v) = \#\{j \mid y_j<y\} \quad \text{for } v=(x,y).
\end{equation}
We note that $0\leq h_N(v)\leq 3N$.

The limit shape describes the (normalized) limiting behavior of $\mathbb{E}[h_N(v_N)]$ as $N\to\infty$, where $v_N$ is a vertex that depends on $N$ in the following way. We let
\begin{equation}\label{eq:vN} 
    v_N = (3N(1+\xi_{1,N}),3N(1+\xi_{2,N}))+v
\end{equation}
with $v$ fixed, $(\xi_{1,N},\xi_{2,N})\in\mathcal{H}$, and $N\xi_{1,N},N\xi_{2,N}\in\mathbb{Z}$. For a sequence $(\xi_{1,N},\xi_{2,N})$ that converges to a point $(\xi_1,\xi_2)\in\mathcal{H}$, one can take the limit of the normalized expectation
\begin{equation}\label{eq:EhN_vN_limit}
    \lim_{N\to\infty}\frac{1}{N}\mathbb{E}[h_N(v_N)] = H(\xi_1,\xi_2),
\end{equation}
and the resulting function $H$ is called the limit shape.

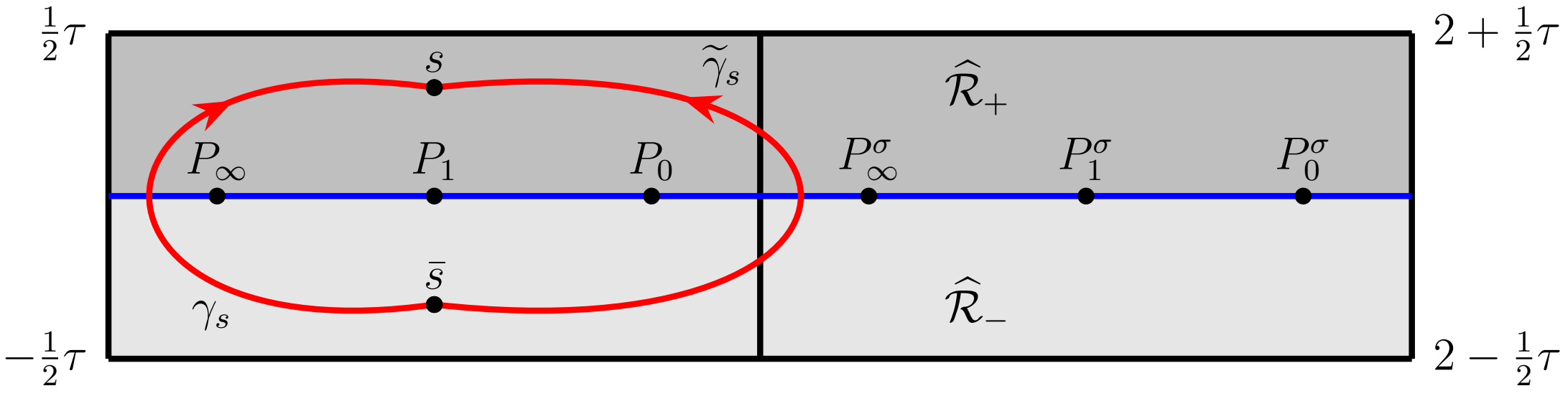
\begin{figure}[t]
	\centering
	\begin{tikzpicture}[scale=0.5](15,10)(0,0)
        \filldraw[gray!50!white] (-6,0) --(18,0) --(18,3) --(-6,3) --(-6,0);
        \filldraw[gray!20!white] (-6,0) --(18,0) --(18,-3) --(-6,-3) --(-6,0);

        \draw (-6.1,3) node[left] {$\frac{1}{2}\tau$};
        \draw (-6.1,-3) node[left] {$-\frac{1}{2}\tau$};
        \draw (18.1,3) node[right] {$2+\frac{1}{2}\tau$};
        \draw (18.1,-3) node[right] {$2-\frac{1}{2}\tau$};

		\draw[ultra thick] (-6,-3)--(-6,3) (6,-3)--(6,3) (18,-3)--(18,3);

		\draw[ultra thick,blue] (-6,0)--(18,0);
		\draw[ultra thick] (-6,3)--(18,3) (-6,-3)--(18,-3);

        \draw (10,2) node {$\widehat{\mathcal{R}}_+$};
        \draw (10,-2) node {$\widehat{\mathcal{R}}_-$};

		\filldraw (-4,0) circle (4pt);
		\filldraw (0,0) circle (4pt);
		\filldraw (4,0) circle (4pt);
        \draw (-4,0) node[above] {$P_\infty$};
        \draw (0,0) node[above] {$P_1$};
        \draw (4,0) node[above] {$P_0$};

		\filldraw (8,0) circle (4pt);
		\filldraw (12,0) circle (4pt);
		\filldraw (16,0) circle (4pt);
		\draw (8,0) node[above] {$P_\infty^\sigma$};
		\draw (12,0) node[above] {$P_1^\sigma$};
		\draw (16,0) node[above] {$P_0^\sigma$};

		\draw[postaction={decorate,decoration={markings,mark=at position 0.7 with {\arrow[rotate=-3]{Stealth}}}},ultra thick,red] (0,-2) .. controls (9,-3) and (9,3) .. (0,2);
		\draw[postaction={decorate,decoration={markings,mark=at position 0.7 with {\arrow[rotate=5]{Stealth}}}},ultra thick,red] (0,-2) .. controls (-7,-3) and (-7,3) .. (0,2);
		\filldraw (0,-2) circle (4pt);
		\filldraw (0,2) circle (4pt);
		\draw (5.3,2.4) node {$\widetilde{\gamma}_s$};
		\draw (-4.1,-2.2) node {$\gamma_s$};

		\filldraw (0,2) circle (4pt);
		\filldraw (0,-2) circle (4pt);
		\draw (0,2) node[above] {$s$};
		\draw (0,-2) node[above] {$\bar{s}$};
	\end{tikzpicture}
	\caption{The complex torus with the integration contours (red) from $\bar{s}\in\mathcal{R}_-$ to $s\in\mathcal{R}_+$ that occur in Theorem~\ref{thm:limit_shape}. The contour $\gamma_s$ intersects the unbounded oval (blue) in between $P_0^\sigma$ and $P_\infty$, and $\widetilde{\gamma}_s$ intersects it between $P_0$ and $P_\infty^\sigma$. There are no other intersections with the real ovals, and the contours are taken to be symmetric with respect to complex conjugation.}
    \label{fig:double_torus_EhN_contours}
\end{figure}

\begin{theorem}\label{thm:limit_shape}
    Suppose that $(\xi_{1,N},\xi_{2,N})\in\mathcal{H}$ with $N\xi_{1,N},N\xi_{2,N}\in\mathbb{Z}$ for all $N\in\mathbb{N}$. Then, the following hold.
	\begin{enumerate}
    	\item[\rm{(a)}] For $v_N$ as in \eqref{eq:vN} with $v=(0,0)$ and for each $\delta\in(0,1)$, one has
        \begin{multline}\label{eq:EhN_vN}
            \mathbb{E}[h_N(v_N)] = 3N - \Tr\left[\frac{1}{(2\pi i)^2}\oint_{(1-\delta)\mathbb{T}} dz_2\oint_\mathbb{T} \frac{dz_1}{z_1-z_2}\frac{W(z_2)^{N(1-\xi_{1,N})}}{z_2^{N(1-\xi_{2,N})}}\right. \\
            \left.\times R_N(z_1,z_2) \frac{W(z_1)^{N(1+\xi_{1,N})}}{z_1^{N(1+\xi_{2,N})}}\right].
    	\end{multline}
        \item[\rm{(b)}] For $(\xi_{1,N},\xi_{2,N}) \to (\xi_1,\xi_2)\in\mathfrak{L}$ and $v_N$ as in \eqref{eq:vN}, the limit \eqref{eq:EhN_vN_limit} holds with the limit shape $H$ defined as follows. Let $s=\Omega(\xi_1,\xi_2)$, with the KO map $\Omega$ as in \eqref{eq:KO_map}, and let $\gamma_s$ be a contour in $\widehat{\mathcal{R}}$ from $\bar{s}$ to $s$ as indicated in Figure~\ref{fig:double_torus_EhN_contours}. Then,
    	\begin{equation}\label{eq:limit_shape_1}
            H(\xi_1,\xi_2) = \frac{1}{2\pi i}\int_{\gamma_s}\omega_{(\xi_1,\xi_2)},
    	\end{equation}
        and the gradient of $H$ is given by
        \begin{equation}\label{eq:limit_shape_gradient}
            \nabla H(\xi_1,\xi_2) = \left(-\frac{1}{2\pi i}\int_{\gamma_s}d\log\lambda, \frac{1}{2\pi i}\int_{\gamma_s}d\log z\right)\!.
    	\end{equation}
    \end{enumerate}
\end{theorem}

Alternatively, \eqref{eq:limit_shape_1} can be written as
\begin{equation}\label{eq:limit_shape_2}
    H(\xi_1,\xi_2) = 3 + \frac{1}{2\pi i}\int_{\widetilde{\gamma}_s}\omega_{(\xi_1,\xi_2)},
\end{equation}
where $\widetilde{\gamma}_s$ is the other contour shown in Figure~\ref{fig:double_torus_EhN_contours}. This follows by Cauchy's theorem as the residues of $\omega_{(\xi_1,\xi_2)}$ at $P_\infty$, $P_1$, and $P_0$ add up to $-3$.

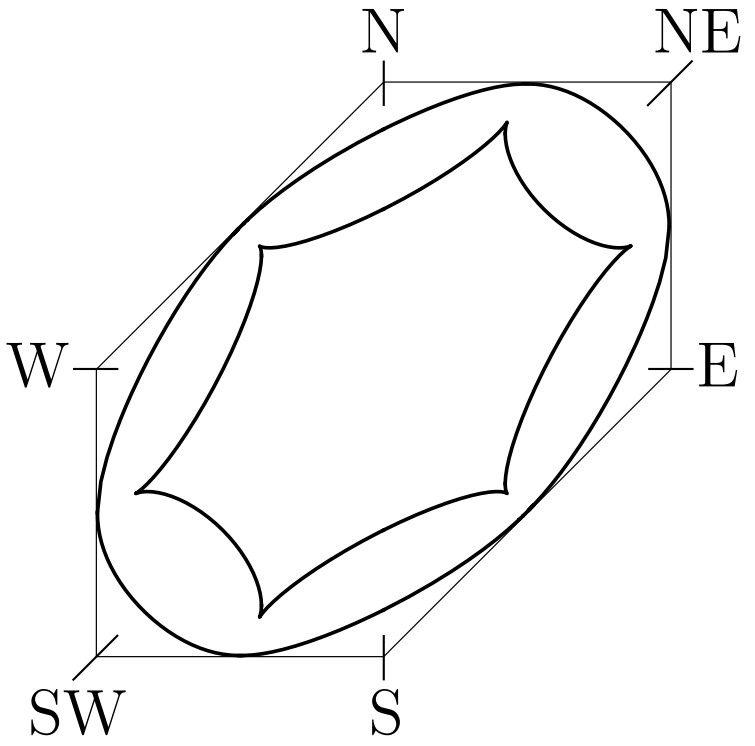
\begin{figure}[ht]
    \centering
    \begin{overpic}[width=0.3\linewidth]{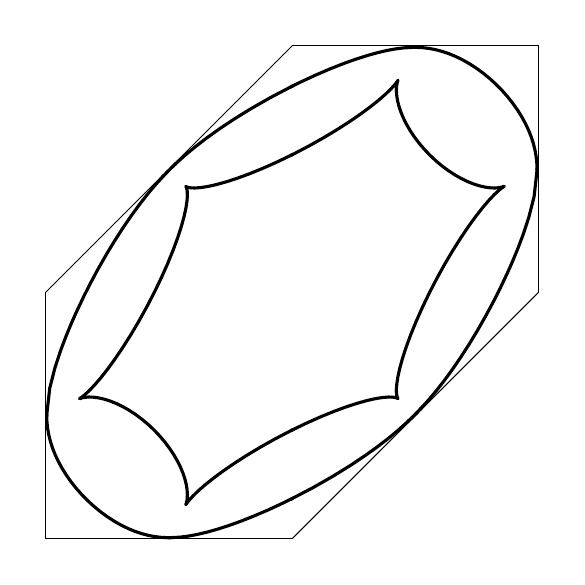}
        \put(46.5,96.5){N}
        \put(46.8,88.5){
        \begin{tikzpicture}
            \draw[-] (0,0) -- (0,-0.3);
        \end{tikzpicture}}
        \put(89.5,96.5){NE}
        \put(85.45,88.5){
        \begin{tikzpicture}
            \draw[-] (0,0) -- (0.3,0.3);
        \end{tikzpicture}}
        \put(95.9,47.5){E}
        \put(85.6,49.9){
        \begin{tikzpicture}
            \draw[-] (0,0) -- (-0.3,0);
        \end{tikzpicture}}
        \put(47.8,-3.5){S}
        \put(46.8,4.2){
        \begin{tikzpicture}
            \draw[-] (0,0) -- (0,-0.3);
        \end{tikzpicture}}
        \put(-2.2,-3.5){SW}
        \put(1.15,4.2){
        \begin{tikzpicture}
            \draw[-] (0,0) -- (-0.3,-0.3);
        \end{tikzpicture}}
        \put(-5.5,47.5){W}
        \put(1.2,49.9){
        \begin{tikzpicture}
            \draw[-] (0,0) -- (0.3,0);
        \end{tikzpicture}}
	\end{overpic}
    \vspace{0.1cm}
	\caption{The three asymptotic phases (frozen, rough, and smooth) and the arctic curves. In particular, the figure shows the labeling of the six connected components of the frozen phase used in \eqref{eq:limit_shape_frozen_phase}.}
    \label{fig:frozen_phases}
\end{figure}

\begin{remark}\label{rem:limit_shape}
	For $(\xi_1,\xi_2)\in\mathfrak{G}$ (smooth phase), one can show that the limit shape is given by
	\begin{equation}\label{eq:limit_shape_smooth_phase}
        H(\xi_1,\xi_2) = \frac{1}{2\pi i}\oint_\gamma \omega_{(\xi_1,\xi_2)},
    \end{equation}
    where $\gamma$ is the resulting cycle obtained by taking the contour $\gamma_s$ in \eqref{eq:limit_shape_smooth_phase} and letting $s=\Omega(\xi_1,\xi_2)$ approach the bounded oval. In fact, one may take any cycle homotopic to $\mathbb{T}_1\cup\mathbb{T}_2$ that intersects the bounded oval at an arbitrary point and the unbounded oval in between $P_0^\sigma$ and $P_\infty$. Thus, the contour $\gamma$ can be taken independently of the point $(\xi_1,\xi_2) \in \mathfrak{S}$ and, in particular, the gradient of \eqref{eq:limit_shape_smooth_phase}, see also \eqref{eq:limit_shape_gradient}, is constant throughout the smooth phase.

    Similarly, an analogue of \eqref{eq:limit_shape_smooth_phase} can be shown for $(\xi_1,\xi_2)\in\mathfrak{S}$ (frozen phase). The contour is then replaced by one of six null-homotopic cycles that intersect the unbounded oval twice; once in between $P_0^\sigma$ and $P_\infty$, and once at an arbitrary point on one of the six connected components of $\widehat{\mathcal{O}}_u\!\setminus\!\widehat{\mathcal{P}}$. Thus, the gradient of the limit shape is piecewise constant throughout the frozen phase. Moreover, it also follows that the frozen phase consists of six connected components, and using residue calculus together with the residues listed in Table~\ref{table:poles_omega_xi1_xi2}, one finds
    \begin{equation}\label{eq:limit_shape_frozen_phase}
        H(\xi_1,\xi_2) = 
        \begin{cases}
            3 & \text{in the north (N) frozen phase}, \\
            3\xi_2 & \text{in the north-east (NE) frozen phase}, \\
            3-3\xi_1+3\xi_2 & \text{in the east (E) frozen phase}, \\
            0 & \text{in the south (S) frozen phase}, \\
            3+3\xi_2 & \text{in the south-west (SW) frozen phase}, \\
            3\xi_2-3\xi_1 & \text{in the west (W) frozen phase}.
    	\end{cases}
    \end{equation}
    These six frozen phases are shown in Figure~\ref{fig:frozen_phases}.
\end{remark}

\subsection{Local scaling limit}
The local scaling limit of the correlation kernel describes the limiting behavior of $K_N(v_{1,N},v_{2,N})$ as $N\to\infty$, where $v_{1,N}$ and $v_{2,N}$ are vertices that depend on $N$ as in \eqref{eq:vN}. More precisely, we let
\begin{align}
    v_{1,N} &= (3N(1+\xi_{1,N}), 3N(1+\xi_{2,N})) + v_1, \label{eq:v1N} \\
    v_{2,N} &= (3N(1+\xi_{1,N}), 3N(1+\xi_{2,N})) + v_2 \label{eq:v2N}
\end{align}
with $v_1$ and $v_2$ fixed, $(\xi_{1,N},\xi_{2,N})\in\mathcal{H}$, and $N\xi_{1,N},N\xi_{2,N}\in\mathbb{Z}$. We consider a point $(\xi_1,\xi_2)\in\mathfrak{L}$ (rough phase) and suppose that $(\xi_{1,N},\xi_{2,N}) \to (\xi_1,\xi_2)$. Then, we prove that $K_N(v_{1,N},v_{2,N})$ has a limit as $N\to\infty$. The limiting kernel $K_{(\xi_1,\xi_2)}(v_1,v_2)$ is the kernel of a determinantal point process on $\mathbb{Z}^2$ that is also $3\times3$-doubly periodic.

To describe this limiting kernel, we need $W$ and $T_{0,j}$ as in \eqref{eq:Wz}--\eqref{eq:T0jz}, as well as the spectral projection $\Pi$ associated with $W(z)$, which is also known as the Riesz projector. For $p=(z,\lambda)\in\mathcal{R}$, it is given by 
\begin{equation}\label{eq:Riesz_projector}
    \Pi(p) = \frac{1}{2\pi i}\oint_{C_\lambda}\left(\mu I - W(z)\right)^{-1} d\mu \quad \text{for } p = (z,\lambda)\in\mathcal{R}\setminus\{P_\infty\},
\end{equation}
where $C_\lambda$ is a small circle around $\lambda$ that does not enclose any of the other eigenvalues of $W(z)$. If $W(z)$ is diagonalizable, which is the case for every $z\in\mathbb{C}\setminus\{0,z_{\min},z_{\max}\}$, then $\Pi(p)$ is a projection matrix that projects for $p=(z,\lambda)$ onto the eigenspace corresponding to the eigenvalue $\lambda$, and the eigenvectors corresponding to the other eigenvalues are in its kernel.

The spectral projection $\Pi$ is analytic in $\mathcal{R}\setminus\{P_\infty,P_0,B_1,B_3\}$ with simple poles at the branch points $B_1$ and $B_3$, see \eqref{eq:B1_B3}. The differential $dz/z$ has simple zeros at these branch points and only (simple) poles at $P_\infty$ and $P_0$. Thus, $\Pi\,dz/z$ is a matrix-valued meromorphic differential on $\mathcal{R}$ with only poles at $P_\infty$ and $P_0$.

As part of our result, we show that $K_{(\xi_1,\xi_2)}$ is given by a contour integral involving, among other factors, the differential $\Pi\,dz/z$. The dependence on the point $(\xi_1,\xi_2)$ only appears in the integration contour.

We denote the projection map by
\begin{equation}\label{eq:proj_map}
    \pi \colon \widehat{\mathcal{R}} \to \mathcal{R}.
\end{equation}

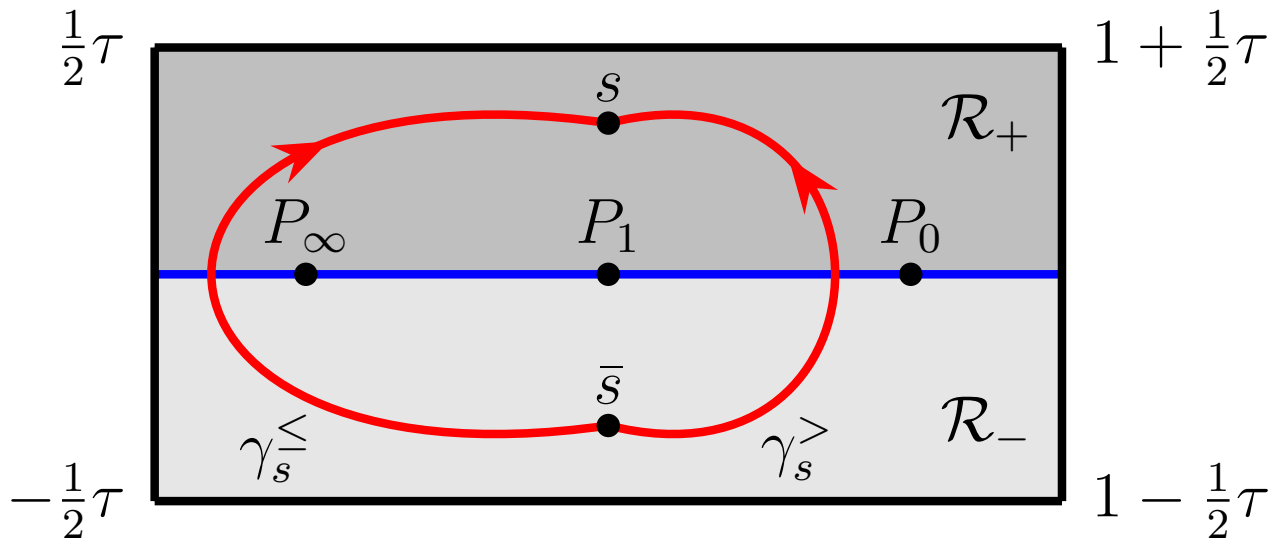
\begin{figure}[t]
	\centering
	\begin{tikzpicture}[scale=0.5](15,10)(0,0)
		\filldraw[gray!50!white] (-6,0) --(6,0) --(6,3) --(-6,3) --(-6,0);
		\filldraw[gray!20!white] (-6,0) --(6,0) --(6,-3) --(-6,-3) --(-6,0);

        \draw (-6.1,3) node[left] {$\frac{1}{2}\tau$};
        \draw (-6.1,-3) node[left] {$-\frac{1}{2}\tau$};
        \draw (6.1,3) node[right] {$1+\frac{1}{2}\tau$};
        \draw (6.1,-3) node[right] {$1-\frac{1}{2}\tau$};

        \draw (5,2) node {$\mathcal{R}_+$};
        \draw (5,-2) node {$\mathcal{R}_-$};

		\draw[ultra thick,blue] (-6,0)--(6,0);
		\draw[ultra thick] (-6,3)--(6,3) (-6,-3)--(6,-3);

        \draw[ultra thick] (-6,-3)--(-6,3) (6,-3)--(6,3);

		\filldraw (-4,0) circle (4pt);
		\filldraw (0,0) circle (4pt);
		\filldraw (4,0) circle (4pt);
		\draw (-4,0) node[above] {$P_\infty$};
		\draw (0,0) node[above] {$P_1$};
		\draw (4,0) node[above] {$P_0$};

        \draw[postaction={decorate, decoration={markings,mark=at position 0.7 with {\arrow[rotate=-9]{Stealth}}}},ultra thick,red] (0,-2) .. controls (4,-3) and (4,3) .. (0,2);
		\draw[postaction={decorate, decoration={markings,mark=at position 0.7 with {\arrow[rotate=4]{Stealth}}}},ultra thick,red] (0,-2) .. controls (-7,-3) and (-7,3) .. (0,2);
		\draw (-4.4,-2.3) node {$\gamma_s^\leq$};
		\draw (2.5,-2.3) node {$\gamma_s^>$};

		\filldraw (0,2) circle (4pt);
		\filldraw (0,-2) circle (4pt);
		\draw (0,2) node[above] {$s$};
		\draw (0,-2) node[above] {$\bar{s}$};
	\end{tikzpicture}
    \caption{The complex torus with the integration contours (red) from $\bar{s}\in\mathcal{R}_-$ to $s\in\mathcal{R}_+$ that are used the limiting kernel \eqref{eq:limiting_kernel} in Theorem~\ref{thm:local_scaling_limit}. The contour $\gamma_s^>$ intersects the unbounded oval (blue) in the part between $P_\infty$ and $P_0$ that contains $P_1$ (it could intersect on either side of $P_1$), and $\gamma_s^\leq$ intersects the unbounded oval in the complementary part between $P_0$ and $P_\infty$. The contours are in $\mathcal{R}$, where they are symmetric with respect to complex conjugation, and there are no other intersections with the real ovals in $\mathcal{R}$.}
    \label{fig:torus_KN_contours}
\end{figure}

\begin{theorem}\label{thm:local_scaling_limit}
    Suppose that $(\xi_{1,N},\xi_{2,N})\in\mathcal{H} \to (\xi_1,\xi_2)\in\mathfrak{L}$ with $N\xi_{1,N},N\xi_{2,N}\in\mathbb{Z}$, and let $v_{1,N}$ and $v_{2,N}$ be as in \eqref{eq:v1N}--\eqref{eq:v2N}. Then,
    \begin{equation}\label{eq:KN_v1N_v2N_limit}
        \lim_{N\to\infty} K_N(v_{1,N},v_{2,N}) = K_{(\xi_1,\xi_2)}(v_1,v_2),
    \end{equation}
    where the limiting kernel $K_{(\xi_1,\xi_2)} \colon \mathbb{Z}^2 \times \mathbb{Z}^2 \to \mathbb{R}$ is defined as follows. Let $s=\Omega(\xi_1,\xi_2)$, with the KO map $\Omega$ as in \eqref{eq:KO_map}, and the projection map \eqref{eq:proj_map}. Let $\gamma_s^\leq$ and $\gamma_s^>$ be any contours from $\bar{s}$ to $s$ in $\mathcal{R}$ as indicated in Figure~\ref{fig:torus_KN_contours}. Then, for
    \begin{equation}\label{eq:v1_v2}
        v_1 = (3x_1 + j_1, 3y_1 + k_1), \quad v_2 = (3x_2+j_2, 3y_2 + k_2),
    \end{equation}
    with $x_1, y_1, x_2, y_2\in\mathbb{Z}$ and $j_1, k_1, j_2, k_2\in\{0,1,2\}$,
    \begin{equation}\label{eq:limiting_kernel}
    	K_{(\xi_1,\xi_2)}(v_1,v_2) =
        \begin{cases}
            \!\raisebox{-12pt}{$\left[\dfrac{1}{2\pi i}{\displaystyle\int_{\gamma_s^\leq}}\dfrac{\lambda^{x_1-x_2}}{z^{y_1-y_2}}T_{0,j_2}^{-1}(z)\Pi\,T_{0,j_1}(z)\dfrac{dz}{z}\right]_{k_2,k_1}
             \text{if } 3x_1+j_1\leq 3x_2+j_2,$} & \!\!\!\! \\
            \!\raisebox{-20pt}{$\left[\dfrac{1}{2\pi i}{\displaystyle\int_{\gamma_s^>}}\dfrac{\lambda^{x_1-x_2}}{z^{y_1-y_2}}T_{0,j_2}^{-1}(z)\Pi\,T_{0,j_1}(z) \dfrac{dz}{z}\right]_{k_2,k_1} \text{if } 3x_1+j_1>3x_2+j_2,$} & \!\!\!\!
    	\end{cases}
    \end{equation}
    where we recall that $T_{0,j}$ is given by \eqref{eq:T0jz} and $\Pi$ is the
    spectral projection \eqref{eq:Riesz_projector}.
\end{theorem}
In case $3x_1 + j_1 > 3x_2 + j_2$, we have that $x_1 \geq x_2$, and therefore $\lambda^{x_1-x_2}$ is analytic at $P_1$. Then, the integrand in \eqref{eq:limiting_kernel} has no pole at $P_1$, and we can move $\gamma_s^>$ to the other side of $P_1$ (if we wish), as also indicated in the caption of Figure~\ref{fig:torus_KN_contours}.


\subsection{Organization of the rest of the paper}
The remaining part of the paper is divided into three sections.

Section~\ref{sec:four_zeros_KO_map} is devoted to the proofs of Theorems~\ref{thm:four_zeros_on_bounded_oval},~\ref{thm:KO_map}, and~\ref{thm:six_cusps}. It also sets the stage for the steepest descent analysis developed in the following sections. In particular, we first relate the $g$-function associated with the equilibrium measure with the meromorphic differential $\omega_{(\xi_1,\xi_2)}$. We also study the integral of this differential, which is called the action function, and plays a central role in the steepest descent analysis.

The main novelty of our work lies in Section~\ref{sec:local_scaling_limit}, where we carry out the steepest descent analysis of the correlation kernel. Ultimately, this allows us to prove the local scaling limits of the correlation kernel as stated in Theorem~\ref{thm:local_scaling_limit}.

We find that the asymptotic behavior of the double contour integral that represents the correlation kernel is determined by the action function. Thus, we are first led to study the critical level sets of this action function. We then apply the first two transformations arising in RH analysis to the double contour integral and interpret the resulting integrals on the Riemann surface. Using the jump condition of the RH problem, we split the integrals and then show the existence of the suitable curves that locally follow the paths of steepest descent and ascent. Finally, we deform the integration contours to these suitable contours, and in this deformation process, we pick up a residue contribution that gives the leading-order contribution to the correlation kernel.

In Section~\ref{sec:limit_shape}, we adapt the steepest descent analysis developed in Section~\ref{sec:local_scaling_limit} to establish Theorem~\ref{thm:limit_shape}.


\subsection{Outlook: transitions to two smooth phases}\label{sec:outlook}
Our main results, Theorems~\ref{thm:limit_shape} and~\ref{thm:local_scaling_limit}, are concerned with the limit shape and local scaling limits of the particular model. We believe that the techniques that are developed in this paper will help to establish other phenomena, such as Airy-like behavior for fluctuations around the arctic curves, and the Gaussian free field fluctuations (possibly with a discrete component as in \cite{BN25+}) for the height function in the rough phase.

As mentioned in the introduction, we also would like to extend the analysis to periodic weightings with more general weights, and to hexagons of arbitrary dimensions.

We expect interesting features to occur if we use the same $3\times3$-doubly periodic weights as in the paper, but change the dimensions of the hexagon, as encoded by the scaling parameters $B$ and $C$ in \eqref{eq:hexagon_dimensions}. Since the weights are the same, the spectral curve remains the same Harnack curve. However, the double cover of the associated Riemann surface may be different, as it depends on the support of the equilibrium measure in an appropriate external field that now depends on $B$ and $C$.

For general $B$ and $C$, the symmetries that allow us to determine the relevant equilibrium measure explicitly are lost. However, for the case
\begin{equation}\label{eq:BC_assumption_extended}
	B = C,
\end{equation}
we retain one symmetry in the model that ensures that the relevant equilibrium measure is still supported on $\mathbb{T}_1\cup\mathbb{T}_2\cup\mathbb{T}_3$. In particular, its support will have the $S$-property, see also Section~\ref{sec:equilibrium_measure}. This property is essential for the analysis of the Riemann--Hilbert problem, see \cite{Kui25}, and is also used at several points throughout this paper. In \cite[Remark~1.8]{Kui25}, it was already argued that the support of this equilibrium measure need not be the cycle $\mathbb{T}_1\cup\mathbb{T}_2$, and in such a case the geometry of the problem changes.

\begin{figure}[t]
	\centering
    \hspace{0.2cm}
    \includegraphics[scale=0.5]{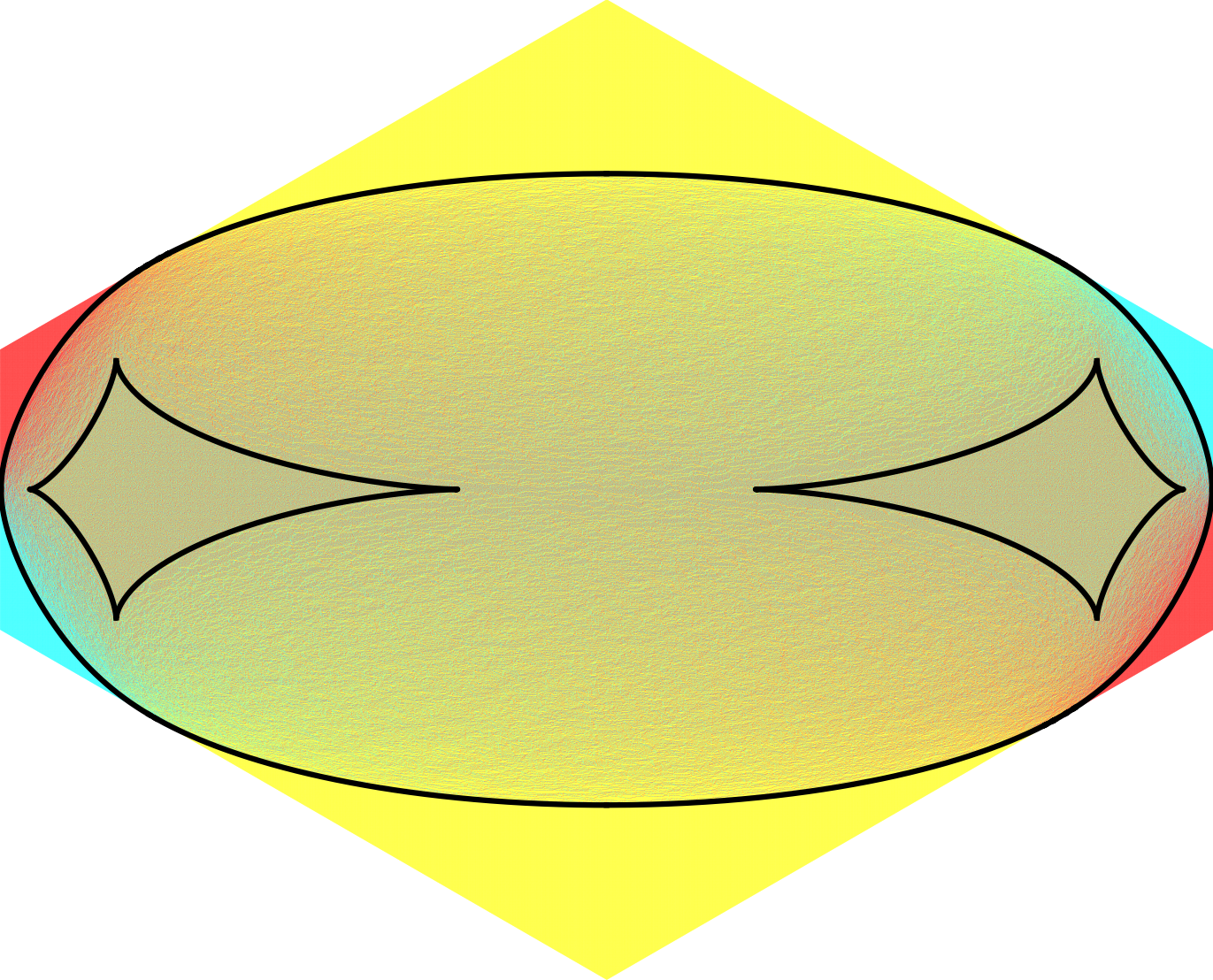}
    \hspace{2cm}
    \includegraphics[scale=0.5]{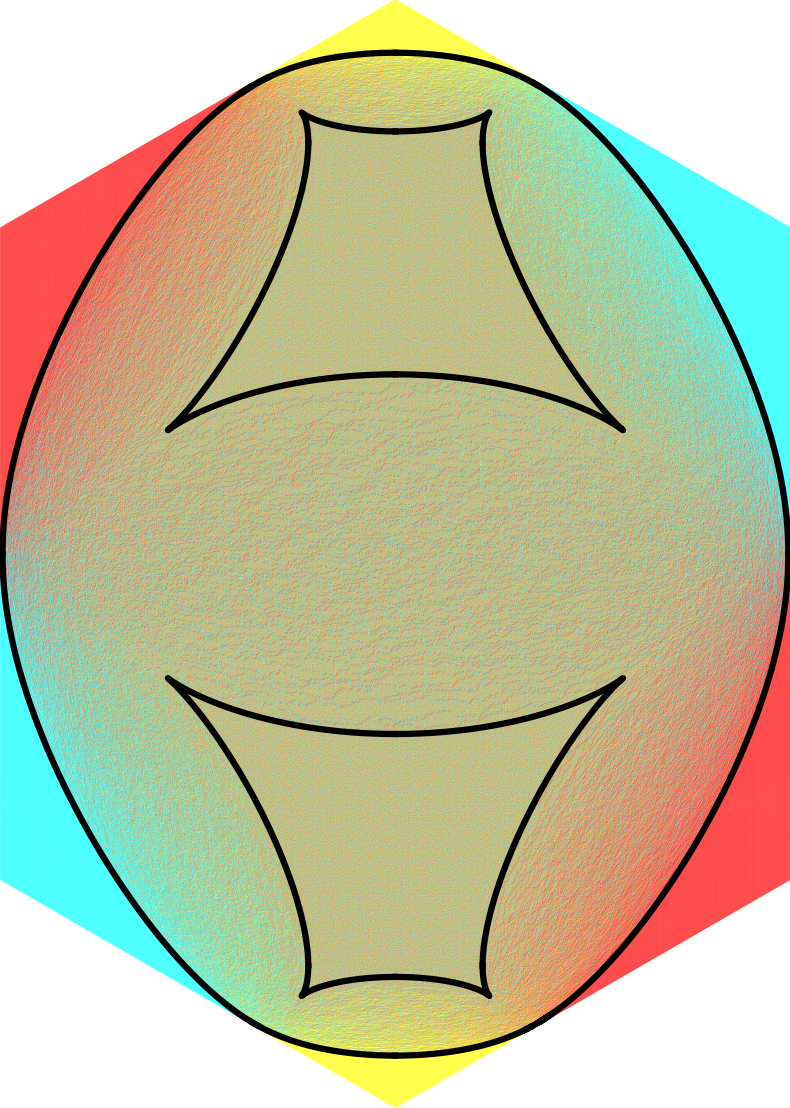}
    \hspace{0.5cm}
    \caption{Large-size simulations of  random tilings of a hexagon of size $3N\times3BN\times3BN$ for $B>B_{\mathrm{crit},2}$ (left subfigure) and $B<B_{\mathrm{crit},1}$ (right subfigure) with the same $3\times3$-doubly periodic weights as in this paper. The figure shows the arctic curves (the boundaries between the phases). Two smooth phases appear in both subfigures, each with four cusp points. As $B\to B_{\mathrm{crit},2}$, we observe that the two smooth phases in the left subfigure merge at a cusp point, while they merge along a line segment connecting two cusp points as $B \to B_{\mathrm{crit},1}$ in the right subfigure.}
    \label{fig:tiling_simulations_2}
\end{figure}

There are two transitions, at two critical parameters (depending on $\beta$, see \eqref{eq:spectral_curve})
\[
    0<B_{\mathrm{crit},1}<1<B_{\mathrm{crit},2}<\infty
\]
that characterize the support of the equilibrium measure as follows.
\begin{itemize}
    \item For $B=C\in(B_{\mathrm{crit},1}, B_{\mathrm{crit},2})$, the support is the full cycle $\mathbb{T}_1\cup\mathbb{T}_2$.
    \item For $B=C\in(B_{\mathrm{crit},2},\infty)$, the support is a subarc of the cycle $\mathbb{T}_1\cup\mathbb{T}_2$.
    \item For $B=C\in(0,B_{\mathrm{crit},1})$, the support consists of the union of the cycle $\mathbb{T}_1\cup\mathbb{T}_2$ with a subarc of $\mathbb{T}_3$.
\end{itemize}
In all cases, the equilibrium measure has a positive density away from the end points of its support, where it has square-root vanishing.

The first case is expected to be similar to the case $B=C=1$ considered in this paper, in the sense that there is one smooth phase with six cusp points. The double cover is then the same as the one used in this paper.

The other two cases are expected to lead to two smooth phases, each with four cusp points. This is supported by the outcome of numerical simulations as shown in Figure~\ref{fig:tiling_simulations_2}. In these two cases, the double cover will depend on the parameter $B$ through the support of the equilibrium measure. The double cover is ramified and has genus $2$. Further numerical simulations around the critical values $B_{\mathrm{crit},1}$ and $B_{\mathrm{crit},2}$ indicate that the two transitions from one to two smooth phases are of a different nature. We hope to be able to report on this intriguing phenomenon in future work.


\section{Proofs of Theorems~\ref{thm:four_zeros_on_bounded_oval} and~\ref{thm:KO_map}}\label{sec:four_zeros_KO_map}
We first recall results from \cite{Kui25} concerning the equilibrium measure and the $g$-function that are needed for our analysis. Then, we complement these results by showing that the differential $\omega$ as defined in Definition~\ref{def:omega} can be expressed in terms of the differential of the $g$-function.


\subsection{Equilibrium measure}\label{sec:equilibrium_measure}
The equilibrium measure is defined as the minimizer of the bipolar Green's energy in the external field $\Re V$ on an $S$-curve, with 
\begin{equation}\label{eq:ext_field}
    V = 2\log z-2\log\lambda \quad \text{on } \mathcal{R}.
\end{equation}

Note that the logarithmic energy and the logarithmic potential from potential theory in the complex plane should be adapted to the Riemann surface $\mathcal{R}$ in order to study potential theory on it. In our case, this is done by introducing the bipolar Green's kernel on $\mathcal{R}$ with a pole at $P_\infty$.

\begin{definition}[{\cite[Definition~1.6]{Kui25}}]\label{def:bipolar_greens_kernel}
	A function
    \begin{equation}
        G_{P_\infty} \colon \{(p,q)\in(\mathcal{R}\setminus\{P_\infty\})^2 \mid p \neq q\} \to \mathbb{R}
    \end{equation}
    is called a bipolar Green's kernel (on $\mathcal{R}$ with a pole at $P_\infty$) if for each $q\in\mathcal{R}\setminus\{P_\infty\}$:
	\begin{enumerate}
        \item[\rm{(a)}] the function $p \mapsto G_{P_\infty}(p,q)$ is harmonic on $\mathcal{R}\setminus\{P_\infty, q\}$;
		\item[\rm{(b)}] for a local coordinate $z_q$ at $q$, one has
		\[
		G_{P_\infty}(p,q) = -\log|z_q(p)| + \mathcal{O}(1) \quad \text{as } p \to q;
		\]
		\item[\rm{(c)}] for a local coordinate $z_\infty$ at $P_\infty$, one has
		\[
		G_{P_\infty}(p,q) = \log|z_\infty(p)| + \mathcal{O}(1) \quad \text{as } p \to P_\infty;
		\]
		\item[\rm{(d)}] $G_{P_\infty}(p,q) = G_{P_\infty}(q,p)$ for all $p, q\in\mathcal{R}$ with $p \neq q$.
	\end{enumerate}
	The properties (a)--(d) define a function uniquely up to an additive constant. We impose the additional condition (\cite[Lemma~4.1]{Kui25}):
	\begin{enumerate}
		\item[\rm{(e)}] for a local coordinate $z_\infty$ at $P_\infty$, one has
		\begin{equation}\label{eq:bipolar_G_at_Pinfty}
		      G_{P_\infty}(p,q) = \log|z_\infty(p)| + o(1) \quad \text{as } p \to P_\infty.
		\end{equation}
	\end{enumerate}
	In the sequel, we will reserve the notation $G_{P_\infty}(p,q)$ for the unique bipolar Green's kernel that also satisfies the condition (e).
\end{definition}

For a measure $\mu$ on $\mathcal{R}$, the bipolar Green's energy of $\mu$ in external field $\Re V$, with $V$ as in \eqref{eq:ext_field}, is defined as
\begin{equation}\label{eq:EV_mu}
    E_V(\mu) = 3\iint G_{P_\infty}(p,q)d\mu(p)d\mu(q) +\int \Re V d\mu.
\end{equation}

Let $\Gamma$ be a finite union of closed contours in $\mathcal{R}\setminus\{P_\infty,P_0\}$. Then, there exists a unique probability measure on $\Gamma$ that minimizes \eqref{eq:EV_mu} among all probability measures on $\Gamma$. This is the equilibrium measure in the external field of $\Gamma$, and we denote it by $\mu_\mathrm{eq}$. The equilibrium measure is characterized by the Euler--Lagrange variational equations:
\begin{equation}\label{eq:EL_conditions_1}
    6\int G_{P_\infty}(p,q) d\mu_\mathrm{eq}(q) + \Re V
    \begin{cases}
        = 2\ell & \text{on } \supp(\mu_\mathrm{eq}), \\
        \geq 2\ell & \text{on } \Gamma,
	\end{cases}
\end{equation}
which hold for a certain constant $\ell$.

The system of contours $\Gamma$ has the $S$-property in the external field if
\begin{multline}\label{eq:S_property_1}
    \frac{\partial}{\partial n_+}\left(6\int G_{P_\infty}(p,q) d\mu_\mathrm{eq}(q) + \Re V \right) = \frac{\partial}{\partial n_-}\left(6\int G_{P_\infty}(p,q) d\mu_\mathrm{eq}(q) + \Re V\right) \\
	\quad \text{on } \supp(\mu_\mathrm{eq}),
\end{multline}
where $\frac{\partial}{\partial n_\pm}$ denote the derivatives in the two normal directions to the contour. We emphasize that the equality of the normal derivatives only needs to hold on the support of the equilibrium measure.

For our purposes, we want a system of contours $\Gamma$ that goes around the origin once on each sheet of the Riemann surface, and that also has the $S$-property in the external field $\Re V$.

The equilibrium problem in the case that we consider was solved in \cite{Kui25}, where it was found that one may take
\[
    \Gamma = z^{-1}(\mathbb{T}) = \mathbb{T}_1\cup\mathbb{T}_2\cup\mathbb{T}_3.
\]
We recall that $\mathbb{T}_j$ is the unit circle on the $j$th sheet of the Riemann surface, and that $z^{-1}(\mathbb{T})$ consists of two disjoint closed contours, namely $\mathbb{T}_1\cup\mathbb{T}_2$ and $\mathbb{T}_3$. Then, it was shown that
\[
    \supp(\mu_\mathrm{eq}) = \mathbb{T}_1\cup\mathbb{T}_2,
\]
and that $\mu_\mathrm{eq}$ is absolutely continuous with respect to Lebesgue measure  with a density that is positive on $\mathbb{T}_1\cup\mathbb{T}_2$. Moreover, a strict inequality holds in the variational conditions \eqref{eq:EL_conditions_1} on $\Gamma\setminus\supp(\mu_\mathrm{eq})=\mathbb{T}_3$. Thus, we obtain the following from \eqref{eq:ext_field} and \eqref{eq:EL_conditions_1}:
\begin{equation}\label{eq:EL_conditions_2} 
	3\int G_{P_\infty}(p,q)d\mu_\mathrm{eq}(q)+\log|z(p)|-\log|\lambda(p)|
    \begin{cases}
        = \ell & \text{on } \mathbb{T}_1\cup\mathbb{T}_2, \\
        > \ell & \text{on } \mathbb{T}_3.
    \end{cases}
\end{equation}

We observe that the equality in \eqref{eq:EL_conditions_2} is essential to the geometry of the double cover $\widehat{\mathcal{R}}$, as we constructed it by gluing two copies of $\mathcal{R}$ along the support of the equilibrium measure. The inequality in \eqref{eq:EL_conditions_2} plays a less prominent role in this paper. However, we emphasize that it is indispensable to the steepest descent analysis of the RH problem in \cite{Kui25}.


\subsection{Meromorphic differential $\omega$ and the $g$-function}\label{sec:g_function}
The left-hand side of \eqref{eq:EL_conditions_2} defines a function that is harmonic on $\mathcal{R}\setminus(\mathbb{T}_1\cup\mathbb{T}_2\cup\mathcal{P})$. By \cite[Lemma~6.3 and Definition~6.4]{Kui25}, there exists an analytic function $g$ on $\mathcal{R}\setminus(\mathbb{T}_1\cup\mathbb{T}_2\cup z^{-1}(\overline{\mathbb{R}}))$, with $\overline{\mathbb{R}} = \mathbb{R}\cup\{\infty\}$, such that
\begin{equation}\label{eq:Reg_p_1}
    \Re g(p) = -3\int G_{P_\infty}(p,q)d\mu_\mathrm{eq}(q).
\end{equation} 
In fact, $\Re g$ extends harmonically to $\mathcal{R}\setminus(\mathbb{T}_1\cup\mathbb{T}_2\cup\{P_\infty\})$ with 
\begin{equation}\label{eq:Reg_at_Pinfty}
    \Re g(p) = -3 \log|z_\infty(p)| + o(1) \quad \text{as } p \to P_\infty,
\end{equation}
see \eqref{eq:bipolar_G_at_Pinfty}. This $g$-function plays a crucial role in the steepest descent analysis of the RH problem in \cite{Kui25}.

From \eqref{eq:EL_conditions_2} and \eqref{eq:Reg_p_1}, it follows that $g_+ + g_- + 2\log\lambda - 2\log z$ has constant real part on $\mathbb{T}_1\cup\mathbb{T}_2$. The $S$-property \eqref{eq:S_property_1} implies by means of the Cauchy--Riemann equations that the imaginary part is constant as well on $\mathbb{T}_1\cup\mathbb{T}_2$. Thus, there is a constant $\widehat{\ell}$ with $\Re\widehat{\ell} = \ell$ such that
\begin{equation}\label{eq:S_property_2}
    g_+ + g_- + 2\log\lambda - 2\log z = -2\widehat{\ell} - \pi i \quad \text{on } \mathbb{T}_1\cup\mathbb{T}_2,
\end{equation}
see also \cite[Section~6]{Kui25}. The identity \eqref{eq:S_property_2} assumes a continuous choice of the branches of the logarithms along $\mathbb{T}_1\cup\mathbb{T}_2$. If we do not specify the branches of $\log z$ and $\log\lambda$, then \eqref{eq:S_property_2} should be considered as an identity modulo $2\pi i$.

\begin{proposition}\label{prop:omega_as_dg_dlogl_dlogz}
    The differential $\omega$ as defined in Definition~\ref{def:omega} is given by
    \begin{equation}\label{eq:omega_as_dg_dlogl_dlogz}
		\omega =
		\begin{cases}
            \hspace{0.32cm} dg + d\log\lambda - d\log z & \text{on } \mathcal{R}, \\
			-dg - d\log\lambda + d\log z & \text{on } \mathcal{R}^\sigma.
		\end{cases}
	\end{equation}
\end{proposition}

\begin{proof}
    Recall from Definition~\ref{def:omega} that $\omega$ is the unique meromorphic differential on $\widehat{\mathcal{R}}$ with (simple) poles and residues as listed in Table~\ref{table:poles_omega}, and with the additional requirement that all its periods are purely imaginary. We proceed to show that the differential defined by the right-hand side of \eqref{eq:omega_as_dg_dlogl_dlogz} satisfies these requirements.

    Since the $g$-function is analytic in $\mathcal{R}\setminus(\mathbb{T}_1\cup\mathbb{T}_2\cup z^{-1}(\overline{\mathbb{R}}))$, it is clear that \eqref{eq:omega_as_dg_dlogl_dlogz} is analytic in that set. The $g$-function does not extend analytically across $z^{-1}(\overline{\mathbb{R}})$, as it has different limits from the two sides. However, these limits differ by purely imaginary constants, see \cite[Lemma~6.3 and Definition~6.4]{Kui25}, which is enough to conclude that \eqref{eq:omega_as_dg_dlogl_dlogz} extends to a meromorphic differential on $\mathcal{R}\setminus(\mathbb{T}_1\cup\mathbb{T}_2)$.

    There are simple poles at $P_1$ and $P_0$ with respective residues $3$ and $-3$ due to the terms $d\log\lambda$ and $-d\log z$, as $\lambda$ has triple zero at $P_1$ and $z$ a triple zero at $P_0$. All three terms in \eqref{eq:omega_as_dg_dlogl_dlogz} have a simple pole at $P_\infty$. The residue of $dg$ at $P_\infty$ is $-3$ because of the behavior \eqref{eq:Reg_at_Pinfty} of the $g$-function at $P_\infty$. The remaining two terms in \eqref{eq:omega_as_dg_dlogl_dlogz} have residues $-3$ and $3$ at $P_\infty$, and hence $\omega$ has residues $-3$ at $P_\infty$.

    From these considerations, we also conclude that \eqref{eq:omega_as_dg_dlogl_dlogz} is meromorphic on $\mathcal{R}^\sigma\setminus(\mathbb{T}_1^\sigma\cup\mathbb{T}_2^\sigma)$ with simple poles $P_\infty^\sigma$, $P_1^\sigma$, and $P_0^\sigma$ with respective residues $3$, $-3$, and $3$. To prove that \eqref{eq:omega_as_dg_dlogl_dlogz} defines a meromorphic differential on $\widehat{\mathcal{R}}$, we have to show that the expressions in the two cases of \eqref{eq:omega_as_dg_dlogl_dlogz} have the same limiting values on $\mathbb{T}_1\cup\mathbb{T}_2$.
	By \eqref{eq:EL_conditions_2} and \eqref{eq:Reg_p_1}, we have
    \begin{equation}\label{eq:Reg_p_2}
        \Re\!\left(g+\log\lambda-\log z\right) = -\ell \quad \text{on } \mathbb{T}_1\cup\mathbb{T}_2.
    \end{equation}
    Thus, the tangential derivative of $\Re(g+\log\lambda-\log z)$ vanishes along $\mathbb{T}_1\cup\mathbb{T}_2$. By the Cauchy--Riemann equations, it follows that
	\begin{equation}\label{eq:Reg_p_3}
	   \frac{\partial}{\partial n_+}\Im(g+\log\lambda-\log z) = \frac{\partial}{\partial n_-}\Im(g+\log\lambda-\log z) = 0
    \end{equation}
	on $\mathbb{T}_1\cup\mathbb{T}_2$. Recall that the $S$-property \eqref{eq:S_property_1} states that
	\begin{equation}\label{eq:Reg_p_4}
        \frac{\partial}{\partial n_+} \Re(g+\log\lambda-\log z) = \frac{\partial}{\partial n_-}\Re(g+\log\lambda-\log z)
    \end{equation}
    on $\mathbb{T}_1\cup\mathbb{T}_2$. Combining \eqref{eq:Reg_p_3} and \eqref{eq:Reg_p_4} yields
    \begin{equation}\label{eq:Reg_p_5}
        \frac{\partial}{\partial n_+}(g+\log\lambda-\log z) = \frac{\partial}{\partial n_-}(g+ \log\lambda-\log z)\quad \text{on } \mathbb{T}_1\cup\mathbb{T}_2.
	\end{equation}
	
    Since the function $g+\log\lambda-\log z$ is analytic on both sides of $\mathbb{T}_1\cup\mathbb{T}_2$, we get from \eqref{eq:Reg_p_5} that $d(g+\log\lambda-\log z)_+ = -d(g+\log\lambda-\log z)_-$ on $\mathbb{T}_1\cup\mathbb{T}_2$. This indeed shows that $-g-\log\lambda+\log z$ is the meromorphic continuation of $g+\log\lambda-\log z$ across $\mathbb{T}_1\cup\mathbb{T}_2$.

    Finally, since $\log z$ and $\log\lambda$ have single-valued real parts, the differentials $d\log z$ and $d\log\lambda$ have purely imaginary periods. Also, $\Re g$ is well-defined and single-valued on $\mathcal{R}\setminus(\mathbb{T}_1\cup\mathbb{T}_2)$. It follows that the integral of \eqref{eq:omega_as_dg_dlogl_dlogz} along any closed contour that is entirely contained in $\mathcal{R}$ or in $\mathcal{R}^\sigma$ and avoids the poles of \eqref{eq:omega_as_dg_dlogl_dlogz} is purely imaginary. By symmetry, it is clear that the integral of \eqref{eq:omega_as_dg_dlogl_dlogz} along the bounded oval $\widehat{\mathcal{O}}_b$ equals $0$. This is sufficient to conclude that all periods are purely imaginary. This completes the proof.
\end{proof}


\subsection{Symmetries of the double cover}\label{sec:symmetries}
In the proofs that follow, we will make use of certain symmetries of the spectral curve and their extensions to the double cover.

As noted in \cite[(2.6)]{Kui25}, the spectral curve \eqref{eq:spectral_curve} admits the symmetries
\begin{equation}\label{eq:sigmas}
	(z,\lambda) \mapsto (-\lambda,-z), \quad (z,\lambda) \mapsto \bigl(\tfrac{1}{z},\tfrac{\lambda}{z}\bigr), \quad (z, \lambda) \mapsto \bigl(-\tfrac{z}{\lambda},\tfrac{1}{\lambda}\bigr),
\end{equation}
which are all holomorphic involutions on $\mathcal{R}$. These symmetries extend to the double cover $\widehat{\mathcal{R}}$. The first symmetry in \eqref{eq:sigmas} on $\mathcal{R}$ interchanges $P_0$ to $P_1$ and leaves $P_\infty$ invariant. We extend it to $\widehat{\mathcal{R}}$ in such a way that it still leaves $P_\infty$ invariant, but then to obtain an automorphism, we are forced to let it map $P_0$ to $P_1^\sigma$ and $P_1$ to $P_0^\sigma$. The action is more clear on the double complex torus (see also Figure~\ref{fig:double_torus}), where it becomes
\begin{equation}\label{eq:sigma_1}
    \sigma_1 \colon \widehat{\mathcal{R}} \to \widehat{\mathcal{R}}, \quad u \mapsto \tfrac{1}{3} - u \quad (\mod 2\mathbb{Z} + \tau \mathbb{Z}). 
\end{equation}
Similarly, the second and third symmetries in \eqref{eq:sigmas} correspond to
\begin{align}
	\sigma_2 \colon \widehat{\mathcal{R}} \to \widehat{\mathcal{R}}, \quad u \mapsto 1 - u \quad (\mod 2\mathbb{Z} + \tau \mathbb{Z}), \label{eq:sigma_2} \\
    \sigma_3 \colon \widehat{\mathcal{R}} \to \widehat{\mathcal{R}}, \quad u \mapsto \tfrac{5}{3} - u \quad (\mod 2\mathbb{Z} + \tau \mathbb{Z}), \label{eq:sigma_3}
\end{align}
and the involution \eqref{eq:sigma} that interchanges $\mathcal{R}$ to $\mathcal{R}^\sigma$ leads to
\begin{equation}\label{eq:sigma_0}
    \sigma \colon \widehat{\mathcal{R}} \to \widehat{\mathcal{R}}, \quad u \mapsto 1 + u \quad (\mod 2\mathbb{Z} + \tau \mathbb{Z}).
\end{equation}

The involutions \eqref{eq:sigma_1}--\eqref{eq:sigma_0} generate the dihedral group $D_6$, that is, the symmetry group of a regular hexagon. This group is also generated by 
\begin{equation}\label{eq:rho_def} 
	\rho = \sigma_1\circ\sigma_2\circ\sigma \colon \widehat{\mathcal{R}} \to \widehat{\mathcal{R}}, \quad u \mapsto \tfrac{1}{3} + u \quad (\mod 2\mathbb{Z} + \tau \mathbb{Z}),
\end{equation}
together with one of the $\sigma_j$'s with $j=1,2,3$. Then, $\rho$ acts as a rotation of the hexagon, and $\sigma_j$ acts as a reflection.

\begin{figure}[t]
    \centering
    \begin{tikzpicture}[scale=0.5](15,10)(0,0)
        \filldraw[gray!50!white] (-6,0) --(18,0) --(18,3) --(-6,3) --(-6,0);
        \filldraw[gray!20!white] (-6,0) --(18,0) --(18,-3) --(-6,-3) --(-6,0);

		\draw[ultra thick] (-6,0)--(18,0);
		\draw[ultra thick] (-6,3)--(18,3) (-6,-3)--(18,-3);

        \draw (-6.3,3) node[left] {$\frac{1}{2}\tau$};
        \draw (-6.3,-3) node[left] {$-\frac{1}{2}\tau$};
        \draw (18.3,3) node[right] {$2+\frac{1}{2}\tau$};
        \draw (18.3,-3) node[right] {$2-\frac{1}{2}\tau$};

        \draw[ultra thick,red] (-6,-3)--(-6,3) (6,-3)--(6,3);
        \draw[ultra thick,dashed,red] (-2,-2.85)--(-2,3) (2,-2.85)--(2,3);
        \draw[ultra thick,dashed,blue] (-4,-2.85)--(-4,3) (0,-2.85)--(0,3) (4,-2.85)--(4,3);
        \draw (-6,-1.5) node[right] {$\Gamma_0$};
        \draw (-2,-1.5) node[right] {$\Gamma_1$};
        \draw (2,-1.5) node[right] {$\Gamma_2$};
        \draw (-4,-1.5) node[right] {$\Upsilon_0$};
        \draw (0,-1.5) node[right] {$\Upsilon_1$};
        \draw (4,-1.5) node[right] {$\Upsilon_2$};
        \draw (18,-1.5) node[right] {$\Gamma_0$};

        \draw[ultra thick,red] (18,-3)--(18,3);
        \draw[ultra thick,dashed,red] (10,-2.85)--(10,3) (14,-2.85)--(14,3);
        \draw[ultra thick,dashed,blue] (8,-2.85)--(8,3) (12,-2.85)--(12,3) (16,-2.85)--(16,3);
        \draw (6,-1.5) node[right] {$\Gamma_0^\sigma$};
        \draw (10,-1.5) node[right] {$\Gamma_1^\sigma$};
        \draw (14,-1.5) node[right] {$\Gamma_2^\sigma$};
        \draw (8,-1.5) node[right] {$\Upsilon_0^\sigma$};
        \draw (12,-1.5) node[right] {$\Upsilon_1^\sigma$};
        \draw (16,-1.5) node[right] {$\Upsilon_2^\sigma$};

        \filldraw (-6,3) circle (4pt);
        \filldraw (-2,3) circle (4pt);
        \filldraw (2,3) circle (4pt);
        \filldraw (6,3) circle (4pt);
        \filldraw (-6,-3) circle (4pt);
        \filldraw (-2,-3) circle (4pt);
        \filldraw (2,-3) circle (4pt);
        \filldraw (6,-3) circle (4pt);
        \draw (-5.9,3.075) node[above] {$R_0$};
        \draw (-1.9,3.075) node[above] {$R_1$};
        \draw (2.1,3.075) node[above] {$R_2$};
        \draw (18.1,3.075) node[above] {$R_0$};
        \draw (-5.9,-3.075) node[below] {$R_0$};
        \draw (-1.9,-3.075) node[below] {$R_1$};
        \draw (2.1,-3.075) node[below] {$R_2$};
        \draw (18.1,-3.075) node[below] {$R_0$};

        \filldraw (-4,3) circle (4pt);
        \filldraw (0,3) circle (4pt);
        \filldraw (4,3) circle (4pt);
        \filldraw (-4,-3) circle (4pt);
        \filldraw (0,-3) circle (4pt);
        \filldraw (4,-3) circle (4pt);
        \draw (-3.9,3.075) node[above] {$Z_0$};
        \draw (0.1,3.075) node[above] {$Z_1$};
        \draw (4.1,3.075) node[above] {$Z_2$};
        \draw (-3.9,-3.075) node[below] {$Z_0$};
        \draw (0.1,-3.075) node[below] {$Z_1$};
        \draw (4.1,-3.075) node[below] {$Z_2$};

        \filldraw (10,3) circle (4pt);
        \filldraw (14,3) circle (4pt);
        \filldraw (18,3) circle (4pt);
        \filldraw (10,-3) circle (4pt);
        \filldraw (14,-3) circle (4pt);
        \filldraw (18,-3) circle (4pt);
        \draw (6.1,3) node[above] {$R_0^\sigma$};
        \draw (10.1,3) node[above] {$R_1^\sigma$};
        \draw (14.1,3) node[above] {$R_2^\sigma$};
        \draw (6.1,-3) node[below] {$R_0^\sigma$};
        \draw (10.1,-3) node[below] {$R_1^\sigma$};
        \draw (14.1,-3) node[below] {$R_2^\sigma$};

        \filldraw (8,3) circle (4pt);
        \filldraw (12,3) circle (4pt);
        \filldraw (16,3) circle (4pt);
        \filldraw (8,-3) circle (4pt);
        \filldraw (12,-3) circle (4pt);
        \filldraw (16,-3) circle (4pt);
        \draw (8.1,3) node[above] {$Z_0^\sigma$};
        \draw (12.1,3) node[above] {$Z_1^\sigma$};
        \draw (16.1,3) node[above] {$Z_2^\sigma$};
        \draw (8.1,-3) node[below] {$Z_0^\sigma$};
        \draw (12.1,-3) node[below] {$Z_1^\sigma$};
        \draw (16.1,-3) node[below] {$Z_2^\sigma$};
    \end{tikzpicture}
    \caption{The double complex torus with the cycles $\Gamma_j$, $\Gamma_j^\sigma$ (red), and $\Upsilon_j$, $\Upsilon_j^\sigma$ (blue) for $j=0,1,2$, together with their intersections $R_j$, $R_j^\sigma$, and $Z_j$, $Z_j^\sigma$ with the bounded oval.}
    \label{fig:torus_cycles}
\end{figure}
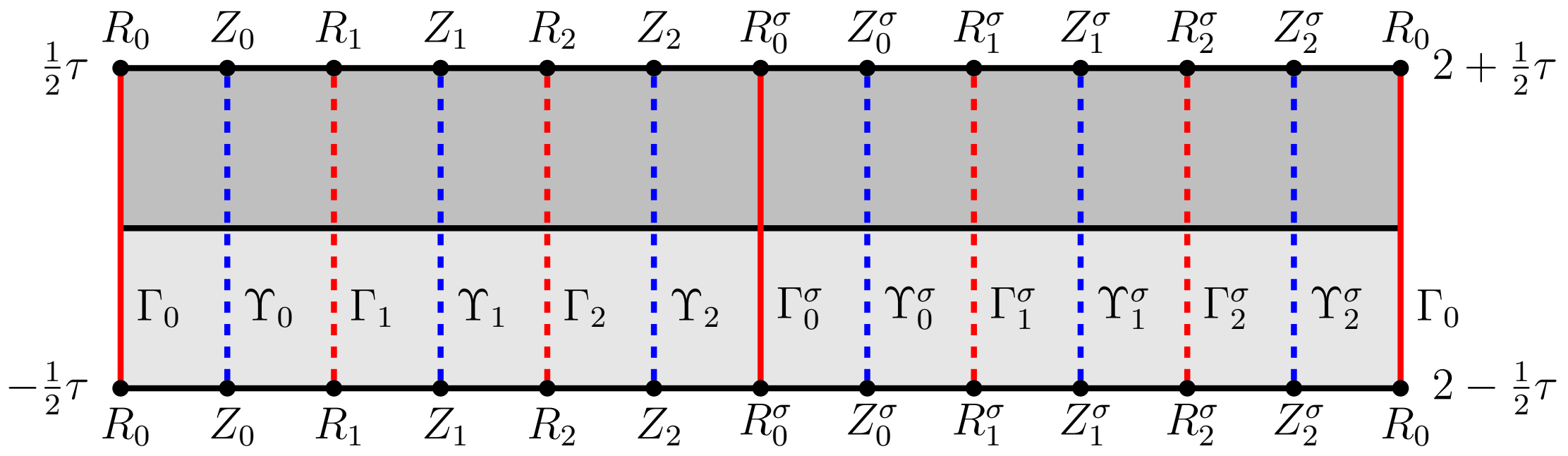

We also introduce a number of cycles in $\widehat{\mathcal{R}}$ that play a role in what follows. These cycles are again easiest to describe on the double complex torus:
\begin{equation}\label{eq:Gammas_def}
	\begin{aligned}
        \Gamma_0 \colon \quad & u\in[0,\tau], & \qquad \Gamma_0^\sigma \colon \quad & u\in1+[0,\tau], \\
		\Gamma_1 \coloneqq \; \rho\,(\Gamma_0) \colon \quad & u\in\tfrac{1}{3}+[0,\tau], & \qquad \Gamma_1^\sigma \colon \quad & u\in\tfrac{4}{3}+[0,\tau], \\
		\Gamma_2 \coloneqq \rho^2(\Gamma_0) \colon \quad & u\in\tfrac{2}{3}+[0,\tau], & \qquad \Gamma_2^\sigma \colon \quad & u\in\tfrac{5}{3}+[0,\tau],
	\end{aligned}
\end{equation}
and
\begin{equation}\label{eq:Upsilons_def}
	\begin{aligned}
        \Upsilon_0 \colon \quad & u\in\tfrac{1}{6}+[0,\tau], & \qquad \Upsilon_0^\sigma \colon \quad & u\in\tfrac{7}{6}+[0,\tau], \\
		\Upsilon_1 \coloneqq \; \rho\,(\Upsilon_0) \colon \quad & u\in\tfrac{1}{2}+[0,\tau], & \qquad \Upsilon_1^\sigma \colon \quad & u\in\tfrac{3}{2}+[0,\tau], \\
		\Upsilon_2 \coloneqq \rho^2(\Upsilon_0) \colon \quad & u\in\tfrac{5}{6}+[0,\tau], & \qquad \Upsilon_2^\sigma \colon \quad & u\in\tfrac{11}{6}+[0,\tau].
	\end{aligned}
\end{equation}
Figure~\ref{fig:torus_cycles} shows these cycles on the double complex torus. In the $(z,\lambda)$-coordinates, these cycles on $\mathcal{R}$ are characterized by:
\begin{equation}\label{eq:Gamma_Upsilons_char}
	\begin{aligned}
        \Gamma_0 \colon \quad & |z|=1<|\lambda|, & \qquad \Upsilon_0 \colon \quad & |z|=|\lambda|>1, \\
		\Gamma_1 \colon \quad & |\lambda|=1<|z|, & \qquad \Upsilon_1 \colon \quad & |z|=1>|\lambda|, \\
		\Gamma_2 \colon \quad & |z|=|\lambda|<1, & \qquad \Upsilon_2 \colon \quad & |\lambda|=1>|z|.
	\end{aligned}
\end{equation}
Note that $\Gamma_0 = \mathbb{T}_1\cup\mathbb{T}_2$ and $\Upsilon_1 = \mathbb{T}_3$.

\subsection{Zeros of the meromorphic differential $\omega$}
Our first use of the symmetries of the double cover discussed in Section~\ref{sec:symmetries} is to localize the zeros of the meromorphic differential $\omega$ given in \eqref{eq:omega_as_dg_dlogl_dlogz}.

\begin{proposition}\label{prop:omega_and_dg_2}
	The meromorphic differential $\omega$ has six simple zeros that are all on the bounded oval, which are denoted by $Z_0$, $Z_1$, and $Z_2$ on $\mathcal{R}$, and $Z_0^\sigma$, $Z_1^\sigma$, and $Z_2^\sigma$ on $\mathcal{R}^\sigma$, as shown in Figure~\ref{fig:torus_cycles}.
\end{proposition}

\begin{proof}
    Since $\mathcal{\widehat{R}}$ is of genus one, a meromorphic differential has the same number of zeros and poles, counted with multiplicity. Thus, there are six zeros of $\omega$, as it has six (simple) poles.
    
    Recall that $\omega$ is characterized by the location and residues of its poles, as well as by the property that it has purely imaginary periods. Together with the symmetries in the pole and residue structure, this will allow us to locate the zeros of $\omega$. In particular, for $\sigma_1$ and $\sigma_2$ from \eqref{eq:sigma_2}--\eqref{eq:sigma_3}, we find
    \begin{equation}\label{eq:omega_sigmas}
        \sigma_*(\omega) = -\omega, \quad \sigma_{1,*}(\omega) = \omega, \quad \sigma_{2,*}(\omega) = \omega
    \end{equation}
    by comparing the poles and their residues.

    Recall from \eqref{eq:rho_def} that the rotation symmetry $\rho$ is given by $\sigma_1\circ\sigma_2\circ\sigma$. It now follows from \eqref{eq:omega_sigmas} that $\rho_*(\omega) = -\omega$. Hence, any zero $Z$ of $\omega$ gives rise to five more zeros: $Z + \frac{j}{3}$ $(\mod 2\mathbb{Z}+\tau\mathbb{Z})$ for $j=1,\ldots,5$. We may assume that $0\leq\Re Z<\frac{1}{3}$.

    Composing $\sigma_2$ and $\sigma$ as in \eqref{eq:sigma_2} and \eqref{eq:sigma_0} yields the reflection symmetry $\sigma_2\circ\sigma \colon u \mapsto -u$ $(\mod 2\mathbb{Z} + \tau \mathbb{Z})$. Using \eqref{eq:omega_sigmas} again, we have $(\sigma_2\circ\sigma)_*(\omega) = -\omega$. Thus, $-Z$ is a zero of $\omega$ as well. It follows that either $-Z \equiv Z$ or $-Z \equiv Z + \frac{5}{3}$ $(\mod 2\mathbb{Z} + \tau \mathbb{Z})$. Together with $0\leq\Re Z<\frac{1}{3}$, this gives us the following four possibilities:
    \[
        Z \equiv 0, \quad Z \equiv \frac{1}{2}\tau, \quad Z \equiv \frac{1}{6}, \quad \text{or} \quad Z \equiv \frac{1}{6} + \frac{1}{2}\tau.
    \]

	According to Abel's theorem, the sum of the zeros equals the sum of the poles. The sum of the poles is
	\[
        \frac{1}{6} + \frac{1}{2} + \frac{5}{6} + \frac{7}{6} + \frac{3}{2} + \frac{11}{6} \equiv 0  \quad (\mod 2\mathbb{Z} + \tau \mathbb{Z}),
    \]
	and the zeros add up to
	\[
        \sum_{j=0}^5\left(Z + \frac{j}{3}\right) \equiv 6Z + 1 \quad (\mod 2\mathbb{Z} + \tau \mathbb{Z}).
	\]
	Therefore, we should have $6Z \equiv 1$. This means that only $Z \equiv \frac{1}{6}$ and $Z \equiv \frac{1}{6} + \frac{1}{2}\tau$ are the two remaining possibilities.

    We may exclude $Z\equiv \frac{1}{6}$, since in this case $Z$ would coincide with the pole $P_\infty$. Thus, $Z\equiv \frac{1}{6}+\frac{1}{2}\tau$, which implies that all the zeros of $\omega$ are given by
    \[
        Z_0 = Z, \quad Z_1 = Z+\frac{1}{3}, \quad Z_2 = Z+\frac{2}{3} \quad \text{on } \mathcal{R},
    \]
    and
    \[
        Z_0^\sigma = Z+1, \quad Z_1^\sigma = Z+\frac{4}{3}, \quad Z_2^\sigma = Z+\frac{5}{3} \quad \text{on } \mathcal{R}^\sigma.
    \]
    We observe that these zeros are equally spaced on the bounded oval, as shown in Figure~\ref{fig:torus_cycles}.
\end{proof}

\begin{remark}
	The zero $Z_1$ is the intersection of $\Upsilon_1$ with the bounded oval. In the $(z,\lambda)$-coordinates, it is given by $Z_1 = (1,\lambda_3(1))$. The remaining zeros on $\mathcal{R}$ are the images of $Z_1$ under the symmetries as in \eqref{eq:sigmas}. That is, $Z_0 = (-\lambda_3(1)^{-1},-\lambda_3(1)^{-1})$ and $Z_2 = (-\lambda_3(1),-1)$.
\end{remark}


\subsection{Action function and its saddle points}
The action function is a fundamental ingredient in the steepest descent analysis. It is defined as follows.

\begin{definition}\label{def:Phi_xi1_xi2}
	For $(\xi_1,\xi_2)\in\mathcal{H}$, the action function is given by
	\begin{equation}\label{eq:Phi_xi1_xi2_def}
		\Phi_{(\xi_1,\xi_2)} =
		\begin{cases}
            \hspace{0.32cm} g + (1-\xi_1)\log\lambda - (1-\xi_2)\log z + \widehat{\ell} & \text{on } \mathcal{R}, \\
			-g - (1+\xi_1)\log\lambda + (1+\xi_2)\log z - \widehat{\ell} + \pi i & \text{on } \mathcal{R}^\sigma,
		\end{cases}
	\end{equation}
    with $\widehat{\ell}$ as in \eqref{eq:S_property_2}. The special case corresponding to $(\xi_1,\xi_2) = (0,0)$ is denoted by
	\begin{equation}\label{eq:Phi_00_def}
		\Phi = \Phi_{(0,0)}.
	\end{equation}
\end{definition}

The definition~\eqref{eq:Phi_xi1_xi2_def} depends on the choice of branches of the logarithms. For definiteness, we choose the principal branches on the upper and lower half-planes of each sheet $\mathcal{R}_j$ and $\mathcal{R}^\sigma_j$ for $j=1,2,3$. Then, $\Phi_{(\xi_1,\xi_2)}$ has jump discontinuities whenever $z$ is real, that is on $z^{-1}(\overline{\mathbb{R}})$. There is an additional discontinuity across the cycles $\Gamma_0$ and $\Gamma_0^\sigma$, which connect $\mathcal{R}$ to its copy $\mathcal{R}^\sigma$, which is due to the appearance of the $g$-function in \eqref{eq:Phi_xi1_xi2_def}.

\begin{lemma}\label{lem:Phi_analytic}
	For $(\xi_1,\xi_2)\in\mathcal{H}$, the function $\Phi_{(\xi_1,\xi_2)}$ is analytic in $\widehat{\mathcal{R}}\setminus(\Gamma_0\cup\Gamma_0^\sigma$ $\cup\,z^{-1}(\overline{\mathbb{R}}))$ with piecewise constant jump discontinuities. Moreover, the jump discontinuities of $\Phi = \Phi_{(0,0)}$ are integer multiples of $2\pi i$, except on the bounded oval, where its jump discontinuity is an odd integer multiple of $\pi i$.
\end{lemma}

\begin{proof}
    The analyticity of $\Re\Phi$ was already addressed in the discussion preceding the lemma, hence we restrict ourselves to the case $(\xi_1,\xi_2) = 0$.

	The discontinuity on $\Gamma_0\cup\Gamma_0^\sigma$ comes from that of the $g$-function. It follows from the detailed formulas in \cite[Lemma~6.3 and Definition~6.4]{Kui25} that the jump discontinuities are integer multiples of $2\pi i$ (with possibly different multiples on the different sheets).
	
	The discontinuity on $z^{-1}(\overline{\mathbb{R}})$ comes from the chosen branches of the logarithms, as discussed above, as well as the discontinuity of $g$. The jump discontinuities of the logarithms are clearly multiples of $2\pi i$. The same holds for the jump discontinuities of $g$ on $z^{-1}(\overline{\mathbb{R}})$, as follows once again from the precise formulas in \cite{Kui25}.
\end{proof}

Comparing \eqref{eq:Phi_xi1_xi2_def} with \eqref{eq:omega_as_dg_dlogl_dlogz}, one can observe that
\begin{equation}\label{eq:omega_xi1_xi2_2}
    d\Phi_{(\xi_1,\xi_2)} =	\omega_{(\xi_1,\xi_2)} = \omega - \xi_1d\log\lambda + \xi_2d\log z,
\end{equation}
which is a meromorphic differential on $\widehat{\mathcal{R}}$. Thus, saddle points of $\Phi_{(\xi_1,\xi_2)}$ are precisely the zeros of $\omega_{(\xi_1,\xi_2)}$.

Before we can localize these zeros, more properties of $\Phi_{(\xi_1,\xi_2)}$ are needed. In particular, we again use the symmetries from Section~\ref{sec:symmetries} to determine the behavior of $\Re\Phi$ on the cycles as in \eqref{eq:Gammas_def}--\eqref{eq:Upsilons_def}.

\begin{lemma}\label{lem:Phi_xi1_xi2_rho}
	For the rotation symmetry $\rho$ as in \eqref{eq:rho_def}, one has
    \begin{align} 
        \Re\Phi_{(\xi_1,\xi_2)}\circ\rho &= -\Re\Phi_{(\xi_2,\xi_2-\xi_1)}, \label{eq:Phi_xi1_xi2_rho} \\
        \Re\Phi_{(\xi_1,\xi_2)}\circ\rho^{-1} &= -\Re\Phi_{(\xi_1-\xi_2,\xi_1)}. \label{eq:Phi_xi1_xi2_rho_inv}
    \end{align}
\end{lemma}

\begin{proof}
    By definition~\eqref{eq:rho_def}, $\rho$ is the composition $\sigma_1\circ\sigma_2\circ\sigma$. Hence, \eqref{eq:Phi_xi1_xi2_rho}--\eqref{eq:Phi_xi1_xi2_rho_inv} follow once we establish the following:
    \begin{equation}\label{eq:Phi_xi1_xi2_sigmas}
        \begin{aligned}
            \Re\Phi_{(\xi_1,\xi_2)}\circ\sigma\,\, &= -\Re\Phi_{(-\xi_1,-\xi_2)}, \\
            \Re\Phi_{(\xi_1,\xi_2)}\circ\sigma_1 &= \Re\Phi_{(-\xi_2,-\xi_1)}, \\
            \Re\Phi_{(\xi_1,\xi_2)}\circ\sigma_2 &= \Re\Phi_{(\xi_1,\xi_1-\xi_2)}.
        \end{aligned}
    \end{equation}

    The first identity in \eqref{eq:Phi_xi1_xi2_sigmas} is an immediate consequence of the definition~\eqref{eq:Phi_xi1_xi2_def} of $\Phi_{(\xi_1,\xi_2)}$. For the second identity, we recall from \eqref{eq:omega_sigmas} that $\sigma_{1,*}(\omega) = \omega$. Together with $\sigma_{1,*}(d\log z) = d\log\lambda$ and $\sigma_{1,*}(d\log\lambda) = d\log z$, we see that $\sigma_{1,*}(\omega_{(\xi_1,\xi_2)}) = \omega_{(-\xi_2,-\xi_1)}$. Or equivalently, $d(\Phi_{(\xi_1,\xi_2)}\circ\sigma_1) = d\Phi_{(-\xi_2,-\xi_1)}$, which implies the second identity in \eqref{eq:Phi_xi1_xi2_sigmas} up to a constant. This constant can then be shown to be zero by comparing the values at the fixed point $Z_0$ of $\sigma_1$.

    The proof of the third identity in \eqref{eq:Phi_xi1_xi2_sigmas} is similar to that of the second.
\end{proof}

\begin{lemma}\label{lem:Phi_00_on_Gammas_Upsilons} For $\Phi = \Phi_{(0,0)}$ the following hold.
	\begin{enumerate}
        \item[\rm{(a)}] $\Re\Phi = 0$ on $\Gamma_j\cup\Gamma_j^\sigma$ for $j=0,1,2$.
        \item[\rm{(b)}] $\Re\Phi$ is positive and strictly increasing on $\Upsilon_0\cap\mathcal{R}_+$, $\Upsilon_2\cap\mathcal{R}_+$, and $\Upsilon_1^\sigma\cap\mathcal{R}_+^\sigma$.
        \item[\rm{(c)}] $\Re\Phi$ is negative and strictly decreasing on $\Upsilon_1\cap\mathcal{R}_+$, $\Upsilon_0^\sigma\cap\mathcal{R}_+^\sigma$, and $\Upsilon_2^\sigma\cap\mathcal{R}_+^\sigma$.
    \end{enumerate}
    Here, the curves are oriented from the bounded oval to the unbounded oval.
\end{lemma}

\begin{proof}
    (a) Recalling the definition~\eqref{eq:Phi_xi1_xi2_def} of $\Re\Phi = \Re\Phi_{(0,0)}$ and $\Gamma_0 = \mathbb{T}_1\cup\mathbb{T}_2$, one sees that part~(a) for $j=0$ is an immediately consequence of \eqref{eq:EL_conditions_2} and \eqref{eq:Reg_p_1}. The cases $j=1,2$ are obtained from this by using the identity $\Re\Phi\circ\rho\,^j = (-1)^j\Re\Phi$, which follows from \eqref{eq:Phi_xi1_xi2_rho}, recalling that $\Gamma_j = \rho\,^j(\Gamma_0)$.

    \medskip

    (b) \& (c) The function $\Re\Phi$ is harmonic on the double cover, except at the points in $\widehat{\mathcal{P}}$, see \eqref{eq:calP_hat}, where there is a logarithmic behavior tending to $+\infty$ (at $P_\infty$, $P_0$, and $P_1^\sigma$) or $-\infty$ (at the other three points $P_\infty^\sigma$, $P_0^\sigma$, and $P_1$). Each connected component of $\widehat{\mathcal{R}}\setminus \bigcup_{j=0}^2(\Gamma_j\cup\Gamma_j^\sigma)$ contains one of the points in $\widehat{\mathcal{P}}$, and $\Re\Phi = 0$ on its boundary by part (a). The maximum principle of harmonic functions then implies that $\Re\Phi$ has a constant (non-zero) sign on each component. This proves the statements about the sign of $\Re\Phi$ in parts (b) and (c).    
    
    The critical points of $\Re\Phi$ are precisely the zeros of the meromorphic differential $\partial\Re\Phi = \frac{1}{2}\omega$. By Proposition~\ref{prop:omega_and_dg_2}, the critical points are at $Z_j$ and $Z_j^\sigma$ for $j=0,1,2$, as shown in Figure~\ref{fig:torus_cycles}, and these are the only critical points.
    By \eqref{eq:Phi_xi1_xi2_sigmas}, we have $\Re\Phi\circ\sigma_2 = \Re\Phi$. Combining this with $\Re\Phi(\bar{p}) = \Re\Phi(p)$ yields the symmetry
    \begin{equation}\label{eq:Phi_sigma2_bar}
        \Re\Phi(\sigma_2(\bar{p})) = \Re\Phi(p), \quad p\in\widehat{\mathcal{R}}.
    \end{equation}
    From the symmetry \eqref{eq:Phi_sigma2_bar} and noting that $\Upsilon_1 = \mathbb{T}_3$ is fixed under the map $p \mapsto \sigma_2(\bar{p})$, the Cauchy--Riemann equations imply that $\Im\Phi$ is constant along $\Upsilon_1$. Since $Z_1$ is the only critical point on $\Upsilon_1\cap\mathcal{R}_+$, it follows that $\Re\Phi$ is strictly monotone along this curve. It is then strictly decreasing as $\Re\Phi$ tends to $-\infty$ at $P_1$.

    The remaining cases are obtained from the first by applying the symmetries \eqref{eq:Phi_xi1_xi2_rho}--\eqref{eq:Phi_xi1_xi2_rho_inv}. Note that the appearance of the minus signs in \eqref{eq:Phi_xi1_xi2_rho}--\eqref{eq:Phi_xi1_xi2_rho_inv} leads to the alternating pattern of decrease and increase.
\end{proof}


\subsection{Proof of Theorem~\ref{thm:four_zeros_on_bounded_oval}}
\begin{proof}
    (a) Since the meromorphic differential $\omega_{(\xi_1,\xi_2)}$ has six simple poles and the genus of $\widehat{\mathcal{R}}$ is one, it follows that $\omega_{(\xi_1,\xi_2)}$ has six zeros. The quotient
    \begin{equation}\label{eq:f_xi1_xi2}
        f_{(\xi_1,\xi_2)} = \frac{\omega_{(\xi_1,\xi_2)}}{d\log z},
    \end{equation}
    is a meromorphic function on $\widehat{\mathcal{R}}$ that is  real-valued on the real ovals, due to the symmetry
    \begin{equation}\label{eq:f_xi1_xi2_bar}
        f_{(\xi_1,\xi_2)}(\bar{p}) = \overline{f_{(\xi_1,\xi_2)}(p)}, \quad p\in\widehat{\mathcal{R}}.
    \end{equation}
    It then follows by \eqref{eq:omega_xi1_xi2_2} that $\Phi_{(\xi_1,\xi_2)}$ has a constant imaginary part on the bounded oval $\widehat{\mathcal{O}}_b$. The constant imaginary part is actually a different one on the two sides of $\widehat{\mathcal{O}}_b$, see Lemma~\ref{lem:Phi_analytic}. Critical points of $\Re\Phi$ on the bounded oval thus occur at a zero of $\omega_{(\xi_1,\xi_2)} =  d\Phi_{(\xi_1,\xi_2)}$.	By compactness, there is a global maximum and a global minimum on $\widehat{\mathcal{O}}_b$. Hence, there are at least two zeros of $\omega_{(\xi_1,\xi_2)}$ on $\widehat{\mathcal{O}}_b$. From these considerations, we also conclude that there are an even number of zeros on $\widehat{\mathcal{O}}_b$.

    For $j=0,1,2$, we use $R_j$ and $R_j^\sigma$ to denote the intersection points of the cycles $\Gamma_j$ and $\Gamma_j^\sigma$ with the bounded oval $\widehat{\mathcal{O}}_b$, see also Figure~\ref{fig:torus_cycles}. By Lemma~\ref{lem:Phi_00_on_Gammas_Upsilons}(a), we have that $\Re\Phi$ vanishes at these points. It then follows from \eqref{eq:Phi_xi1_xi2_def} and \eqref{eq:Phi_00_def} that 
    \begin{equation}\label{eq:Phi_xi1_xi2_Rjs}
        \Re\Phi_{(\xi_1,\xi_2)}(R_j) = \Re\Phi_{(\xi_1,\xi_2)}(R_j^\sigma), \quad j=0,1,2.
    \end{equation}
    By Rolle's theorem, there is a local maximum or local minimum in the part from $R_0$ to $R_0^\sigma$ along the bounded oval (going from left to right in Figure~\ref{fig:torus_cycles}), and also in the complementary part from $R_0^\sigma$ to $R_0$. Recall that these critical points are zeros of $\omega_{(\xi_1,\xi_2)}$, and their multiplicity is odd.

    Suppose, to get a contradiction, that these are the only zeros of $\omega_{(\xi_1,\xi_2)}$ on $\widehat{\mathcal{O}}_b$ with odd multiplicity. Then, $\Re\Phi_{(\xi_1,\xi_2)}$ has one global minimum and one global maximum on the bounded oval, and no other local extrema. Because of \eqref{eq:Phi_xi1_xi2_Rjs} with $j=0$, the global extrema are not at $R_0$ or $R_0^\sigma$. Let us consider the case that $\Re\Phi_{(\xi_1,\xi_2)}$ has its global maximum in the open part between $R_0$ and $R_0^\sigma$ and its global minimum in the complementary part.  The other case is treated similarly. Then,
    \[
        \Re\Phi_{(\xi_1,\xi_2)}(p)>\Re\Phi_{(\xi_1,\xi_2)}(R_0)
    \]
    on the part of the bounded oval from $R_0$ to $R_0^\sigma$, and
    \[
        \Re\Phi_{(\xi_1,\xi_2)}(p)<\Re\Phi_{(\xi_1,\xi_2)}(R_0)
    \]
    on the complementary part. Evaluating these inequalities at $R_1$ and $R_1^\sigma$ gives $\Re\Phi_{(\xi_1,\xi_2)}(R_1)>\Re\Phi_{(\xi_1,\xi_2)}(R_1^\sigma)$. However, this contradicts the equality \eqref{eq:Phi_xi1_xi2_Rjs} for $j=1$.

    Therefore, there are at least three zeros of $\omega_{(\xi_1,\xi_2)}$ on the bounded oval with odd multiplicity. The total number of zeros (counted with multiplicities) is even, and it follows that there are at least four zeros on the unbounded oval with odd multiplicity.

    \medskip

    (b) This follows directly from the symmetry \eqref{eq:f_xi1_xi2_bar}.

    \medskip

    (c) The statement follows immediately from parts (a) and (b) when not all zeros lie on the bounded oval.

    In case that all zeros lie on the bounded oval, we distinguish between two cases. In the first case, there are at least four simple zeros on the bounded oval. Then, the remaining two zeros could either be simple as well, or form one double zero. Second, there are at most three simple zeros on the bounded oval. Since there are at least four zeros with odd multiplicity, and the total number of zeros is six, there must be a triple zero in this case, along with three simple zeros. In both cases, we indeed observe that there is at most one zero with higher multiplicity.

    \medskip

    (d) Suppose that $(\xi_1,\xi_2)\in\mathbb{R}^2\setminus\mathcal{H}$. If none of the residues in  Table~\ref{table:poles_omega_xi1_xi2} vanishes, then there are consecutive poles with residues of the same sign. There is a zero on the unbounded oval in between two such poles. If one of the residues in Table~\ref{table:poles_omega_xi1_xi2} vanishes, then according to the convention used in Remark~\ref{rem:double_zero}, the corresponding point is treated as a zero. In both cases, there is a zero on the unbounded oval. The other remaining zero will then also be on the unbounded oval.

    The cases in which a double zero occurs on the unbounded oval can be identified explicitly. Otherwise, there are two simple zeros. By the pervious discussion, a double zero, in the sense of Remark~\ref{rem:double_zero}, can only occur when one of the residues Table~\ref{table:poles_omega_xi1_xi2} vanishes, and it must occur at one of the points in $\widehat{\mathcal{P}}$. Moreover, we will show that this happens if and only if $(\xi_1,\xi_2)\in\mathfrak{T}$, with $\mathfrak{T}$ as in \eqref{eq:turning_points}. Table~\ref{table:zero_at_P} lists the corresponding point $P(\xi_1,\xi_2)$ for each of the six cases.
	\begin{table}[ht]
		\centering
		\renewcommand{\arraystretch}{1.25}
		\setlength{\tabcolsep}{5pt}
		\begin{tabular}{c | c c c c c c}
		\hline
		$(\xi_1,\xi_2)$ & $(\frac{1}{2},-\frac{1}{2})$ & $(1,\frac{1}{2})$ & $(\frac{1}{2},1)$ & $(-\frac{1}{2}, \frac{1}{2})$ & $(-1,-\frac{1}{2})$ & $(-\frac{1}{2},-1)$ \\ 
		\hline
		$P(\xi_1,\xi_2)$ & $P_\infty$ & $P_1$ & $P_0$ & $P_\infty^\sigma$ & $P_1^\sigma$ & $P_0^\sigma$ \\
		\hline
		\end{tabular}
		\renewcommand{\arraystretch}{1}
        \caption{Points $(\xi_1,\xi_2)$ for which $\omega_{(\xi_1,\xi_2)}$ has a genuine zero (and a double zero in the sense of Remark~\ref{rem:double_zero}) at one of the points in $\widehat{\mathcal{P}}$.}
		\label{table:zero_at_P}
	\end{table}

    One of the residues in Table~\ref{table:poles_omega_xi1_xi2} is zero if and only if $\xi_1=\pm1,\xi_2=\pm1$, or $ \xi_1-\xi_2=\pm1$. Let us consider $\xi_2=1$, which corresponds to a zero residue at $P_0$. Then, $\omega_{(\xi_1,1)}=\omega-\xi_1d\log\lambda+d\log z$ has a genuine zero at $P_0$ if and only if
    \begin{equation}\label{eq:xi1_at_P0}
        \xi_1 = \left(\frac{\omega+d\log z}{d\log\lambda}\right)\!(P_0).
    \end{equation}
    Thus, there is a unique $\xi_1$ such that $\omega_{(\xi_1,1)}$ has a double zero, which then occurs at $P(\xi_1,1) = P_0$. One might evaluate $\xi_1$ directly from \eqref{eq:xi1_at_P0}. However, it is easier to look at the residues. The respective residues are as in Table~\ref{table:poles_omega_P0}, which follow immediately from Table~\ref{table:poles_omega_xi1_xi2}. 
    \begin{table}[ht]
    	\centering
    	\renewcommand{\arraystretch}{1.25}
    	\setlength{\tabcolsep}{5pt}
    	\begin{tabular}{c | c c c c c c}
    	\hline
    	Pole & $P_\infty$ & $P_1$ & $P_0$ & $P_\infty^\sigma$ & $P_1^\sigma$ & $P_0^\sigma$ \\
    	\hline
    	Residue & $-6 + 3 \xi_1$ & $3-3\xi_1$ & $0$ & $3\xi_1$ & $-3-3\xi_1$ & $6$ \\
    	\hline
    	\end{tabular}
    	\renewcommand{\arraystretch}{1}
    	\caption{Poles and residues of $\omega_{(\xi_1,1)}$.}
        \label{table:poles_omega_P0}
	\end{table}
    The residues are symmetric around $P_0$ if and only if $\xi_1=\frac{1}{2}$. This symmetry forces a genuine zero at $P_0$, which then implies that $\xi_1=\frac{1}{2}$. The remaining cases follow in a similar fashion. 
\end{proof}


\subsection{Proof of Theorem~\ref{thm:KO_map}(a)}\label{sec:KO_map_a}
The proof of Theorem~\ref{thm:KO_map}(a) depends on certain properties of $d\log\lambda/d\log z$, which is a meromorphic function on $\mathcal{R}$. It also extends to a meromorphic function on $\widehat{\mathcal{R}}$. The poles of $d\log z$ on $\mathcal{R}$ are at $P_\infty$ and $P_0$ with residues $\pm 3$, and the poles of $d\log\lambda$ are at $P_\infty$ and $P_1$, also with residues $\pm 3$, as can be readily verified in local coordinates. The poles and residues repeat on the second copy $\mathcal{R}^\sigma$.

\begin{figure}[t]
	\centering
	\begin{tikzpicture}[scale=0.5](15,10)(0,0)
		\filldraw[gray!50!white] (-6,0) --(6,0) --(6,3) --(-6,3) --(-6,0);
		\filldraw[gray!20!white] (-6,0) --(6,0) --(6,-3) --(-6,-3) --(-6,0);

        \draw (-6.1,3) node[left] {$\frac{1}{2}\tau$};
        \draw (-6.1,-3) node[left] {$-\frac{1}{2}\tau$};
        \draw (6.1,3) node[right] {$1+\frac{1}{2}\tau$};
        \draw (6.1,-3) node[right] {$1-\frac{1}{2}\tau$};

        \draw (-4,1) node {$\mathcal{R}_+$};
        \draw (-4,-1) node {$\mathcal{R}_-$};

		\draw[ultra thick] (-6,0)--(6,0);
		\draw[ultra thick] (-6,3)--(6,3) (-6,-3)--(6,-3);

        \draw[ultra thick] (-6,-3)--(-6,3) (6,-3)--(6,3);

		\filldraw (0,0) circle (4pt);
		\filldraw (4,0) circle (4pt);
		\draw (0,0) node[above] {$P_1$};
		\draw (4,0) node[above] {$P_0$};

        \filldraw (-5.25,3) circle (4pt);
        \filldraw (-2.75,3) circle (4pt);
        \filldraw (1.25,3) circle (4pt);
        \filldraw (2.75,3) circle (4pt);
        \filldraw (-5.25,-3) circle (4pt);
        \filldraw (-2.75,-3) circle (4pt);
        \filldraw (1.25,-3) circle (4pt);
        \filldraw (2.75,-3) circle (4pt);
        \draw (-5.15,3) node[above] {$B_0$};
        \draw (-2.65,3) node[above] {$B_1$};
        \draw (1.35,3) node[above] {$B_2$};
        \draw (2.85,3) node[above] {$B_3$};
        
        \draw (-5.15,-3) node[below] {$B_0$};
        \draw (-2.65,-3) node[below] {$B_1$};
        \draw (1.35,-3) node[below] {$B_2$};
        \draw (2.85,-3) node[below] {$B_3$};

        \draw (1,2) node {$\Im(d\log\lambda/d\log z) > 0$};
        \draw (1,-2) node {$\Im(d\log\lambda/d\log z) < 0$};
	\end{tikzpicture}
    \caption{The complex torus with the sign pattern of $\Im(d\log\lambda/d\log z)$, together with the zeros $B_0$, $B_2$, $P_0$ and poles $B_1$, $B_3$, $P_1$ of $d\log\lambda/d\log z$ as in Lemma~\ref{lem:Im_dlogl_dlogz}.}
    \label{fig:torus_sign_pattern}
\end{figure}
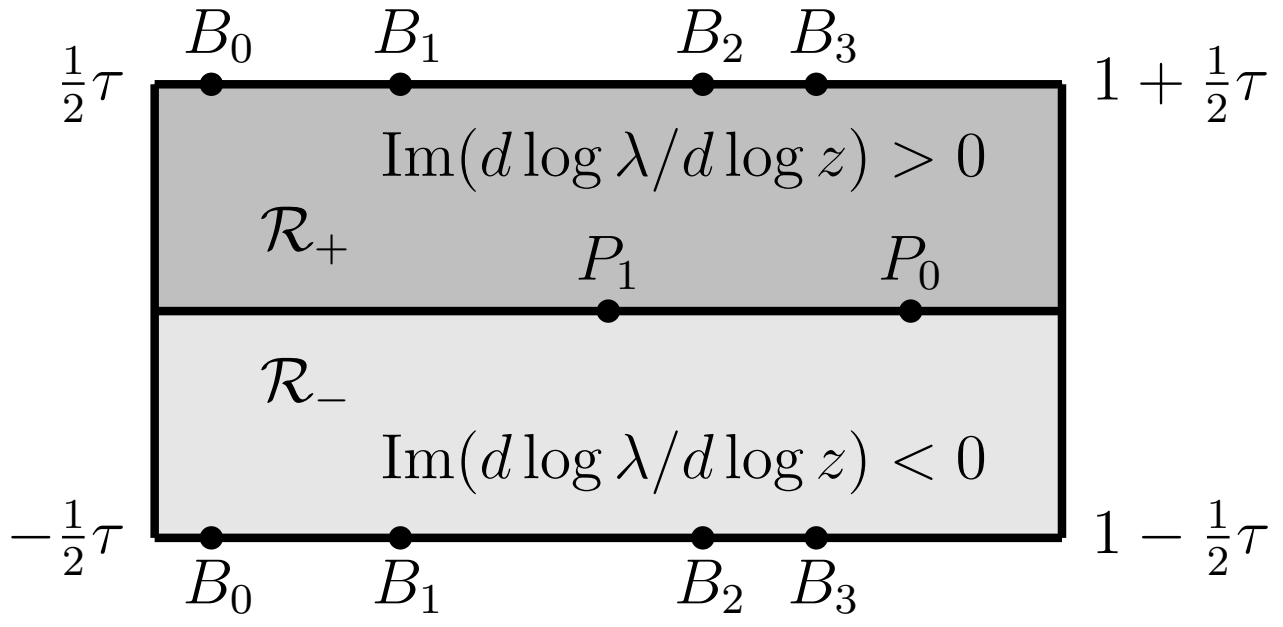

The zeros of $d\log z$ and $d\log\lambda$ are at the points where $z$ and $\lambda$ attain their extremal values on the bounded oval. Recall that $z$ attains its extremal values at the branch points $B_1$ and $B_3$, given explicitly in \eqref{eq:B1_B3}. The points $B_0$ and $B_2$ where $\lambda$ attains its extremal values are obtained by applying the first symmetry in \eqref{eq:sigmas} to $B_1$ and $B_3$, respectively. Thus, we have
\begin{equation}\label{eq:B0_B2}
    B_0 = (-\lambda_{\max},-z_{\max}), \quad B_2 = (-\lambda_{\min},-z_{\min}). 
\end{equation}
Figure~\ref{fig:torus_sign_pattern} shows the points $B_j$, for $j=0,1,2,3$, on the complex torus.

\begin{lemma}\label{lem:Im_dlogl_dlogz}
    The function $d\log\lambda/d\log z$ is meromorphic on $\mathcal{R}$ with simple zeros at $B_0$, $B_2$, and $P_0$, and simple poles at $B_1$, $B_3$, and $P_1$. The function is real-valued on the real ovals, and
	\begin{equation}\label{eq:sign_im_dlogl_dlogz}
        \Im\!\left(\dfrac{d\log\lambda}{d\log z}\right)>0 \quad \text{on } \mathcal{R}_+, \qquad \Im\!\left(\dfrac{d\log\lambda}{d\log z}\right)<0 \quad \text{on } \mathcal{R}_-,
    \end{equation}
    see also Figure~\ref{fig:torus_sign_pattern}.

    Furthermore, $d\log\lambda/d\log z$ strictly increases, when we go around the unbounded oval starting from the pole at $P_1$, and the function strictly decreases on any of the two intervals determined by the poles $B_1$ and $B_3$ on the bounded oval.
\end{lemma}

\begin{proof}
    The statement about zeros and poles of $d\log\lambda/d\log z$ is an immediate consequence of the discussion on the zeros and poles of $d\log\lambda$ and $d\log z$ before the statement of the lemma. The simple pole $P_\infty$ of both $d\log z$ and $d\log\lambda$ disappears in the quotient. It follows that $d\log\lambda/d\log z$ has degree~$3$.

    Moreover, the function is real-valued on the real ovals, and its simple zeros and poles interlace on both ovals. Due to interlacing, every real value is attained on the two parts of the bounded oval between the poles at $B_1$ and $B_3$. Every real value is also attained on the unbounded oval due to the simple pole at $P_1$. Thus, every real value is attained at least three times on the real ovals of $\mathcal{R}$. Since the function has degree~$3$, we have accounted for all real values, and we conclude that
    \[
        \Im\!\left(\dfrac{d\log\lambda}{d\log z}\right) \neq 0 \quad \text{on } \mathcal{R}\setminus(\mathcal{O}_b\cup\mathcal{O}_u).
    \]

    To determine the sign, we note that $d\log\lambda$ has a simple pole at $P_1$ with positive residue and simple pole at $P_0$ with negative residues. Writing $u$ for the natural coordinate on the complex torus, it follows that $d\log\lambda/du$ is negative on the unbounded oval between $P_\infty$ and $P_1$, and positive in between $P_1$ and $P_0$. Similarly, $d\log z/du$ is negative on the unbounded oval between $P_\infty$ and $P_0$. Around $P_1$, there is then a sign change from positive to negative in the quotient $d\log\lambda/d\log z$, which implies that it is strictly increasing on the unbounded oval. By the Cauchy--Riemann equations, the sign of the imaginary part is positive when moving from the unbounded oval into $\mathcal{R}_+$, and negative when moving into $\mathcal{R}_-$. This proves \eqref{eq:sign_im_dlogl_dlogz}.
     
    The above argument already established the strict increase of $d\log\lambda/d\log z$ along the unbounded oval. In a similar way, we find that it strictly decreases on the two intervals between the two poles on the bounded oval.    
\end{proof}

\begin{remark}\label{rem:increase}
    The proof also gives that the strict increase/decrease in Lemma~\ref{lem:Im_dlogl_dlogz} comes from the fact that the derivative $(d\log\lambda/d\log z)'$ in some local coordinate is positive on $\mathcal{O}_u\setminus\{P_1\}$ and negative on $\mathcal{O}_b\setminus\{B_1,B_3\}$. The same strict inequalities hold on the real ovals of $\mathcal{R}^\sigma$.
\end{remark}

Now, we can give the proof of Theorem~\ref{thm:KO_map}(a).

\begin{proof}[Proof of Theorem~\ref{thm:KO_map}(a)]
    Recall that the KO map \eqref{eq:KO_map} maps $(\xi_1,\xi_2)\in\mathfrak{L}$ to the unique $s = \Omega(\xi_1,\xi_2)\in\widehat{\mathcal{R}}_+$ that is a zero of $\omega_{(\xi_1,\xi_2)}$. Dividing the right-hand side of \eqref{eq:omega_xi1_xi2_2} by $d\log z$ and evaluating at $s$, we arrive at
	\begin{equation}\label{eq:saddle_point_equation_1}
		\xi_1\!\left(\dfrac{d\log\lambda}{d\log z}\right)\!(s) - \xi_2 = \left(\dfrac{\omega}{d\log z}\right)\!(s),
	\end{equation}
	Taking the real and imaginary parts of \eqref{eq:saddle_point_equation_1} and noting that $\xi_1$ and $\xi_2$ are real, we find the linear system of equations
    \begin{equation}\label{eq:saddle_point_equation_2}
        \begin{pmatrix} \Re\!\left(\dfrac{d\log\lambda}{d\log z}\right)\!(s) & -1 \vspace{0.1cm} \\ \Im\!\left(\dfrac{d\log\lambda}{d\log z}\right)\!(s) & 0 \end{pmatrix}\!\!\begin{pmatrix} \raisebox{-0.25cm}{$\xi_1$} \\[0.9cm] \raisebox{0.25cm}{$\xi_2$} \end{pmatrix} = \begin{pmatrix} \Re\!\left(\dfrac{\omega}{d\log z}\right)\!(s) \vspace{0.1cm} \\ \Im\!\left(\dfrac{\omega}{d\log z}\right)\!(s) \end{pmatrix}.
    \end{equation}

    For a given $s\in\widehat{\mathcal{R}}_+$, the linear system \eqref{eq:saddle_point_equation_2} has a unique solution, since $\Im(d\log\lambda/d\log z)(s)$ is non-zero by Lemma~\ref{lem:Im_dlogl_dlogz}. Thus, we can uniquely recover $(\xi_1,\xi_2)$ from $s$, and it follows that $\Omega$ is injective.

    To prove surjectivity, we need to see that the map $s \in\widehat{\mathcal{R}}_+ \mapsto (\xi_1(s),\xi_2(s))$ with $(\xi_1(s),\xi_2(s))$ the solution of \eqref{eq:saddle_point_equation_2} maps $\widehat{\mathcal{R}}_+$ into $\mathfrak{L}$. Explicitly, one has
	\begin{align}
		\xi_1(s) &= \frac{\Im\!\left(\dfrac{\omega}{d\log z}\right)\!(s)}{\Im\!\left(\dfrac{d\log\lambda}{d\log z}\right)\!(s)}, \label{eq:KO_map_inv_xi1} \\
		\xi_2(s) &= \Re\!\left(\dfrac{d\log\lambda}{d\log z}\right)\!(s)\frac{\Im\!\left(\dfrac{\omega}{d\log z}\right)\!(s)}{\Im\!\left(\dfrac{d\log\lambda}{d\log z}\right)\!(s)} - \Re\!\left(\frac{\omega}{d\log z}\right)\!(s) \nonumber \\
		&= -\frac{\Im\!\left(\dfrac{\omega}{d\log\lambda}\right)\!(s)}{\Im\!\left(\dfrac{d\log z}{d\log\lambda}\right)\!(s)}. \label{eq:KO_map_inv_xi2}
	\end{align}
    Then, writing $(\xi_1,\xi_2)=(\xi_1(s),\xi_2(s))$
    for a given $s \in \widehat{\mathcal{R}}_+$, we see that $\omega_{(\xi_1,\xi_2)}$ has $s$ as a zero that is away from the real ovals. Moreover, $\overline{s}$ is also a zero, and there are four more zeros on the bounded oval by Theorem~\ref{thm:four_zeros_on_bounded_oval}(a). Since these are all six zeros, there are no zeros on the unbounded oval. Hence, $(\xi_1,\xi_2)\in\mathcal{H}$ by Theorem~\ref{thm:four_zeros_on_bounded_oval}(d) and then $(\xi_1,\xi_2)\in\mathfrak{L}$ by Definition~\ref{def:phases}. 

    The smoothness of $\Omega$ is clear, while the smoothness of its inverse follows from its explicit form \eqref{eq:KO_map_inv_xi1}--\eqref{eq:KO_map_inv_xi2}. 
\end{proof}


\subsection{Proof of Theorem~\ref{thm:KO_map}(b)}\label{sec:KO_map_b}
For the proof of part (b), we introduce
\begin{align}
    \mathfrak{A} &= \mathfrak{A}_\mathrm{out}\cup\mathfrak{A}_\mathrm{in} \nonumber \\
    &= \{(\xi_1,\xi_2)\in\mathbb{R}^2 \mid \omega_{(\xi_1,\xi_2)} \text{ has a zero with higher multiplicity}\}. \label{eq:frakA}
\end{align}
It follows from Theorem~\ref{thm:four_zeros_on_bounded_oval}(c) that for each $(\xi_1,\xi_2)\in\mathfrak{A}$, there exists a unique zero $s_0$ of $\omega_{(\xi_1,\xi_2)}$ with higher multiplicity. This zero necessarily lies on one of the real ovals. Thus, there is a well-defined map
\begin{equation}\label{eq:Omega_frakA}
    \Omega_\mathfrak{A} \colon \mathfrak{A} \to \widehat{\mathcal{O}}_b \cup\widehat{\mathcal{O}}_u, \quad (\xi_1,\xi_2) \mapsto s_0.
\end{equation}
Here, zeros with higher multiplicity are understood in the sense of Remark~\ref{rem:double_zero}. In particular, if $(\xi_1,\xi_2)\in\mathfrak{T}$, then $\omega_{(\xi_1,\xi_2)}$ has a genuine zero (and a double zero in the sense of Remark~\ref{rem:double_zero}) at $s_0=P(\xi_1,\xi_2)$, see Table~\ref{table:zero_at_P}, by Theorem~\ref{thm:four_zeros_on_bounded_oval}(d). Therefore, $(\xi_1,\xi_2)$ is included in $\mathfrak{A}$ and $\Omega_\mathfrak{A}(\xi_1,\xi_2)=s_0$.

The next lemma provides the main step towards the proof of Theorem~\ref{thm:KO_map}.

\begin{lemma}\label{lem:boundary_rough_phase}
    The boundary of the rough phase $\mathfrak{L}$ is given by
    \begin{equation}\label{eq:boundary_rough_phase_2}
        \partial\mathfrak{L}=\mathfrak{A},
    \end{equation}
    and \eqref{eq:Omega_frakA} maps \eqref{eq:boundary_rough_phase_2} bijectively onto $\widehat{\mathcal{O}}_b \cup\widehat{\mathcal{O}}_u$.
\end{lemma}

\begin{proof}
    At the boundary $\partial\mathfrak{L}$, a transition takes place from only zeros on the real ovals to one complex-conjugate pair of zeros away from the real ovals. Such a transition only takes place when at least two zeros coalesce, either on the unbounded oval (corresponding to the boundary between $\mathfrak{L}$ and $\mathfrak{S}$) or on the bounded oval (corresponding to the boundary between $\mathfrak{L}$ and $\mathfrak{G}$). Hence, $\partial\mathfrak{L}\subset\mathfrak{A}$.

    For the converse, we let $(\xi_1,\xi_2)\in\mathfrak{A}$. We first suppose that $s_0 = \Omega_\mathfrak{A}(\xi_1,\xi_2)\notin\widehat{\mathcal{P}}$. From \eqref{eq:omega_xi1_xi2_1}, we then find, after dividing by $d\log z$ and evaluating at $s_0$,
	\begin{align*} 
		\left(\frac{\omega}{d\log z}\right)\!(s_0)-\xi_1\left(\frac{d\log\lambda}{d\log z}\right)\!(s_0)+\xi_2 &= 0, \\
		\left(\frac{\omega}{d\log z}\right)'\!\!(s_0)-\xi_1\left(\frac{d\log\lambda}{d\log z}\right)'\!\!(s_0) &= 0,
	\end{align*}
    where $'$ denotes differentiation in some local coordinate around $s_0$. The derivative $(d\log\lambda/d\log z)'$ is non-zero on the real ovals, due to Lemma~\ref{lem:Im_dlogl_dlogz} and Remark~\ref{rem:increase}. Then, solving for $(\xi_1,\xi_2)$ gives
	\begin{align}
        \xi_1 = \xi_1(s_0) &= \frac{\left(\dfrac{\omega}{d\log z}\right)'\!\!(s_0)}{\left(\dfrac{d\log\lambda}{d\log z}\right)'\!\!(s_0)}, \label{eq:KO_map_boundary_1} \\
        \xi_2 = \xi_2(s_0) &= -\left(\dfrac{\omega}{d\log z}\right)\!(s_0)+\dfrac{\left(\dfrac{\omega}{d\log z}\right)'\!\!(s_0)}{\left(\dfrac{d\log\lambda}{d\log z}\right)'\!\!(s_0)}\left(\dfrac{d\log\lambda}{d\log z}\right)\!(s_0) \nonumber \\
        &= -\frac{\left(\dfrac{\omega}{d\log\lambda}\right)'\!\!(s_0)}{\left(\dfrac{d\log z}{d\log\lambda}\right)'\!\!(s_0)}. \label{eq:KO_map_boundary_2}
	\end{align}
    One can verify that \eqref{eq:KO_map_boundary_1}--\eqref{eq:KO_map_boundary_2} are the limiting expressions of \eqref{eq:KO_map_inv_xi1}--\eqref{eq:KO_map_inv_xi2} as $s \in\widehat{\mathcal{R}}_+$ tends to $s_0\in\bigl(\widehat{\mathcal{O}}_b \cup\widehat{\mathcal{O}}_u\bigr)\!\setminus\!\widehat{\mathcal{P}}$. Since $(\xi_1(s),\xi_2(s))\in\mathfrak{L}$ for $s \in\widehat{\mathcal{R}}_+$ by Theorem~\ref{thm:four_zeros_on_bounded_oval}, it then follows that $(\xi_1,\xi_2)\in\partial\mathfrak{L}$.

    Second, we suppose that $s_0=\Omega_\mathfrak{A}(\xi_1,\xi_2)\in\widehat{\mathcal{P}}$. Then, we can obtain $s_0$ as a limit of $s_0'$ with $s_0'\in\widehat{\mathcal{O}}_u\setminus\widehat{\mathcal{P}}$. For each $s_0'$, the formulas \eqref{eq:KO_map_boundary_1}--\eqref{eq:KO_map_boundary_2} give a unique $(\xi_1',\xi_2')\in\partial\mathfrak{L}$ with $s_0' = \Omega_\mathfrak{A}(\xi_1',\xi_2')$. Letting $s_0' \to s_0 \in\widehat{\mathcal{P}}$, we then find (by a compactness argument and passing to a subsequence) that $(\xi_1',\xi_2')\to(\xi_1,\xi_2)$ and thus $(\xi_1,\xi_2)\in\partial\mathfrak{L}$. This completes the proof that $\partial\mathfrak{L}=\mathfrak{A}$.

    Along the way we also found the inverse of $\Omega_\mathfrak{A}$. It is given by \eqref{eq:KO_map_boundary_1}--\eqref{eq:KO_map_boundary_2} in case $s_0\in\bigl(\widehat{\mathcal{O}}_b \cup\widehat{\mathcal{O}}_u\bigr)\setminus\widehat{\mathcal{P}}$, while for $s_0$, it follows from Table~\ref{table:zero_at_P}. It follows that $\Omega_\mathfrak{A}$ is a bijection, which completes the proof.
\end{proof}

The proof of Theorem~\ref{thm:KO_map}(b) is now straightforward.

\begin{proof}[Proof of Theorem~\ref{thm:KO_map}(b)]
    The extension of the KO map \eqref{eq:KO_map} to the boundary $\partial\mathfrak{L}$ is given by
    \begin{equation}\label{eq:KO_map_ext_2}
        \Omega \colon \overline{\mathfrak{L}} \to \widehat{\mathcal{R}}_+ \cup\widehat{\mathcal{O}}_b \cup\widehat{\mathcal{O}}_u, \quad (\xi_1,\xi_2) \mapsto
        \begin{cases}
            \Omega(\xi_1,\xi_2) & \text{if } (\xi_1,\xi_2)\in\mathfrak{L}, \\
            \Omega_\mathfrak{A}(\xi_1,\xi_2) & \text{if } (\xi_1,\xi_2)\in\partial\mathfrak{L}=\mathfrak{A}.
        \end{cases}
    \end{equation}
    It is immediate that \eqref{eq:KO_map_ext_2} is continuous. The continuity of its inverse is a corollary of the proof of Lemma~\ref{lem:boundary_rough_phase}.

    Finally, recall that the values of \eqref{eq:KO_map_ext_2} at the turning points $\mathfrak{T}$ as in \eqref{eq:turning_points} follow immediately from Table~\ref{table:zero_at_P}.
\end{proof}


\subsection{The KO map and the double amoeba}
In this section, we examine the action of the KO map \eqref{eq:KO_map} more closely.

\begin{figure}[t]
	\centering
	\begin{tikzpicture}[scale=0.5](15,10)(0,0)
		\filldraw[gray!50!white] (-6,0) --(18,0) --(18,3) --(-6,3) --(-6,0);
		\filldraw[gray!20!white] (-6,0) --(18,0) --(18,-3) --(-6,-3) --(-6,0);

		\draw (-6.1,3) node[left] {$\frac{1}{2}\tau$};
		\draw (-6.1,-3) node[left] {$-\frac{1}{2}\tau$};
		\draw (18.1,3) node[right] {$2+\frac{1}{2}\tau$};
		\draw (18.1,-3) node[right] {$2-\frac{1}{2}\tau$};

		\draw[ultra thick] (-6,0)--(18,0);
		\draw[ultra thick] (-6,3)--(18,3) (-6,-3)--(18,-3);

		\draw (-4,1) node {$\mathcal{R}_+$};
		\draw (-4,-1) node {$\mathcal{R}_-$};
		\draw (8,1) node {$\mathcal{R}^\sigma_+$};
		\draw (8,-1) node {$\mathcal{R}^\sigma_-$};

		\filldraw (0,0) circle (4pt);
		\filldraw (12,0) circle (4pt);
		\draw (0,0) node[above] {$P_1$};
		\draw (12,0) node[above] {$P_1^\sigma$};

		\draw[ultra thick] (-6,-3)--(-6,3) (6,-3)--(6,3) (18,-3)--(18,3);

		\filldraw (-2.75,3) circle (4pt);
		\filldraw (2.75,3) circle (4pt);
		\filldraw (-2.75,-3) circle (4pt);
		\filldraw (2.75,-3) circle (4pt);
		\draw (-2.65,3.075) node[above] {$B_1$};
		\draw (2.85,3.075) node[above] {$B_3$}; 
		\draw (-2.65,-3.075) node[below] {$B_1$};
		\draw (2.85,-3.075) node[below] {$B_3$};

		\filldraw (-4,3) circle (4pt);
		\filldraw (0,3) circle (4pt);
		\filldraw (4,3) circle (4pt);
		\filldraw (-4,-3) circle (4pt);
		\filldraw (0,-3) circle (4pt);
		\filldraw (4,-3) circle (4pt);
		\draw (-3.9,3.075) node[above]  {$Z_0$};
		\draw (0.1,3.075) node[above] {$Z_1$};
		\draw (4.1,3.075) node[above] {$Z_2$}; 
		\draw (-3.9,-3.075) node[below]  {$Z_0$};
		\draw (0.1,-3.075) node[below] {$Z_1$};
		\draw (4.1,-3.075) node[below] {$Z_2$};

		\filldraw (9.25,3) circle (4pt);
		\filldraw (14.75,3) circle (4pt);
		\filldraw (9.25,-3) circle (4pt);
		\filldraw (14.75,-3) circle (4pt);
		\draw (9.35,3) node[above] {$B_1^\sigma$};
		\draw (14.85,3) node[above] {$B_3^\sigma$};
		\draw (9.35,-3) node[below] {$B_1^\sigma$};
		\draw (14.85,-3) node[below] {$B_3^\sigma$};

		\filldraw (8,3) circle (4pt);
		\filldraw (12,3) circle (4pt);
		\filldraw (16,3) circle (4pt);
		\filldraw (8,-3) circle (4pt);
		\filldraw (12,-3) circle (4pt);
		\filldraw (16,-3) circle (4pt);
		\draw (8.1,3) node[above] {$Z_0^\sigma$};
		\draw (12.1,3) node[above] {$Z_1^\sigma$};
		\draw (16.1,3) node[above] {$Z_2^\sigma$};
		\draw (8.1,-3) node[below] {$Z_0^\sigma$};
		\draw (12.1,-3) node[below] {$Z_1^\sigma$};
		\draw (16.1,-3) node[below] {$Z_2^\sigma$};

		\draw (1,2) node {$\Im(\omega/d\log z) > 0$};
		\draw (13,2) node {$\Im(\omega/d\log z) < 0$};
		\draw (1,-2) node {$\Im(\omega/d\log z) < 0$};
		\draw (13,-2) node {$\Im(\omega/d\log z) > 0$};

		\draw (-6,-1) node[left] {$\Gamma_0$};
		\draw (6,-1) node[left] {$\Gamma_0^\sigma$};
		\draw (18,-1) node[right] {$\Gamma_0$};
	\end{tikzpicture}
    \caption{The double complex torus with the sign pattern of $\Im(\omega/d\log z)$, together with the zeros $Z_j$, $Z_j^\sigma$, $j=0,1,2$, and poles $B_1, B_3, P_1, B_1^\sigma, B_3^\sigma, P_1^\sigma$ of $\omega/d\log z$ as given in Lemma~\ref{lem:Im_omega_dlogz}.}
	\label{fig:double_torus_sign_pattern}
\end{figure}
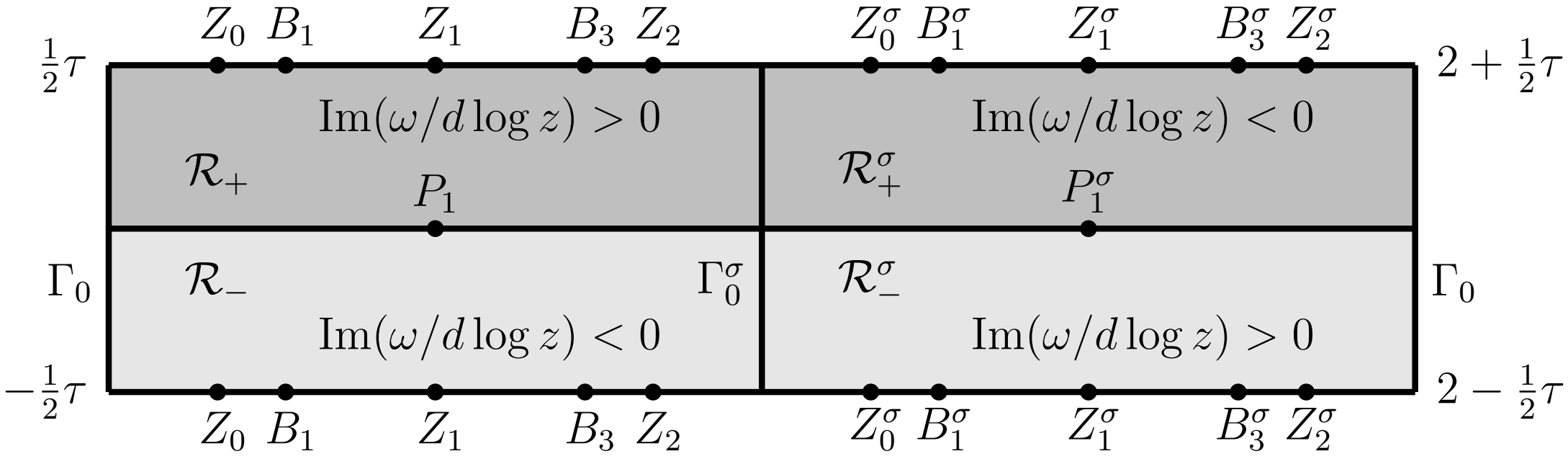

\begin{lemma}\label{lem:Im_omega_dlogz}
    Let $f$ be the function	
	\begin{equation}\label{eq:f_00}
		f = \dfrac{\omega}{d\log z}.
	\end{equation}
    Then, $f$ is meromorphic on $\widehat{\mathcal{R}}$ of degree six, with simple zeros at $Z_j$ and $Z_j^\sigma$ for $j=0,1,2$, and simple poles at $B_1$, $B_3$, $B_1^\sigma$, $B_3^\sigma$, $P_1$, and $P_1^\sigma$. Furthermore, this function is real-valued on the real ovals and on $\Gamma_0\cup\Gamma_0^\sigma$, and
	\[
    	\begin{aligned}
    		&\Im f>0 \quad \text{on } \mathcal{R}_+, \qquad &\Im f<0 \quad \text{on } \mathcal{R}^\sigma_+, \\
    		&\Im f<0 \quad \text{on } \mathcal{R}_-, \qquad &\Im f>0 \quad \text{on } \mathcal{R}_-^\sigma,
    	\end{aligned}
	\]
	see also Figure~\ref{fig:double_torus_sign_pattern}.
\end{lemma}

\begin{proof}
    The statement about zeros and poles of $f$ is an immediate consequence of Proposition~\ref{prop:omega_and_dg_2} and the discussion on the zeros and poles of $d\log z$ before Lemma~\ref{lem:Im_dlogl_dlogz}, see also Figure~\ref{fig:double_torus_sign_pattern} for the zeros and poles on the double complex torus. It follows that $f$ has degree~$6$.
	
    Note that $f = f_{(0,0)}$ with $f_{(\xi_1,\xi_2)}$ as in \eqref{eq:f_xi1_xi2}, and recall that this function is real-valued on the real ovals, which was established in the proof of Theorem~\ref{thm:four_zeros_on_bounded_oval} using the symmetry \eqref{eq:f_xi1_xi2_bar}. By \eqref{eq:omega_sigmas}, we have $(\sigma_2\circ\sigma)_*(\omega) = -\omega$. Combining this with $(\sigma_2\circ\sigma)_*(\log z) = -\log z$ and \eqref{eq:f_xi1_xi2_bar} yields
	\begin{equation}\label{eq:f_00_sigmas_bar}
		f(\sigma_2\circ\sigma(\bar{p})) = \overline{f(p)}, \quad p\in\widehat{\mathcal{R}}.
	\end{equation}
    From the symmetry \eqref{eq:f_00_sigmas_bar}, it follows that $f$ is also real-valued on $\Gamma_0\cup\Gamma_0^\sigma$, as this set is fixed under the map $p \mapsto \sigma_2\circ\sigma(\bar{p})$.
	
	Following the proof of Lemma~\ref{lem:Im_dlogl_dlogz}, one can show that
	\[
        \Im f \neq 0 \quad \text{on } \widehat{\mathcal{R}}\setminus(\widehat{\mathcal{O}}_b\cup\widehat{\mathcal{O}}_u\cup\Gamma_0\cup\Gamma_0^\sigma)
	\]
	by accounting for all real values, which occur precisely six times since $f$ has degree~$6$. The corresponding sign pattern of $\Im f$ can likewise be determined by following the same proof.
\end{proof}

Recall that \eqref{eq:KO_map_inv_xi1} gives $\xi_1$ as a function of $s \in\widehat{\mathcal{R}}_+$. Combining this formula with Lemmas~\ref{lem:Im_dlogl_dlogz} and~\ref{lem:Im_omega_dlogz}, we find that $\xi_1 > 0$ for $s \in\mathcal{R}_+\setminus(\Gamma_0 \cup\Gamma_0^\sigma)$ and $\xi_1 < 0$ for $s \in\mathcal{R}^\sigma_+\setminus(\Gamma_0 \cup\Gamma_0^\sigma)$.

\begin{figure}[t]
	\centering
	\hspace{0.6cm}
	\begin{overpic}[width=0.4\linewidth]{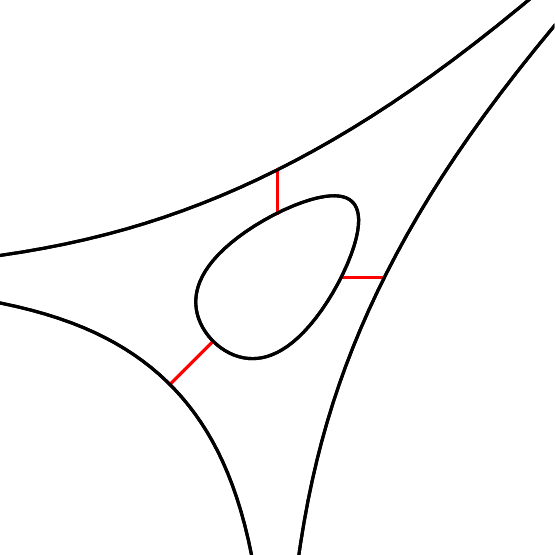}
		\put(-9,47.5){$P_0$}
		\put(46.5,-7){$P_1$}
		\put(99,100){$P_\infty$}
		\put(56.5,83){$\Gamma_0$}
		\put(47.6,66){
		\begin{tikzpicture}
			\draw[-] (0,0) -- (-0.6,-0.9);
		\end{tikzpicture}}
		\put(82,39){$\Gamma_1$}
		\put(63,41.6){
		\begin{tikzpicture}
				\draw[-] (0,0) -- (-0.9,0.5);
		\end{tikzpicture}}
		\put(12,29){$\Gamma_2$}
		\put(17.5,31.3){
			\begin{tikzpicture}
				\draw[-] (0,0) -- (-0.9,-0.2);
		\end{tikzpicture}}
		\put(1,47.5){$0<\xi_1<\xi_2$}
		\put(63,64.5){\rotatebox{45}{$\xi_2<0<\xi_1$}}
		\put(46.5,34){\rotatebox{-90}{$0<\xi_2<\xi_1$}}
	\end{overpic}
	\hfill
	\begin{overpic}[width=0.4\linewidth]{figures/level_set_0.pdf}
		\put(-9.5,47.5){$P_0^\sigma$}
		\put(46,-7){$P_1^\sigma$}
		\put(99,100){$P_\infty^\sigma$}
        \put(56.5,83){$\Gamma_0^\sigma$}
		\put(47.6,66){
		\begin{tikzpicture}
			\draw[-] (0,0) -- (-0.6,-0.9);
		\end{tikzpicture}}
		\put(82,39){$\Gamma_1^\sigma$}
		\put(63,41.6){
		\begin{tikzpicture}
				\draw[-] (0,0) -- (-0.9,0.5);
		\end{tikzpicture}}
		\put(12,29){$\Gamma_2^\sigma$}
		\put(17.5,31.3){
			\begin{tikzpicture}
				\draw[-] (0,0) -- (-0.9,-0.2);
		\end{tikzpicture}}
		\put(1,47.5){$\xi_2<\xi_1<0$}
		\put(63,64.5){\rotatebox{45}{$\xi_1<0<\xi_2$}}
		\put(46.5,34){\rotatebox{-90}{$\xi_1<\xi_2<0$}}
	\end{overpic}
	\hspace{0.6cm}
	\vspace{0.4cm}
	\caption{The double amoeba with the sign pattern of $\xi_1 = \xi_1(s)$ and $\xi_2 = \xi_2(s)$ for $s\in\widehat{\mathcal{R}}_+$, see also Lemma~\ref{lem:xi1_xi2_conditions}. The figure also shows the cycles $\Gamma_j$ and $\Gamma_j^\sigma$ (red) for $j=0,1,2$.}
	\label{fig:xi1_xi2_conditions}
\end{figure}
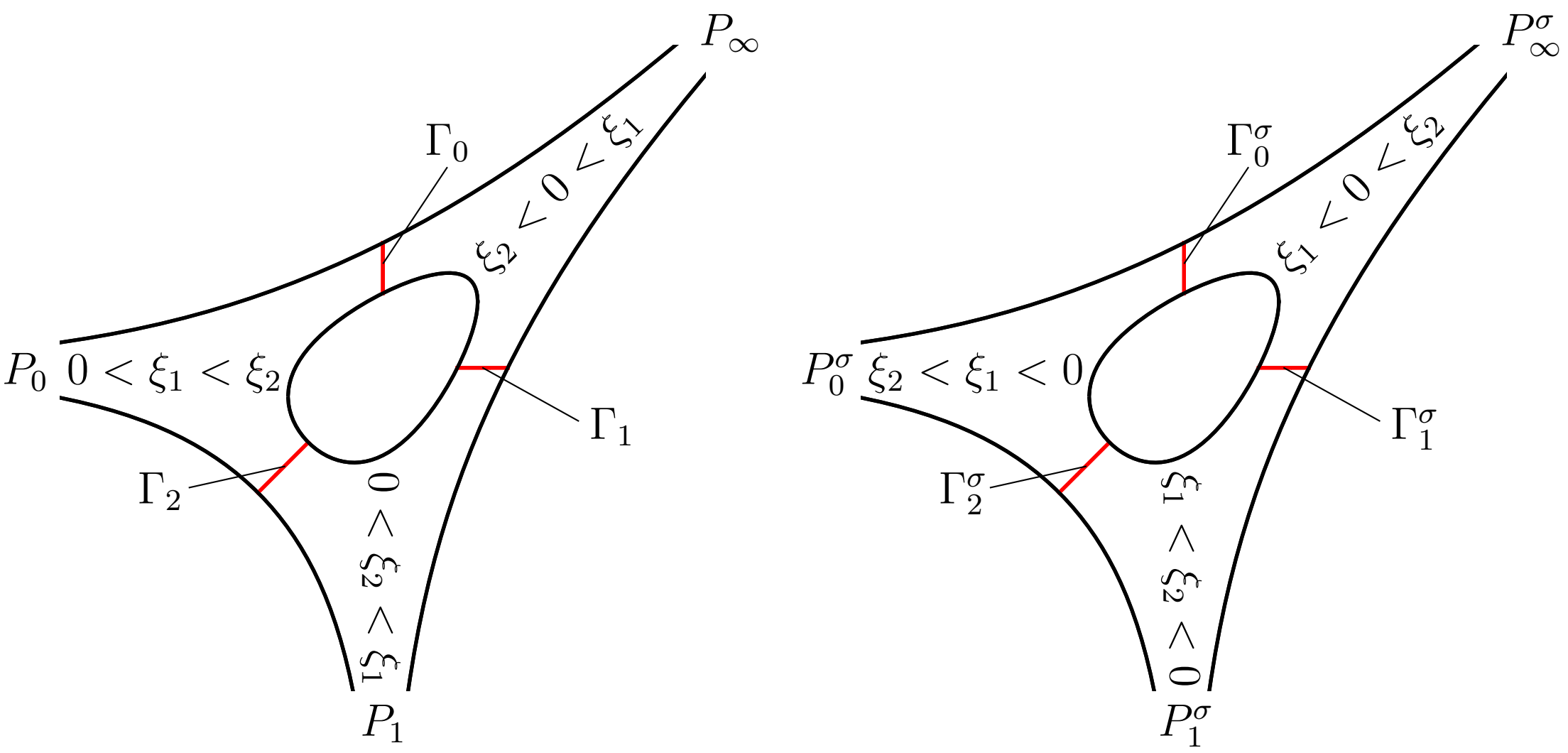

Using the symmetries from Section~\ref{sec:symmetries}, we can also identify the regions determined by the sign of $\xi_2$ and of $\xi_1-\xi_2$. The results are listed in the next lemma, and they are illustrated in Figure~\ref{fig:xi1_xi2_conditions}, where $\widehat{\mathcal{R}}_+$ is visualized as the interior of the double amoeba.

The results of Lemma~\ref{lem:xi1_xi2_conditions} will be used in the following section. 

\begin{lemma}\label{lem:xi1_xi2_conditions}
	Let $(\xi_1,\xi_2)\in\mathfrak{L}$ and write $s = \Omega(\xi_1,\xi_2)$. Then, the following hold:
	\begin{align}\label{eq:xi1_xi2_conditions}
		\begin{aligned}
			s &\in\Gamma_0\cup\Gamma_0^\sigma \iff \xi_1=0, \\
			s &\in\Gamma_1\cup\Gamma_1^\sigma \iff \xi_2=0, \\
			s &\in\Gamma_2\cup\Gamma_2^\sigma \iff \xi_1=\xi_2,
		\end{aligned}
		\qquad
		\begin{aligned}
			s &\in\mathcal{R}_+\setminus(\Gamma_0\cup\Gamma_0^\sigma) \iff \xi_1>0, \\
			\rho^{-1}(s) &\in\mathcal{R}_+\setminus(\Gamma_0\cup\Gamma_0^\sigma) \iff \xi_2>0, \\
			\rho\,(s) &\in\mathcal{R}_+\setminus(\Gamma_0\cup\Gamma_0^\sigma) \iff \xi_1>\xi_2,
		\end{aligned}
	\end{align}
	see also Figure~\ref{fig:xi1_xi2_conditions}.
\end{lemma}

\begin{proof}
    Note that the first line in \eqref{eq:xi1_xi2_conditions} was established in the discussion leading up to this lemma.
	
	It follows from the proof of Lemma~\ref{lem:Phi_xi1_xi2_rho} that
	\begin{equation}\label{eq:omega_xi1_xi2_rho}
		\rho_*(\omega_{(\xi_1,\xi_2)}) = -\omega_{(\xi_1-\xi_2,\xi_1)}, \qquad \rho^{-1}_*(\omega_{(\xi_1,\xi_2)}) = -\omega_{(\xi_2,\xi_2-\xi_1)}.
	\end{equation}
	Since $\rho\,(s),\rho^{-1}(s)\in\widehat{\mathcal{R}}_+$, we have
	\begin{equation}\label{eq:KO_map_rho}
		\rho\,(s) = \Omega(\xi_1-\xi_2,\xi_1), \qquad \rho^{-1}(s) = \Omega(\xi_2,\xi_2-\xi_1).
	\end{equation}
	The second and third lines of \eqref{eq:xi1_xi2_conditions} now follow from the first by the identities in \eqref{eq:KO_map_rho}, recalling that $\rho\,(\Gamma_0\cup\Gamma_0^\sigma) = \Gamma_1\cup\Gamma_1^\sigma$ and $\rho^{-1}(\Gamma_0\cup\Gamma_0^\sigma) = \Gamma_2\cup\Gamma_2^\sigma$.
\end{proof}

\subsection{Proof of Theorem~\ref{thm:six_cusps}}\label{sec:six_cusps}
Part (a) of Theorem~\ref{thm:six_cusps} is contained in Lemma~\ref{lem:boundary_rough_phase}. We begin with the proofs of parts (b) and (c). The proof of part (d) relies on the parallelity property, which we discuss in Lemma~\ref{lem:parallelity_property} below.


\subsubsection{Proof of Theorem~\ref{thm:six_cusps}(b)} 
\begin{proof}
    Let $(\eta_1,\eta_2)\in\partial\mathfrak{L}$ and suppose that $\omega_{(\eta_1,\eta_2)}$ has a double zero at $s_0$. To prove that $\partial\mathfrak{L}$ is smooth at $(\eta_1,\eta_2)$, we use the Weierstrass preparation theorem and the implicit function theorem.
	
    We assume that $s_0$ is not a zero of $d\log z$. Otherwise, $s_0$ is not a zero of $d\log\lambda$, and the argument will be similar. Let $s$ be a local coordinate around $s_0$. Then,
	\begin{equation}\label{eq:F_xi}
		F(s,\xi_1,\xi_2) \coloneqq \frac{\omega(s)-\xi_1d\log\lambda(s)+\xi_2d\log z(s)}{d\log z(s)}. 
	\end{equation}
    is an analytic function near $(s_0, \eta_1,\eta_2)$, that for $(\xi_1,\xi_2)=(\eta_1,\eta_2)$ has a double zero at $s_0$. By the Weierstrass preparation theorem, see e.g.~\cite[Theorem~2.1]{Ebe07}, there is a unique factorization  
	\begin{equation}\label{eq:WPT_1} 
		F = UP \qquad \text{near } (s_0,\eta_1,\eta_2),
    \end{equation}
	with $U$ analytic and non-zero, and 
	\begin{equation}\label{eq:WPT_2} 
		P(s,\xi_1,\xi_2) = (s-s_0)^2+a(\xi_1,\xi_2)(s-s_0)+b(\xi_1,\xi_2),
	\end{equation}
    where $a$ and $b$ are analytic functions that vanish at $(\eta_1,\eta_2)$. Because of the special form \eqref{eq:F_xi} one has $\partial F/\partial\xi_2 = 1$. Then, from \eqref{eq:WPT_1}--\eqref{eq:WPT_2}, one finds
	\[
        1 = \frac{\partial F}{\partial\xi_2}(s_0, \eta_1, \eta_2) = u(s_0, \eta_1,\eta_2)\frac{\partial b}{\partial\xi_2}(\eta_1,\eta_2),
    \]
    and thus $\partial b/\partial\xi_2\neq0$ at $(\eta_1,\eta_2)$. The discriminant of \eqref{eq:WPT_2} is $\Delta = a^2 - 4b$ and it satisfies
	\[
        \frac{\partial\Delta}{\partial\xi_2}(\eta_1,\eta_2) = -4 \frac{\partial b}{\partial\xi_2}(\eta_1,\eta_2) \neq 0,
    \]
    as $a^2$ vanishes quadratically at $(\eta_1,\eta_2)$. Then, by the analytic version of the implicit function theorem, there is a unique analytic function $\varphi$ defined near $\eta_1$ with $\varphi(\eta_1) = \eta_2$ and  
	\begin{equation}\label{eq:WPT_3}
		\Delta(\xi_1,\varphi(\xi_1)) = 0.
	\end{equation}
    This implies that $P$ and thus $F$ has a double zero for the parameters $\xi_1$ and $\xi_2 = \varphi(\xi_1)$. From conjugation symmetry, it will follow that $\xi_2$ is real for real $\xi_1$ near $\eta_1$. Hence, $(\xi_1,\varphi(\xi_1)) \in \partial\mathfrak{L}$ by Lemma~\ref{lem:boundary_rough_phase} if $\xi_1$ is real and close enough to $\eta_1$. This gives a smooth parametrization of $\partial\mathfrak{L}$ near $(\eta_1,\eta_2)$.	
\end{proof}


\subsubsection{Proof of Theorem~\ref{thm:six_cusps}(c)}
\begin{proof}
    Suppose $\omega_{(\eta_1,\eta_2)}$ has a triple zero at $s_0$. As in the proof of part (b), we consider the function $F$ from \eqref{eq:F_xi} which now has a triple zero at $s_0$ for $(\xi_1,\xi_2) = (\eta_1,\eta_2)$. In this case, the Weierstrass preparation theorem \cite[Theorem~2.1]{Ebe07} tells us that \eqref{eq:WPT_1} holds with $U$ analytic and non-zero, and a cubic
	\begin{equation}\label{eq:WPT_4}
		P(s,\xi_1,\xi_2) = (s-s_0)^3 + a(\xi_1,\xi_2)(s-s_0)^2 + b(\xi_1,\xi_2)(s-s_0) + c(\xi_1,\xi_2)
    \end{equation}
	with $a,b$, and $c$ analytic functions that vanish at $(\eta_1,\eta_2)$. We have
	\begin{align*} 
        \frac{\partial F}{\partial\xi_j} &= \frac{\partial U}{\partial\xi_j} P + U \frac{\partial P}{\partial\xi_j}, \\
        \frac{\partial^2 F}{\partial s \partial\xi_j} &= \frac{\partial^2 U}{\partial s \partial\xi_j} P + \frac{\partial U}{\partial\xi_j} \frac{\partial P}{\partial s} + \frac{\partial U}{\partial s} \frac{\partial P}{\partial\xi_j} + U \frac{\partial^2 P}{\partial s \partial\xi_j}
    \end{align*}
    for $j=1,2$. Evaluating these at $(s_0, \eta_1, \eta_2)$ and taking note of \eqref{eq:WPT_4}, we obtain
	\[
	   \frac{\partial F}{\partial\xi_j} = U \frac{\partial c}{\partial\xi_j}, \quad \frac{\partial^2 F}{\partial s \partial\xi_j}  = \frac{\partial U}{\partial s} \frac{\partial c}{\partial\xi_j} + U \frac{\partial b}{\partial\xi_j} \qquad \text{at } (s_0, \eta_1,\eta_2). 
	\]
	In particular,
	\[
        U\det\begin{pmatrix} \dfrac{\partial b}{\partial\xi_1} & \dfrac{\partial b}{\partial\xi_2} \\[10pt] \dfrac{\partial c}{\partial\xi_1} & \dfrac{\partial c}{\partial\xi_2}\end{pmatrix} = \det \begin{pmatrix} \dfrac{\partial^2 F}{\partial s \partial\xi_1} & \dfrac{\partial^2 F}{\partial s \partial\xi_2} \\[10pt] \dfrac{\partial F}{\partial\xi_1} & \dfrac{\partial F}{\partial\xi_2}\end{pmatrix}\qquad \text{at } (s_0,\eta_1,\eta_2).
    \]
    From the special form \eqref{eq:WPT_1} of $F$, we obtain
    \[
        \frac{\partial^2 F}{\partial s \partial\xi_1} = -\frac{d}{ds}\!\left(\frac{d\log\lambda}{d\log z}\right), \quad \frac{\partial F}{\partial\xi_2} = 1, \quad \frac{\partial^2 F}{\partial s \partial\xi_2} = 0,
    \]
    and therefore
	\[
        U\det\begin{pmatrix} \dfrac{\partial b}{\partial\xi_1} & \dfrac{\partial b}{\partial\xi_2} \\[10pt] \dfrac{\partial c}{\partial\xi_1} & \dfrac{\partial c}{\partial\xi_2} \end{pmatrix} = -\frac{d}{ds}\!\left(\frac{d\log\lambda}{d\log z} \right) \qquad \text{at } (s_0,\eta_1,\eta_2).
    \]
    The right-hand side is non-zero by Lemma~\ref{lem:Im_dlogl_dlogz} and Remark~\ref{rem:increase}. We conclude that
	\begin{equation}\label{eq:WPT_5} 
        \det\begin{pmatrix} \dfrac{\partial b}{\partial\xi_1} & \dfrac{\partial b}{\partial\xi_2} \\[10pt] \dfrac{\partial c}{\partial\xi_1} & \dfrac{\partial c}{\partial\xi_2} \end{pmatrix} \neq 0 \quad \text{at } (\eta_1,\eta_2).
    \end{equation}

    The change of variables $s = s_0 + z - \frac{1}{3}a$ transforms \eqref{eq:WPT_4} into
	\begin{equation}\label{eq:WPT_6} 
		P(s,\xi_1,\xi_2) = z^3+u(\xi_1,\xi_2)z+v(\xi_1,\xi_2),
    \end{equation}
    with $u=-\frac{1}{3}a^2+b$ and $v=\frac{2}{27}a^3-\frac{1}{3}ab+c$. Then, $\partial u/\partial\xi_j = \partial b/\partial\xi_j$ and $\partial v/\partial\xi_j = \partial c/\partial\xi_j$ at $(\eta_1,\eta_2)$, and from \eqref{eq:WPT_5}, we get
	\begin{equation}\label{eq:WPT_7} 
        \det\begin{pmatrix}\dfrac{\partial u}{\partial\xi_1} & \dfrac{\partial u}{\partial\xi_2} \\[10pt] \dfrac{\partial v}{\partial\xi_1} & \dfrac{\partial v}{\partial\xi_2} \end{pmatrix} \neq 0 \quad \text{at } (\eta_1,\eta_2).
    \end{equation}
    The inverse function theorem then tells us that, locally, the map $(\xi_1,\xi_2) \mapsto (u,v)$ has an analytic inverse $(u,v) \mapsto (\xi_1,\xi_2)$, while mapping $(0,0)$ to $(\eta_1,\eta_2)$.

    The form \eqref{eq:WPT_6} is known to be the generic form for the unfolding of a cusp singularity, see e.g.~\cite[Section~2.4]{Mon21}, where this is explained in detail. The cubic \eqref{eq:WPT_6} has a multiple root if and only if $(u,v)$ satisfies the equation $\bigl(\frac{u}{3}\bigr)^3 + \bigl(\frac{v}{2}\bigr)^2 = 0$, which is a cusp in the $u$-$v$ plane. The inverse map $(u,v) \mapsto (\xi_1,\xi_2)$ maps it to a cusp in $\xi_1$-$\xi_2$ plane. This proves part (c).
\end{proof}


\subsubsection{Parallelity property}
To prepare for the proof of part (d) of Theorem~\ref{thm:six_cusps} we need a lemma. The statement in Lemma~\ref{lem:parallelity_property}(c) is called the parallelity property.

\begin{lemma}\label{lem:parallelity_property}
    Let $(\eta_1,\eta_2)\in\partial\mathfrak{L}$ and write $s_0 = \Omega(\eta_1,\eta_2)$.
	\begin{enumerate}
    \item[\rm{(a)}] For $(\xi_1,\xi_2)\in\mathbb{R}^2$, $\omega_{(\xi_1,\xi_2)}$ has a zero at $s_0$ if and only if $(\xi_1,\xi_2)$ is on the line $L$ through $(\eta_1,\eta_2)$ with slope $(d\log\lambda/d\log z)(s_0)$.
       
    \item[\rm{(b)}] If $\partial\mathfrak{L}$ is smooth at $(\eta_1,\eta_2)$, then the line $L$ from part (a) is tangent to $\partial\mathfrak{L}$ at $(\eta_1,\eta_2)$, and $\partial\mathfrak{L}$ is locally convex near $(\eta_1,\eta_2)$.
     
    \item[\rm{(c)}] If $(\eta_1,\eta_2)\in\partial\mathfrak{L}\setminus\mathfrak{T}$ (so that $s_0\notin\widehat{\mathcal{P}}$), then the tangent line to $\partial\mathcal{A}$ at $\Log s_0 \in\partial\mathcal{A}$ has the same slope as $L$. That is, the tangent line to $\partial\mathcal{A}$ at $\Log s_0$ is parallel to the tangent line at $\partial\mathfrak{L}$ at $(\eta_1,\eta_2)$.
     
    \item[\rm{(d)}] If $(\eta_1,\eta_2)\in\mathfrak{T}$, then the tangent line to $\partial\mathfrak{L}$ at $(\eta_1,\eta_2)$ is also tangent to the boundary $\partial\mathcal{H}$ of the hexagon. The tangent lines are given explicitly in Table~\ref{table:tangent_lines}.
	\end{enumerate}
	\begin{table}[ht]
		\centering
		\renewcommand{\arraystretch}{1.25}
		\setlength{\tabcolsep}{5pt}
		\begin{tabular}{c | c c c c c c}
		\hline
		$(\eta_1,\eta_2)$ & $(\frac{1}{2},-\frac{1}{2})$ & $(1,\frac{1}{2})$ & $(\frac{1}{2},1)$ & $(-\frac{1}{2}, \frac{1}{2})$ & $(-1,-\frac{1}{2})$ & $(-\frac{1}{2},-1)$ \\ 
		\hline
		tangent & $\eta_2=\eta_1-1$ & $\eta_1=1$ & $\eta_2=1$ & $\eta_1=\eta_2+1$ & $\eta_1=-1$ & $\eta_2=-1$ \\
		\hline
		\end{tabular}
		\renewcommand{\arraystretch}{1}
        \caption{The tangent lines to $\partial\mathfrak{L}$ at the turning points $(\eta_1,\eta_2)\in\mathfrak{T}$.}
		\label{table:tangent_lines}
	\end{table} 
\end{lemma}

\begin{proof}
    (a) Note that $s_0$ is a zero of $\omega_{(\eta_1,\eta_2)} = \omega-\eta_1d\log\lambda+\eta_2d\log z$. Then, $s_0$ is a zero of $\omega_{(\xi_1,\xi_2)}$ if and only if $s_0$ is also a zero of
    \[
        \omega_{(\xi_1,\xi_2)} - \omega_{(\eta_1,\eta_2)} = -(\xi_1-\eta_1) d\log\lambda + (\xi_2-\eta_2)d\log z.
    \]
    This occurs if and only if
    \[
        \frac{\xi_2-\eta_2}{\xi_1-\eta_1} = \left(\frac{d\log\lambda}{d\log z}\right)\!(s_0),
    \]
    with appropriate modifications if $d\log\lambda/d\log z$ has a pole at $s_0$. This proves part~(a). 
     
    \medskip

    (b) From parts (b) and (c) of Theorem~\ref{thm:six_cusps}, that we already proved, we conclude that $\partial\mathfrak{L}$ is smooth at $(\eta_1,\eta_2)$ if and only if $\omega_{(\eta_1,\eta_2)}$ has a double zero at $s_0$. We assume that this is the case.

    We have that $s_0\in\widehat{\mathcal{O}}_b \cup\widehat{\mathcal{O}}_u$. First, we suppose that $s_0\in\widehat{\mathcal{O}}_u$, so that $(\eta_1,\eta_2) \in \mathfrak{A}_\mathrm{out}$. Then, for any $(\xi_1,\xi_2)\in L$, we have by part (a) that $\omega_{(\xi_1,\xi_2)}$ has the zero $s_0$ on the unbounded oval. Thus, $(\xi_1,\xi_2)$ cannot belong to the rough region $\mathfrak{L}$, and we conclude
    \begin{equation}\label{eq:L_tangent}
        L \cap \mathfrak{L} = \emptyset.
    \end{equation}

    The second case is that $s_0 \in\widehat{\mathcal{O}}_b$, and $(\eta_1,\eta_2)$ is on the inner arctic curve. We are in a situation where $s_0$ is a double zero of  $\omega_{(\eta_1,\eta_2)}$. Due to Theorem~\ref{thm:four_zeros_on_bounded_oval}(a,c), there are four additional simple zeros on the bounded oval $\widehat{\mathcal{O}}_b$. This is an open condition, and therefore these four zeros remain distinct and simple if instead of $(\eta_1,\eta_2)$, we take $(\xi_1,\xi_2)\in L$ in a small neighborhood $U$ around $(\eta_1,\eta_2)$. Note that we always keep a zero at $s_0$. Thus, for $(\xi_1,\xi_2)\in L \cap U$, there are at least five zeros on the bounded oval, and $(\xi_1,\xi_2)$ is not in the rough phase  $\mathfrak{L}$. Hence,
    \begin{equation}\label{eq:LU_tangent}
        (L\cap U)\cap\mathfrak{L} = \emptyset
    \end{equation}
    for some neighborhood $U$ of $(\eta_1,\eta_2)$. 

    Since $\partial\mathfrak{L}$ is smooth at $(\eta_1,\eta_2)$, any line $L$ through $(\eta_1,\eta_2)$ that satisfies \eqref{eq:L_tangent} or \eqref{eq:LU_tangent} has to be the tangent line to $\partial\mathfrak{L}$. In addition, it follows that $\partial\mathfrak{L}$ is locally convex near $(\eta_1,\eta_2)$.

    \medskip

    (c) Writing $s_0 = (z,\lambda)$, we have $z,\lambda\in\mathbb{R}\setminus\{0\}$ because $s_0\notin\widehat{\mathcal{P}}$. It follows that $\Log s_0 \in\partial\mathcal{A}$. We introduce $s_0(t) = (z(t),\lambda(t))$ for $t \in (-\varepsilon, \varepsilon)$, with $s_0(0)=s_0$, such that
    \[
        \Log s_0(t) = \bigl(\log|z(t)|,\log|\lambda(t)|\bigr), \quad t\in(-\varepsilon,\varepsilon),
    \]
    is a smooth parametrization of $\partial\mathcal{A}$ near $\Log s_0$. Then,
    \[
        \frac{d}{dt}\Log s_0(t) = \left(\frac{d}{dt}\log|z(t)|,\frac{d}{dt}\log|\lambda(t)|\right)
    \]
    is a tangent vector to $\partial\mathcal{A}$ at $s_0(t)$. The slope of the tangent line is given by
    \[
        \frac{\frac{d}{dt}\log|\lambda(t)|}{\frac{d}{dt}\log|z(t)|}.
    \]
    Taking $t=0$, this slope equals $(d\log\lambda/d\log z)(s_0)$, which coincides with the slope of $L$ by part (a).

    \medskip

    (d) Using \eqref{eq:spectral_curve}, we evaluate the slope $(d\log\lambda/d\log z)(P_\infty)=1$ at $P_\infty$, which results in the tangent line as given in Table~\ref{table:tangent_lines}. By Lemma~\ref{lem:Im_dlogl_dlogz}, $d\log\lambda/d\log z$ has a simple zero at $P_0$ and a simple pole at $P_1$. It follows that the respective tangent lines at these points are also as given in Table~\ref{table:tangent_lines}. The same applies to the tangent lines at the remaining points $P_\infty^\sigma$, $P_1^\sigma$, and $P_0^\sigma$.
\end{proof}


\subsubsection{Proof of Theorem~\ref{thm:six_cusps}(d)}

\begin{proof} 
    The homeomorphism $\Log\circ\,\Omega$ maps the boundary of the smooth phase $\partial\mathfrak{G}$ to the bounded component of the boundary of the double amoeba, which we denote by $(\partial\widehat{\mathcal{A}})_b$. This map reverses the orientation, as can be seen from e.g.~Table~\ref{table:zero_at_P}. The top row in Table~\ref{table:zero_at_P} gives the turning points and they are listed in counterclockwise orientation if we go around the outer arctic curve. The second row gives the corresponding points in $\widehat{\mathcal{P}}$, that when mapped to the double amoeba, follow the clockwise orientation along the outer boundary of the amoeba. 

    We follow the inverse map $(\Log\circ\,\Omega)^{-1}$ and consider its effect on the orientation of outward normal vectors to  $(\partial\widehat{\mathcal{A}})_b$. Since it is a combination of two smooth closed curves, there is an increase in the angle by $4\pi$ if we go once along $(\partial\widehat{\mathcal{A}})_b$ in the positive direction. By the parallelity property, i.e., Lemma~\ref{lem:parallelity_property}(c), this implies that there is a combined decrease by $4\pi$ in the angle of the outward normal vector along the smooth parts of $\mathfrak{A}_\mathrm{in}=\partial\mathfrak{G}$. It is a decrease because the orientation is reversed by the inverse map.

    At a cusp point of $\partial\mathfrak{G}$, the angle increases by $\pi$. If there are $c$ cusp points, then the total increase of the angle along the full boundary of $\mathfrak{G}$ is equal to $c\pi-4\pi$. Since $\partial\mathfrak{G}$ is a simple closed curve, the total increase is $2\pi$, and therefore $c = 6$, as claimed in part (d).
\end{proof}


\section{Proof of Theorem~\ref{thm:local_scaling_limit}}\label{sec:local_scaling_limit}
The main difficulty in the proof of Theorem~\ref{thm:local_scaling_limit} lies in the behavior of the double contour integral in the formula~\eqref{eq:KN_v1_v2} for the correlation kernel $K_N$ as $N\to\infty$. We denote this double contour integral by $\mathsf{K}_N(v_{1,N},v_{2,N})$, with $v_{1,N}$ and $v_{2,N}$ as in \eqref{eq:v1N}--\eqref{eq:v2N}, and it is given by
\begin{multline}\label{eq:KN_matrix}
	\mathsf{K}_N(v_{1,N},v_{2,N}) = \frac{1}{(2\pi i)^2}\oint_\mathbb{T}\oint_\mathbb{T}T_{0,j_2}^{-1}(z_2) \frac{W(z_2)^{N(1-\xi_{1,N})-x_2}}{z_2^{N(1-\xi_{2,N})-y_2}} \\
	\times R_N(z_2,z_1)\frac{W(z_1)^{N(1+\xi_{1,N})+x_1}}{z_1^{N(1+\xi_{2,N})+y_1}}T_{0,j_1}(z_1)\frac{dz_1dz_2}{z_1}.
\end{multline}
Recall that $R_N$ is given in terms of the solution $\mathbf{Y}_N$ of the RH problem as in \eqref{eq:RN_formula}, and note that $N\xi_{1,N},N\xi_{2,N},x_1,x_2,y_1,y_2\in\mathbb{Z}$ and $j_1,j_2\in\{0,1,2\}$.

Instead of using the asymptotics of $\mathbf{Y}_N$ directly, we apply the first two transformation $\mathbf{Y}_N \mapsto \mathbf{X}_N \mapsto \mathbf{T}_N$ from \cite{Kui25} to the double contour integral.

More precisely, we apply $\mathbf{Y}_N \mapsto \mathbf{X}_N$ and then move the contours of integration from the unit circle $\mathbb{T}$ in the complex plane to the contours $z^{-1}(\mathbb{T}) = \mathbb{T}_1\cup\mathbb{T}_2\cup\mathbb{T}_3$ on the Riemann surface. One of the integrals is then deformed to a certain contour $\gamma_B$ on the Riemann surface. The other contour integral splits as a sum of two integrals, and together these are deformed to contours that we call $\gamma_A$ and $\gamma_C$ in Section~\ref{sec:deform_p2}. The result is stated in Lemma~\ref{lem:KN_split}, which gives \eqref{eq:KN_matrix} as a sum of three terms.

In particular, the term \eqref{eq:KN_sN} will be the leading term as $N\to\infty$. The other two terms will tend to zero as $N\to\infty$. In order to prove this, the contours $\gamma_A$ and $\gamma_B$ have to pass through  $s_N = \Omega(\xi_{1,N},\xi_{2,N})$ and $\overline{s}_N$, which are the saddle points of the action function $\Phi_{(\xi_{1,N},\xi_{2,N})}$ from Definition~\ref{def:Phi_xi1_xi2}. Locally, near the saddle points, the contour $\gamma_B$ will follow the path of steepest descent, and $\gamma_A$ follows the path of steepest ascent.

To show that the appropriate contour deformations can indeed be carried out, we first study the critical level sets of $\Re\Phi_{(\xi_1,\xi_2)}$ passing through the saddle points. We will use the double amoeba to visualize these level sets.


\subsection{Critical level sets}\label{sec:critical_level_sets}
The critical level set of $\Re\Phi_{(\xi_1,\xi_2)}$ is defined as follows.

\begin{definition}\label{def:critical_level_sets}
    For $(\xi_1,\xi_2)\in\mathfrak{L}$, set $C_{(\xi_1,\xi_2)} = \Re\Phi_{(\xi_1,\xi_2)}(\Omega(\xi_1,\xi_2))$ and define the critical level set as
    \begin{equation}\label{eq:critical_level_set_on_Rhat}
        \widehat{\mathcal{N}}_{(\xi_1,\xi_2)} = \{p\in\widehat{\mathcal{R}} \mid \Re\Phi_{(\xi_1,\xi_2)}(p) = C_{(\xi_1,\xi_2)}\}.
    \end{equation}
    The two regions that are separated by $\widehat{\mathcal{N}}_{(\xi_1,\xi_2)}$ are denoted by
	\begin{equation}\label{eq:critical_level_set_regions_on_Rhat}
		\begin{aligned}
        	\widehat{\mathcal{N}}^>_{(\xi_1,\xi_2)} &= \{p\in\widehat{\mathcal{R}} \mid \Re\Phi_{(\xi_1,\xi_2)}(p) > C_{(\xi_1,\xi_2)}\}, \\
        	\widehat{\mathcal{N}}^<_{(\xi_1,\xi_2)} &= \{p\in\widehat{\mathcal{R}} \mid \Re\Phi_{(\xi_1,\xi_2)}(p) < C_{(\xi_1,\xi_2)}\}. 
		\end{aligned}
	\end{equation}
\end{definition}

The restrictions to $\mathcal{R}$ are denoted without the use of the decoration $\,\widehat{\cdot}\,$. For example,
\begin{equation}\label{eq:critical_level_set_on_R}
    \mathcal{N}_{(\xi_1,\xi_2)} = \{p\in\mathcal{R} \mid \Re\Phi_{(\xi_1,\xi_2)}(p) = C_{(\xi_1,\xi_2)}\}.
\end{equation}
This is well-defined since $\Re\Phi_{(\xi_1,\xi_2)}$ is a continuous (but no longer harmonic) function on $\mathcal{R}$.

Our goal is to understand the global structure of the critical level set $\widehat{\mathcal{N}}_{(\xi_1,\xi_2)}$ for $(\xi_1,\xi_2)\in\mathfrak{L}$ (rough phase). We first restrict ourselves to the case $0\leq\xi_2\leq\xi_1$. Then, $\xi_1 > 0$ and by \eqref{eq:xi1_xi2_conditions}, the saddle point $s = \Omega(\xi_1,\xi_2)$ is in $\mathcal{R}$.

In case of strict inequalities $0 < \xi_2 < \xi_1$, the saddle point is in the connected component of $\mathcal{R}_+ \setminus(\Gamma_1\cup\Gamma_2)$ that contains $P_1$. If $\xi_2 =0$, then $s \in \Gamma_1$, and if $\xi_1 = \xi_2$, then $s \in \Gamma_2$. This all follows from \eqref{eq:xi1_xi2_conditions}, but it is easily seen in Figure~\ref{fig:xi1_xi2_conditions}, where the amoeba on the left has a region with $0 < \xi_2 < \xi_1$.

\begin{figure}[t]
    \hspace{0.3cm}
	\begin{overpic}[width=0.3\linewidth]{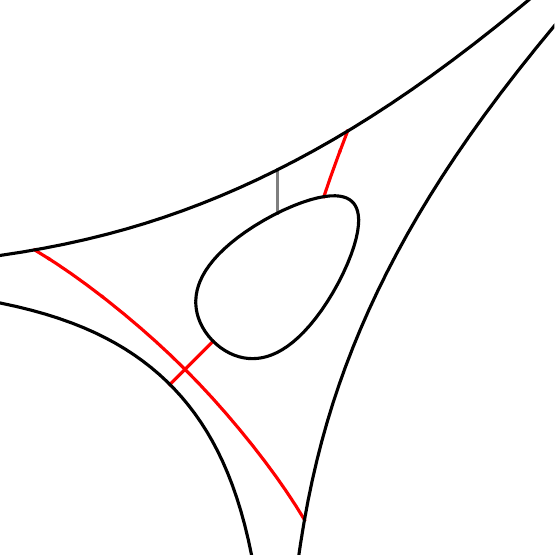}	
		\put(-11,47){$P_0$}
		\put(46,-8){$P_1$}
	\end{overpic}
    \hfill
	\begin{overpic}[width=0.3\linewidth]{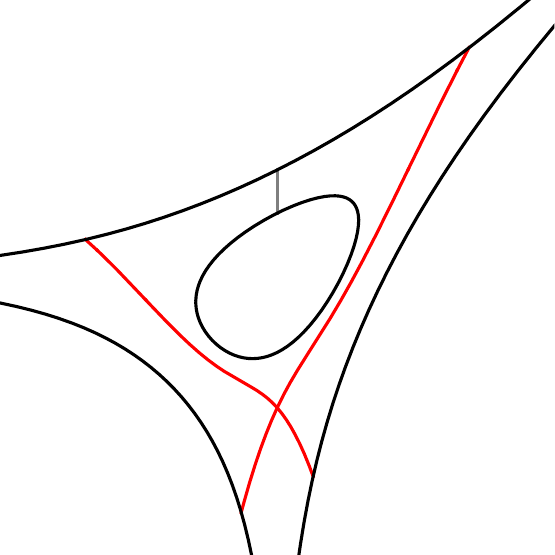}	
        \put(-10,99){$P_\infty$}
		\put(-11,47){$P_0$}
		\put(46,-8){$P_1$}
	\end{overpic}
    \hfill
	\begin{overpic}[width=0.3\linewidth]{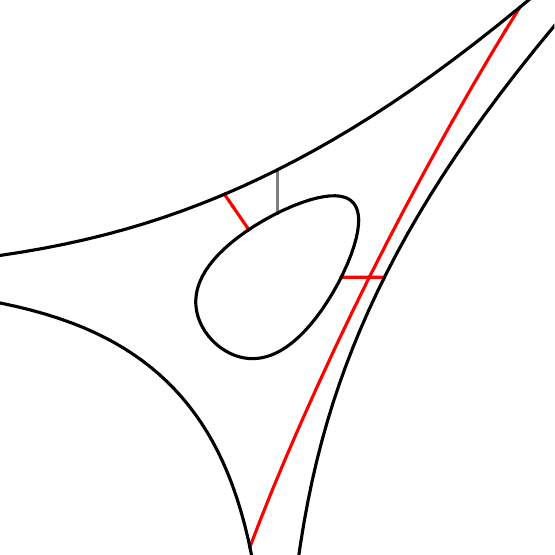}	
        \put(-10,99){$P_\infty$}
		\put(-11,47){$P_0$}
		\put(46,-8){$P_1$}
		\put(98,99){$P_\infty$}
	\end{overpic}
    \hspace{0.2cm}
    \vspace{0.25cm}
	\caption{Level set $\mathcal{N}_{(\xi_1,\xi_2)}$ (red) for various $(\xi_1,\xi_2)\in\mathfrak{L}$ with $0\leq\xi_1\leq\xi_2$.}
    \label{fig:critical_level_sets_1}
\end{figure}
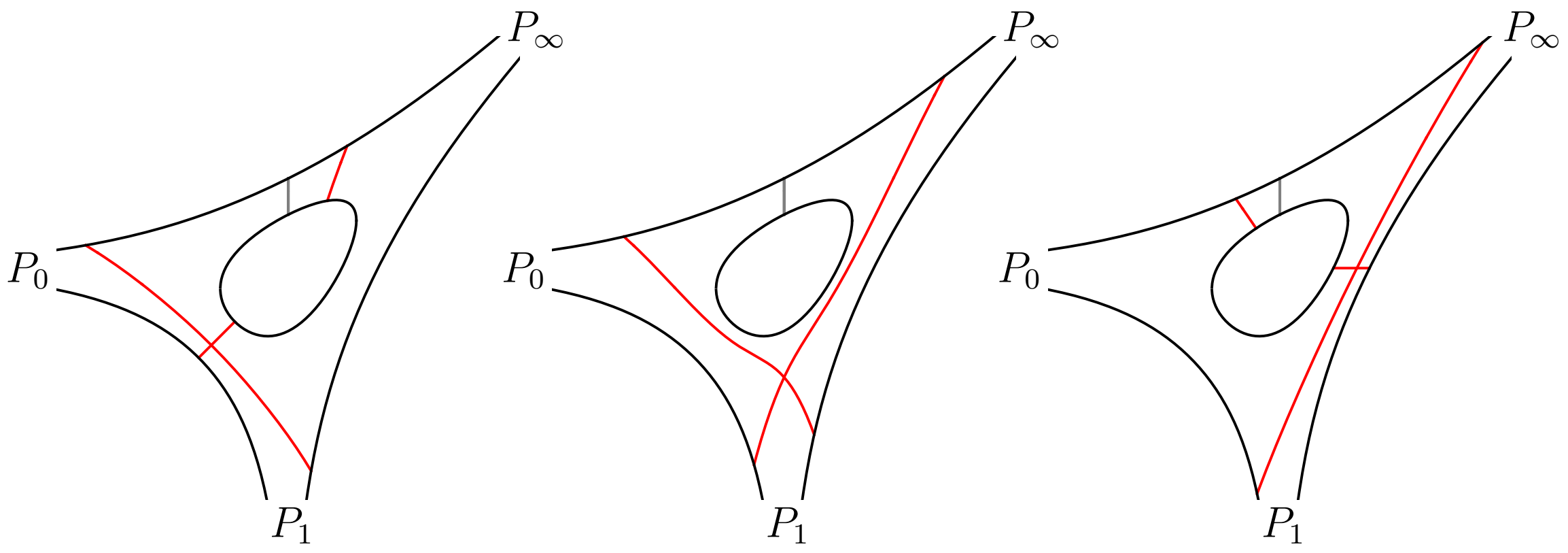

We will prove that, when $0\leq\xi_2\leq\xi_1$, the level set $\mathcal{N}_{(\xi_1,\xi_2)}$ admits one of three topologically distinct configurations. These three distinct configurations are shown in Figure~\ref{fig:critical_level_sets_1}. Throughout this section, we visualize the Riemann surface $\mathcal{R}$ via its amoeba.

Locally, near the saddle point $s$, the level set $\widehat{\mathcal{N}}_{(\xi_1,\xi_2)}$ looks like two crossing lines. Thus, there exist four level curves that emanate from $s$. We show that in this case all these curves stay within $\mathcal{R}$, and that three or four of them pass through different components of the unbounded oval $\widehat{\mathcal{O}}_u$. Here, a component of $\widehat{\mathcal{O}}_u$ refers to one of the six connected components of $\widehat{\mathcal{O}}_u\setminus\widehat{\mathcal{P}}$. Furthermore, at most one of these curves can pass through the bounded oval $\mathcal{O}_b$.

In order to be able to prove this, we first need to establish some monotonicity properties of $\Re\Phi_{(\xi_1,\xi_2)}$.

\begin{lemma}\label{lem:Phi_xi1_xi2_on_Gammas}
    Let $(\xi_1,\xi_2)\in\mathfrak{L}$. Then, the following hold.
	\begin{enumerate}
    	\item[\rm{(a)}] $\Re\Phi_{(\xi_1,\xi_2)}$ is strictly monotone along each component of the unbounded oval $\widehat{\mathcal{O}}_u$, tending to $+\infty$ and $-\infty$ at the two respective ends of the component.
        \item[\rm{(b)}] On the cycles $\Gamma_j$ and $\Gamma_j^\sigma$, for $j=0,1,2$, one has
    	\begin{equation}\label{eq:Phi_on_Gammas}
    		\Re\Phi_{(\xi_1,\xi_2)}(p) = 
    		\begin{cases}
    			-\xi_1\log|\lambda(p)| & \text{on } \Gamma_0\cup\Gamma_0^\sigma, \\
    			\xi_2\log|z(p)| & \text{on } \Gamma_1\cup\Gamma_1^\sigma, \\
    			(\xi_2-\xi_1)\log|z(p)| & \text{on } \Gamma_2\cup\Gamma_2^\sigma.
    		\end{cases}
    	\end{equation}
    For each $j=0,1,2$, $\Re\Phi_{(\xi_1,\xi_2)}$ is either identically zero on $\Gamma_j\cup\Gamma_j^\sigma$, or it has constant sign and is strictly monotone along $\Gamma_j\cap\mathcal{R}_+$ and $\Gamma_j^\sigma\cap\mathcal{R}^\sigma_+$.
    \end{enumerate}
\end{lemma}

\begin{proof}
    (a) Note that $\Re\Phi_{(\xi_1,\xi_2)}$ tends to $+\infty$ and $-\infty$ on each component of $\widehat{\mathcal{O}}_u$. Suppose that it is not monotone along one of these components. Then, $\Re\Phi_{(\xi_1,\xi_2)}$ attains a local maximum and a local minimum on that component, and this critical point is a zero of $\omega_{(\xi_1,\xi_2)}$ due to \eqref{eq:omega_xi1_xi2_2}. However, this contradicts the assumption $(\xi_1,\xi_2)\in\mathfrak{L}$, as there are no zeros on the unbounded oval $\widehat{\mathcal{O}}_u$ in this case. This concludes part~(a).

    \medskip

    (b) By Lemma~\ref{lem:Phi_00_on_Gammas_Upsilons}(a), we have $\Re\Phi = 0$ on $\Gamma_0\cup\Gamma_0^\sigma$. This implies by \eqref{eq:Phi_xi1_xi2_def} that $\Re\Phi_{(\xi_1,\xi_2)}(p) = -\xi_1\log|\lambda(p)|$, since $\log|z(p)| = 0$ on $\Gamma_0\cup\Gamma_0^\sigma$. The remaining cases in \eqref{eq:Phi_on_Gammas} follow similarly.

    The monotonicity of $\Re\Phi_{(\xi_1,\xi_2)}$ along the cycles $\Gamma_j$ and $\Gamma_j^\sigma$ within $\mathcal{R}_+$ or $\mathcal{R}_+^\sigma$ follows immediately from \eqref{eq:Phi_on_Gammas}.
\end{proof}

The monotonicity properties from Lemma~\ref{lem:Phi_xi1_xi2_on_Gammas} have two immediate consequences for the structure of $\widehat{\mathcal{N}}_{(\xi_1,\xi_2)}$.

\begin{corollary}\label{cor:Nhat_xi1_xi2_and_Gammas}
    For $(\xi_1,\xi_2)\in\mathfrak{L}$, the following hold.
	\begin{enumerate}
        \item[\rm{(a)}] The level set $\widehat{\mathcal{N}}_{(\xi_1,\xi_2)}$ intersects each component of the unbounded oval $\widehat{\mathcal{O}}_u$ at exactly one point.
        \item[\rm{(b)}] For each $j=0,1,2$, the level set $\widehat{\mathcal{N}}_{(\xi_1,\xi_2)}$ either intersects each of $\Gamma_j \cap \mathcal{R}_+$ and $\Gamma_j^\sigma \cap \mathcal{R}_+^\sigma$ in at most one point, or it contains the entire union $\Gamma_j\cup\Gamma_j^\sigma$.
    \end{enumerate}
\end{corollary}

\begin{proof}
    This is immediate from Lemma~\ref{lem:Phi_xi1_xi2_on_Gammas}.
\end{proof}

\begin{remark}\label{rem:Phi_xi1_xi2_on_Gammas}
    The monotonicity statement in Lemma~\ref{lem:Phi_xi1_xi2_on_Gammas}(b) can be made more precise for the case $0<\xi_2<\xi_1$.
    
    In this case, $\Re\Phi_{(\xi_1,\xi_2)}$ is positive on $\Gamma_1$ and $\Gamma_2$, and negative on $\Gamma_0$. Indeed, this follows immediately from \eqref{eq:Phi_on_Gammas}, since $\log|\lambda(p)|>0$ on $\Gamma_0$, $\log|z(p)|>0$ on $\Gamma_1$, and $\log|z(p)|<0$ on $\Gamma_2$. Moreover, $\Re\Phi_{(\xi_1,\xi_2)}$ is increasing along $\Gamma_1\cap\mathcal{R}_+$ and $\Gamma_2\cap\mathcal{R}_+$, and decreasing along $\Gamma_0\cap\mathcal{R}_+$, where the curves are oriented from the bounded oval to the unbounded oval. These monotonicity properties again follow immediately from \eqref{eq:Phi_on_Gammas}.

    Different patterns of signs and monotonicity arise for the cases $0<\xi_1<\xi_2$ and $\xi_2<0<\xi_1$.
\end{remark}

We are now in a position to describe the different topological structures of $\mathcal{N}_{(\xi_1,\xi_2)}$.

\begin{proposition}\label{prop:critical_level_sets}
    Let $(\xi_1,\xi_2)\in\mathfrak{L}$ with $0\leq\xi_2\leq\xi_1$. Then, the level set $\mathcal{N}_{(\xi_1,\xi_2)}$ is topologically equivalent to one of the subfigures in Figure~\ref{fig:critical_level_sets_1}. If $\xi_1=\xi_2>0$, then $\mathcal{N}_{(\xi_1,\xi_2)}$ is the level set shown in the left subfigure, and if $0=\xi_2<\xi_1$, it is shown in the right subfigure.
\end{proposition}

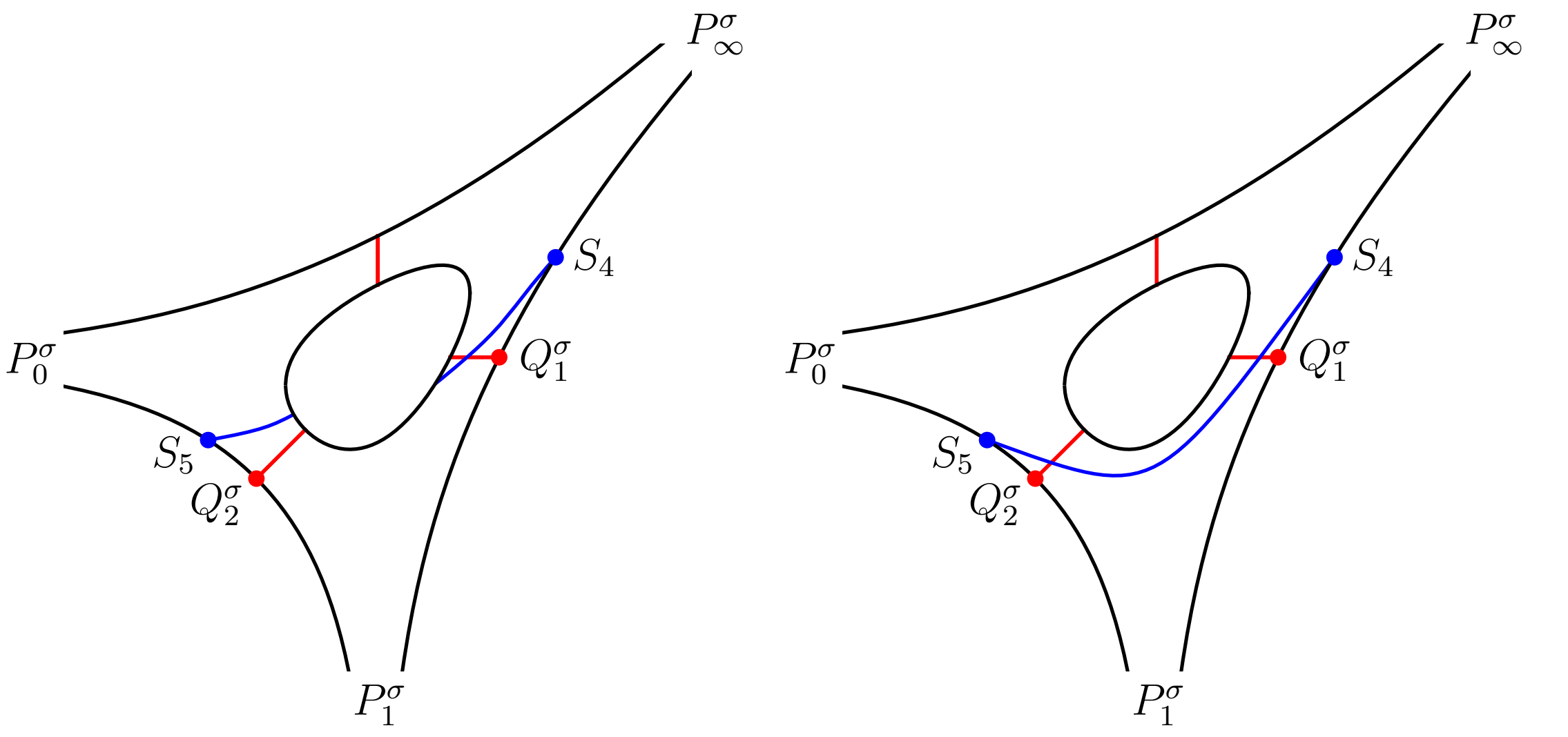
\begin{figure}[t]
    \centering
    \hspace{0.6cm}
	\begin{overpic}[width=0.4\linewidth]{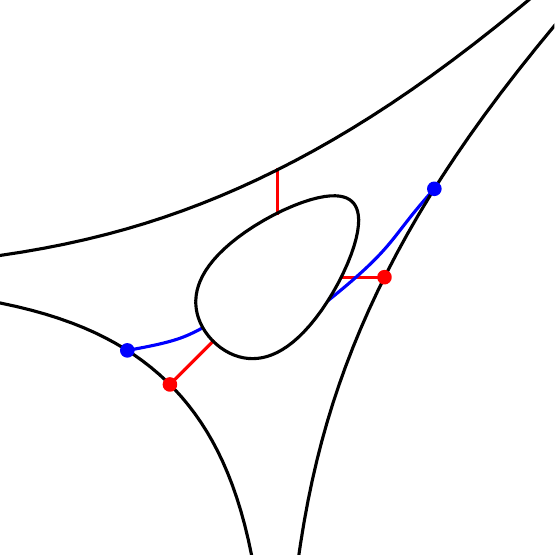}
		\put(-9.5,47.5){$P_0^\sigma$}
		\put(46,-7){$P_1^\sigma$}
		\put(99,100){$P_\infty^\sigma$}
		\put(72.5,48){$Q_1^\sigma$}
		\put(20,25){$Q_2^\sigma$}
		\put(81,64){$S_4$}
		\put(14,32.5){$S_5$}
    \end{overpic}
    \hfill
	\begin{overpic}[width=0.4\linewidth]{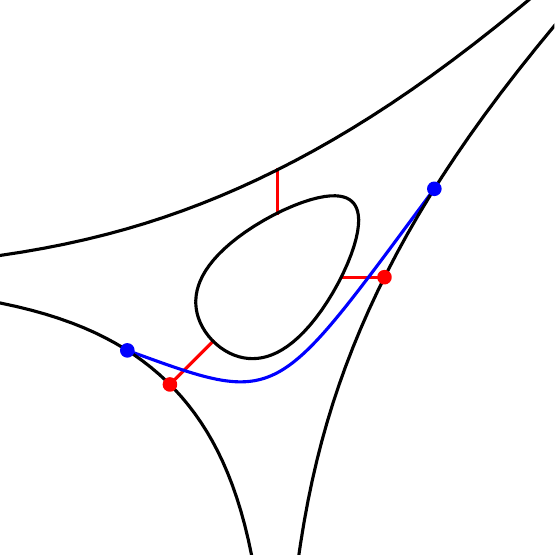}
		\put(-9.5,47.5){$P_0^\sigma$}
		\put(46,-7){$P_1^\sigma$}
		\put(99,100){$P_\infty^\sigma$}
		\put(72.5,48){$Q_1^\sigma$}
		\put(20,25){$Q_2^\sigma$}
		\put(81,64){$S_4$}
		\put(14,32.5){$S_5$}
    \end{overpic}
    \hspace{0.6cm}
    \vspace{0.4cm}
    \caption{
    The part of the level set $\widehat{\mathcal{N}}_{(\xi_1,\xi_2)}$
    that is in $\mathcal{R}^\sigma$ is shown
    in the amoeba $\mathcal{A}^\sigma$
    for the case $(\xi_1,\xi_2) \in \mathfrak{L}$ with  $0 \leq \xi_2 \leq \xi_1$.
    There are two possible topological
    situations as shown in the left and right subfigure. This illustrates Lemma~\ref{lem:Nsigma}.
	\label{fig:critical_level_sets_3}}
\end{figure}

\begin{proof}
    Let us first assume that $0<\xi_2<\xi_1$. Recall that by \eqref{eq:xi1_xi2_conditions}, the saddle point $s = \Omega(\xi_1,\xi_2)$ then belongs to the connected component of $\mathcal{R}_+\setminus(\Gamma_1\cup\Gamma_2)$ that contains $P_1$. The boundary of this region of $\mathcal{R}_+$ consists of five parts, see also Figure~\ref{fig:critical_level_sets_2} for its image in the amoeba. For $j=0,1,2$, we denote by $Q_j$ and $R_j$ the points of intersection of $\Gamma_j$ with the unbounded oval $\mathcal{O}_u$ and the bounded oval $\mathcal{O}_b$, respectively. The points $R_j$ are the same as the ones that were used before in Figure~\ref{fig:torus_cycles}. Figure~\ref{fig:critical_level_sets_2} shows their images on the boundary of the amoeba.

    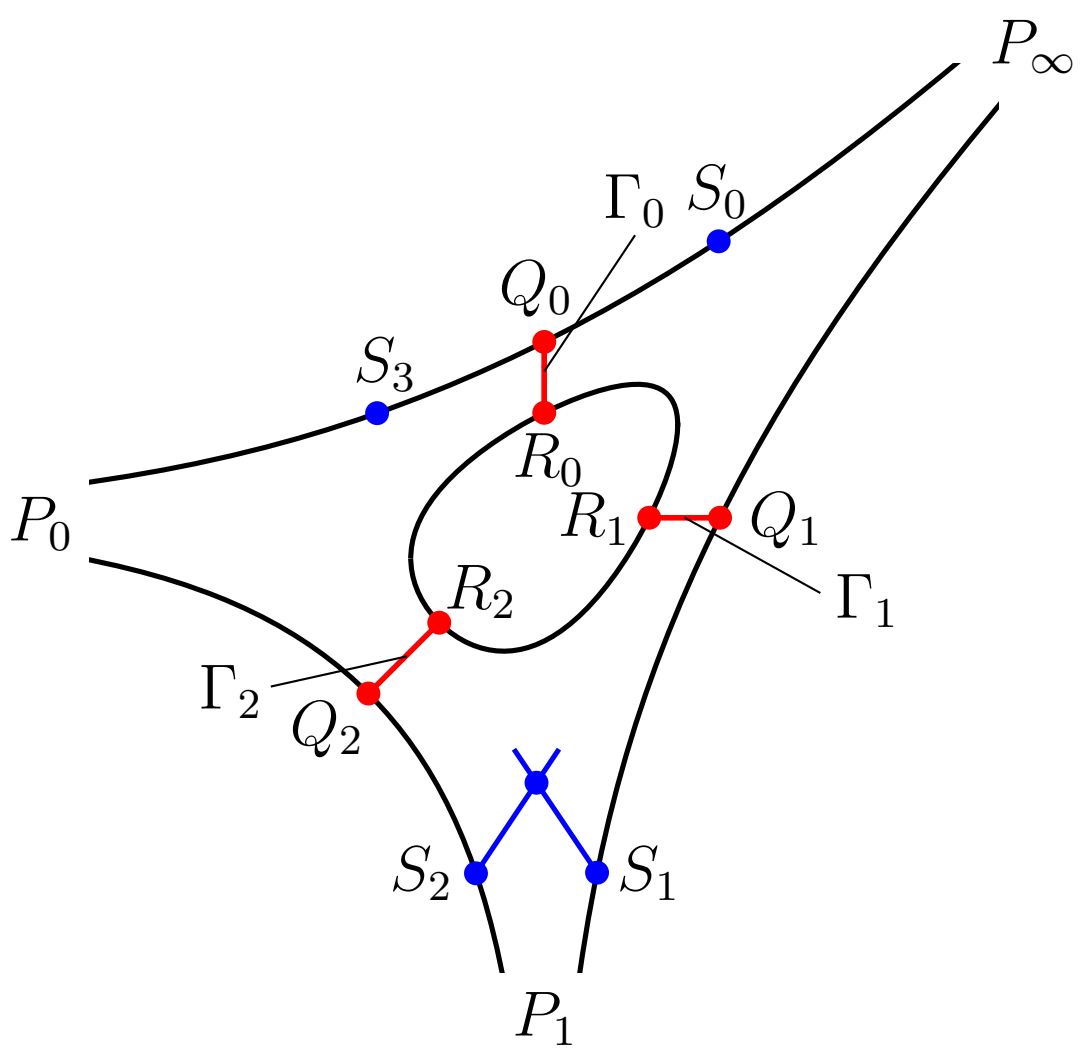
\begin{figure}[t]
    	\centering
    	\begin{overpic}[width=0.4\linewidth]{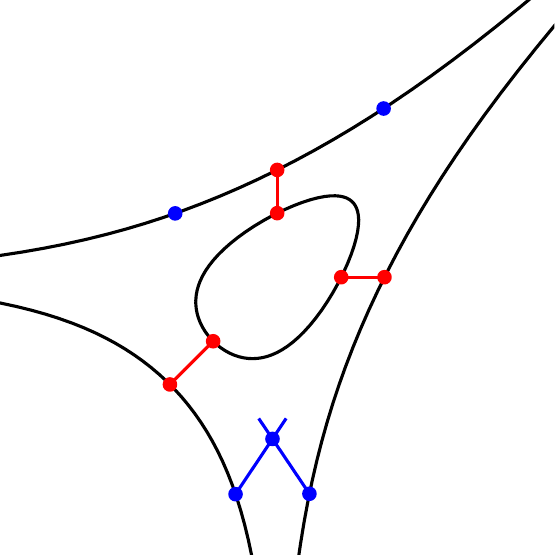}
    		\put(-9,47.5){$P_0$}
    		\put(46.5,-7){$P_1$}
    		\put(99,100){$P_\infty$}
    		\put(56.5,83){$\Gamma_0$}
    		\put(47.6,66){
    			\begin{tikzpicture}
    				\draw[-] (0,0) -- (-0.6,-0.9);
    		\end{tikzpicture}}
    		\put(82,39){$\Gamma_1$}
    		\put(63,41.6){
    			\begin{tikzpicture}
    				\draw[-] (0,0) -- (-0.9,0.5);
    		\end{tikzpicture}}
    		\put(12,29){$\Gamma_2$}
    		\put(17.5,31.3){
    			\begin{tikzpicture}
    				\draw[-] (0,0) -- (-0.9,-0.2);
    		\end{tikzpicture}}
    		\put(45,73.5){$Q_0$}
    		\put(72.5,48){$Q_1$}
    		\put(22,25){$Q_2$}
    		\put(46.5,54.5){$R_0$}
    		\put(51.5,48){$R_1$}
    		\put(39,40){$R_2$}
    		\put(65.5,84){$S_0$}
    		\put(58,9){$S_1$}
    		\put(33,9){$S_2$}
    		\put(29,65){$S_3$}
    	\end{overpic}
    	\vspace{0.4cm}
    	\caption{The amoeba with the saddle point $s$ and the points $S_0$, $S_1$, $S_2$, and $S_3$ (blue) on the unbounded oval that are used in the proof of Proposition~\ref{prop:critical_level_sets}. The figure also shows the cycles $\Gamma_j$ (red), together with their intersections $Q_j$ and $R_j$ with the unbounded and bounded ovals, for $j=0,1,2$.}
    	\label{fig:critical_level_sets_2}
    \end{figure}

    The five parts of the boundary are:
    \begin{description}[font=\normalfont,itemsep=0pt,labelwidth=52pt,labelsep=0pt]
        \item[Part I:] the contour $\Gamma_1\cap\mathcal{R}_+$;
        \item[Part II:] the part of the unbounded oval between $Q_1$ and $P_1$;
        \item[Part III:] the part of the unbounded oval between $P_1$ and $Q_2$;
        \item[Part IV:] the contour $\Gamma_2\cap\mathcal{R}_+$;
        \item[Part V:] the part of the bounded oval between $R_2$ and $R_1$.
    \end{description}

    Each of the four level curves that emanate from $s$ must reach one of the five parts. Indeed, if such a level curve does not reach the boundary, then it encloses a bounded domain in $\mathcal{R}$ that does not contain any poles of $\Re\Phi_{(\xi_1,\xi_2)}$. Since $\Re\Phi_{(\xi_1,\xi_2)}$ is constant on the boundary of this domain, this then violates the maximum principle for harmonic functions, as strict inequality holds in the interior. A similar argument shows that Part~V can be reached by a level curve at most once. This is also the case for Parts~I--IV, but this is a consequence of the strict monotonicity of $\Re\Phi_{(\xi_1,\xi_2)}$ along these parts, see Lemma~\ref{lem:Phi_xi1_xi2_on_Gammas}(b) and Remark~\ref{rem:Phi_xi1_xi2_on_Gammas}.

    It follows that one part is not reached by a level curve, while the other four parts are reached exactly once. In particular, at least one of $\Gamma_1$ and $\Gamma_2$ is reached by a level curve. Since $\Re\Phi_{(\xi_1,\xi_2)}(p) = C_{(\xi_1,\xi_2)}$ for $p \in \mathcal{N}_{(\xi_1,\xi_2)}$, we have by \eqref{eq:Phi_on_Gammas} that either
	\[
        C_{(\xi_1,\xi_2)} = \xi_2\log|z(p)| \quad \text{for some } p\in\Gamma_1,
    \]
    which is positive because $\log|z(p)|>0$ on $\Gamma_1$ and $\xi_2>0$, or
	\[
        C_{(\xi_1,\xi_2)} = (\xi_2-\xi_1)\log|z(p)| \quad \text{for some } p\in\Gamma_2.
    \]
    which is positive as well because $\log|z(p)|<0$ on $\Gamma_2$ and $\xi_2-\xi_1<0$. Thus, in both cases, we see that
    \begin{equation}\label{eq:C_xi1_xi_2}
        C_{(\xi_1,\xi_2)}>0.
    \end{equation}
    As noted in Remark~\ref{rem:Phi_xi1_xi2_on_Gammas}, $\Re\Phi_{(\xi_1,\xi_2)}(p)<0$ on $\Gamma_0\cup\Gamma_0^\sigma$. Together with $\Re\Phi_{(\xi_1,\xi_2)} = C_{(\xi_1,\xi_2)}>0$ on $\mathcal{N}_{(\xi_1,\xi_2)}$, we find that the four level curves cannot come to  $\Gamma_0\cup\Gamma_0^\sigma$. Hence, these curves stay within $\mathcal{R}$.

    Recall from Corollary~\ref{cor:Nhat_xi1_xi2_and_Gammas}(a) that $\mathcal{N}_{(\xi_1,\xi_2)}$ intersects each component of the unbounded oval exactly once. These points will be denoted by $S_j$ in what follows, and their images under the amoeba map are in the Figures~\ref{fig:critical_level_sets_2} and~\ref{fig:critical_level_sets_3}. We will show that these Figures give the accurate picture, with four points $S_0$, $S_1$, $S_2$ and $S_3$ on $\mathcal{R}$ and two point $S_4$ and $S_5$ on $\mathcal{R}^\sigma$. 
    
    Since $\Re\Phi_{(\xi_1,\xi_2)}$ tends to $+\infty$ at both $P_\infty$ and $P_0$ and $\Re\Phi_{(\xi_1,\xi_2)}(Q_0)<0<C_{(\xi_1,\xi_2)}$, it follows that $\mathcal{N}_{(\xi_1,\xi_2)}$ intersects the unbounded oval between $P_\infty$ and $Q_0$ at a unique point $S_0$, as well as the unbounded oval between $P_0$ and $Q_0$ at $S_3$, as shown in Figure~\ref{fig:critical_level_sets_2}.

	Next, we claim that 
	\begin{equation}\label{eq:Phi_at_Q12}
		\Re\Phi_{(\xi_1,\xi_2)}(Q_j) \geq C_{(\xi_1,\xi_2)} \quad \text{for } j=1,2.
	\end{equation}
    Suppose, to get a contradiction, that \eqref{eq:Phi_at_Q12} is not true for $j=1$. By Remark~\ref{rem:Phi_xi1_xi2_on_Gammas}, $\Re\Phi_{(\xi_1,\xi_2)}$ decreases along $\Gamma_1$ as one moves from $Q_1$ towards $R_1$. It also decreases along the unbounded oval from $Q_1$ towards $P_1$. Then, it follows that both Part~I and Part~II are entirely contained in $\mathcal{N}_{(\xi_1,\xi_2)}^<$. However, this contradicts the fact that exactly one part is not reached by a level curve. Thus, \eqref{eq:Phi_at_Q12} holds for $j=1$, and a similar reasoning shows it holds for $j=2$. The same reasoning applies to $Q_2$.
	
    Since $\Re\Phi_{(\xi_1,\xi_2)}$ tends to $-\infty$ at $P_1$, it follows from \eqref{eq:Phi_at_Q12} that $\mathcal{N}_{(\xi_1,\xi_2)}$ intersects the unbounded oval between $P_1$ and $Q_j$ at a point $S_j$ for $j=1,2$. The images of these points in the amoeba are also shown in Figure~\ref{fig:critical_level_sets_2}.

    There is one level curve emanating from $s$ that reaches Part~II at $S_1$. Indeed, otherwise, Part~I and Parts~III--V would be reached by such level curves, and then the level curve through $S_1$ has nowhere to go. Similarly, there is another level curve emanating from $s$ that reaches Part~III at $S_2$.

    There is at least one other level curve emanating from $s$ that does not come to the bounded oval. This curve reaches either Part~I or Part~IV, and then it necessarily continues to either $S_0$ or $S_3$, as the curve cannot pass through $\Gamma_0$ or the bounded oval.

    The remaining level curve emanating from $s$ either reaches the other one of $S_0$ and $S_3$, or it comes to the bounded oval. In the latter case, the level curve does not necessarily reach Part~V, as it could first pass through Part~I or Part~IV and only come to the bounded oval after that.

    We now distinguish two cases depending on the behavior of the level curves passing through $S_0$ and $S_3$. As argued above, at least one of these points is reached by a level curve emanating from $s$. The three cases are now as follows (see also Figure~\ref{fig:critical_level_sets_1}):
    \begin{enumerate}
        \item[\rm{(i)}] Only $S_0$ is reached by a level curve emanating from $s$, another such level curve passes through the bounded oval, and the level curve emanating from $S_3$ also passes through the bounded oval. This is the situation from the left subfigure of Figure~\ref{fig:critical_level_sets_1}.
        \item[\rm{(ii)}] Both $S_0$ and $S_3$ are reached by a level curve emanating from $s$. This is what happens in the middle subfigure of Figure~\ref{fig:critical_level_sets_1}.
        \item[\rm{(iii)}] A case analogous to case (i), obtained by interchanging the roles of $S_0$ and $S_3$ as in the right subfigure).
    \end{enumerate}
    The cases (i)--(iii) exhaust all possible configurations. This completes the proof for the case $0<\xi_2<\xi_1$

    The remaining cases $\xi_1=\xi_2>0$ and $0=\xi_2<\xi_1$ are similar, as they are limiting cases from the ones  considered above. In case $\xi_2=0$, we have $s\in\Gamma_1$ and the full cycle $\Gamma_1$ belongs to $\mathcal{N}_{(\xi_1,\xi_2)}$ by \eqref{eq:Phi_on_Gammas}. This leads to the right subfigure of Figure~\ref{fig:critical_level_sets_1}. For $\xi_1=\xi_2$, we similarly have that $s \in \Gamma_2 \subset \mathcal{N}_{(\xi_1,\xi_2)}$, and we obtain the left subfigure.
\end{proof}	

We complete the study of $\widehat{\mathcal{N}}_{(\xi_1,\xi_2)}$, for the case $0 \leq \xi_2 \leq \xi_1$, by describing its part on $\mathcal{R}^\sigma$. The result is illustrated in Figure~\ref{fig:critical_level_sets_3}.

\begin{lemma}\label{lem:Nsigma} 
    Let $(\xi_1,\xi_2)\in\mathfrak{L}$ with $0\leq\xi_2\leq\xi_1$. Then, the level set $\mathcal{N}_{(\xi_1,\xi_2)}^\sigma$ is topologically equivalent to one of the subfigures in Figure~\ref{fig:critical_level_sets_3}.
\end{lemma}

\begin{proof}
    Recall from Corollary~\ref{cor:Nhat_xi1_xi2_and_Gammas}(a) that $\widehat{\mathcal{N}}_{(\xi_1,\xi_2)}$ intersects each of the six components of the unbounded oval exactly once. Four of these intersection points are in $\mathcal{R}$ and they were called $S_0$, $S_1$, $S_2$ and $S_3$ in the proof of Proposition~\ref{prop:critical_level_sets}, see also Figure~\ref{fig:critical_level_sets_2}. That leaves two intersection points $S_4$ and $S_5$ on $\mathcal{R}^\sigma$, with $S_4$ on the component between $P_\infty^\sigma$ and $P_1^\sigma$, and $S_5$ on the component between $P_1^\sigma$ and $P_0^\sigma$.
    
    By \eqref{eq:Phi_on_Gammas} and \eqref{eq:Phi_at_Q12}, we have
    \begin{equation}\label{eq:Phi_at_Q12_sigma}
        \Re\Phi_{(\xi_1,\xi_2)}(Q_j^\sigma) = \Re\Phi_{(\xi_1,\xi_2)}(Q_j) \geq C_{(\xi_1,\xi_2)} \quad \text{for } j=1,2.
    \end{equation}
    Since $\Re\Phi_{(\xi_1,\xi_2)}$ tends $-\infty$ at both $P_\infty^\sigma$ and $P_0^\sigma$, it follows that $S_4$ is on the unbounded oval in between $P_\infty^\sigma$ and $Q_1^\sigma$, and $S_5$ is in between $P_0^\sigma$ and $Q_2^\sigma$, as shown in Figure~\ref{fig:critical_level_sets_3}. Note that equality with $Q_1^\sigma$ or $Q_2^\sigma$ is possible.

    The level lines emanating from $S_4$ and $S_5$ can either both come to  the bounded oval, (left subfigure of Figure~\ref{fig:critical_level_sets_3}), or they connect to each other while passing through both $\Gamma_1^\sigma$ and $\Gamma_2^\sigma$ (right subfigure). These are the only possibilities, since the level lines cannot reach any other components of the unbounded oval or $\Gamma_0\cup\Gamma_0^\sigma$.
\end{proof}

The structure of $\widehat{\mathcal{N}}_{(\xi_1,\xi_2)}$ for the case $0\leq\xi_2\leq\xi_1$ is now well established by Proposition~\ref{prop:critical_level_sets} and Lemma~\ref{lem:Nsigma}. The remaining cases follow by applying the symmetries of the double cover discussed in Section~\ref{sec:symmetries}.

\begin{lemma}\label{lem:critical_level_sets}
   With $\rho$ the rotation symmetry as in \eqref{eq:rho_def}, one has
    \begin{align}
    	\widehat{\mathcal{N}}_{(\xi_1,\xi_2)} &= \rho\bigl(\widehat{\mathcal{N}}_{(\xi_2,\xi_2-\xi_1)}\bigr), \label{eq:critical_level_set_rho} \\
        \widehat{\mathcal{N}}_{(\xi_1,\xi_2)} &= \rho^{-1}\bigl(\widehat{\mathcal{N}}_{(\xi_1-\xi_2,\xi_1)}\bigr), \label{eq:critical_level_set_rho_inv}
    \end{align}
    and similarly, with the double cover involution $\sigma$ from \eqref{eq:sigma}, one has
    \begin{equation}\label{eq:critical_level_set_sigma}
        \widehat{\mathcal{N}}_{(\xi_1,\xi_2)} = \sigma\bigl(\widehat{\mathcal{N}}_{(-\xi_1,-\xi_2)}\bigr).
    \end{equation}
\end{lemma}

\begin{proof}
    We note that \eqref{eq:critical_level_set_rho_inv}--\eqref{eq:critical_level_set_rho} follow immediately from \eqref{eq:Phi_xi1_xi2_rho}--\eqref{eq:Phi_xi1_xi2_rho_inv}, and \eqref{eq:critical_level_set_sigma} follows from the first identity in \eqref{eq:Phi_xi1_xi2_sigmas}.
\end{proof}

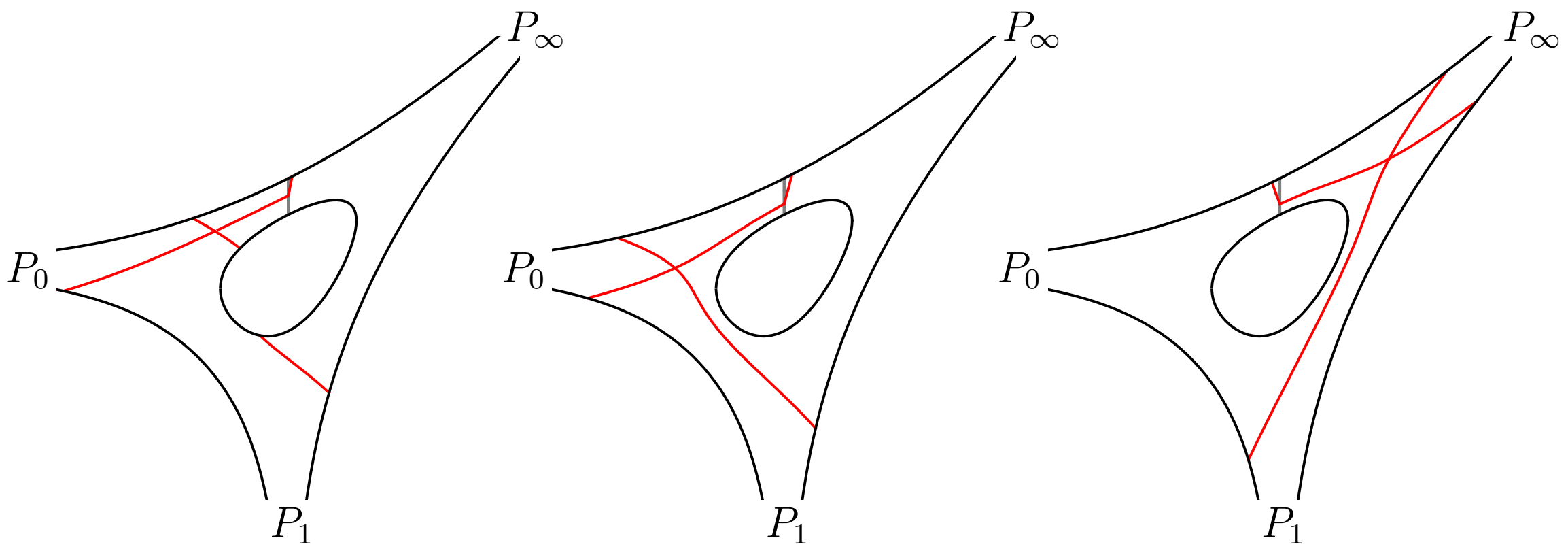
\begin{figure}[t]
    \hspace{0.3cm}
	\begin{overpic}[width=0.3\linewidth]{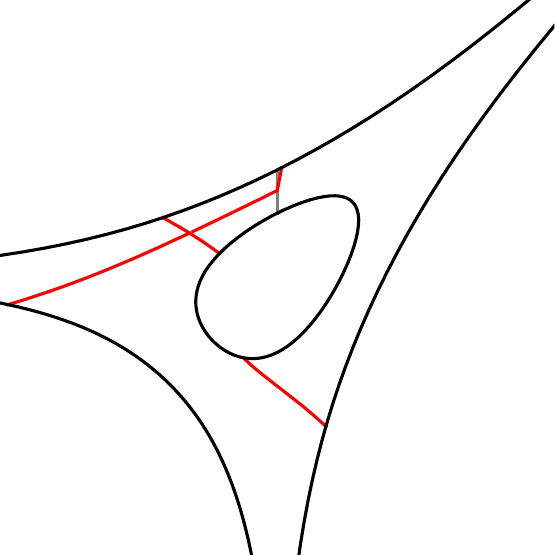}	
		\put(-11,47){$P_0$}
		\put(46,-8){$P_1$}
	\end{overpic}
    \hfill
	\begin{overpic}[width=0.3\linewidth]{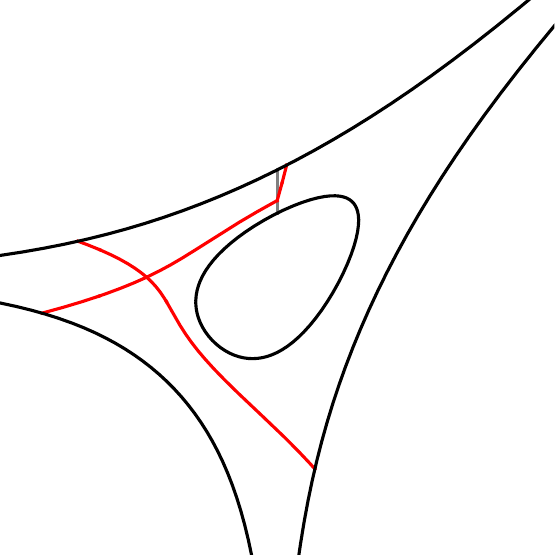}	
        \put(-10,99){$P_\infty$}
		\put(-11,47){$P_0$}
		\put(46,-8){$P_1$}
	\end{overpic}
    \hfill
	\begin{overpic}[width=0.3\linewidth]{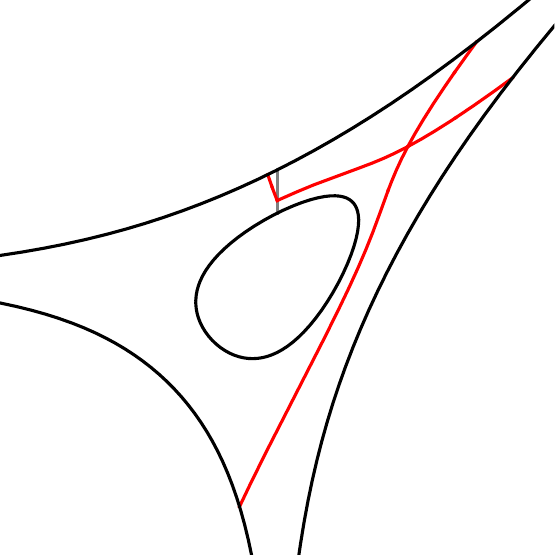}	
        \put(-10,99){$P_\infty$}
		\put(-11,47){$P_0$}
		\put(46,-8){$P_1$}
		\put(98,99){$P_\infty$}
	\end{overpic}
    \hspace{0.2cm}
    \vspace{0.25cm}
    \caption{Level sets $\mathcal{N}_{(\xi_1,\xi_2)}$ for various $(\xi_1,\xi_2)\in\mathfrak{L}$ with $0\leq\xi_1<\xi_2$ (left and middle subfigure) and $\xi_2<0 \leq\xi_1$ (right subfigure).}
    \label{fig:critical_level_sets_4}
\end{figure}

Recall from \eqref{eq:xi1_xi2_conditions} that $\xi_1\geq0$, which corresponds to either $s\in\mathcal{R}_+$ (in case of strict inequality) or $s\in\Gamma_0\cup\Gamma_0^\sigma$ (if $\xi_1=0$). Due to \eqref{eq:critical_level_set_sigma}, the level set $\widehat{\mathcal{N}}_{(\xi_1,\xi_2)}$ for $\xi_1<0$ is simply obtained from the case $\xi_1>0$ by interchanging the roles of $\mathcal{R}$ and $\mathcal{R}^\sigma$. We are therefore left to consider the cases $0\leq\xi_1<\xi_2$ and $\xi_2<0\leq\xi_1$.

Suppose that $0<\xi_1<\xi_2$. Since then $0<\xi_2-\xi_1<\xi_2$, the level sets $\mathcal{N}_{(\xi_2,\xi_2-\xi_1)}$ and $\mathcal{N}_{(\xi_2,\xi_2-\xi_1)}^\sigma$ are topologically equivalent to one of the subfigures Figures~\ref{fig:critical_level_sets_1} and~\ref{fig:critical_level_sets_3}, respectively. By \eqref{eq:critical_level_set_rho}, applying the rotation symmetry $\rho$ yields $\widehat{\mathcal{N}}_{(\xi_1,\xi_2)}$. In particular, the level set $\mathcal{N}_{(\xi_1,\xi_2)}$ is topologically equivalent to either the left and middle subfigures of Figures~\ref{fig:critical_level_sets_4}, or the left subfigure of Figures~\ref{fig:critical_level_sets_1}.

Similarly, the level set $\widehat{\mathcal{N}}_{(\xi_1,\xi_2)}$ for $\xi_2<0<\xi_1$ can be determined. The right subfigure of Figure~\ref{fig:critical_level_sets_4} shows one example of $\mathcal{N}_{(\xi_1,\xi_2)}$.

Finally, the case $\xi_1=0$ is special in our analysis and will be treated separately in what follows. The level set $\widehat{\mathcal{N}}_{(\xi_1,\xi_2)}$ in this case is given in Section~\ref{sec:remaining_cases}, see Figure~\ref{fig:deform_BD} below.


\subsection{First transformation $\mathbf{Y}_N \mapsto \mathbf{X}_N$}\label{sec:first_transformation}
The first transformation $\mathbf{Y}_N \mapsto \mathbf{X}_N$ involves the matrix $E$ consisting of eigenvectors of $W$, see \cite[Definition~5.1]{Kui25} for the precise definition of $E$. We recall from \cite[Lemma~5.2]{Kui25} that the matrix-valued function $E$ is analytic and invertible on $\mathbb{C}\setminus\mathbb{R}$, and it satisfies the jump condition $E_+ = E_-J_E$ with
\begin{align}\label{eq:jump_E}
	J_E =
	\begin{cases}
		\sigma_{01} \coloneqq \left(\begin{smallmatrix} 0 & 1 & 0 \\ 1 & 0 & 0 \\ 0 & 0 & 1 \end{smallmatrix}\right) & \text{on } (-\infty,0], \vspace{0.1cm} \\
		\sigma_{12} \coloneqq \left(\begin{smallmatrix} 1 & 0 & 0 \\ 0 & 0 & 1 \\ 0 & 1 & 0 \end{smallmatrix}\right) & \text{on } (0,z_{\min})\cup(z_{\max},\infty), \\
        I_3 & \text{on } [z_{\min},z_{\max}].
	\end{cases} 
\end{align}
Note that $\sigma_{01}$ and $\sigma_{12}$ are permutation matrices. We keep the convention that rows and columns of $3\times3$ matrices are labeled by $0$, $1$, and $2$. The matrix-valued function, see \cite[Definition~5.5]{Kui25},
\begin{equation}\label{eq:XN}
	\mathbf{X}_N = \mathbf{Y}_N \begin{pmatrix} E & 0_3 \\ 0_3 & E \end{pmatrix}
\end{equation}
is then analytic in $\mathbb{C}\setminus(\mathbb{R}\cup\mathbb{T})$, and it satisfies the jump condition $\mathbf{X}_{N,+} = \mathbf{X}_{N,-}J_{\mathbf{X}_N}$ with
\begin{equation}\label{eq:jump_X} 
	J_{\mathbf{X}_N}(z) =
	\begin{cases}
		\begin{pmatrix} J_E(z) & 0_3 \\ 0_3 & J_E(z) \end{pmatrix}, & z\in\mathbb{R}, \\
		\begin{pmatrix} I_3 & z^{-2N} \Lambda(z)^{2N} \\ 0_3 & I_3 \end{pmatrix}, & z\in\mathbb{T},
	\end{cases}
\end{equation}
see also \cite[RH problem~5.6]{Kui25}. Here, $\Lambda(z) = \diag(\lambda_1(z),\lambda_2(z),\lambda_3(z))$ comes from the spectral decomposition $W(z) = E(z)\Lambda(z)E^{-1}(z)$ of $W(z)$.

\begin{lemma}\label{lem:KN_XN}
	After applying the transformation $\mathbf{Y}_N \mapsto \mathbf{X}_N$, one obtains that \eqref{eq:KN_matrix} can be written as
	\begin{multline}\label{eq:KN_XN} 
		\mathsf{K}_N(v_{1,N},v_{2,N}) = \frac{1}{(2\pi i)^2}\oint_\mathbb{T} dz_2 \oint_\mathbb{T} \frac{dz_1}{z_1(z_1-z_2)} T_{0,j_2}^{-1}(z_2) E(z_2) \frac{\Lambda(z_2)^{N(1-\xi_{1,N})-x_2}}{z_2^{N(1-\xi_{2,N})-y_2}} \\
		\times\begin{pmatrix} 0_3 & I_3 \end{pmatrix} \mathbf{X}_N^{-1}(z_2) \mathbf{X}_N(z_1) \begin{pmatrix} I_3 \\ 0_3 \end{pmatrix} \frac{\Lambda(z_1)^{N(1+\xi_{1,N})+x_1}}{z_1^{N(1+\xi_{2,N})+y_1}} E^{-1}(z_1)T_{0,j_1}(z_1).
	\end{multline}
\end{lemma}

\begin{proof}
	The integrand in \eqref{eq:KN_XN} is the result of substituting \eqref{eq:RN_formula} and \eqref{eq:XN} into \eqref{eq:KN_v1_v2}, using that $W(z) = E(z)\Lambda(z)E^{-1}(z)$.
\end{proof}

Despite the factor $(z_1-z_2)^{-1}$ in \eqref{eq:KN_XN}, there is no singularity in the integrand for $z_1=z_2$. This is because the factor
\[
    \begin{pmatrix} 0_3 & I_3 \end{pmatrix} \mathbf{X}_N^{-1}(z_2)\mathbf{X}_N(z_1) \begin{pmatrix} I_3 \\ 0_3 \end{pmatrix}
\]
vanishes for $z_1=z_2$.


\subsection{Integrals on the Riemann surface}\label{sec:integrals_Riemann_surface}
We rewrite \eqref{eq:KN_XN} by putting the integrals on the Riemann surface $\mathcal{R}$. We use $p$ to denote a generic point on $\mathcal{R}$, and $z=z(p)$ and $\lambda=\lambda(p)$ are meromorphic functions of $p$. We also define the piecewise constant matrix-valued function defined on $\mathcal{R}$ by
\begin{equation}\label{eq:Xi_def}
    \Xi(p) =
    \begin{cases}
        \Xi_1 \coloneqq \left(\begin{smallmatrix} 1 & 0 & 0 \\ 0 & 0 & 0 \\ 0 & 0 & 0 \end{smallmatrix}\right) & \text{on } \mathcal{R}_1, \\[5pt]
        \Xi_2 \coloneqq \left(\begin{smallmatrix} 0 & 0 & 0 \\ 0 & 1 & 0 \\ 0 & 0 & 0 \end{smallmatrix}\right) & \text{on } \mathcal{R}_2, \\[5pt]
 		\Xi_3 \coloneqq \left(\begin{smallmatrix} 0 & 0 & 0 \\ 0 & 0 & 0 \\ 0 & 0 & 1 \end{smallmatrix}\right) & \text{on } \mathcal{R}_3.
    \end{cases}
\end{equation}
Then, $\Xi$ satisfies the following jump condition on the branch cuts
\begin{equation}\label{eq:Xi_jump}
    \begin{aligned}
        \Xi_+(p) &= \sigma_{01} \Xi_-(p) \sigma_{01} \quad \text{between } \mathcal{R}_1 \text{ and } \mathcal{R}_2, \\
        \Xi_+(p) &= \sigma_{12} \Xi_-(p) \sigma_{12} \quad \text{between } \mathcal{R}_2 \text{ and } \mathcal{R}_3,
    \end{aligned}
\end{equation}
where $\sigma_{01}$ and $\sigma_{12}$ are the permutation matrices from \eqref{eq:jump_E}.

Recalling $N\in\mathbb{N}$ and $N\xi_{1,N},N\xi_{2,N}\in\mathbb{Z}$, we introduce
\begin{align} 
	\mathbf{B}_N(p) &= \mathbf{Y}_N^{-t}(z(p)) \begin{pmatrix} 0_3 \\ I_3 \end{pmatrix} \left(E(z(p)) \Xi(p) E^{-1}(z(p))\right)^t\frac{\lambda(p)^{N(1-\xi_{1,N})}}{z(p)^{N(1-\xi_{2,N})}} \nonumber \\
	&= \mathbf{X}_N^{-t}(z(p)) \begin{pmatrix} 0_3 \\ I_3 \end{pmatrix} \Xi(p)  E^{t}(z(p))  \frac{\lambda(p)^{N(1-\xi_{1,N})}}{z(p)^{N(1-\xi_{2,N})}}, \label{eq:BN} \\
	\mathbf{C}_N(p) &= \mathbf{Y}_N(z(p)) \begin{pmatrix} I_3 \\ 0_3 \end{pmatrix} E(z(p)) \Xi(p) E^{-1}(z(p))\frac{\lambda(p)^{N(1+\xi_{1,N})}}{z(p)^{N(1+\xi_{2,N})}} \nonumber \\
	&= \mathbf{X}_N(z(p)) \begin{pmatrix} I_3 \\ 0_3 \end{pmatrix} \Xi(p)E^{-1}(z(p))\frac{\lambda(p)^{N(1+\xi_{1,N})}}{z(p)^{N(1+\xi_{2,N})}}. \label{eq:CN}
\end{align}
Note that $\mathbf{B}_N$ and $\mathbf{C}_N$ are $6\times3$ matrix-valued functions on the Riemann surface. The second identities in \eqref{eq:BN} and \eqref{eq:CN} are immediate from \eqref{eq:XN} and \eqref{eq:Xi_def}.

\begin{lemma}\label{lem:BNCN_poles}
    For $(\xi_{1,N},\xi_{2,N})\in\mathcal{H}$ with $N\xi_{1,N},N\xi_{2,N}\in\mathbb{Z}$, the following hold.
    \begin{enumerate}
        \item[\rm{(a)}] Both $\mathbf{B}_N$ and $\mathbf{C}_N$ are meromorphic on $\mathcal{R}$ with higher-order poles at $P_\infty$ and $P_0$, a higher-order zero at $P_1$, and at most simple poles at the branch points $B_1$ and $B_3$, see \eqref{eq:B1_B3}.
        \item[\rm{(b)}] One has
        \begin{multline}\label{eq:KN_XN_on_R} 
        	\mathsf{K}_N(v_{1,N},v_{2,N}) = \frac{1}{(2\pi i)^2}\oint\limits_{z^{-1}(\mathbb{T})} dz(p_2) \oint\limits_{z^{-1}(\mathbb{T})}  \frac{dz(p_1)}{z(p_1)(z(p_1)-z(p_2))} \\
            \times\frac{\lambda(p_2)^{-x_2}}{z(p_2)^{-y_2}}\frac{\lambda(p_1)^{x_1}}{z(p_1)^{y_1}}T_{0,j_2}^{-1}(z(p_2))\mathbf{B}^t_N(p_2)\mathbf{C}_N(p_1) T_{0,j_1}(z(p_1)).
        \end{multline}
    \end{enumerate}
\end{lemma}

\begin{proof}
    (a) The matrices $E$ and $E^{-1}$ are built out of meromorphic functions $e_k$ and $\widetilde{e}_l$ on $\mathcal{R}$ in the sense that
    \[
        E(z) = \begin{pmatrix} e_k(z, \lambda_j(z)) \end{pmatrix}_{k,j=1}^3, \quad E^{-1}(z) = \begin{pmatrix}\widetilde{e}_l(z,\lambda_j(z)) \end{pmatrix}_{j,k=1}^3
    \]
    see \cite[Definition~5.1 and Lemma~5.4]{Kui25}. It follows that
    \begin{equation}\label{eq:EXi_Einv}
        E(z(p)) \Xi(p) E^{-1}(z(p)) = \begin{pmatrix} e_k(p)\widetilde{e}_l(p) \end{pmatrix}_{k,l=1}^3.
    \end{equation}
    Each entry of \eqref{eq:EXi_Einv} is thus a meromorphic function on $\mathcal{R}$. The poles are easy to determine from \cite[Table~1]{Kui25}. In particular, there are simple poles at $B_1$ and $B_3$, and at most double poles at $P_\infty$ and $P_0$, which are the only possible poles.
    
    The scalar factors in \eqref{eq:BN} and \eqref{eq:CN} are clearly meromorphic on $\mathcal{R}$ with poles at $P_0$ and possibly at $P_\infty$. Since
    \[
        \mathbf{Y}_N^{-t}(z) \begin{pmatrix} 0_3 \\ I_3 \end{pmatrix}, \quad \mathbf{Y}_N(z) \begin{pmatrix} I_3 \\ 0_3 \end{pmatrix}
    \]
    in the first lines of \eqref{eq:BN} and \eqref{eq:CN} consist of polynomials in $z$, the corresponding first factors in \eqref{eq:BN} and \eqref{eq:CN} contribute to give a pole at $P_\infty$, but no other poles on the Riemann surface.
    
    \medskip
    
    (b)	The contour $z^{-1}(\mathbb{T})$ consists of the three circles $\mathbb{T}_1$, $\mathbb{T}_2$, and $\mathbb{T}_3$ on the Riemann surface. Suppose $p_1\in\mathbb{T}_k$ and $p_2\in\mathbb{T}_l$. Writing $z_1=z(p_1)$ and $z_2=z(p_2)$, we have $\lambda(p_1)=\lambda_k(z_1)$ and $\lambda(p_2)=\lambda_l(z_2)$. Hence, by \eqref{eq:BN},
    \begin{equation}\label{eq:BN_p2}
        \mathbf{B}_N(p_2) = \mathbf{X}_N^{-t}(z_2) \begin{pmatrix} 0_3 \\ I_3 \end{pmatrix} \Xi_l E^t(z_2)\frac{\lambda_l(z_2)^{N(1-\xi_{1,N})}}{z_2^{N(1-\xi_{2,N})}}, \quad p_2\in\mathcal{R}_l,
    \end{equation}
    for $l=1,2,3$, with $\Xi_l$ as in \eqref{eq:Xi_def}. Similarly, by \eqref{eq:CN},
    \begin{equation}\label{eq:CN_p1}
        \mathbf{C}_N(p_1) = \mathbf{X}_N(z_1) \begin{pmatrix} I_3 \\ 0_3 \end{pmatrix}\Xi_k E^{-1}(z_1)\frac{\lambda_k(z_1)^{N(1+\xi_{1,N})}}{z_1^{N(1+\xi_{2,N})}}, \quad p_1\in\mathcal{R}_k,
    \end{equation}
    for $k=1,2,3$. Using \eqref{eq:BN_p2} and \eqref{eq:CN_p1}, we find that the right-hand side of \eqref{eq:KN_XN_on_R} is equal to
    \begin{multline}\label{eq:KN_XN_on_R2}
    	\frac{1}{(2\pi i)^2}\sum_{k=1}^3 \sum_{l=1}^3\oint_\mathbb{T}\oint_\mathbb{T}\frac{dz_2 dz_1}{z_1(z_1-z_2)}\frac{\lambda_l(z_2)^{N(1-\xi_{1,N})-x_2}}{z_2^{N(1-\xi_{2,N})-y_2}}\frac{\lambda_k(z_1)^{N(1+\xi_{1,N})+x_1}}{z_1^{N(1+\xi_{2,N})+y_1}} \\
    	\times T_{0,j_2}^{-1}(z_2) E(z_2) \Xi_l \begin{pmatrix} 0_3 & I_3 \end{pmatrix} \mathbf{X}_N^{-1}(z_2) \mathbf{X}_N(z_1) \begin{pmatrix} I_3 \\ 0_3 \end{pmatrix} \Xi_k  E^{-1}(z_1)T_{0,j_1}(z_1).
    \end{multline}
    We observe  
    \begin{align*} 
        \sum_{k=1}^3 \lambda_k(z_1)^{N(1+\xi_{1,N})+x_1} \Xi_k &= \Lambda(z_1)^{N(1+\xi_{1,N})+x_1}, \\
        \sum_{l=1}^3 \lambda_l(z_2)^{N(1-\xi_{1,N})-x_2} \Xi_l &= \Lambda(z_2)^{N(1-\xi_{1,N})-x_2},
    \end{align*} 
    and applying this to \eqref{eq:KN_XN_on_R2}, we find that \eqref{eq:KN_XN_on_R2} is equal to the right-hand side of \eqref{eq:KN_XN}, and the identity \eqref{eq:KN_XN_on_R} follows.
\end{proof}

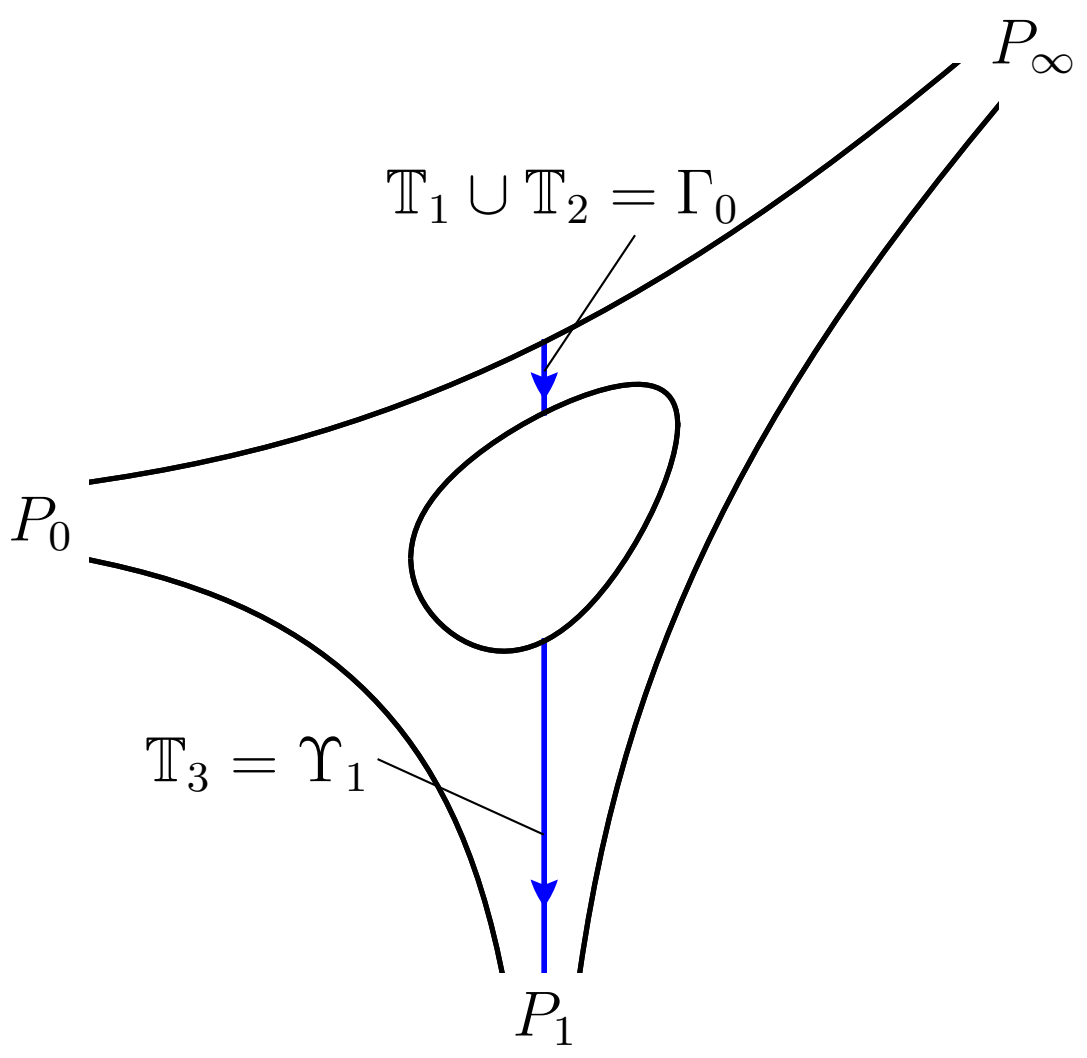
\begin{figure}[t]
    \centering
	\begin{overpic}[width=0.4\linewidth]{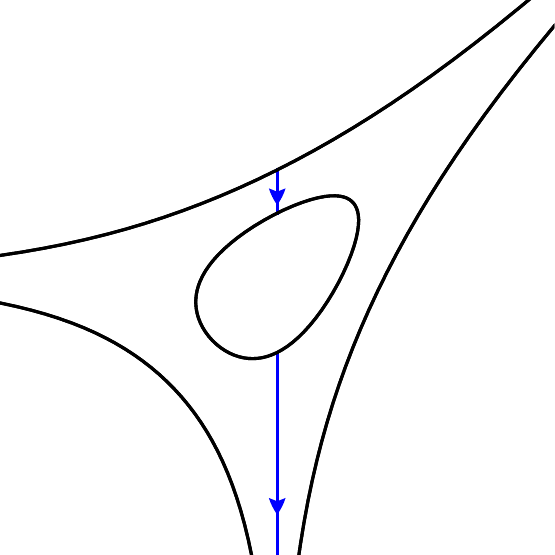}
		\put(-9,47.5){$P_0$}
		\put(46.5,-7){$P_1$}
		\put(99,100){$P_\infty$}
        \put(32.5,83.5){$\mathbb{T}_1\cup\mathbb{T}_2=\Gamma_0$}
		\put(47.6,66){
		\begin{tikzpicture}
			\draw[-] (0,0) -- (-0.6,-0.9);
		\end{tikzpicture}}
        \put(6,21){$\mathbb{T}_3=\Upsilon_1$}
        \put(29.25,15){
        \begin{tikzpicture}
        	\draw[-] (0,0) -- (1.1,-0.5);
        \end{tikzpicture}}
    \end{overpic}
    \vspace{0.4cm}
    \caption{The contours $z^{-1}(\mathbb{T})=\mathbb{T}_1\cup\mathbb{T}_2\cup\mathbb{T}_3=\Gamma_0\cup\Upsilon_1$ on the amoeba.}
    \label{fig:zT_contour}
\end{figure}

\begin{remark}\label{rem:B1_B3_poles}
    While $\mathbf{B}_N$ and $\mathbf{C}_N$ may have simple poles at $B_1$ and $B_3$, the integrals in \eqref{eq:KN_XN_on_R} do not have singularities at these points. This comes from the fact that $dz(p_1)$ and $dz(p_2)$ have simple zeros at $B_1$ and $B_3$, and therefore the meromorphic differentials $\mathbf{B}_N^t(p_2)dz(p_2)$ and $\mathbf{C}_N(p_1)dz(p_1)$ are analytic at $B_1$ and $B_3$.

    The same applies to $\mathbf{A}_N(p_1)dz(p_1)$, where $\mathbf{A}_N$ will be introduced in \eqref{eq:AN} below.
\end{remark}

We visualize the contour $z^{-1}(\mathbb{T})$ on the amoeba as in Figure~\ref{fig:zT_contour}. In what follows, we identify $\mathcal{R}_+$ with the interior of the amoeba. Recall that the two are homeomorphic via the amoeba map. Contours that are shown on the amoeba are taken to be symmetric with respect to complex conjugation. We recall that $z^{-1}(\mathbb{T})$ consists of two disjoint closed contours, namely $\mathbb{T}_1\cup\mathbb{T}_2$ and $\mathbb{T}_3$. Previously, these contours were denoted by $\Gamma_0$ and $\Upsilon_1$, respectively, see \eqref{eq:Gammas_def}--\eqref{eq:Gamma_Upsilons_char}. In what follows, we retain the notation $\mathbb{T}_1\cup\mathbb{T}_2$ and $\mathbb{T}_3$. Each circle $\mathbb{T}_j$ is equipped with the positive (counterclockwise) orientation. On the amoeba, the contour $z^{-1}(\mathbb{T})$ is oriented downwards, as illustrated in Figure~\ref{fig:zT_contour}.


\subsection{Deformation of the $p_2$-integral for $\xi_1>0$}\label{sec:deform_p2}
By analyticity, the contours of integration in \eqref{eq:KN_XN_on_R} may be deformed away from $z^{-1}(\mathbb{T})$. We want to choose the deformed contours in a way that allows us to identify the leading-order contribution to \eqref{eq:KN_XN_on_R} as $N\to\infty$.

The precise deformation depends on $(\xi_{1,N},\xi_{2,N})\in\mathfrak{L}$ and, in particular, on the level set $\mathcal{N}_{(\xi_{1,N},\xi_{2,N})}$, which was studied in Section~\ref{sec:critical_level_sets}. We assume first that $\xi_1>0$, and that $N$ is sufficiently large so that also $\xi_{1,N}>0$. The possible configurations of $\mathcal{N}_{(\xi_{1,N},\xi_{2,N})}$, are shown in Figures~\ref{fig:critical_level_sets_1} and~\ref{fig:critical_level_sets_4}. In our discussion, we focus on the middle subfigure of Figure~\ref{fig:critical_level_sets_1} and the left and right subfigures of Figure~\ref{fig:critical_level_sets_4}, but the reasoning for the other cases is similar.

\begin{figure}[t]
    \hspace{0.3cm}
	\begin{overpic}[width=0.3\linewidth]{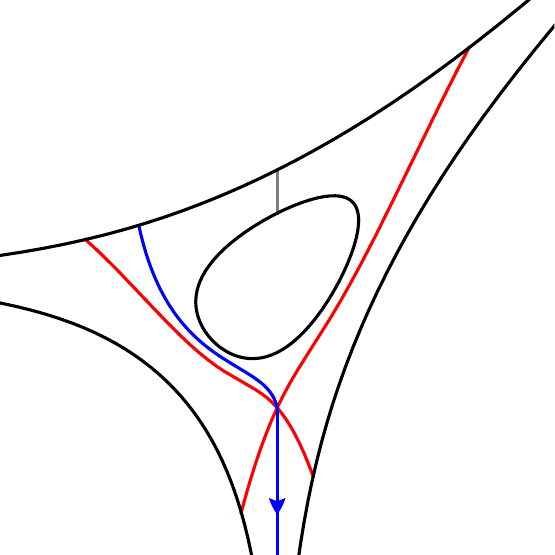}	
		\put(-11,47){$P_0$}
		\put(46,-8){$P_1$}
        \put(22.75,55){
        \begin{tikzpicture}
            \draw[-] (0,0) -- (-0.3,-0.7);
        \end{tikzpicture}}
		\put(28,74){$\gamma_B$}
	\end{overpic}
    \hfill
	\begin{overpic}[width=0.3\linewidth]{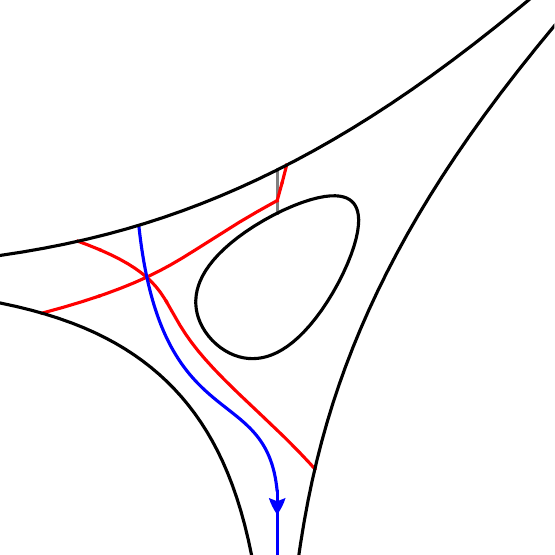}	
        \put(-10,99){$P_\infty$}
		\put(-11,47){$P_0$}
		\put(46,-8){$P_1$}
        \put(22.25,55){
        \begin{tikzpicture}
            \draw[-] (0,0) -- (-0.3,-0.7);
        \end{tikzpicture}}
		\put(27.5,74){$\gamma_B$}
	\end{overpic}
    \hfill
	\begin{overpic}[width=0.3\linewidth]{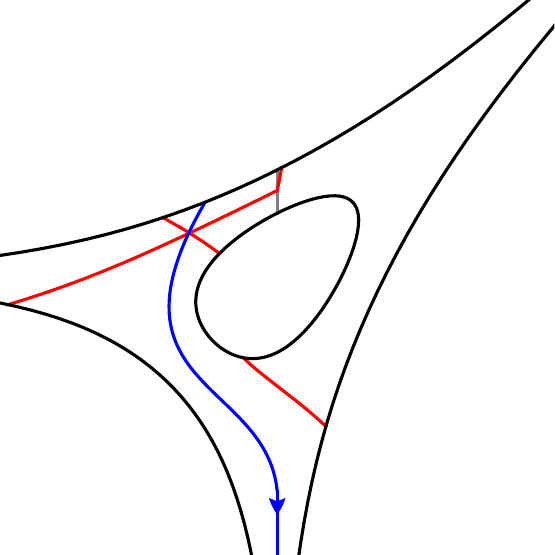}	
        \put(-10,99){$P_\infty$}
		\put(-11,47){$P_0$}
		\put(46,-8){$P_1$}
		\put(98,99){$P_\infty$}
        \put(25.25,60){
        \begin{tikzpicture}
            \draw[-] (0,0) -- (0.3,-0.7);
        \end{tikzpicture}}
		\put(23.5,79){$\gamma_B$}
	\end{overpic}
    \hspace{0.2cm}
    \vspace{0.25cm}
    \caption{The $p_2$-integral is deformed from $z^{-1}(\mathbb{T})$ to $\gamma_B$. The contour $\gamma_B$ lies in $\mathcal{N}^<_{(\xi_{1,N},\xi_{2,N})}$, except for the two saddle points $s_N$ and $\bar{s}_N$. The figure shows $\gamma_B$ on the amoeba for various cases with $\xi_{1,N}>0$.}
    \label{fig:deform_B}
\end{figure}
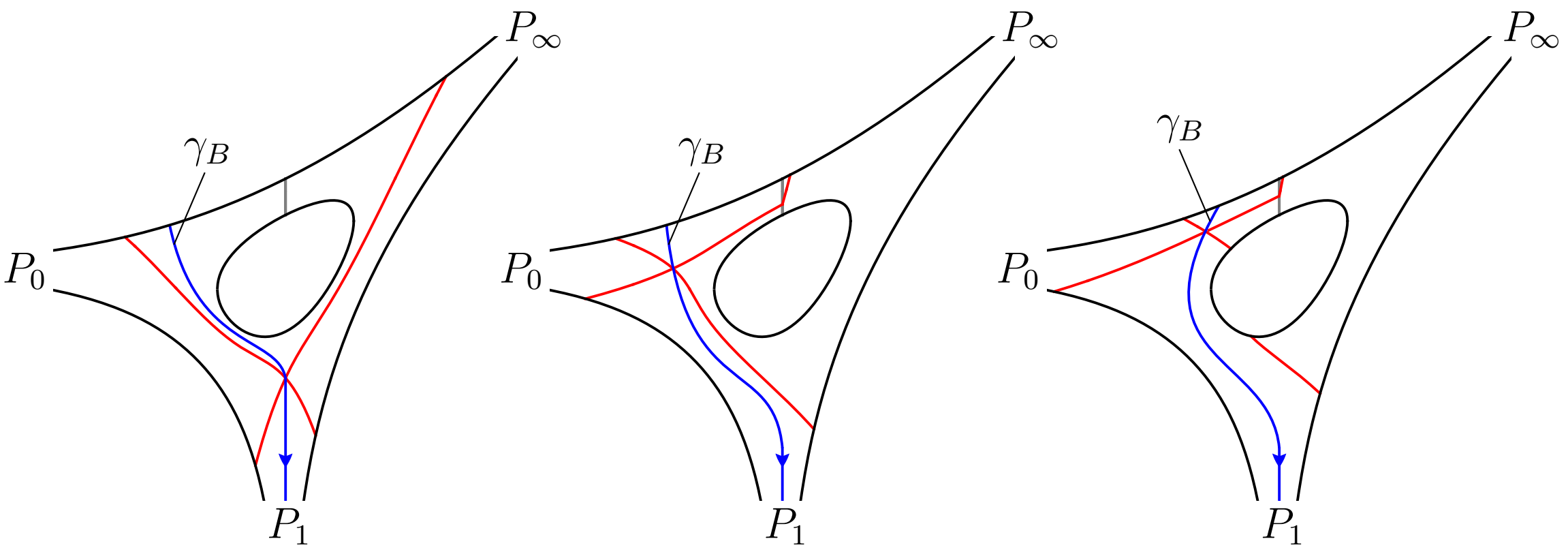

We deform $z^{-1}(\mathbb{T})$ for the $p_2$-integral in \eqref{eq:KN_XN_on_R} to one closed contour $\gamma_B$ in $\mathcal{R}$ that passes through $P_1 \in \mathcal{N}^<_{(\xi_{1,N}, \xi_{2,N})}$,  the saddle points 
\[
    s_N = \Omega(\xi_{1,N},\xi_{2,N})
\]
and $\bar{s}_N$, and that is otherwise in $\mathcal{N}^<_{(\xi_{1,N},\xi_{2,N})}$. Locally, near the saddle points, we let $\gamma_B$ follow the paths of steepest descent, and we also take $\gamma_B$ to be symmetric under complex conjugation. Inspection of the various cases shows that it is always possible to deform to $\gamma_B$ in this way, provided that $\xi_{1,N}>0$. A number of representative cases are given in Figure~\ref{fig:deform_B}. The figure shows the image of $\gamma_B$ on $\mathcal{R}_+$ on the amoeba, which then determines the entire curve due to the symmetry under complex conjugation. In some cases, $\gamma_B$ may also be chosen to pass around the other side of the bounded oval, while still going from $\bar{s}_N$ to $s_N$ within $\mathcal{N}^<_{(\xi_{1,N},\xi_{2,N})}$. Note that $\gamma_B$ depends on $N$, although we will not indicate the $N$-dependence in the notation.

Because of the analyticity properties of $\mathbf{B}_N$, see Lemma~\ref{lem:BNCN_poles} and Remark~\ref{rem:B1_B3_poles}, as well as of the other factors containing $p_2$ in \eqref{eq:KN_XN_on_R}, we are allowed to deform to the contour $\gamma_B$ and obtain by Cauchy's theorem:
\begin{multline}\label{eq:KN_BNCN} 
	\mathsf{K}_N(v_{1,N},v_{2,N}) = \frac{1}{(2\pi i)^2}\oint\limits_{\gamma_B} dz(p_2) \oint\limits_{z^{-1}(\mathbb{T})}  \frac{dz(p_1)}{z(p_1)(z(p_1)-z(p_2))} \\ \times
	\frac{\lambda(p_2)^{-x_2}}{z(p_2)^{-y_2}}
	\frac{\lambda(p_1)^{x_1}}{z(p_1)^{y_1}} 
	T_{0,j_2}^{-1}(z(p_2)) \mathbf{B}^t_N(p_2)
	\mathbf{C}_N(p_1) T_{0,j_1}(z(p_1)).
\end{multline}

We cannot deform the $p_1$-integral in a similar fashion.

\subsection{Splitting and deformation of the $p_1$-integral for $\xi_1>0$}\label{sec:split_p1}
To handle the $p_1$-integral, we introduce another $6\times3$ matrix-valued function
\begin{equation}\label{eq:AN}
	\mathbf{A}_N(p) = \mathbf{X}_N(z(p)) \begin{pmatrix} 0_3 \\ I_3 \end{pmatrix} \Xi(p) E^{-1}(z(p))\frac{\lambda(p)^{-N(1-\xi_{1,N})}}{z(p)^{-N(1-\xi_{2,N})}}.
\end{equation}

\begin{lemma}\label{lem:AN_poles}\leavevmode
    \begin{enumerate}
        \item[\rm{(a)}] $\mathbf{A}_N$ is meromorphic on $\mathcal{R}\setminus z^{-1}(\mathbb{T})$ with higher-order poles at $P_1$ and $P_\infty$, a higher-order zero at $P_0$, and at most simple poles at $B_1$ and $B_3$.
        \item[\rm{(b)}] One has
        \begin{equation}\label{eq:CN_split}
            \mathbf{C}_N(p) = \mathbf{A}_{N,+}(p) - \mathbf{A}_{N,-}(p), \quad p\in z^{-1}(\mathbb{T}).
        \end{equation}
    \end{enumerate}
\end{lemma}

\begin{proof}
    (a) This is similar to the proof of Lemma~\ref{lem:BNCN_poles}, except that there is no analyticity across $z^{-1}(\mathbb{T})$. The zero at $P_0$ and the pole at $P_1$ come from the different signs in the exponents of $z(p)$ and $\lambda(p)$ in \eqref{eq:AN}, when compared to \eqref{eq:BN} and \eqref{eq:CN}.

    \medskip

    (b) We use the jump condition $\mathbf{X}_{N,+}(z) = \mathbf{X}_{N,-}(z) J_{\mathbf{X}_N}(z)$  with $J_{\mathbf{X}_N}(z)$ as in \eqref{eq:jump_X} for $z \in \mathbb{T}$. From this, it follows that
    \[
        \mathbf{X}(z) \begin{pmatrix} I_3 \\ 0_3 \end{pmatrix} \frac{\Lambda(z)^{2N}}{z^{2N}} = \left(\mathbf{X}_{N,+}(z) - \mathbf{X}_{N,-}(z)\right) \begin{pmatrix} 0_3 \\ I_3 \end{pmatrix}, \quad z\in\mathbb{T}.
    \]
    Inserting this into \eqref{eq:CN} and taking note of \eqref{eq:AN}, we obtain \eqref{eq:CN_split}. Here, we also used that $\Lambda(z(p))^{2N}\Xi(p)=\Xi(p)\lambda(p)^{2N}$.
\end{proof}

We apply the formula~\eqref{eq:CN_split} to only a part of $z^{-1}(\mathbb{T})$, and then use the resulting expression in \eqref{eq:KN_XN_on_R}. Which part of $z^{-1}(\mathbb{T})$ depends on the topological structure of the level set $\mathcal{N}_{(\xi_{1,N},\xi_{2,N})}$. We split the $p_1$-integral according to the two terms that come from \eqref{eq:CN_split}. The term involving $\mathbf{X}_{N,+}$ is deformed inwards, towards $P_0$ (which means to the left in the amoeba), and the term involving $\mathbf{X}_{N,-}$ is deformed outwards
(to the right in the amoeba).

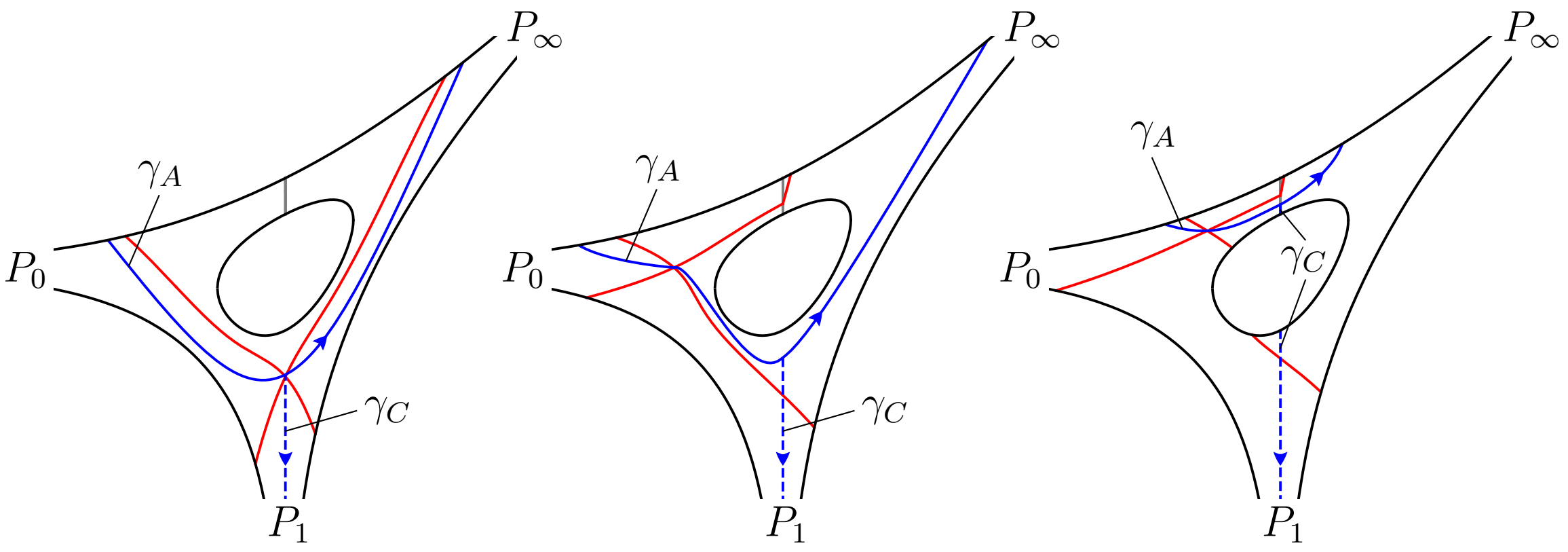
\begin{figure}[t]
    \hspace{0.25cm}
	\begin{overpic}[width=0.3\linewidth]{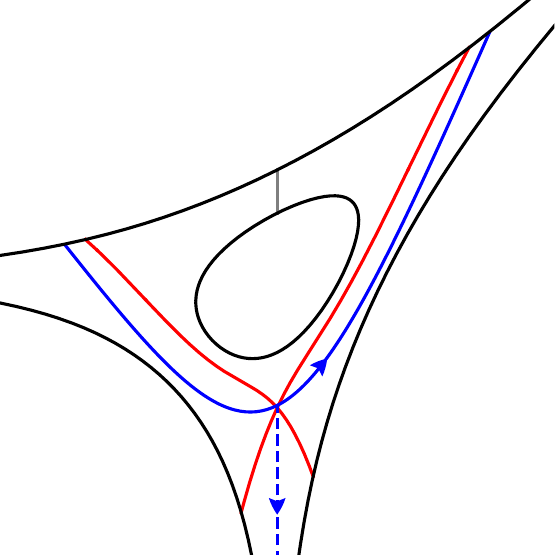}	
		\put(-11,47){$P_0$}
		\put(46,-8){$P_1$}
        \put(12.85,50){
        \begin{tikzpicture}
            \draw[-] (0,0) -- (-0.3,-0.7);
        \end{tikzpicture}}
		\put(18,69){$\gamma_A$}
        \put(46.75,14.5){
        \begin{tikzpicture}
        	\draw[-] (0,0) -- (-0.7,-0.2);
        \end{tikzpicture}}
        \put(67,18){$\gamma_C$}
	\end{overpic}
    \hfill
	\begin{overpic}[width=0.3\linewidth]{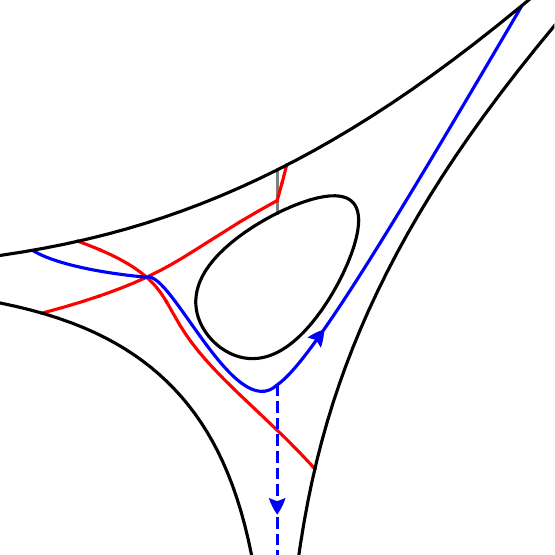}	
        \put(-10,99){$P_\infty$}
		\put(-11,47){$P_0$}
		\put(46,-8){$P_1$}
        \put(12.85,51.25){
        \begin{tikzpicture}
            \draw[-] (0,0) -- (-0.3,-0.7);
        \end{tikzpicture}}
		\put(18,70.25){$\gamma_A$}
        \put(46.75,14.5){
        \begin{tikzpicture}
        	\draw[-] (0,0) -- (-0.7,-0.2);
        \end{tikzpicture}}
        \put(67,18){$\gamma_C$}
	\end{overpic}
    \hfill
	\begin{overpic}[width=0.3\linewidth]{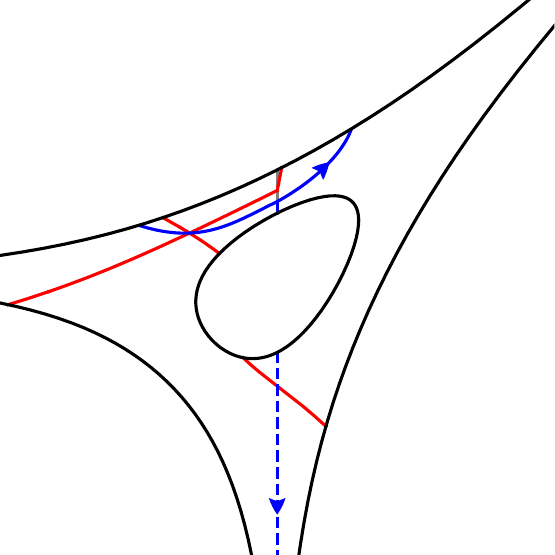}	
        \put(-10,99){$P_\infty$}
		\put(-11,47){$P_0$}
		\put(46,-8){$P_1$}
		\put(98,99){$P_\infty$}
        \put(19,58.15){
        \begin{tikzpicture}
            \draw[-] (0,0) -- (0.3,-0.7);
        \end{tikzpicture}}
		\put(17.5,77.5){$\gamma_A$}
        \put(46.7,55.75){
            \begin{tikzpicture}
        		\draw[-] (0,0) -- (-0.25,0.3);
            \end{tikzpicture}}
        \put(46.7,32.2){
            \begin{tikzpicture}
        		\draw[-] (0,0) -- (-0.25,-0.68);
            \end{tikzpicture}}
		\put(50,50.5){$\gamma_C$}
	\end{overpic}
    \hspace{0.2cm}
    \vspace{0.25cm}
	\caption{After splitting part of the $p_1$-integral, part of $z^{-1}(\mathbb{T})$ is deformed to $\gamma_A$, and the remainder is denoted by $\gamma_C$. The contour $\gamma_A$ lies in $\mathcal{N}^>_{(\xi_{1,N},\xi_{2,N})}$, except for the two saddle points $s_N$ and $\bar{s}_N$. The figure shows $\gamma_A$ on the amoeba for various cases with $\xi_{1,N}>0$. In the third subfigure, $\gamma_C$ consists of two pieces.}
    \label{fig:deform_AC}
\end{figure}

Part of the $p_1$-integral can then be deformed to a closed contour $\gamma_A$ as shown in Figure~\ref{fig:deform_AC} for a number of cases. The contour $\gamma_A$ passes through the saddle points $s_N = \Omega(\xi_{1,N},\xi_{2,N})$ and $\bar{s}_N$, while otherwise remaining in the region $\mathcal{N}^>_{(\xi_{1,N},\xi_{2,N})}$. Locally, near the saddle points, $\gamma_A$ follows the paths of steepest ascent, and $\gamma_A$ is also symmetric under complex conjugation. The remaining part of $z^{-1}(\mathbb{T})$ that is not deformed is called $\gamma_C$. In the last subfigure, $\gamma_C$ is equal to the union of $\mathbb{T}_3$ with a small subset of $\mathbb{T}_1\cup\mathbb{T}_2$, and in the other two subfigures it is a strict subset of $\mathbb{T}_3$.

While deforming to $\gamma_A$, we must cross part of $\gamma_B$. We denote the part that is crossed by $\gamma_{s_N}$. It is precisely the part of $\gamma_B$ that goes from $\bar{s}_N$ to $s_N$ and intersects the bounded oval between $P_\infty$ and $P_0$. In particular, it does not contain $P_1$. By analyticity, we have the freedom to deform the contour in an integral such as \eqref{eq:KN_sN}, as long as it starts at $\bar{s}_N$, ends at $s_N$, and passes through the unbounded oval in between $P_\infty$ and $P_0$.

In this deformation process, we pick up a residue contribution that will give the leading-order contribution to $\mathsf{K}_N(v_{1,N},v_{2,N})$ as $N\to\infty$. To explain the appearance of the Riesz projector $\Pi$, as defined in \eqref{eq:Riesz_projector}, in this leading-order contribution, we make the following observation. We first recall the spectral decomposition $W(z) = E(z)\Lambda(z)E^{-1}(z)$ of $W(z)$ that was used in the first transformation $\mathbf{Y}_N \mapsto \mathbf{X}_N$. Then, 
\begin{equation}\label{eq:Riesz_projector_E_Xi_E_inv}
    \Pi(p) = E(z)\Xi(p)E^{-1}(z) \quad \text{for } p=(z,\lambda),
\end{equation}
with $\Xi$ as in \eqref{eq:Xi_def}.

\begin{lemma}\label{lem:KN_split}
    After applying the above contour deformations, for $(\xi_1,\xi_2)\in\mathfrak{L}$ with $\xi_1>0$, one obtains
    \begin{equation}\label{eq:KN_ANBNCN}
    	\mathsf{K}_N(v_{1,N},v_{2,N}) = \mathsf{K}_N^{(B,C)}(v_{1,N},v_{2,N}) + \mathsf{K}_N^{(B,A)}(v_{1,N},v_{2,N}) + \mathsf{K}_N^{(s_N)}(v_{1,N},v_{2,N}),
    \end{equation}
    with
    \begin{align}
    	\mathsf{K}_N^{(B,C)}(v_{1,N},v_{2,N}) &= \frac{1}{(2\pi i)^2}\oint\limits_{\gamma_B} dz(p_2)\oint\limits_{\gamma_C}\frac{dz(p_1)}{z(p_1)(z(p_1)-z(p_2))} \nonumber \\
    	&\times\frac{\lambda(p_2)^{-x_2}}{z(p_2)^{-y_2}}\frac{\lambda(p_1)^{x_1}}{z(p_1)^{y_1}} T_{0,j_2}^{-1}(z(p_2))\mathbf{B}^t_N(p_2)\mathbf{C}_N(p_1) T_{0,j_1}(z(p_1)), \label{eq:KN_BC} \\
    	\mathsf{K}_N^{(B,A)}(v_{1,N},v_{2,N}) &= \frac{1}{(2\pi i)^2}\oint\limits_{\gamma_B} dz(p_2) \oint\limits_{\gamma_A} \frac{dz(p_1)}{z(p_1)(z(p_1)-z(p_2))} \nonumber \\
        &\times\frac{\lambda(p_2)^{-x_2}}{z(p_2)^{-y_2}}\frac{\lambda(p_1)^{x_1}}{z(p_1)^{y_1}} T_{0,j_2}^{-1}(z(p_2)) \mathbf{B}^t_N(p_2)\mathbf{A}_N(p_1) T_{0,j_1}(z(p_1)), \label{eq:KN_BA} \\
    	\mathsf{K}_N^{(s_N)}(v_{1,N},v_{2,N}) &= \frac{1}{2\pi i}\int_{\gamma_{s_N}}\frac{\lambda(p)^{x_1-x_2}}{z(p)^{y_1-y_2}}T_{0,j_2}^{-1}(z(p))\Pi(p)T_{0,j_1}(z(p_1))\frac{dz(p)}{z(p)}. \label{eq:KN_sN}
	\end{align}
\end{lemma}

\begin{proof}
    The first term \eqref{eq:KN_BC} is obtained by restricting the $p_1$-integral in \eqref{eq:KN_BNCN} to $\gamma_C$. 
    
    In the remaining part of \eqref{eq:KN_BNCN}, the $p_1$-integral is over $z^{-1}(\mathbb{T}) \setminus \gamma_C$. To this part, we apply \eqref{eq:CN_split} to obtain an integral with $A_{N,+}(p_1)$ and an integral with $-A_{N,-}(p_1)$. We invert the orientation of the latter integral, and then deform the two $p_1$-integrals away from $z^{-1}(\mathbb{T})$. Together, they recombine to give an integral over $\gamma_A$, see Figure~\ref{fig:deform_AC}, and we arrive at the second term \eqref{eq:KN_BA}.

    In the process of deforming $z^{-1}(\mathbb{T}) \setminus \gamma_C$ to $\gamma_A$, we have crossed the part $\gamma_{s_N}$ of $\gamma_B$. For each $p\in\gamma_{s_N}$, we pick up the following residue contribution:
    \[
        \frac{\lambda(p)^{x_1-x_2}}{z(p)^{y_1-y_2}}T_{0,j_2}^{-1}(z(p)) \mathbf{B}_N^t(p)\mathbf{A}_N(p) T_{0,j_1}(z(p)) \frac{1}{z(p)}.
    \]
    The product $\mathbf{B}_N^t(p)\mathbf{A}_N(p)$ simplifies considerably. Using \eqref{eq:BN} and \eqref{eq:AN}, we find 
    \[
        \mathbf{B}_N^t(p)\mathbf{A}_N(p) = E(z(p)) \Xi(p)E^{-1}(z(p)) = \Pi(p),
    \]
    where the second equality is precisely the identity \eqref{eq:Riesz_projector_E_Xi_E_inv}. Here, we also used that $\Xi(p)^2 = \Xi(p)$, which is clear from \eqref{eq:Xi_def}. These residues combine into an integral over $\gamma_{s_N}$, and this gives rise to the third term \eqref{eq:KN_sN}.
\end{proof}

The single integral \eqref{eq:KN_sN} depends only on $N$ via the saddle point $s_N$. Then, it is easy to calculate the limit as $N\to\infty$.

\begin{lemma}\label{lem:KN_limit}
    For $(\xi_1,\xi_2) \in \mathfrak{L}$ with $\xi_1 > 0$, one has
    \begin{equation}\label{eq:KN_sN_limit}
    	\lim_{N\to\infty}\mathsf{K}_N^{(s_N)}(v_{1,N},v_{2,N}) = \frac{1}{2\pi i}\int_{\gamma_s^\leq}\frac{\lambda(p)^{x_1-x_2}}{z(p)^{y_1-y_2}}T_{0,j_2}^{-1}(z(p))\Pi(p)T_{0,j_1}(z(p_1))\frac{dz(p)}{z(p)},
    \end{equation}
    where $\gamma_s^\leq$ is as in Theorem~\ref{thm:local_scaling_limit}, see also Figure~\ref{fig:torus_KN_contours}.
\end{lemma}

\begin{proof}
    As $N\to\infty$, we have $s_N \to s = \Omega(\xi_1,\xi_2)$ and $\gamma_{s_N} \to \gamma_s$. By analyticity, the contour $\gamma_s$ can be deformed to any contour $\gamma_s^\leq$ from Theorem~\ref{thm:local_scaling_limit}.
\end{proof}

The limit \eqref{eq:KN_sN_limit} will give the leading-order contribution to $\mathsf{K}_N(v_{1,N},v_{2,N})$ as $N\to\infty$, as the two terms \eqref{eq:KN_BC} and \eqref{eq:KN_BA} with the double contour integrals tend to $0$ as $N\to\infty$. This is shown in Lemma~\ref{lem:KN_small}. In order to prove this, we apply the second transformation of the RH problem, and use the outcome of the steepest descent analysis of the RH problem from \cite{Kui25}.


\subsection{Second transformation $\mathbf{X}_N \mapsto \mathbf{T}_N$}\label{sec:second_transformation}
For the second transformation of the RH problem, we follow \cite[Definition~6.5]{Kui25}, however, with a modification. If we would define $\mathbf{T}_N$ exactly as in \cite[Definition~6.5]{Kui25}, then by \cite[Lemma~6.7]{Kui25}, it would asymptotically behave as 
\[
    \left(I_6 + \mathcal{O}(z^{-1})\right) \begin{pmatrix} L^N E(z) & 0_3 \\ 0_3 & L^{-N} E(z) \end{pmatrix} \quad \text{as } z \to \infty,
\]
for a certain unit lower triangular matrix $L$ of size $3 \times 3$. It is convenient for us to modify the definition of $\mathbf{T}_N$ by including the factor $\bigl(\!\begin{smallmatrix} L^{-N} & 0_3 \\ 0_3 & L^N \end{smallmatrix}\!\bigr)$, which has the effect\vspace{-0.1cm} that $\mathbf{T}_N\bigl(\!\begin{smallmatrix} E^{-1} & 0_3 \\ 0_3 & E^{-1} \end{smallmatrix}\!\bigr)$ remains bounded at infinity. 

We thus define
\begin{multline}\label{eq:TN_def}
	\mathbf{T}_N = \begin{pmatrix} e^{N \widehat{\ell}} L^{-N} & 0_3 \\ 0_3 & e^{-N \widehat{\ell}} L^N \end{pmatrix} \mathbf{X}_N \diag\left(e^{-Ng_1}, e^{-Ng_2}, e^{-Ng_3}, e^{Ng_1}, e^{Ng_2}, e^{Ng_3}\right) \\
	\times \begin{pmatrix} e^{-N \widehat{\ell}} I_3 & 0_3 \\ 0_3 & e^{N\widehat{\ell}} I_3 \end{pmatrix}.
\end{multline} 
Here, $g_1$, $g_2$, and $g_3$ denote the restrictions of the $g$-function, see Section~\ref{sec:g_function}, to the three sheets of the Riemann surface and $\widehat{\ell}$ is the constant from \eqref{eq:S_property_2}.

We recall the definition~\eqref{eq:Phi_xi1_xi2_def} of $\Phi_{(\xi_1,\xi_2)}$, and use it with parameters $(\xi_{1,N},\xi_{2,N})$ to define
\begin{align}
	\Phi_N :\!&= \Phi_{(\xi_{1,N},\xi_{2,N})}\label{eq:PhiN_def} \\ 
	&=
    \begin{cases}
        g+(1-\xi_{1,N})\log\lambda - (1-\xi_{2,N})\log z + \widehat{\ell}, & \text{on } \mathcal{R}, \\	
    	-g-(1+\xi_{1,N})\log\lambda + (1 +\xi_{2,N})\log z - \widehat{\ell} + \pi i, & \text{on } \mathcal{R}^\sigma.
	\end{cases} \nonumber
\end{align}
Using \eqref{eq:TN_def} and \eqref{eq:PhiN_def} in \eqref{eq:BN}--\eqref{eq:CN} and \eqref{eq:AN}, we obtain
\begin{align}
	\mathbf{A}_N(p) &= \begin{pmatrix} e^{-N\widehat{\ell}} L^{-N} & 0_3 \\ 0_3 & e^{N\widehat{\ell}} L^N \end{pmatrix} \mathbf{T}_N(z(p)) \begin{pmatrix} 0_3 \\ I_3 \end{pmatrix}\Xi(p)E^{-1}(z(p))e^{-N\Phi_N(p)}, \label{eq:AN_TN} \\
	\mathbf{B}_N(p) &= \begin{pmatrix} e^{N\widehat{\ell}} (L^t)^{-N} & 0_3 \\ 0_3 & e^{-N \widehat{\ell}} (L^t)^N \end{pmatrix} \mathbf{T}_N^{-t}(z(p)) \begin{pmatrix} 0_3 \\ I_3 \end{pmatrix}\Xi(p)E^{t}(z(p))e^{N\Phi_N(p)}, \label{eq:BN_TN} \\
	\mathbf{C}_N(p) &= (-1)^N\begin{pmatrix} e^{-N\widehat{\ell}} L^{-N} & 0_3 \\ 0_3 & e^{N\widehat{\ell}} L^N \end{pmatrix} \mathbf{T}_N(z(p)) \begin{pmatrix} I_3 \\ 0_3 \end{pmatrix}\Xi(p)E^{-1}(z(p))e^{-N\Phi_N(\sigma(p))}\label{eq:CN_TN}
\end{align}
for $p\in\mathcal{R}$. We emphasize that $\Phi_N$ in \eqref{eq:CN_TN} is not evaluated at $p$, but at $\sigma(p)$, which is on the second copy $\mathcal{R}^\sigma$.
 
The prefactors that contain $L$ and $L^t$ cancel out in the products $\mathbf{B}_N^t(p_2)\mathbf{A}_N(p_1)$ and $\mathbf{B}_N^t(p_2)\mathbf{C}_N(p_1)$ that we have in \eqref{eq:KN_BC} and \eqref{eq:KN_BA}. These products have the scalar factors $e^{N\Phi_N(p_2)} e^{-N\Phi_N(p_1)}$ and $e^{N\Phi_N(p_2)} e^{-N\Phi_N(\sigma(p_1))}$, respectively, that we need to control in the large-$N$ limit. We do this in the proof of Lemma~\ref{lem:KN_small} below. The other factors remain bounded as $N\to\infty$, as we show in the following lemma.

\begin{lemma}\label{lem:TN_bounded}\leavevmode
    \begin{enumerate}
        \item[\rm{(a)}] The matrix-valued functions
        \[
            \mathbf{T}_N\begin{pmatrix} E^{-1} & 0_3 \\ 0_3 & E^{-1} \end{pmatrix} \quad \text{and} \quad \begin{pmatrix} E & 0_3 \\ 0_3 & E \end{pmatrix}\mathbf{T}_N^{-1}
        \]
        are uniformly bounded for $z$ in compact subsets of $\mathbb{C}\setminus\mathbb{T}$ and large $N$.
        \item[\rm{(b)}] The product
        \[
            E(z(p_2))\begin{pmatrix} 0_3 & I_3 \end{pmatrix} \mathbf{T}_N^{-1}(z(p_2)) \mathbf{T}_N(z(p_1))\begin{pmatrix} 0_3 \\ I_3 \end{pmatrix}E^{-1}(z(p_1))
        \]
        is uniformly bounded for $p_1\in\gamma_A\cup\gamma_C$, $p_2\in\gamma_B$, and large $N$.
        \item[\rm{(c)}] The product
        \[
            \frac{1}{z(p_1)-z(p_2)}E(z(p))\begin{pmatrix} 0_3 & I_3 \end{pmatrix} \mathbf{T}_N^{-1}(z(p_2)) \mathbf{T}_N(z(p_1))\begin{pmatrix} I_3 \\ 0_3 \end{pmatrix}E^{-1}(z(p))
        \]
        is uniformly bounded for $p_1\in\gamma_C$, $p_2\in\gamma_B$, and large $N$. 
    \end{enumerate}
\end{lemma}

\begin{proof}
    (a) This is a consequence of the steepest descent analysis of the RH problem from \cite{Kui25}. Let us now elaborate on this.

    The third transformation $\mathbf{T}_N \mapsto \mathbf{S}_N$ from \cite{Kui25} is based on the opening of an annular lens $\{z \mid 1 - \eta<|z|<1 + \eta\}$ around $\mathbb{T}$ for sufficiently small $\eta>0$ satisfying $1 \pm \eta\in(z_{\min},z_{\max})$. 

    Then, according to \cite[Definition~7.1]{Kui25}, we define $\mathbf{S}_N = \mathbf{T}_N$ outside the annular lens, and inside we have
    \begin{equation}\label{eq:SN}
        \mathbf{S}_N = \mathbf{T}_N \begin{pmatrix} I_3 & 0_3 \\ \pm(-1)^N \diag(e^{2N\varphi_1},e^{2N\varphi_2},0) & I_3 \end{pmatrix}
    \end{equation}
    with $\varphi_k(z) = -\widehat{g}_k(z) + \log z - \log\lambda_k(z) - \tfrac{1}{2}\widehat{\ell}$, see \cite[(6.20)--(6.23)]{Kui25}. Here, we take $+$ for $1 - \eta<|z|<1$ and $-$ for $1<|z|<1 + \eta$. Definition~7.1 in \cite{Kui25} includes the prefactor $\bigl(\!\begin{smallmatrix} L^{-N} & 0_3 \\ 0_3 & L^N \end{smallmatrix}\!\bigr)$. This prefactor is in our definition \eqref{eq:TN_def} of $\mathbf{T}_N$, which means that $\mathbf{S}_N$ defined in \eqref{eq:SN} agrees with the definition of $S$ given in \cite{Kui25}.

    It follows from \cite[Lemma~7.3]{Kui25} that for sufficiently small $\eta>0$, we have $\Re\varphi_1(z)\leq0$ and $\Re\varphi_2(z)\leq 0$ for all $1 - \eta<|z|<1 + \eta$. Thus, to prove part (a), it is enough to prove that 
    \begin{equation}\label{eq:SNE}
        \widehat{\mathbf{S}}_N = \mathbf{S}_N\begin{pmatrix} E^{-1} & 0_3 \\ 0_3 & E^{-1} \end{pmatrix}
    \end{equation} 
    and its inverse are uniformly bounded as $N\to\infty$. 

    The fourth and final transformation $\mathbf{S}_N \mapsto \mathbf{R}_N$ is given by $\mathbf{R}_N = \mathbf{S}_N\mathbf{M}_N^{-1}$, where $\mathbf{M}_N$ is the global parametrix, which is the solution to the model problem \cite[RH problem~8.1]{Kui25}. We are interested in $\widehat{\mathbf{M}}_N$, defined analogously to \eqref{eq:SNE}, which solves a related RH problem \cite[RH problem~8.2]{Kui25}. The benefit of $\widehat{\mathbf{M}}_N$ is that it is normalized at infinity and that $\det \widehat{\mathbf{M}}_N \equiv 1$, so that its inverse is always well-defined. Furthermore, $\widehat{\mathbf{M}}_N$ and $\widehat{\mathbf{M}}_N^{-1}$ are uniformly bounded in $\overline{\mathbb{C}}\setminus(\mathbb{T}\cup[z_{\min},z_{\max}])$. We note that $\widehat{\mathbf{S}}_N\widehat{\mathbf{M}}_N^{-1} = \mathbf{R}_N$. An application of the small-norm theorem yields
    \begin{equation}
        \mathbf{R}_N(z) = I_6 + \mathcal{O}\left(\frac{e^{-cN}}{1 + |z|}\right)
    \end{equation}
    for some $c>0$, which holds uniformly in compact subsets of $\mathbb{C}\setminus(\mathbb{T}\cup(1\pm\eta)\mathbb{T})$, see \cite[Corollary~9.3]{Kui25}. It follows that $\widehat{\mathbf{S}}_N = \mathbf{R}_N\widehat{\mathbf{M}}_N$ and its inverse are uniformly bounded as $N\to\infty$, which completes the proof of part (a).

    \medskip

    (b) This is immediate from part~(a).
    
    \medskip

    (c) This also follows from part~(a), as
    \[
        \begin{pmatrix} 0_3 & I_3 \end{pmatrix} \mathbf{T}_N^{-1}(z(p_2))\mathbf{T}_N(z(p_1)) \begin{pmatrix} I_3 \\ 0_3 \end{pmatrix}
    \]
    vanishes when $z(p_1)=z(p_2)$.
\end{proof}


\subsection{Estimation of \eqref{eq:KN_BC} and \eqref{eq:KN_BA} for $\xi_1>0$}\label{sec:estimates}
Now, we show that the two double contour integrals in \eqref{eq:KN_ANBNCN} tend to $0$ as $N\to\infty$.

\begin{lemma}\label{lem:KN_small}
    For $(\xi_1,\xi_2)\in\mathfrak{L}$ with $\xi_1>0$, the two double contour integrals $\mathsf{K}_N^{(B,C)}(v_{1,N},v_{2,N})$ and $\mathsf{K}_N^{(B,A)}(v_{1,N},v_{2,N})$, given by \eqref{eq:KN_BC} and \eqref{eq:KN_BA}, tend to $0$ as $N\to\infty$.
\end{lemma}

\begin{proof}
	By \eqref{eq:AN_TN} and \eqref{eq:BN_TN}, we have
    \begin{multline}\label{eq:BNt_AN_est}
        \mathbf{B}_N^t(p_2) \mathbf{A}_N(p_1) = e^{N (\Phi_N(p_2) - \Phi_N(p_1))} \\
		\times E(z(p_2)) \Xi(p_2) \begin{pmatrix} 0_3 & I_3 \end{pmatrix} \mathbf{T}_N^{-1}(z(p_2)) 
		\mathbf{T}_N(z(p_1)) \begin{pmatrix} 0_3 \\ I_3 \end{pmatrix} \Xi(p_1) E^{-1}(z(p_1)).
    \end{multline}
    The matrix-valued factors remain bounded for $p_1\in\gamma_A$ and $p_2\in\gamma_B$, as $N\to\infty$, due to Lemma~\ref{lem:TN_bounded}(b).

    For $p_1\in\gamma_A \setminus \{s_N, \bar{s}_N\}$ and $p_2\in\gamma_B \setminus \{s_N, \bar{s}_N\}$, we have $p_1\in\mathcal{N}^>_{(\xi_{1,N},\xi_{2,N})}$ and $p_2\in\mathcal{N}^<_{(\xi_{1,N},\xi_{2,N})}$, and thus by \eqref{eq:critical_level_set_regions_on_Rhat},
    \begin{equation}\label{eq:PhiN_p1p2_est}
        \Re\Phi_N(p_2) < C_{(\xi_{1,N},\xi_{2,N})} < \Re\Phi_N(p_1).
    \end{equation}
    It follows that  $|e^{N(\Phi_N(p_2)-\Phi_N(p_1))}|\leq1$ for $p_1 \in \gamma_A$, $p_2 \in \gamma_B$, with equality if and only if $\{p_1,p_2\} \subset \{s_N, \bar{s}_N\}$. Thus, in view of \eqref{eq:BNt_AN_est}, we see that $\mathbf{B}_N^t(p_2) \mathbf{A}_N(p_1)$. Note that $e^{N(\Phi_N(p_2) - \Phi_N(p_1))}$ is exponentially small as $N\to\infty$ if at least one of $p_1$ and $p_2$ is bounded away from the saddle points. 
   
    It follows that the dominant contribution to the double contour integral \eqref{eq:KN_BA} comes from the saddle points. Observe that $(z(p_1)-z(p_2))^{-1}$ becomes unbounded as $p_1$ and $p_2$ approach the same saddle point. With appropriate care it can be shown that the contribution from a small neighborhood of the saddle points is $\mathcal{O}(N^{-1/2})$ as $N\to\infty$, see e.g.~\cite[\S7.5.1]{CDKL20}. The contribution away from the saddle points is exponentially small, and it follows that the double contour integral \eqref{eq:KN_BA} tends to $0$ as $N\to\infty$.

    We now turn to \eqref{eq:KN_BC}. Using \eqref{eq:BN_TN} and \eqref{eq:CN_TN}, we get
    \begin{multline*}
        \mathbf{B}_N^t(p_2) \mathbf{C}_N(p_1) = (-1)^Ne^{N (\Phi_N(p_2) - \Phi_N(\sigma(p_1)))} \\
		\times E(z(p_2)) \Xi(p_2) \begin{pmatrix} 0_3 & I_3 \end{pmatrix} \mathbf{T}_N^{-1}(z(p_2)) 
		\mathbf{T}_N(z(p_1)) \begin{pmatrix} I_3 \\ 0_3 \end{pmatrix} \Xi(p_1) E^{-1}(z(p_1)).
    \end{multline*}
    The matrix-valued factors in $\mathbf{B}_N^t(p_2) \mathbf{C}_N(p_1)$ are of order $\mathcal{O}(z(p_1)-z(p_2))$ as $N\to\infty$, uniformly for $p_1\in\gamma_C$ and $p_2\in\gamma_B$, because of Lemma~\ref{lem:TN_bounded}(c). It is then enough to show that there is a constant $c>0$ such that 
    \begin{equation}\label{eq:PhiN_estimate_BNCN}
        \Re\Phi_N(\sigma(p_1)) - \Re\Phi_N(p_2) \geq c, \quad p_1 \in \gamma_C, \, p_2 \in \gamma_B.
    \end{equation}
    Indeed, then the scalar factor $e^{N(\Phi_N(p_2)-\Phi_N(\sigma(p_1)))}$ in $\mathbf{B}_N^t(p_2) \mathbf{C}_N(p_1)$ is exponentially small as $N\to\infty$ uniformly for $p_1 \in \gamma_C$ and $p_2 \in \gamma_B$, and \eqref{eq:KN_BC} tends to $0$ as $N\to\infty$ at an exponential rate.

    To obtain \eqref{eq:PhiN_estimate_BNCN}, we need to estimate $\Re\Phi_N(\sigma(p_1))$ and $\Re\Phi_N(p_2)$. Since $p_2 \in \gamma_B$, and $\gamma_B \setminus \{s_N, \bar{s}_N\}$ lies in $\mathcal{N}_{(\xi_{1,N}, \xi_{2,N})}^<$, we have
    \begin{equation}\label{eq:PhiN_p2_est}
        \Re\Phi_N(p_2) \leq C_{(\xi_{1,N},\xi_{2,N})}, \quad  p_2 \in \gamma_B.
    \end{equation}

    The estimate on $\Re\Phi_N(\sigma(p_1))$ requires a more elaborate argument. We begin by noting that $\log|z(p)| = 0$ for $p \in z^{-1}(\mathbb{T}) = \Gamma_0 \cup \Upsilon_1$, see also Figure~\ref{fig:zT_contour}. Therefore, by \eqref{eq:Phi_xi1_xi2_def},
    \begin{equation}\label{eq:RePhi_on_zT}
        \Re\Phi_N(\sigma(p)) = \Re\Phi(\sigma(p)) - \xi_{1,N}\log|\lambda(p)|, \quad p \in \Gamma_0 \cup \Upsilon_1.
    \end{equation}
    We use the orientation on $(\Gamma_0 \cup \Upsilon_1) \cap \mathcal{R}_+$ as in Figure~\ref{fig:zT_contour}. Then, the second term $-\xi_{1,N} \log|\lambda(p)|$ in \eqref{eq:RePhi_on_zT} is strictly increasing on $(\Gamma_0 \cup \Upsilon_1) \cap \mathcal{R}_+$, as $\xi_{1,N} > 0$. The first term $\Re\Phi(\sigma(p))$ is also increasing, which is due to Lemma~\ref{lem:Phi_00_on_Gammas_Upsilons}(a,b).

    Let $p_*\in\mathcal{R}_+$ be the intersection of $\gamma_A$ with $z^{-1}(\mathbb{T})$. This is also the point of intersection of $\gamma_A$ with $\gamma_C$, see Figure~\ref{fig:deform_AC}. Then, by the increasing behavior of \eqref{eq:RePhi_on_zT}, we have
    \begin{equation}\label{eq:PhiN_sigma_p1_est}
    	\Re\Phi_N(\sigma(p_1)) \geq \Re\Phi_N(\sigma(p_*)), \quad p_1\in\gamma_C.
    \end{equation}
    Because of \eqref{eq:Phi_xi1_xi2_def}, we have
    \begin{equation}\label{eq:PhiN_sigma_pstar_pstar} 
    	\Re\Phi_N(\sigma(p_*)) = \Re\Phi_N(p_*) - 2\Re\Phi(p_*).
    \end{equation}
    Combining \eqref{eq:PhiN_p2_est}, \eqref{eq:PhiN_sigma_p1_est}, and \eqref{eq:PhiN_sigma_pstar_pstar}, we obtain \eqref{eq:PhiN_estimate_BNCN} with constant $c = c_1 + c_2$, where $c_1 = \Re\Phi_N(p_*) - C_{(\xi_{1,N}, \xi_{2,N})}$ and $c_2 = -2\Re\Phi(p_*)$. It remains to show that $c > 0$.

    We have $c_1 \geq 0$, since $p_* \in \gamma_A$, and $\gamma_A \setminus \{s_N, \bar{s}_N\} \subset \mathcal{N}^>_{\xi_{1,N} ,\xi_{2,N}}$. To see that $c_2 \geq 0$ as well, we recall that $p_* \in \Gamma_0 \cup \Upsilon_1$ and note that it then follows from  Lemma~\ref{lem:Phi_00_on_Gammas_Upsilons}(a) in case that $p_* \in \Gamma_0$ and Lemma~\ref{lem:Phi_00_on_Gammas_Upsilons}(c) in case that $p_* \in \Upsilon_1$. Moreover, if $p_* \in \Upsilon_1$, then we in fact obtain that $c_2>0$. Note that $c_1 = 0$ can only occur when $p_* = s_N$, as is the case in the left subfigure of Figure~\ref{fig:deform_AC}. However, then it must be the case that $p_* \in \Upsilon_1$, so that $c_2>0$. We conclude that $c=c_1+c_2>0$, which completes the proof.
\end{proof}

Having established Lemmas~\ref{lem:KN_limit} and~\ref{lem:KN_small}, we are ready to prove Theorem~\ref{thm:local_scaling_limit} for the case $\xi_1>0$.


\subsection{Proof of Theorem~\ref{thm:local_scaling_limit} for $\xi_1>0$}
\begin{proof}
    Combining Lemmas~\ref{lem:KN_limit} and~\ref{lem:KN_small} with \eqref{eq:KN_ANBNCN}, we find that $\mathsf{K}_N(v_{1,N}, v_{2,N})$ has the limit
    \begin{equation}\label{eq:KN_limit}
        \lim_{N\to\infty}\mathsf{K}_N(v_{1,N},v_{2,N}) = \frac{1}{2\pi i}\int_{\gamma_s^\leq}\frac{\lambda^{x_1-x_2}}{z^{y_1-y_2}}T_{0,j_2}^{-1}(z)\Pi T_{0,j_1}(z)\frac{dz}{z},
    \end{equation}
    as $N\to\infty$.

    From \eqref{eq:KN_v1_v2} with $v_{1,N}$ and $v_{2,N}$ as in \eqref{eq:v1N}--\eqref{eq:v2N} and the definition \eqref{eq:KN_matrix}, we find that
    \begin{equation}\label{eq:KN_case1a}
        K_N(v_{1,N},v_{2,N}) = \left[\mathsf{K}_N(v_{1,N},v_{2,N})\right]_{k_2,k_1}
    \end{equation}
    in case $3x_1 + j_1 \leq 3x_2 + j_2$. Then, by \eqref{eq:KN_limit}, we see that \eqref{eq:KN_case1a} has the limit
    \begin{equation}\label{eq:KN_case1b}
        \left[\frac{1}{2\pi i}\int_{\gamma_s^\leq}\frac{\lambda^{x_1-x_2}}{z^{y_1-y_2}}T_{0,j_2}^{-1}(z)\Pi T_{0,j_1}(z)\frac{dz}{z} \right]_{k_2,k_1}
    \end{equation}
    which agrees with \eqref{eq:limiting_kernel} in case $3x_1+j_1\leq3x_2+j_2$.

    In the other case $3x_2+j_2<3x_1+j_1$, we have the additional single contour integral
    \begin{equation}\label{eq:KN_case2a}
        \left[-\frac{1}{2\pi i}\oint_{\mathbb{T}}T_{0,j_2}^{-1}(z) W(z)^{x_1-x_2} T_{0,j_1}(z) \frac{dz}{z^{y_1-y_2+1}}\right]_{k_2,k_1}
    \end{equation}
    in formula \eqref{eq:KN_v1_v2}. This expression does not change with $N$. We move \eqref{eq:KN_case2a} to an integral on $\mathcal{R}$, in the same way as we did it in Section~\ref{sec:integrals_Riemann_surface} for the double integral (it is even simpler). Then, \eqref{eq:KN_case2a} is equal to
    \begin{equation}\label{eq:KN_case2b} 
        \left[-\frac{1}{2\pi i} \oint_{z^{-1}(\mathbb{T})}\frac{\lambda^{x_1-x_2}}{z^{y_1-y_2}}T_{0,j_2}^{-1}(z)\Pi T_{0,j_1}(z) \frac{dz}{z} \right]_{k_2,k_1}.
    \end{equation}
    We still have \eqref{eq:KN_case1b} for the limit of the double contour integral. Observe that the integrals \eqref{eq:KN_case1b} and \eqref{eq:KN_case2b} have the same integrand, and also observe that the integrand does not have a pole at $P_1$ as now we have $x_1 \geq x_2$.

    By analyticity, we can deform $z^{-1}(\mathbb{T})$ within $\mathcal{R} \setminus \{P_0, P_\infty\}$ to a contour passing through $s$, and we can take it to be $\gamma_s^\leq$ followed by $\gamma_s^>$ in the opposite direction, see also Figure~\ref{fig:torus_KN_contours}. Then, the sum of \eqref{eq:KN_case1b} and \eqref{eq:KN_case2b} is
    \begin{equation}\label{eq:KN_case2c} 
        \left[\frac{1}{2\pi i} \oint_{\gamma_s^>}\frac{\lambda^{x_1-x_2}}{z^{y_1-y_2}}T_{0,j_2}^{-1}(z)\Pi T_{0,j_1}(z)\frac{dz}{z}\right]_{k_2,k_1},
    \end{equation}
    which agrees with \eqref{eq:limiting_kernel} in case $3x_1+j_1>3x_2+j_2$.

    This completes the proof of Theorem~\ref{thm:local_scaling_limit} for the case $\xi_1>0$.
\end{proof}


\subsection{Proof of Theorem~\ref{thm:local_scaling_limit} for $\xi_1<0$ and $\xi_1=0$}\label{sec:remaining_cases}
In order to complete the proof of Theorem~\ref{thm:local_scaling_limit}, we establish the limiting behavior \eqref{eq:KN_limit} of the double contour integral \eqref{eq:KN_matrix} as $N\to\infty$ for the cases $\xi_1<0$ and $\xi_1=0$. Then, by applying the same arguments as we did above for the case $\xi_1 >0$, we obtain Theorem~\ref{thm:local_scaling_limit} for the remaining cases.


\subsubsection{The case $\xi_1<0$}
In Sections~\ref{sec:deform_p2}--\ref{sec:estimates}, we restricted ourselves to the case $\xi_1>0$, so that the saddle points $s = \Omega(\xi_1,\xi_2)$ and $\bar{s}$ lie on the first copy $\mathcal{R}$. In the case $\xi_1<0$, the saddle points $s = \Omega(\xi_1,\xi_2)$ and $\bar{s}$ belong to the second copy $\mathcal{R}^\sigma$, which also holds for $s_N = \Omega(\xi_{1,N},\xi_{2,N})$ and $\bar{s}_N$ for $N$ sufficiently large so that $\xi_{1,N}<0$. We recall \eqref{eq:critical_level_set_sigma}, which states that $\widehat{\mathcal{N}}_{(\xi_1,\xi_2)} = \sigma(\widehat{\mathcal{N}}_{(-\xi_1,-\xi_2)})$. It follows that the case $\xi_1<0$ is simply obtained from the case $\xi_1>0$ by interchanging the roles of $\mathcal{R}$ and $\mathcal{R}^\sigma$.

More precisely, we first deform the $p_1$-integral (instead of the $p_2$-integral) similar to Section~\ref{sec:deform_p2}. Then, we split and deform the $p_2$-integral following Section~\ref{sec:split_p1}. For this splitting step, we use the following analogue of Lemma~\ref{lem:AN_poles}(b). Defining the $6\times3$ matrix-valued function
\begin{equation}\label{eq:DN}
	\mathbf{D}_N(p) = \mathbf{X}_N^{-t}(z(p)) \begin{pmatrix} I_3 \\ 0_3 \end{pmatrix} \Xi(p) E^{t}(z(p))  \frac{\lambda(p)^{-N(1+\xi_{1,N})}}{z(p)^{-N(1+\xi_{2,N})}},
\end{equation}
we have
\begin{equation}\label{eq:DN_split}
    \mathbf{D}_N(p) = \mathbf{B}_{N,+}(p) - \mathbf{B}_{N,-}(p), \quad p\in z^{-1}(\mathbb{T}).
\end{equation}
The remainder of the analysis follows the arguments in Sections~\ref{sec:deform_p2}--\ref{sec:estimates}, and ultimately shows that \eqref{eq:KN_limit} remains valid for $\xi_1<0$.


\subsubsection{The special case $\xi_1=0$}
Let us consider the special case $\xi_1 = 0$, which corresponds to $s\in\Gamma_0\cup\Gamma_0^\sigma$. The steepest descent analysis in this case is still largely the same. However, the splitting step described in Section~\ref{sec:split_p1} must be applied to both $\mathbf{B}_N$ and $\mathbf{C}_N$, see also \eqref{eq:DN_split}.

We first deform the integration contour for the $p_2$-integral to $z^{-1}(\delta\mathbb{T})$ for some small $\delta>0$. Then, we perform the splitting and deformation of the $p_1$-integral as described in Section~\ref{sec:split_p1}. In the process of deforming $z^{-1}(\delta\mathbb{T})$ back to $z^{-1}(\mathbb{T})$, we pick up residues along the contour $\gamma_{s_N}$. After this deformation, we arrive at \eqref{eq:KN_ANBNCN} with $\delta\mathbb{T}$ substituted for $\gamma_B$. 

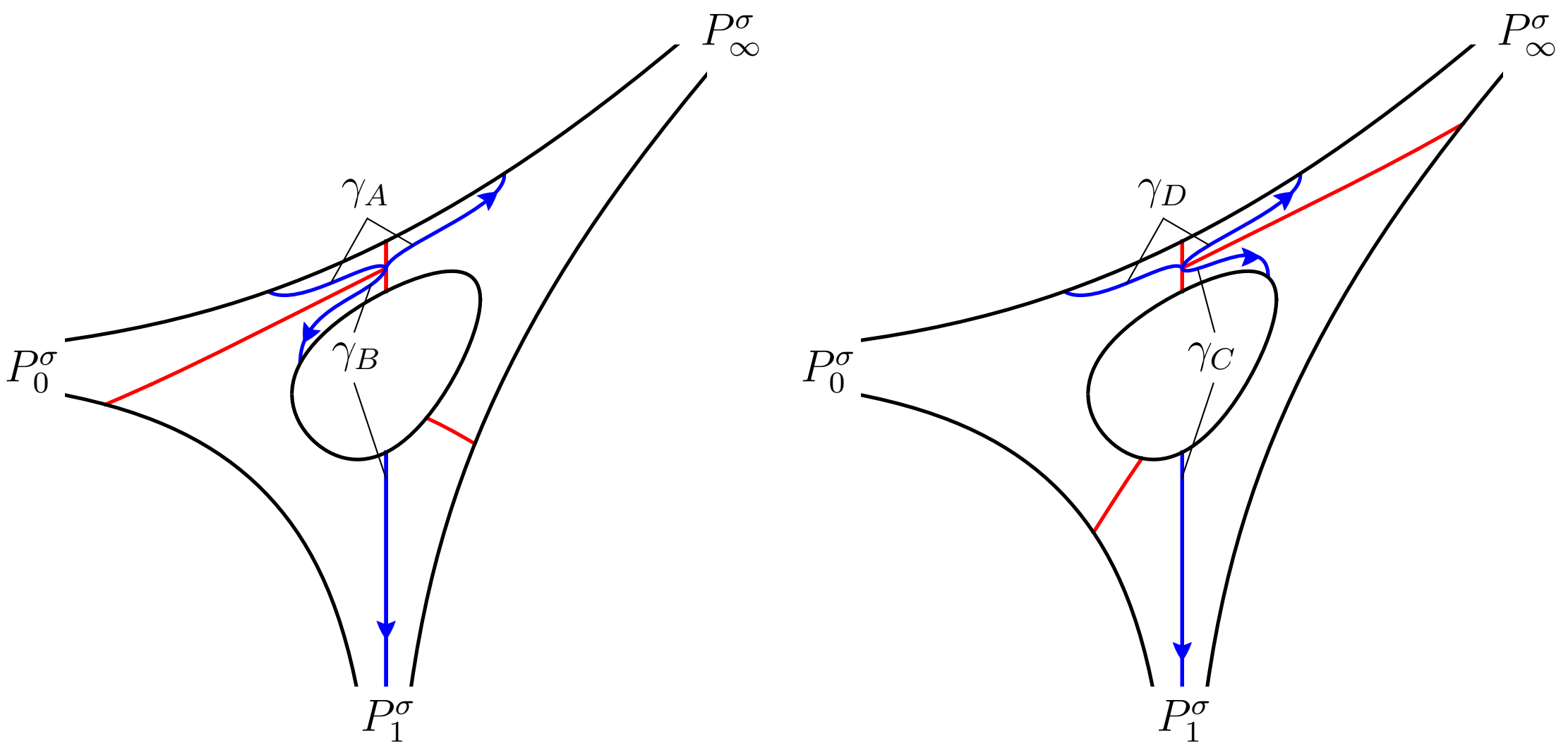
\begin{figure}[t]
    \centering
    \hspace{0.6cm}
	\begin{overpic}[width=0.4\linewidth]{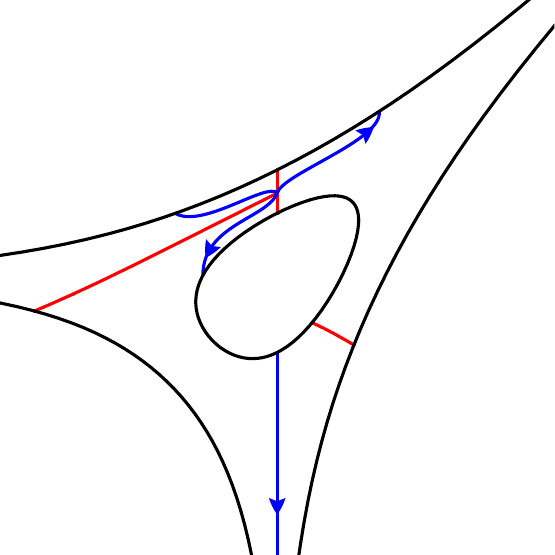}
		\put(-9.5,47.5){$P_0^\sigma$}
		\put(46,-7){$P_1^\sigma$}
		\put(99,100){$P_\infty^\sigma$}
        \put(44.72,68.6){
        \begin{tikzpicture}
        	\draw[-] (0,0) -- (0.43,-0.25);
        \end{tikzpicture}}
        \put(38.9,62.6){
        \begin{tikzpicture}
        	\draw[-] (0,0) -- (-0.35,-0.62);
        \end{tikzpicture}}
		\put(43,76){$\gamma_A$}
        \put(42.6,55){
        \begin{tikzpicture}
            \draw[-] (0,0) -- (0.16,0.45);
        \end{tikzpicture}}
        \put(42.6,32.2){
        \begin{tikzpicture}
            \draw[-] (0,0) -- (0.3,-0.9);
        \end{tikzpicture}}
		\put(41.5,50.5){$\gamma_B$}
    \end{overpic}
    \hfill
	\begin{overpic}[width=0.4\linewidth]{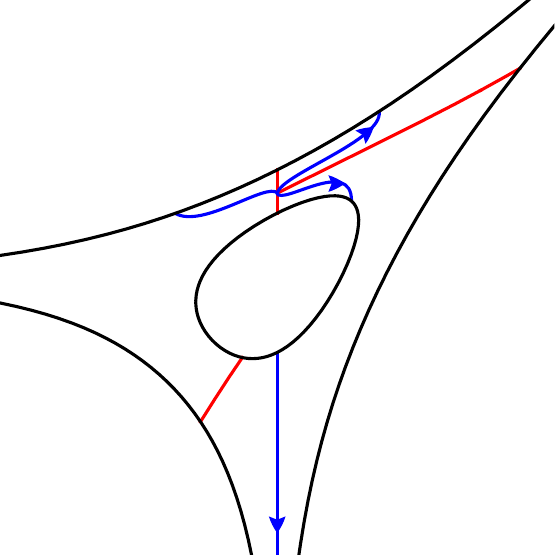}
		\put(-9.5,47.5){$P_0^\sigma$}
		\put(46,-7){$P_1^\sigma$}
		\put(99,100){$P_\infty^\sigma$}
        \put(44.72,68.6){
        \begin{tikzpicture}
        	\draw[-] (0,0) -- (0.43,-0.25);
        \end{tikzpicture}}
        \put(38.9,62.6){
        \begin{tikzpicture}
        	\draw[-] (0,0) -- (-0.35,-0.62);
        \end{tikzpicture}}
		\put(43,76){$\gamma_D$}
        \put(50,55){
        \begin{tikzpicture}
            \draw[-] (0,0) -- (-0.16,0.6);
        \end{tikzpicture}}
        \put(47.5,32.2){
        \begin{tikzpicture}
            \draw[-] (0,0) -- (-0.3,-0.9);
        \end{tikzpicture}}
		\put(50.8,50.5){$\gamma_C$}
    \end{overpic}
    \hspace{0.6cm}
    \vspace{0.4cm}
	\caption{After splitting, parts of the $p_1$- and $p_2$-integrals are deformed to $\gamma_A$ and $\gamma_D$, respectively. The left subfigure shows the contours within $\mathcal{R}$, which are $\gamma_A$ and $\gamma_B$, and the remaining contours $\gamma_C$ and $\gamma_D$, which are within $\mathcal{R}^\sigma$, are shown in the right subfigure.}
    \label{fig:deform_BD}
\end{figure}

We now use \eqref{eq:DN_split} to split the $p_1$-integral, which results in
\begin{align}
	\mathsf{K}_N(v_{1,N},v_{2,N}) &= \mathsf{K}_N^{(B,C)}(v_{1,N},v_{2,N}) + \mathsf{K}_N^{(B,A)}(v_{1,N},v_{2,N}) \nonumber \\
    &+ \mathsf{K}_N^{(D,C)}(v_{1,N},v_{2,N}) + \mathsf{K}_N^{(D,A)}(v_{1,N},v_{2,N}) + \mathsf{K}_N^{(s_N)}(v_{1,N},v_{2,N}). \label{eq:KN_ANBNCNDN}
\end{align}
Here, $\mathsf{K}_N^{(B,C)}(v_{1,N},v_{2,N})$, $\mathsf{K}_N^{(B,A)}(v_{1,N},v_{2,N})$, and $\mathsf{K}_N^{(s_N)}(v_{1,N},v_{2,N})$ are given by \eqref{eq:KN_BC}--\eqref{eq:KN_sN}, and $\mathsf{K}_N^{(D,C)}(v_{1,N},v_{2,N})$ and $\mathsf{K}_N^{(D,A)}(v_{1,N},v_{2,N})$ are defined analogously to \eqref{eq:KN_BC} and \eqref{eq:KN_BA}, with $\mathbf{B}_N$ and $\gamma_B$ replaced by $\mathbf{D}_N$ and $\gamma_D$, respectively. The curves $\gamma_A$, $\gamma_B$, $\gamma_C$, and $\gamma_D$ are illustrated in Figure~\ref{fig:deform_BD}.

The estimation of the double contour integrals in \eqref{eq:KN_ANBNCNDN} is analogous to the analysis in Section~\ref{sec:estimates}. We only make two brief remarks.

In Figure~\ref{fig:deform_BD}, $\gamma_C$ and $\gamma_D$ are drawn as curves in $\mathcal{R}^\sigma$ because we want to interpret them on $\mathcal{R}^\sigma$. The reason for this is the appearance of $\Phi_N(\sigma(p_1))$ and $\Phi_N(\sigma(p_2))$ in the scalar factors. Previously, in Figure~\ref{fig:deform_AC}, the curve $\gamma_C$ was shown as a dashed line within $\mathcal{R}$ for the same reason. 

In addition to the local analysis near the saddle points $s$ and $\bar{s}$, a similar analysis is required near the regular points $\sigma(s)$ and $\sigma(\bar{s})$. A straightforward computation shows that the contributions from these points tend to $0$ as $N\to\infty$ at a rate of $\mathcal{O}(N^{-2})$. Note that this is much faster than the rate $\mathcal{O}(N^{-1/2})$ arising from the saddle points.

This establishes \eqref{eq:KN_limit} for the special case $\xi_1=0$.


\section{Proof of Theorem~\ref{thm:limit_shape}}\label{sec:limit_shape}


\subsection{Proof of Theorem~\ref{thm:limit_shape}(a)}
\begin{proof} 
	The expectation value of the height function at a vertex $v=(x,y)$ is
	\begin{align}
		\mathbb{E}\left[h_N(v)\right] &= \sum_{y'<y}\mathbb{P}\left[\text{particle at } (x,y')\right] \nonumber \\
		&= 3N - \sum_{y'\geq y}\mathbb{P}\left[\text{particle at } (x,y')\right].
	\end{align}
	Because of the determinantal structure this equals
	\begin{equation}\label{eq:EhNv}
		\mathbb{E}\left[h_N(v)\right] = 3N - \sum_{y'\geq y} K_N((x,y'), (x,y')).
	\end{equation}
	Let $v = v_N = (3N(1+\xi_{1,N}), 3N(1+\xi_{2,N}))$ as in \eqref{eq:vN}. Then, by \eqref{eq:KN_v1_v2}, \eqref{eq:KN_matrix}, and \eqref{eq:EhNv}, we have

	\begin{multline}\label{eq:EhN_vN_1}
		\mathbb{E}[h_N(v_N)]
		= 3N - \sum_{y = N(1+\xi_{2,N})}^\infty 
		\Tr \mathsf{K}_N((3N(1+\xi_{1,N}), 3y), (3N(1+\xi_{1,N}),3y)) \\
		= 3N - \sum_{y=N(1+\xi_{2,N})}^\infty
		\Tr \left[\frac{1}{(2\pi i)^2}
		\oint_\mathbb{T}\oint_\mathbb{T} 
		\frac{W(z_2)^{N(1-\xi_{1,N})}}{z_2^{2N-y}} 
		\right. \\
		\left. \times R_N(z_2,z_1) 
		\frac{W(z_1)^{N(1+\xi_{1,N})}}{z_1^y} \frac{dz_1dz_2}{z_1} 
		\right].
	\end{multline} 
	Here, we used \eqref{eq:KN_matrix} with $v_{1,N}=v_{2,N}=(3N(1+\xi_{1,N}),3y)$, replacing $N\xi_{2,N}$ by $y-N$ and setting $x_1$, $x_2$, $y_1$, $y_2$, $j_1$, and $j_2$ to $0$. Note that \eqref{eq:EhN_vN_1} involves an infinite series. However, the terms vanish for $y>2N$ because then we are outside the hexagon, and hence the series reduces to a finite sum.

    To bring the summation inside the double contour integral, we first deform the contour of the $z_2$-integral to the smaller circle $(1-\delta)\mathbb{T}$ with $\delta\in(0,1)$. Then, we use the summation formula of the geometric series
	\[
        \sum_{y=N(1+\xi_{2,N})}^\infty\!\!\left(\frac{z_2}{z_1}\right)^y = \left(\frac{z_2}{z_1}\right) ^{N(1+\xi_{2,N})}\!\frac{z_1}{z_1-z_2} \qquad \text{for } |z_2|=1-\delta<1=|z_1|,
    \]
    and obtain \eqref{eq:EhN_vN}, which concludes the proof.
\end{proof}


\subsection{Strategy for the proof of Theorem~\ref{thm:limit_shape}(b)}
We first observe that it suffices to prove Theorem~\ref{thm:limit_shape}(b) for $v=(0,0)$. Indeed, from the path interpretation it is clear that the height function at two neighboring vertices differs by at most one. Hence,
\[
    \left|h_N(v_N+v)-h_N(v_N)\right| \leq |x|+|y| \quad \text{for } v=(x,y),
\]
and, if the limit \eqref{eq:EhN_vN_limit} holds with $v=(0,0)$, then it also holds for arbitrary $v$.

Therefore, we set $v=(0,0)$, which allows us to use \eqref{eq:EhN_vN_1}. Recall from \eqref{eq:RN_formula} that $R_N$ can be written in terms of the solution $\mathbf{Y}_N$ of the RH problem. We apply the transformation $\mathbf{Y}_N \mapsto \mathbf{X}_N$ from \eqref{eq:XN}, and lift the integrals to the Riemann surface as in Section~\ref{sec:integrals_Riemann_surface} to obtain 
\begin{multline}\label{eq:EhN_vN_with_BC}
    \mathbb{E}[h_N(v_N)] = 3N - \Tr\left[\frac{1}{(2\pi i)^2}\oint_{z^{-1}((1-\delta)\mathbb{T})} dz(p_2) \right. \\
    \left.\oint_{z^{-1}(\mathbb{T})} \frac{dz(p_1)}{(z(p_1)-z(p_2))^2}\mathbf{B}_N^t(p_2) \mathbf{C}_N(p_1)\right]
\end{multline}
for some $\delta\in(0,1)$, where $\mathbf{B}_N$ and $\mathbf{C}_N$ are as given in \eqref{eq:BN} and \eqref{eq:CN}. The formula~\eqref{eq:EhN_vN_with_BC} is the analogue of \eqref{eq:KN_XN_on_R} for the height function.

Our strategy is now to substitute \eqref{eq:EhN_vN_with_BC} for \eqref{eq:KN_XN_on_R} and proceed with the steepest descent analysis developed in Section~\ref{sec:local_scaling_limit}. The most notable difference is that the integrand in \eqref{eq:EhN_vN_with_BC} now has a simple pole at $z(p_1)=z(p_2)$ due to the extra factor of $z(p_1)-z(p_2)$ compared to \eqref{eq:KN_XN_on_R}. The presence of this pole introduces certain complications, as we now explain.

We restrict ourselves to the case $(\xi_1,\xi_2)\in\mathfrak{L}$ (rough phase) with $\xi_1>0$ where $s=\Omega(\xi_1,\xi_2)$ satisfies $|z(s)|<1$. This occurs, for example, when $0<\xi_1<\xi_2$ by \eqref{eq:xi1_xi2_conditions}. The remaining cases can be handled similarly. We take $N$ sufficiently large so that $\xi_{1,N}>0$ and $s_N=\Omega(\xi_{1,N},\xi_{2,N})$ also satisfies $|z(s_N)|<1$.

Then, we deform $z^{-1}((1-\delta)\mathbb{T})$ within $\{p\in\mathcal{R} \mid |z(p)|<1\}$ to a contour passing through the saddle points $s_N$ and $\bar{s}_N$, and lying otherwise in $\mathcal{N}^<_{(\xi_{1,N},\xi_{2,N})}$. Recall from Section~\ref{sec:critical_level_sets} that the level set $\mathcal{N}_{(\xi_{1,N},\xi_{2,N})}$ is topologically equivalent to either the left and middle subfigures of Figures~\ref{fig:critical_level_sets_4}, or the left subfigure of Figures~\ref{fig:critical_level_sets_1}.

The resulting contour is similar to the contour $\gamma_B$ that was used in \eqref{eq:KN_BNCN}, see also Figure~\ref{fig:deform_B}. However, the contour $\gamma_B$ partly overlaps with $z^{-1}(\mathbb{T})$, which must be avoided in the present setting, since the integrand in \eqref{eq:EhN_vN_with_BC} has a pole at $z(p_1)=z(p_2)$. We therefore use a slightly deformed contour $\gamma_{B,\delta}$, so that \eqref{eq:EhN_vN_with_BC} continues to hold when $z^{-1}((1-\delta)\mathbb{T})$ is replaced by $\gamma_{B,\delta}$.

As a next step, we wish to use \eqref{eq:CN_split} to split the $p_1$-integral and thereby obtain the analogue of \eqref{eq:KN_ANBNCN} for the height function. The integration contour $z^{-1}(\mathbb{T})$ then splits into $\gamma_A$ and $\gamma_C$ as in Figure~\ref{fig:deform_AC}. In particular, we would then obtain
\begin{align}\label{eq:EhN_vN_ANBNCN} 
	\mathbb{E}[h_N(v_N)] = 3N - \mathsf{H}_N^{(B,C)}(v_N) - \mathsf{H}_N^{(B,A)}(v_N) - \mathsf{H}_N^{(s_N)}(v_N) \quad \text{(not correct!)}
\end{align}
with
\begin{align}
    \mathsf{H}_N^{(B,C)}(v_N) &= \frac{1}{(2\pi i)^2}\oint_{\gamma_{B,\delta}} dz(p_2)\oint_{\gamma_C}\frac{dz(p_1)}{(z(p_1)-z(p_2))^2}\Tr(\mathbf{B}_N^t(p_2) \mathbf{C}_N(p_1)), \label{eq:HN_BC} \\
    \mathsf{H}_N^{(B,A)}(v_N) &= \frac{1}{(2\pi i)^2}\oint_{\gamma_{B,\delta}} dz(p_2)\oint_{\gamma_A} \frac{dz(p_1)}{(z(p_1)-z(p_2))^2}\Tr(\mathbf{B}_N^t(p_2) \mathbf{A}_N(p_1)), \label{eq:HN_BA} \\
    \mathsf{H}_N^{(s_N)}(v_N) &= \frac{1}{2\pi i}\int_{\gamma_{s_N}}dz(p_2)\Res_{p_1=p_2}\left[\Tr(\mathbf{B}_N^t (p_2)\mathbf{A}_N(p_1))\frac{dz(p_1)}{(z(p_1)-z(p_2))^2}\right]\!. \label{eq:HN_sN}
\end{align}
However, this is not entirely justified, as the double contour integral \eqref{eq:HN_BA} diverges due to the presence of a double pole at $z(p_1)=z(p_2)$.


\subsection{Avoiding the divergent integral}
\begin{figure}[t]
	\centering
	\begin{tikzpicture}
		\def\r{1.8}

		\draw[ultra thick,white,postaction={decorate,decoration={markings,mark=at position 0.55 with {\arrow{<}}}}] (0,{-sqrt(2)*1.2*\r}) -- (0,0);

		\draw[ultra thick,red] (-1.2*\r, -1.2*\r) -- (1.2*\r, 1.2*\r);
		\draw[ultra thick,red] (-1.2*\r, 1.2*\r) -- (1.2*\r, -1.2*\r);

		\foreach \angle in {0, 10, ..., 350} {
            \draw[ultra thick] (\angle-2.5:1.5*\r) arc (\angle-2.5:\angle+2.5:1.5*\r);
		}

		\draw[ultra thick,dotted] (0,0) -- ({sin(-70)*1.5*\r},{cos(-70)*1.5*\r})
		node[midway,above,xshift=-22pt,yshift=8pt] {$\sqrt{\varepsilon}$};
		\draw[ultra thick,dotted] (0,0) -- ({sin(-100)*\r/2},{cos(-100)*\r/2})
		node[midway,below,xshift=-5pt,yshift=0pt] {$\varepsilon$};

		\draw[ultra thick,blue,postaction={decorate,decoration={markings,mark=at position 0.55 with {\arrow{>}}}}] (0,{sqrt(2)*1.2*\r}) -- (0,0) node[midway,right,xshift=-1pt,yshift=6pt] {$\gamma_{B,1}$};

		\draw[ultra thick,blue,postaction={decorate,decoration={markings,mark=at position 0.82 with {\arrow{>}}}}] ({-sqrt(2)*1.2*\r},0) -- (-\r/2,0) arc (180:360:\r/2) -- ({sqrt(2)*1.2*\r},0) node[midway,above,xshift=0pt,yshift=0pt] {$\gamma_{A,1,\varepsilon}$};

		\fill (0,0) circle (2pt);
	\end{tikzpicture}
	\qquad
	\begin{tikzpicture}
		\def\r{1.8}

		\draw[ultra thick,white,postaction={decorate,decoration={markings,mark=at position 0.55 with {\arrow{<}}}}] (0,{sqrt(2)*1.2*\r}) -- (0,0);

		\draw[ultra thick,red] (-1.2*\r, -1.2*\r) -- (1.2*\r, 1.2*\r);
		\draw[ultra thick,red] (-1.2*\r, 1.2*\r) -- (1.2*\r, -1.2*\r);

		\foreach \angle in {0, 10, ..., 350} {
            \draw[ultra thick] (\angle-2.5:1.5*\r) arc (\angle-2.5:\angle+2.5:1.5*\r);
		}

		\draw[ultra thick,dotted] (0,0) -- ({sin(-110)*1.5*\r},{cos(-110)*1.5*\r})
		node[midway,below,xshift=-22pt,yshift=-6pt] {$\sqrt{\varepsilon}$};
		\draw[ultra thick,dotted] (0,0) -- ({sin(-80)*\r/2},{cos(-80)*\r/2})
		node[midway,above,xshift=-5pt,yshift=0pt] {$\varepsilon$};

		\draw[ultra thick,blue,
		postaction={decorate,decoration={markings,mark=at position 0.55 with {\arrow{<}}}}] (0,{-sqrt(2)*1.2*\r}) -- (0,0) node[midway,right,xshift=-1pt,yshift=-6pt] {$\gamma_{B,2}$};

		\draw[ultra thick,blue,postaction={decorate,decoration={markings,mark=at position 0.82 with {\arrow{>}}}}] ({-sqrt(2)*1.2*\r},0) -- (-\r/2,0) arc (180:0:\r/2) -- ({sqrt(2)*1.2*\r},0) node[midway,above,xshift=0pt,yshift=0pt] {$\gamma_{A,2,\varepsilon}$};

		\fill (0,0) circle (2pt);
	\end{tikzpicture}
	\caption{The contours $\gamma_{A,j,\varepsilon}$ and $\gamma_{B,j}$ for $j=1,2$ in a small neighborhood of the saddle point.}
	\label{fig:contours_in_epsilon_neighborhood}
\end{figure}
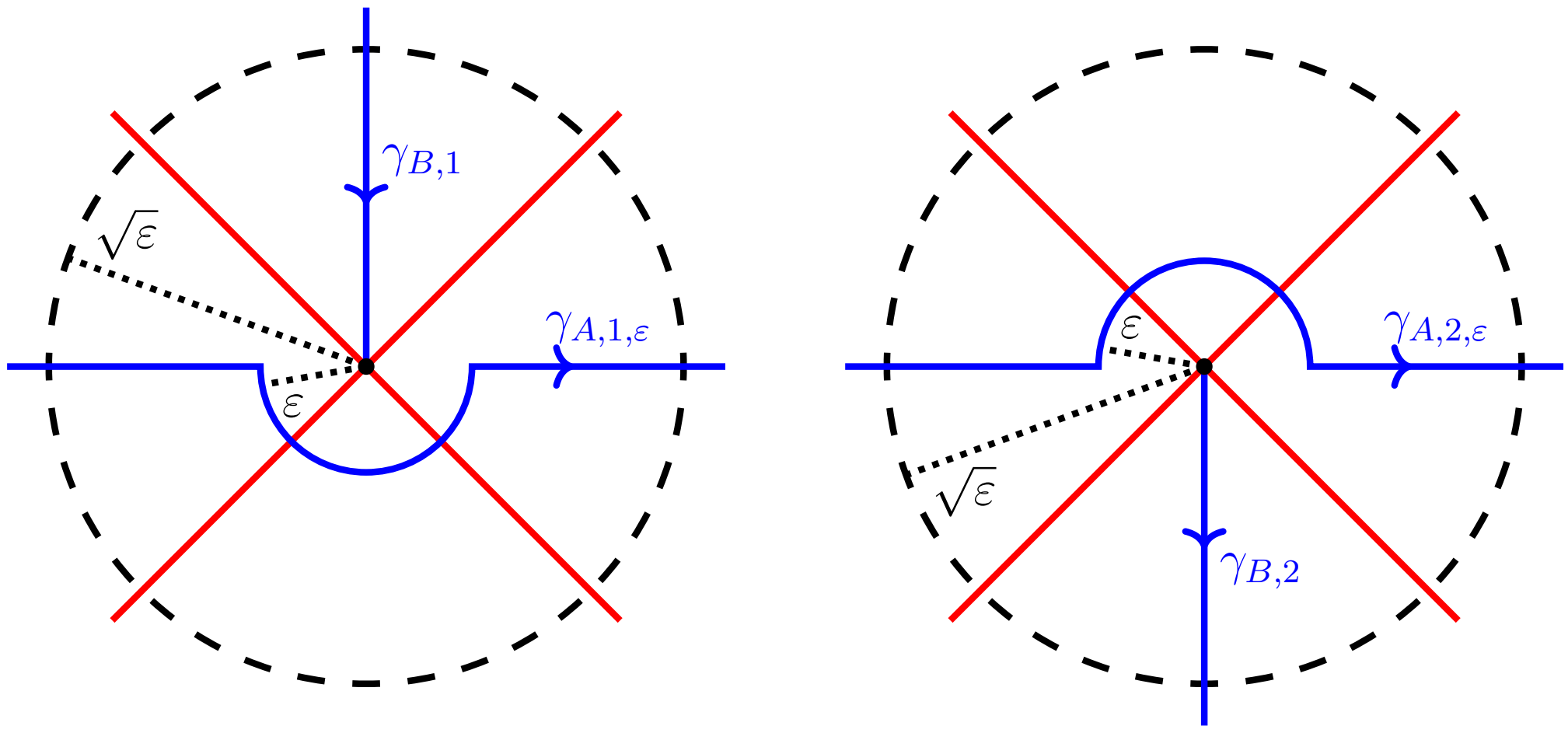

We have to be more careful with the deformation to $\gamma_A$. Note that for this part we may deform $\gamma_{B,\delta}$ back to $\gamma_B$. Before we do the deformation of the $p_1$-integral, we first write the $p_2$-integral over $\gamma_B$ as a sum of two integrals, one over $\gamma_{B,1} = \gamma_{s_N}$ and the other over $\gamma_{B,2} = \gamma_B\setminus\gamma_{s_N}$. Then, we deform the $p_1$-integral to either $\gamma_{A,1,\varepsilon}$ or $\gamma_{A,2,\varepsilon}$ such that $\gamma_{A,j,\varepsilon}$ is disjoint from $\gamma_{B,j}$ for $j=1,2$. More precisely, $\gamma_{A,j,\varepsilon}$ is equal to $\gamma_A$ outside of a small $\varepsilon$-neighborhood of $s_N$, combined with a circular arc of radius $\varepsilon$ around $s_N$ (the metric is with respect to the $z$-variable). The circular arc is different for the two cases $j=1,2$, see Figure~\ref{fig:contours_in_epsilon_neighborhood}. We preserve the conjugation symmetry, and there is a similar construction around $\bar{s}_N$.

In the deformation process, we pick up only the residue contribution \eqref{eq:HN_sN} for $j=1$. There is no residue contribution for $j=2$. The result of this process is the following corrected version of \eqref{eq:EhN_vN_ANBNCN}:
\begin{equation}\label{eq:EhN_vN_ANBNjCN}
	\mathbb{E}[h_N(v_N)] = 3N - \mathsf{H}_N^{(B,C)}(v_N) - \mathsf{H}_{N,1}^{(B,A)}(v_N) - \mathsf{H}_{N,2}^{(B,A)}(v_N) - \mathsf{H}_N^{(s_N)}(v_N),
\end{equation}
where $\mathsf{H}_N^{(B,C)}(v_N)$ and $\mathsf{H}_N^{(s_N)}(v_N)$ are given in \eqref{eq:HN_BC} and \eqref{eq:HN_sN}, and
\begin{equation}\label{eq:HN_BjA}
	\mathsf{H}_{N,j}^{(B,A)} = \frac{1}{(2\pi i)^2}\int_{\gamma_{B,j}} dz(p_2) \oint_{\gamma_{A,j,\varepsilon}} \frac{dz(p_1)}{(z(p_1)-z(p_2))^2} \Tr(\mathbf{B}_N^t(p_2) \mathbf{A}_N(p_1))
\end{equation}
for $j=1,2$.

\begin{remark}\label{rem:HN_BC_small}
    Following the proof of Lemma~\ref{lem:KN_small}, one finds that $\mathsf{H}_N^{(B,C)}(v_N)$ is exponentially small as $N\to\infty$. This is due to the inequality \eqref{eq:PhiN_estimate_BNCN} which remains valid for $p_1\in\gamma_C$ and $p_2\in\gamma_{B,\delta}$. As a consequence, we have that $\mathbf{B}_N^t(p_2)\mathbf{C}_N(p_1)$ is exponentially small as $N\to\infty$. Then, the double contour integral $\mathsf{H}_N^{(B,C)}(v_N)$ tends to $0$ as $N\to\infty$ at an exponential rate.
\end{remark}


\subsection{Estimation of \eqref{eq:HN_BjA}}\label{sec:estimate}
The double contour integrals \eqref{eq:HN_BjA} depend on $\varepsilon$. To show that these terms do not contribute at leading order as $N\to\infty$, we let $\varepsilon$ to depend on $N$ and choose specifically $\varepsilon=N^{-1/2}$. The integrand in \eqref{eq:HN_BjA} is the product of the scalar factor 
\[
    \frac{1}{(z(p_1)-z(p_2))^2}e^{N(\Phi_N(p_2)-\Phi_N(p_1))}
\]
and a matrix-valued part that remains uniformly bounded as $N\to\infty$ due to \eqref{eq:AN_TN}, \eqref{eq:BN_TN}, and Lemma~\ref{lem:TN_bounded}(b).

If $p_1\in\gamma_{A,j,\varepsilon}$ such that $|z(p_1) - z(s_N)|<\sqrt{\varepsilon}$, then
\[
    \Re(\Phi_N(s_N) - \Phi_N(p_1))<c\,\varepsilon^2 = cN^{-1}
\]
for some constant $c>0$ that is independent of $N$. Indeed, if $\varepsilon<|z(p_1) - z(s_N)|<\sqrt{\varepsilon}$, then $\Re(\Phi_N(s_N) - \Phi_N(p_1))<0$, since $p_1$ is on $\gamma_A$, which is the path of steepest ascent from $s_N$. In case that $|z(p_1) - z(s_N)| = \varepsilon$, we use that $s_N$ is a simple saddle point to conclude that $\Re(\Phi_N(s_N) - \Phi_N(p_1))$ is bounded by a constant multiple of the squared distance from $p_1$ to $s_N$. Here, we may take the same constant multiple $c$ for all $N$ sufficiently large, since the second derivatives are uniformly bounded in a neighborhood of $s = \lim_{N\to\infty} s_N$.

Additionally, if $p_2\in\gamma_{B,j}$, then
\[
    \Re(\Phi_N(p_2) - \Phi_N(p_1))<cN^{-1}
\]
because $\Re(\Phi_N(p_2) - \Phi_N(s_N))<0$. Observing that
\[
    |z(p_1) - z(p_2)| \geq \varepsilon = N^{-1/2},
\]
we find that 
\[
    \left|\frac{1}{(z(p_1)-z(p_2))^2}e^{N(\Phi_N(p_2)-\Phi_N(p_1))}\right| \leq N e^c.
\]
The length of $\gamma_{A,j,\varepsilon}$ inside the disk of radius $\sqrt{\varepsilon}$ is of order $N^{-1/4}$. It follows that the contribution of this part to the double contour integral is of order $\mathcal{O}(N^{3/4})$ as $N\to\infty$.

Similarly, if $p_1\in\gamma_{A,j,\varepsilon}$ such that $|z(p_1) - z(s_N)| \geq \sqrt{\varepsilon}$ and $p_2\in\gamma_{B,j}$, then
\[
    \Re(\Phi_N(p_2) - \Phi_N(p_1))<-c\,\varepsilon = -cN^{-1/2}
\]
and
\[
    |z(p_1) - z(p_2)| \geq \sqrt{\varepsilon} = N^{-1/4}.
\]
Therefore,
\[
    \left|\frac{1}{(z(p_1)-z(p_2))^2}e^{N(\Phi_N(p_2)-\Phi_N(p_1))}\right| \leq \sqrt{N}e^{-c\sqrt{N}},
\]
which becomes uniformly small as $N\to\infty$. We therefore conclude that the double contour integrals in \eqref{eq:HN_BjA} are of order $\mathcal{O}(N^{3/4})$. After dividing by $N$, these terms tend to $0$ as $N\to\infty$.


\subsection{Evaluation of \eqref{eq:HN_sN}}
It remains to consider \eqref{eq:HN_sN}. We begin by computing the residue that appears in the integrand. Because of \eqref{eq:BN}, \eqref{eq:AN}, and the cyclic property of the trace, we want to evaluate
\begin{multline}\label{eq:EX_double_pole}
	\Res_{p=p_2}\Tr\left[E^{-1}(z(p)) E(z(p_2))\Xi(p_2) \right. \\
    \times \begin{pmatrix} 0_3 & I_3 \end{pmatrix} \mathbf{X}_N^{-1}(z(p_2))\mathbf{X}_N(z(p)) \begin{pmatrix} 0_3 \\ I_3 \end{pmatrix}\Xi(p)\\
    \times \left.\left(\frac{\lambda(p_2)}{\lambda(p)}\right)^{N(1-\xi_{1,N})}\left(\frac{z(p)}{z(p_2)}\right)^{N(1-\xi_{2,N})}\frac{dz(p)}{(z(p)-z(p_2))^2}\right].
\end{multline}
With the transformation $\mathbf{X}_N \mapsto \mathbf{T}_N$ from \eqref{eq:TN_def}, we obtain that \eqref{eq:EX_double_pole} is equal to
\begin{multline}\label{eq:ET_double_pole}
	\Res_{p=p_2}\Tr\left[E^{-1}(z(p)) E(z(p_2)) \Xi(p_2) \right. \\
    \times \begin{pmatrix} 0_3 & I_3 \end{pmatrix} \mathbf{T}_N^{-1}(z(p_2))\mathbf{T}_N(z(p)) \begin{pmatrix} 0_3 \\ I_3 \end{pmatrix}\Xi(p) \\
    \times\left.e^{N(g(p_2)-g(p))}\left(\frac{\lambda(p_2)}{\lambda(p)}\right)^{N(1-\xi_{1,N})}\left(\frac{z(p)}{z(p_2)}\right)^{N(1-\xi_{2,N})}\frac{dz(p)}{(z(p)-z(p_2))^2}\right].
\end{multline}
The evaluation of \eqref{eq:ET_double_pole} is a straightforward application of the residue formula for double poles. In its basic form, it relies on the observation that for a matrix-valued function $F$ that is invertible and analytic at $z_0$, one has
\begin{align*}
    \Res_{z=z_0}\left[\frac{1}{(z-z_0)^2}F(z)F(z_0)^{-1}\right] &= F'(z_0) F(z_0)^{-1}, \\
	\Res_{z=z_0}\left[\frac{1}{(z-z_0)^2} F(z_0)^{-1}F(z)\right] &= F(z_0)^{-1} F'(z_0).
\end{align*}
If instead of $F(z)F(z_0)^{-1}$ or $F(z_0)^{-1}F(z)$, one has a product of such factors, then each factor contributes separately in an additive way. For example, if both $F$ and $G$ are invertible and analytic at $z_0$, then 
\[
    \Res_{z=z_0}\left[\frac{1}{(z-z_0)^2}F(z)F(z_0)^{-1}G(z_0)^{-1} G(z)\right] = F'(z_0)F(z_0)^{-1} + G(z_0)^{-1}G'(z_0).
\]
Applying these rules to \eqref{eq:ET_double_pole}, we find that \eqref{eq:ET_double_pole}, and hence \eqref{eq:EX_double_pole}, can be expressed as a sum of five terms:
\begin{multline}\label{eq:E_double_pole_5_terms}
	-\Tr\left[\left.E^{-1}(z) E'(z)\right|_{z=z(p_2)} \Xi(p_2)\right] \\
    + \Tr\left[\Xi(p_2) \begin{pmatrix} 0_3 & I_3 \end{pmatrix}\left. \mathbf{T}_N^{-1}(z)\frac{d}{dz} \mathbf{T}_N(z) \right|_{z=z(p_2)} \begin{pmatrix} 0_3 \\ I_3 \end{pmatrix}\right] \\
	- N\frac{dg}{dz}(p_2) - N(1-\xi_{1,N})\frac{1}{\lambda(p_2)}\frac{d\lambda}{dz}(p_2) + N(1-\xi_{2,N})\frac{1}{z(p_2)}.
\end{multline} 

Dividing by $N$ and letting $N\to\infty$, the first term in \eqref{eq:E_double_pole_5_terms} vanishes in the limit, and the last three terms combine in the limit to give
\begin{equation}\label{eq:E_double_pole_limit}
    -\frac{dg}{dz}(p_2) - (1-\xi_1)\frac{1}{\lambda(p_2)}\frac{d\lambda}{dz}(p_2) + (1-\xi_2)\frac{1}{z(p_2)} = \left(\frac{\omega_{(\xi_1,\xi_2)}}{dz}\right)\!(p_2).
\end{equation}
Here, the second equality follows from \eqref{eq:omega_xi1_xi2_1} and \eqref{eq:omega_as_dg_dlogl_dlogz}.

Recall from Lemma~\ref{lem:TN_bounded} that the outcome of the steepest descent analysis of the RH problem is that $\mathbf{T}_N$ and $\mathbf{T}_N^{-1}$ remain uniformly bounded as $N\to\infty$. Then, by Cauchy's theorem, $\frac{d}{dz} \mathbf{T}_N(z)$ also remains uniformly bounded. Thus, the second term in \eqref{eq:E_double_pole_5_terms} remains bounded, and since we divide by $N$ its contribution vanishes in the limit as $N\to\infty$.

We have shown that \eqref{eq:ET_double_pole}, and therefore \eqref{eq:EX_double_pole}, converges to \eqref{eq:E_double_pole_limit} as $N\to\infty$. It follows that
\begin{equation}\label{eq:HN_sN_limit}
	\lim_{N\to\infty}\frac{1}{N}\mathsf{H}_N^{(s_N)}(v_N) = -\frac{1}{2\pi i}\int_{\widetilde{\gamma}_s}\omega_{(\xi_1,\xi_2)},
\end{equation}
where $\widetilde{\gamma}_s$ is a curve from $\bar{s}$ to $s$ satisfying the properties stated in the caption of Figure~\ref{fig:double_torus_EhN_contours}.


\subsection{Completion of the proof of Theorem~\ref{thm:limit_shape}(b)}
\begin{proof}
    Dividing \eqref{eq:EhN_vN_ANBNjCN} by $N$, we find
    \[
    	\frac{1}{N}\mathbb{E}[h_N(v_N)] = 3 - \frac{1}{N}\mathsf{H}_N^{(B,C)}(v_N) - \frac{1}{N}\mathsf{H}_{N,1}^{(B,A)}(v_N) - \frac{1}{N}\mathsf{H}_{N,2}^{(B,A)}(v_N) - \frac{1}{N}\mathsf{H}_N^{(s_N)}(v_N).
    \]
    As argued in Remark~\ref{rem:HN_BC_small} and Section~\ref{sec:estimate}, the second, third, and fourth terms vanish as $N\to\infty$. Hence, the limiting behavior of $\frac{1}{N}\mathbb{E}[h_N(v_N)]$ as $N\to\infty$ is determined entirely by the first and fifth terms, which together converge to the right-hand side of \eqref{eq:limit_shape_2} by \eqref{eq:HN_sN_limit}. This is equivalent to \eqref{eq:limit_shape_1}, which then establishes the first part of Theorem~\ref{thm:limit_shape}(b).

    It remains to determine the gradient of the limit shape $H(\xi_1,\xi_2)$ for $(\xi_1,\xi_2)\in\mathfrak{L}$. We recall that
    \begin{equation}\label{eq:limit_shape_3}
        H(\xi_1,\xi_2) = \frac{1}{2\pi i}\int_{\gamma_s} \omega_{(\xi_1,\xi_2)},
    \end{equation}
    with $s = \Omega(\xi_1,\xi_2)$. The first variation of \eqref{eq:limit_shape_3} due to the dependence of $\gamma_s$ on the parameters vanishes because the differential $\omega_{(\xi_1,\xi_2)}$ has zeros at the endpoints. This implies \eqref{eq:limit_shape_gradient}, since
    \[
        \frac{\partial}{\partial\xi_1} \omega_{(\xi_1,\xi_2)} = -d\log\lambda, \quad \frac{\partial}{\partial\xi_2}\omega_{(\xi_1,\xi_2)} = d\log z,
    \]
    which follows from the precise form~\eqref{eq:omega_xi1_xi2_1} of $\omega_{(\xi_1,\xi_2)}$.
\end{proof}

\section*{Acknowledgments}
We are grateful to Lennart H\"ubner for developing, together with the first-named author, the Python code that allows reliable large-size simulations of random tilings. The code is based on an earlier Matlab code of Christophe Charlier, and we thank him for making this code available to us. We used the Python code to generate Figure~\ref{fig:tiling_simulations_1} and Figure~\ref{fig:tiling_simulations_2} in the paper.

The research of M.~J.~van Horssen is supported by the Research Foundation Flanders (FWO) grant 1139226N, and the research of A.~B.~J.~Kuijlaars is supported by the long term structural funding-Methusalem grant of the Flemish Government.



\begin{thebibliography}{99}
	\bibitem{Agg23} A. Aggarwal,
    Universality for lozenge tiling local statistics,
    Ann. Math. 198 (2023), 881--1012.

	\bibitem{Bai23} E. Bain, 
	Local correlation functions of the two-periodic weighted Aztec diamond in mesoscopic limit,
	J. Math. Phys. 64 (2023), paper 023301, 64pp.

	\bibitem{BCJ18} V. Beffara, S. Chhita, and K. Johansson, 
	Airy point process at the liquid--gas boundary, 
	Ann. Probab. 46 (2018), 2973--3013.

	\bibitem{Ber21} T. Berggren,
    Domino tilings of the Aztec diamond with doubly periodic weightings, 
	Ann. Probab. 49 (2021), 1965--2011.

	\bibitem{BB25} T. Berggren and A. Borodin, 
	Geometry of the doubly periodic Aztec dimer model,
	Comm. Amer. Math. Soc. 5 (2025), 475--570.

	\bibitem{BBG26+} T. Berggren, A. Borodin, and T. George,
    Dimer models on astroidal zig-zag graphs,
    preprint arXiv:2605.03986.

	\bibitem{BB26+} T. Berggren and N. Bradinoff,
	Marked GUE-corners process in doubly periodic dimer models,
    preprint arXiv:2603.27906.

    \bibitem{BD19} T. Berggren and M. Duits, 
    Correlation functions for determinantal processes defined by infinite block Toeplitz minors, 
    Adv. Math. 356 (2019) 106766. 

	\bibitem{BN25+} T. Berggren and M. Nicoletti,
	Gaussian Free field and discrete Gaussians in periodic dimer models,
    preprint arXiv:2502.07241.

	\bibitem{BGK23} M. Bertola, A. Groot, and A.B.J. Kuijlaars,
	Critical measures on higher genus Riemann surfaces,
	Comm. Math. Phys. 404 (2023), 51--95.

	\bibitem{BB26} A.I. Bobenko and N. Bobenko,
    Dimers and M-curves: limit shapes from Riemann surfaces,
    Duke Math. J. 175 (2026), 1861--1930.

	\bibitem{BD23} A. Borodin and M. Duits,
    Biased $2\times2$ periodic Aztec diamond and an elliptic curve,
    Probab. Theory Related Fields 187 (2023), 259--315.

	\bibitem{BCdT23} C. Boutillier, D. Cimasoni, and B. de Tili\`ere, 
	Minimal bipartite dimers and higher genus Harnack curves,
    Probab. Math. Phys. 4 (2023), 151--208.

    \bibitem{BdT27} C. Boutillier and B. de Tili\`ere, 
    Fock's dimer model on the Aztec diamond, 
    preprint arXiv:2405.20284,
    to appear in Ann. Inst. Henri Poincar\'e D.

	\bibitem{CM12} G.A. Cassatella-Contra and M. Ma\~nas,
	Riemann--Hilbert problems, matrix orthogonal polynomials and discrete matrix equations with singularity confinement,
    Stud. Appl. Math. 128 (2012), 252--274.

	\bibitem{Cha21} C. Charlier,
    Doubly periodic lozenge tilings of a hexagon and matrix-valued orthogonal polynomials, 
	Stud. Appl. Math. 146 (2021), 3--80.

	\bibitem{CDKL20} C. Charlier, M. Duits, A.B.J. Kuijlaars, and J. Lenells,
	A periodic hexagon tiling model and non-Hermitian orthogonal polynomials,
	Comm. Math. Phys. 378 (2020), 401--466.

	\bibitem{CDF26+} S. Chhita, D. Dauvergne, and T. Finn,
	The Airy line ensemble at the rough--smooth boundary,
    preprint arXiv:2602.24063.

	\bibitem{CJ16} S. Chhita and K. Johansson,
	Domino statistics of the two-periodic Aztec diamond,
	Adv. Math. 294 (2016), 37--149.

	\bibitem{Dei99} P.A. Deift, 
	Orthogonal polynomials and random matrices: a Riemann--Hilbert approach,
    Courant Lecture Notes in Mathematics, vol. 3, New York University, Courant Institute of Mathematical Sciences, New York; American Mathematical Society, Providence, RI, 1999.

	\bibitem{DKMVZ99} P. Deift, T. Kriecherbauer, K.T.-R. McLaughlin, S. Venakides, and X. Zhou,
    Uniform asymptotics for polynomials orthogonal with respect to varying exponential weights and applications to universality questions in random matrix theory,
	Comm. Pure Appl. Math. 52 (1999), 1335--1425. 

	\bibitem{DZ93} P. Deift and X. Zhou,
	A steepest descent method for oscillatory Riemann--Hilbert problems, asymptotics for the MKdV equation,
    Ann. Math. 137 (1993), 295--368.

	\bibitem{DK21} M. Duits and A.B.J. Kuijlaars, 
	The two-periodic Aztec diamond and matrix-valued orthogonal polynomials, 
	J. Eur. Math. Soc. 23 (2021), 1075--1131.

    \bibitem{Ebe07} W. Ebeling,
    Functions of Several Complex Variables and Their Singularities,
    Amer. Math. Soc., Providence RI, 2007.

	\bibitem{EM98} B. Eynard and M.L. Mehta, 
	Matrices coupled in a chain. I. Eigenvalue correlations,
	J. Phys. A 31 (1998), 4449--4456.

	\bibitem{FIK92} A.S. Fokas, A.R. Its, and A.V. Kitaev,
	The isomonodromy approach to matrix models in 2D quantum gravity,
    Comm. Math. Phys. 147 (1992), 395--430. 

	\bibitem{GV85} I. Gessel and G. Viennot, 
	Binomial determinants, paths, and hook length formulae,
	Adv. Math. 58 (1985), 300--321.

	\bibitem{Gor21} V. Gorin,
    Lectures on random lozenge tilings,
    Cambridge Studies in Advanced Mathematics, vol. 193, Cambridge University Press, Cambridge, 2021.

	\bibitem{GK21} A. Groot and A.B.J. Kuijlaars,
    Matrix-valued orthogonal polynomials related to hexagon tilings, 
	J. Approx. Theory 270 (2021), paper no. 105619, 36pp.

	\bibitem{GIM11} F.A. Gr\"unbaum, M.D. de la Iglesia, and A. Mart\'inez-Finkelshtein,
    Properties of matrix orthogonal polynomials via their Riemann--Hilbert characterization,
	SIGMA 7 (2011), paper 098, 31pp.

	\bibitem{Hua26+} J. Huang,
    Height fluctuations of lozenge tilings of polygons,
    preprint arXiv: 2609.04789.

	\bibitem{Joh02} K. Johansson, Non-intersecting paths, random tilings and random matrices, Probab. Theory  Related Fields 123 (2022), 225--280.

	\bibitem{JM23} K. Johansson and S. Mason,
	Dimer--dimer correlations at the rough-smooth boundary,
	Comm. Math. Phys. 400 (2023), 1255--1315.

	\bibitem{KO06} R. Kenyon and A. Okounkov, 
	Planar dimers and Harnack curves, 
	Duke Math. J. 131 (2006), 499--524.

	\bibitem{KO07} R. Kenyon and A. Okounkov,
	Limit shapes and the complex Burgers equation, 
	Acta Math. 199 (2007), 263--302.

	\bibitem{KOS06} R. Kenyon, A. Okounkov, and S. Sheffield,
	Dimers and amoebae,
    Ann. Math. 163 (2006), 1019--1056.

	\bibitem{Kui25} A.B.J. Kuijlaars,
	Matrix-valued orthogonal polynomials arising from hexagon tilings with $3\times3$-periodic weightings,
	J. Approx. Theory 311 (2025), 54pp.

	\bibitem{KP25} A.B.J. Kuijlaars and M. Piorkowski, 
	Wiener--Hopf factorizations and matrix-valued orthogonal polynomials,
	Probab. Math. Phys. 6 (2025), 547--580. 

	\bibitem{Lin73} B. Lindstr\"om,
	On the vector representations of induced matroids,
	Bull. London Math. Soc. 5 (1973), 85--90.

	\bibitem{Mon21} J. Montaldi,
	Singularities, Bifurcations and Catastrophes,
	Cambridge Univ. Press, Cambridge, 2021.

	\bibitem{Pet14} L. Petrov,
    Asymptotics of random lozenge tilings via Gelfand--Tsetlin schemes,
    Probab. Theory Related Fields 160 (2014), 429--487.

	\bibitem{Pet15} L. Petrov,
    Asymptotics of uniformly random lozenge tilings of polygons. Gaussian free field,
    Ann. Probab. 43 (2015), 1--43.

	\bibitem{Pio26} M. Piorkowski,
	Arctic curves of periodic dimer models and generalized discriminants,
    Selecta Math. 32 (2026), 40pp.
\end{thebibliography}
\end{document}